\documentclass[12pt,a4paper]{amsart}
 \usepackage{url} 
 \usepackage{framed}
\usepackage[T1]{fontenc}
\input xypic
\usepackage{amsmath}
\usepackage{amssymb}
\usepackage[usenames]{color}
\usepackage{amssymb}
\usepackage{latexsym}
 \usepackage{graphicx}
\usepackage{psfrag}
\newcommand{\canberemoved}[1]{}

\newcommand{\ggtext}{}
\newcommand{\ftext}{}
\newcommand{\vtext}{}
\newcommand{\guntext}{}

\newcommand{\mmmmtext}{}
\newcommand{\mmmtext}{}
\newcommand{\mmtext}{}
\newcommand{\mtext}{}
\newcommand{\muutos}{}
\newcommand{\slava}{}

\newcommand{\reftext}{}
\newcommand{\refHOX}[1]{}

\newcommand{\rtext}{}

\newcommand{\rbtext}{}

\def\itext{}
\def\gtext{} 
\def\mattitext{} 
\newcommand{\revtext}{}

\newcommand{\removedEinstein}[1]{}

\def\wflat{\mmmmtext w_0}
\def\hatzeta{\eta}
\def\hatomega{\omega}

\def\nohat{}
\def\vel{}
\def\vex {}
\def\vey {}
\def\vez {}

\newcommand{\noextension}{} 

\newcommand{\arxivpreprint}{} 

\newcommand{\motivation}{} 

\newcommand{\modified}{}

\newcommand{\extension}[1]{} 
\newcommand{\MTEXT}{}

\newcommand{\MATTITEXT} {} 

\newcommand{\observation}[1]{}
\newcommand{\generalizations}[1]{}

\newcommand\underhat{\widetilde} 
\def\g{{\gamma}}
\newcommand\firstpoint{{\mathcal E}} 
\newcommand{\eeta}{\eta}
\newcommand{\bbf}{{\bf f}}
 \def\bs{{\bf s}}

\newcommand{\q}{{\underhat q}}
\newcommand{\bW}{{\underhat W}}
\newcommand\qpoint{\modified{q}}

\newcommand{\Smod}{\rtext S_0}
\newcommand{\ryksi}{r_2}
\newcommand{\rkaksi}{r_1}

\newcommand{\aell}{\ell}
\newcommand{\solw}{w}

\newcommand{\bQ}{Q}
\def\cirt{\overline }
\def\Hsymbol{H}
\def\tsparameter{s}
\newcommand{\symbolP}{v}
\newcommand{\symbolQ}{w}

\newcommand{\cell}{\ell}
\newcommand{\hattuM}{M}
\newcommand{\Sclo}{\mathcal S_{cl}}
\newcommand{\Scle}{\mathcal S_e}

\newcommand{\xxi}{\xi}
\newcommand{\zzeta}{\zeta}

\newcommand{\Ric}{\hbox{Ric}}
\newcommand{\Ein}{\hbox{Ein}}
\def \Box {$\square$}

\newcommand{\newcorrection}{\color{red}} 

\newcommand{\newtext}{}
\newcommand{\tobecheckedtext}{}

\newcommand{\HOX}[1]{}

\newcommand{\hiddenfootnote}[1]{}

\newcommand{\rhoepsilon}{{\rho}}
\newcommand{\bsequence}{{\bf b}}
\def\X{{\mathcal X}}
\def\Y{{\mathcal Y}}
\def\U{{\mathcal U}}
\def\G{{\mathcal G}}

\def\S{{\mathcal S}}
\def\T{{\mathcal T}}
\def\P{{\mathcal P}}
\def\qP {{\mathcal D}}

\def\pointear{{\bf e}}

\renewcommand{\H}{{\mathcal H}} 
\newcommand{\R}{{\mathbb R}} 
\newcommand{\D}{{\cal D}} 
\newcommand{\I}{{\cal I}} 
\newcommand{\be}{{\mathcal C}} 
\newcommand{\E}{{\cal E}} 
\newcommand{\A}{{\cal A}}
\newcommand{\B}{{\cal B}}  
  
\newcommand{\K}{{\cal K}}  

\newcommand{\C}{{\mathbb C}} 
\newcommand{\N}{{\mathbb N}} 
\newcommand{\cal}{\mathcal }

\renewcommand{\L}{{\mathcal L}} 
\newcommand{\V}{{\cal N}} 
\newcommand{\W}{{\cal W}}

\def\hat{\widehat}
\def\tilde{\widetilde}
\def \bfo {\begin {eqnarray*} }
\def \efo {\end {eqnarray*} }
\def \ba {\begin {eqnarray*} }
\def \ea {\end {eqnarray*} }
\def \beq {\begin {eqnarray}}
\def \eeq {\end {eqnarray}}
\def \supp {\hbox{supp}\,}
\def \dim{\hbox{dim}\,}

\def \re {\hbox{Re}\,}
\def \im {\hbox{Im}\,}

\def \WF {\hbox{WF}\,}

\def \dist {\hbox{dist}}
\def\diag{\hbox{diag }}
\def \det {\hbox{det}}

\def\bra{\langle}
\def\cet{\rangle}

\def \e {\varepsilon}
\def \p {\partial}
\def \a {\alpha}
\renewcommand{\b}{\beta}

\def\M{{\mathcal U}}
\def\F{{\mathcal F}}

\def\Z{{\mathbb Z}}

\newtheorem{definition}{Definition}[section] 
\newtheorem{theorem}[definition]{Theorem} 
\newtheorem{lemma}[definition]{Lemma} 
 
\newtheorem{proposition}[definition]{Proposition} 
\newtheorem{corollary}[definition]{Corollary}

\begin{document}
\title[Inverse problems for Lorentzian manifolds] {Inverse problems for Lorentzian manifolds and non-linear hyperbolic equations}
\date{July 28, 2017}
\author{Yaroslav Kurylev, Matti Lassas,  Gunther Uhlmann}
\address{Yaroslav Kurylev, UCL; Matti Lassas,  University of Helsinki;
Gunther Uhlmann, University of Washington,  and 
 University of Helsinki.  {\rm y.kurylev@ucl.ac.uk, Matti.Lassas@helsinki.fi,  
  gunther@math.washington.edu}}

\email{}

\maketitle

{\bf Abstract:} 
{\it  
We study two inverse problems on a globally hyperbolic Lorentzian manifold $(M,g)$. 
 The problems are: 

1. Passive observations in spacetime: Consider
observations in  an open set $V\subset M$.
  The light observation set corresponding to a point source at $q\in M$
is the intersection of {\muutos $V$ and}
 the light-cone emanating from the point $q$.
  Let $W\subset M$ be an unknown open, relatively compact set.
We show that under natural causality conditions,
the family of light observation sets corresponding to point sources at {\muutos points} $q\in W$
determine uniquely 
the conformal type of $W$.

2.  Active measurements in spacetime: We develop a  new method  for inverse problems for non-linear hyperbolic equations 
that utilizes the
non-linearity as a tool. This enables us to solve inverse problems for non-linear equations
for which the corresponding problems for linear equations are still unsolved.
To illustrate this method, we solve an inverse problem for   semilinear wave equations
with quadratic non-linearities. 
We assume that we are given the neighborhood $V$ of the time-like {\rtext path} $\mu$ and  
the source-to-solution operator 
that maps the source
supported on  $V$ to the restriction of the solution of the wave equation {\vtext to} $V$. 
When $M$    is 4-dimensional, we show that these data
determine the topological, differentiable, and conformal structures of the spacetime in  
the maximal set where waves can propagate from $\mu$ and return back
to $\mu$.}



\tableofcontents

\section{Introduction and main results}

We study the question of whether an observer in spacetime can
determine
the structure of the surrounding spacetime by doing measurements near  its
world line. We consider two kinds of problems: inverse problems  for active measurements
and for  passive observations.

For active measurements, we consider the  wave equation
\beq \label{eq: wave-eq general pre}
& &\hspace{-.5cm}  \square_{g}u(x)+a(x)\,u(x)^2=f(x)
\quad\hbox{on }{\rtext M},\hspace{-.5cm}\\ \nonumber
 & &\quad 
u(x)=0,\quad \hbox{outside causal future of }\supp(f),
\eeq
on the Lorentzian manifold $({\rtext M},g)$, 
 a {\gtext future-pointing} time-like {\rtext path} $\hat \mu=\hat \mu([-1,1])\subset {\rtext M}$ and 
an open neighborhood  
$V$ of $\hat \mu$.  {\mattitext The wave equation is considered as a model problem
for which we demonstrate the new techniques {\gtext we develop.
We will consider in a follow up paper the Einstein
equations coupled with scalar fields with applications to general relativity
and} other physical models involving non-linear hyperbolic
equations, see \cite{E-preprint}.} 
 We assume
 that 
we can control sources supported in
$ V$ and measure the physical fields in the same set $V$. Our aim is to 
determine  the conformal class of the metric (or even the metric 
tensor in some cases)
in a suitable {\mattitext larger 
set $J=J^+(p^-)\cap J^-(p^+)$, {\muutos that is the set of the points that are in
the causal future of the point $p^-=\hat \mu (s_-)$ and in the causal past
of the point $p^+=\hat \mu (s_+)$}, {\rtext where $-1<s_-<s_+<1$}, see Fig.\ 1(Left).}
We study the inverse problem for active measurements by considering
the interaction of distorted plane wave packets (Fig.\ 1, Right) reducing the problem to the
problem for passive observations.


{\rtext The new method we introduce in this paper  for inverse problems for non-linear hyperbolic equations 
 {\mmmmtext utilises the
non-linearity as a tool.} This enables us to solve inverse problems for non-linear equations
for which the corresponding problems for linear equations are still unsolved.
Indeed, the existing uniqueness results for linear hyperbolic equations
with vanishing initial data
 are limited to
the time-independent or real-analytic coefficients, {\muutos see e.g.}\ \cite{AKKLT,Bel1,BelKur,Eskin,Eskin2,KKL} since 
these results are based on Tataru's unique 
continuation theorem \cite{Tataru1,Tataru2}. Such unique continuation results 
have   been shown to fail for general metric
tensors which are not analytic in the time variable \cite{Alinhac}.}

{\mmtext 
{\guntext

 \begin{center}

\psfrag{1}{$p^-$}
\psfrag{2}{$p^+$}
\psfrag{4}{$J$}
\psfrag{3}{\hspace{-0.1cm}V}
\psfrag{5}{}
\psfrag{6}{}
\includegraphics[height=5.5cm]{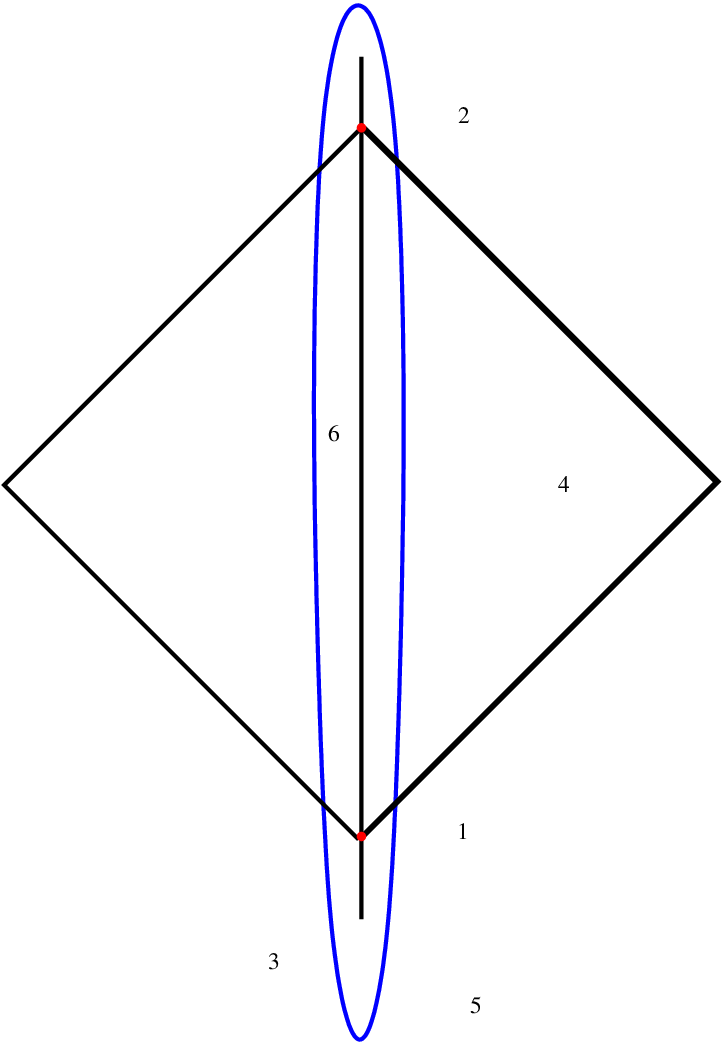}
\psfrag{1}{$y_0$}
\psfrag{2}{}
\psfrag{5}{$\Sigma_0$}
\psfrag{3}{$\Sigma_1$}
\psfrag{4}{$x_0$}
\includegraphics[height=3.5cm]{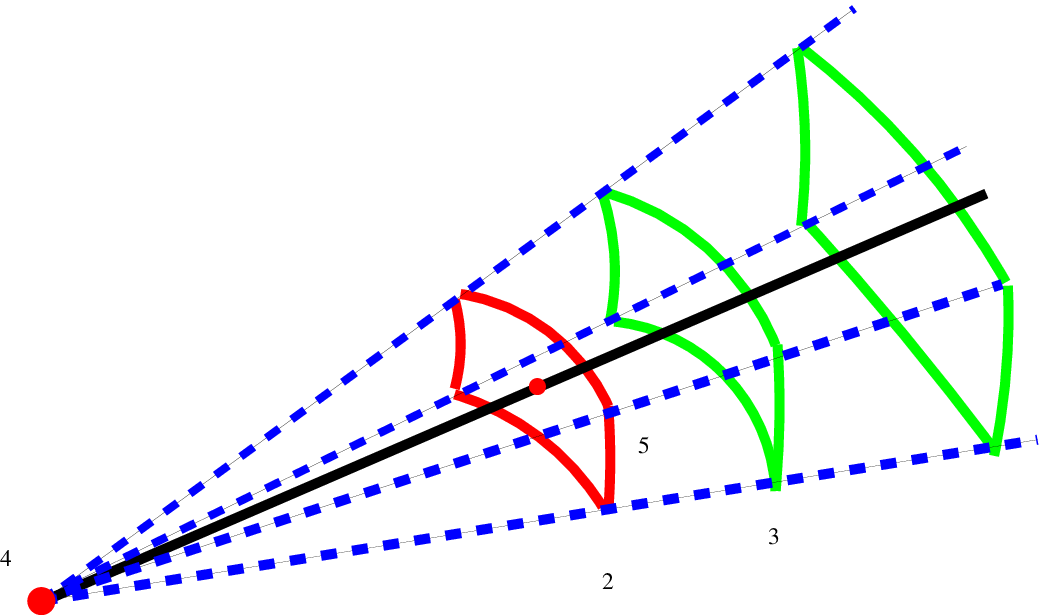}
\end{center}
{\it FIGURE 1. 
{\bf Left:} This is a schematic figure in  $\R^{1+1}$.  
  The black vertical line is {\muutos a time-like path} 
 $\hat \mu$ {\muutos that}  contains the points $p^-$ and $p^+$. The neighborhood $V$ of $\hat \mu$
is marked by a blue curve.
The black ``diamond'' is the set \refHOX{Comment 1: Referee asked us to consider colour figures.}
$J=J(p^-,p^+)=J^+({p^-})\cap J^-({p^+})$.

{\bf Right:} 
This is a schematic figure in the space $\R^{3}$.
It describes the location of a {distorted plane wave (or a piece of a spherical wave)} $u_1$ at different 
time moments. This wave
propagates near a {light-like} geodesic $\gamma_{x_0,\zeta_0}((0,\infty))\subset \R^{1+3}$,
$x_0=(y_0,t_0)$
 and is singular on
 a subset of a light cone emanating {\ftext from $x_0$.} 
 {The black line segment is the projection of $\gamma_{x_0,\zeta_0}\subset \R^{1+3}$  in to $\R^3$.}
  The piece of the distorted plane wave is {\gtext sent} from 
  {\muutos
  the point $y_0$  at the   time $t_0$ and at a later time moment $t_1$ the wave is singular on the red surface
  $\Sigma_0\subset \R^3$.
  At later time moments, it is singular on the green surfaces, like on
 $\Sigma_1$. 
 }%
\smallskip
}

{\ggtext
{\vtext Our method to solve inverse problems for the non-linear wave equation
with active measurements is to apply
global Lorentzian geometry, our results on the inverse problem for passive measurements,
and the results on the 
 non-linear interaction of
non-smooth waves having conormal singularities.}
There are many results on such non-linear interaction,  starting from the studies of
 Bony \cite{Bony}, Melrose and Ritter \cite{MR1,MR2} and Rauch and Reed,  \cite{R-R}.
However, these studies are different {\vtext from the present paper that in these papers it is assumed}  
 that the geometrical setting of the interacting singularities, {\vtext and in particular the locations and types of caustics,} is a priori known. 
 {\vtext In inverse problems  we study waves on
an unknown manifold,
so we do not know the  underlying geometry and, therefore, the location of singularities
of the fields. 
For example, the waves can have caustics that
may even be of an unstable type.
These produce further difficulties in the analysis of the non-linear interaction.
To overcome these difficulties we use methods of the 
 global Lorentzian geometry, results for the passive inverse problem, and the 
 layer-stripping method. These make it possible to reconstruct the accessible part of the Lorentzian manifold step by step.}}}

The inverse problem for passive observations 
 means the reconstruction of a region $W$ of a 
 Lorentzian manifold from light observation sets  $\P_V(q)$ corresponding
 to points $q\in W$. 
 The light observation 
 set
 $\P_V(q)$ is the intersection of set $V$ and the future light cone $\L^+(q)$ emanating from the source  point $q$. 
Physically,  this corresponds to the case of a passive observer, who registers in the set $V$ 
light {\mattitext (or a gravitational wave)} 
coming from a source at $q$.
Due to the existence of  conjugate points (or physically speaking,
gravitational lensing or Einstein rings) such observations can by strongly distorted. 
Under appropriate conditions,
 we first show that  $W$ can be reconstructed as a topological
manifold from these data. After that, we show that the differentiable structure of $W$ and the conformal class of
$g|_W$ can be reconstructed.  


\subsection{Inverse problem for passive observations}\label{subsec: IP for passive}
To formulate the results, we first introduce some definitions.
Let   $(M,g)$ be a $n$-dimensional Lorentzian manifold of signature $(1, n-1),\, n \geq 3.$ 
In this paper we assume that $(M, g)$ is time-oriented so that  we can define
future and past pointing time-like and causal paths.\refHOX{Comment 2. Definition of time-like and causal paths is added}
{\reftext We recall that a smooth path $\mu:(a,b)\to M$ is time-like if
$g(\dot \mu(s),\dot \mu(s))<0$  for all $s\in (a,b)$. Also, $\mu$ is causal if $g(\dot \mu(s),\dot \mu(s))\leq 0$
and $\dot \mu(s)\not =0$ for all for all $s\in (a,b)$.}
 For $ p, q \in M$ we denote  $p\ll q$
if $p\not =q$ and there is a future pointing time-like path from $p$ to $q$. We denote
$p < q$, if $p\not =q$ and there is a future pointing  causal path from $p$ to $q$ and denote $p \leq q$
when either $p=q$ or $p <q$. The chronological future of $p \in M$ is
the set $I^+(p)=\{q\in M«;\ p\ll q\}$ {\muutos and} the causal
future of $p$ is the set  $J^+(p)=\{q\in M;\ p \leq q\}$.
Similarly, we introduce the chronological past, $I^-(p)$, and the causal past,
$J^-(p)$, see  \cite{ONeill}.
Note that $I^{\pm}(p)$ are always open. 
For a set $A\subset M$ we denote $J^\pm(A)=\cup_{p\in A}J^\pm(p)$.
We also denote $J(p,q):=J^+(p)\cap J^-(q)$ and $I(p,q):=I^+(p)\cap I^-(q)$.

By  \cite{Bernal}, a time-orientable Lorentzian manifold $(M,g)$ is globally hyperbolic
if and only if there are no closed causal paths in $M$ and
for all $q_1,q_2\in M$ such that $q_1<q_2$ the set 
$J(q_1,q_2)\subset M$ is compact.
{\revtext  Roughly speaking the last property means
that $M$ has no naked singularities which 
one could reach by moving along a time-like path starting from a point $q^-$ and
ending in a point $q^+$.
{\gtext 
 In particular,}
this condition is needed to make the hyperbolic equations on $(M,g)$ well
posed.}

We assume throughout the paper that $(M,g)$ is globally hyperbolic. 
In this case, $J^{\pm}(p)$
are closed and $\hbox{cl}(I^{\pm}(p))=J^{\pm}(p)$.
\canberemoved{To the knowledge of {the} authors, it is not presently known if there
exist  homeomorphic 4-dimensional, globally hyperbolic, Lorentzian 
manifolds that are not diffeomorphic but  they are known to  exist in higher dimensions, see
\cite{Chernov}. The topology of such manifolds, even with asymptotically flat
vacuum spacetimes can be very complicated, \cite{IMP1,IMP2}.}
%
%

\modified{Let $L_pM=\{\xi\in T_pM\setminus\{0\};\ g(\xi,\xi)=0\}$ {\itext be the set of light-like vectors 
in the tangent space $T_pM$}. Also, 
$L_p^+M\subset L_pM$ and $L^-_pM\subset L_pM$ denote    the future and the past light-like vectors in
$T_pM$.}  

Let $\exp_q: T_qM \to M$
be the exponential map on $(M,g)$.
The geodesic starting at $p$ in the
direction $\xi \in T_pM\setminus\{0\}$ is the curve
 $\g_{p, \xi}(t)=\exp_p(\xi t)$,  $t \geq 0$.  

Let 
$\hat \mu:[-1,1]\to M$ be a $C^\infty$-smooth future pointing time-like {\rbtext path}
and $V\subset M$ be an open {\muutos connected} neighborhood of $\hat \mu([-1,1])$.

\subsubsection{The set of earliest light observations}

{\revtext {\rtext Recall that  $p^\pm=\hat \mu(s_\pm)$ where $-1<s_-<s_+<1$.} 
Next define the light observation sets and consider in particular the case when
 $W \subset I^-(p^+) \setminus J^-(p^-)$ is a relatively compact open set,
 see Fig.\ 2 {\muutos (Right)}.


\begin{center}

\hspace{-4mm}

\psfrag{1} {}
\psfrag{2}{$V$}
\psfrag{3}{\hspace{-2mm}$q$}
\includegraphics[height=4cm]{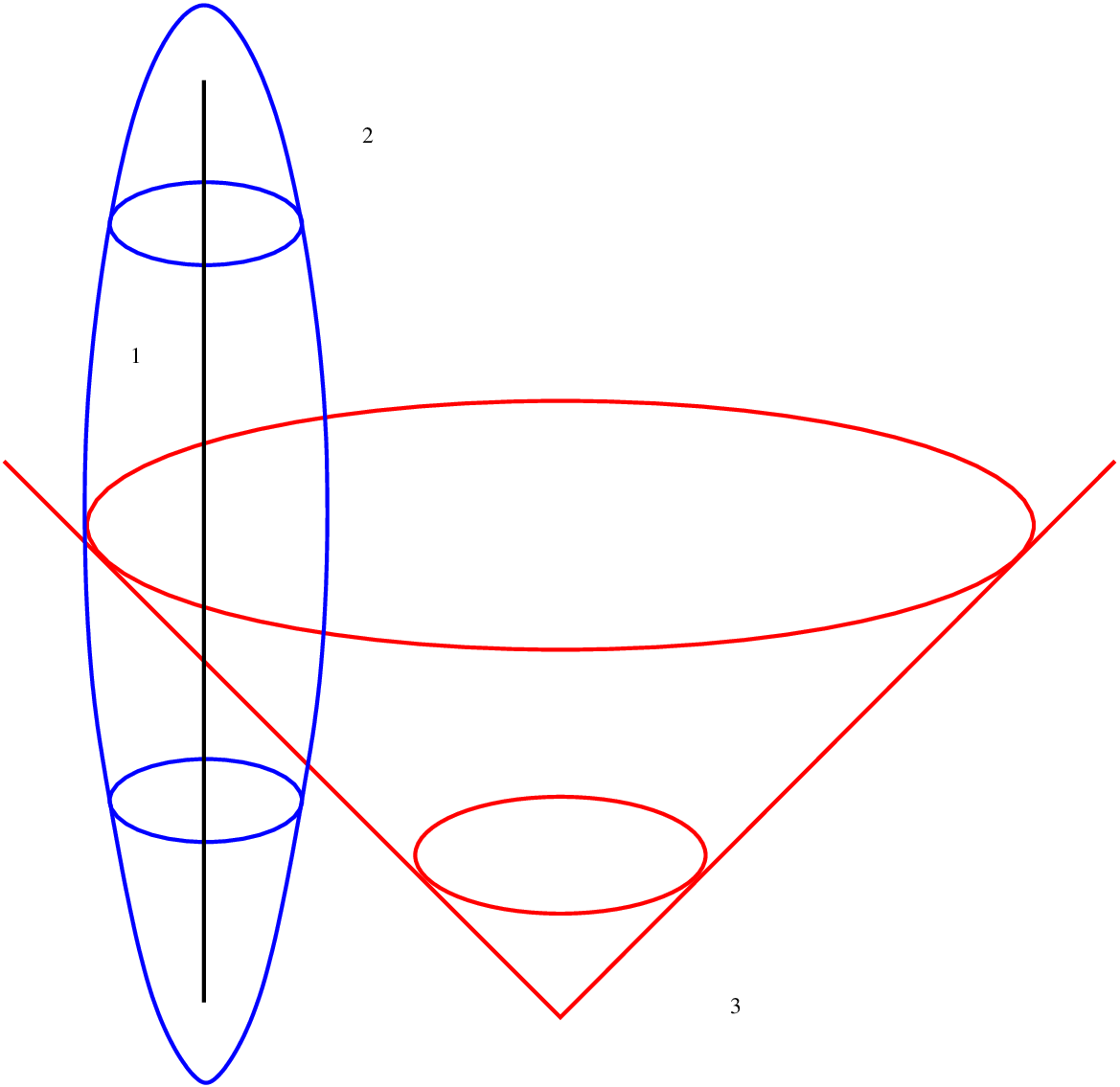}\quad
\psfrag{1}{\hspace{-1.4mm}$W$}
\psfrag{2}{\hspace{-3mm}$q$}
\psfrag{3}{}
\psfrag{6}{\hspace{-2mm}V}
\includegraphics[height=4cm]{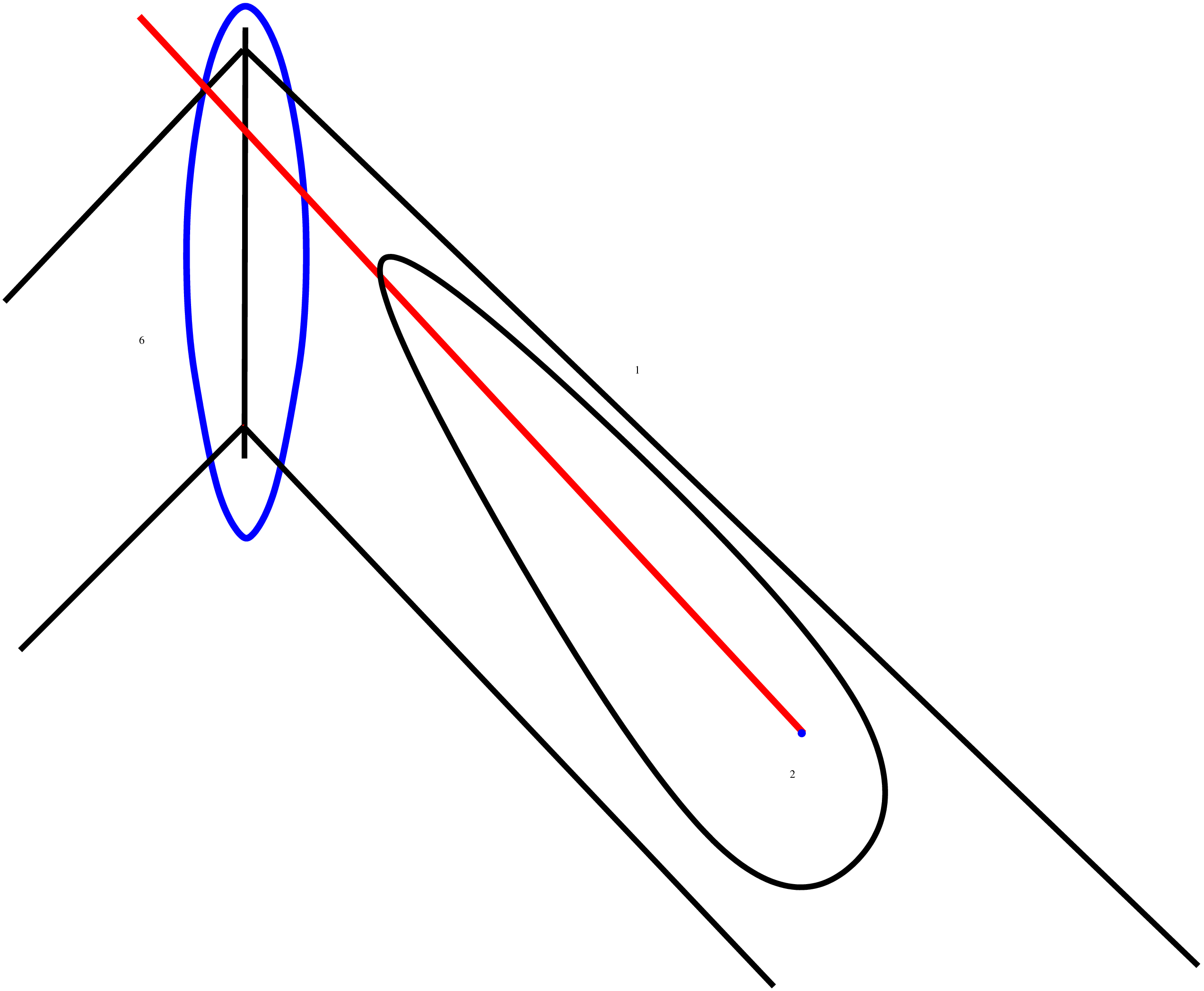}
\end{center}
{\bf Figure 2.}  \refHOX{Comment 3: Earlier Figure 2-Right is moved later.} 
{{\bf Left:} {\it The future light cone $\mathcal L^+(q)$ \modified{from} the point $q$ is  shown
as a red cone. The point $q$ is the tip of the cone. {\itext The set $V$, where observations are done,} is shown in blue.
The light observation set $\mathcal P_V(q)$ with a  point   source at $q$ is the intersection
$\mathcal L^+(q)\cap V$. {\muutos {\bf Right:}}
In Theorem\ \ref{main thm paper part 2}, we consider a set $W\subset I^-(p^+)\setminus J^-(p^-)$.
The boundary of  $W$ is shown in the figure as a black curve. The red line
is a light ray from a point $q\in W$ that is observed in the blue set $V$.
These observations are shown to determine $W$ as a differentiable
manifold and the conformal class of the metric on $W$.}


\begin{definition} \label{def. O_U}

(i) For $q\in M$, let
$$\L^+(q)=\exp_q(L^+_qM)\cup\{q\}=
\{\gamma_{q,\xi}(t)\in M;\ \xi\in L^+_qM,\ t\geq 0\}\subset M
$$ {\itext be  the future directed light-cone emanating from the point $q$.

The light observation set of $q$ in the observation set $V$} is 
\ba
\P_V(q)= {\mathcal L}^+(q) \cap V \in 2^V.
\ea

{\mattitext

(ii) The earliest light observation set of $q\in M$  in $V$ is
\beq\label{eq BB a}
& &\hspace{5mm}\firstpoint_V(q)=\{x\in \P_V(q):\ \hbox{there are no $y\in\P_V(q)$ and}\\
\nonumber
& &\hspace{25mm}\hbox{
future-pointing time-like path 
$\alpha:[0,1]\to V$}\\
\nonumber
& &\hspace{25mm}\hbox{ such that $\a(0)=y$ and $\a(1)=x$}\}\subset V
\eeq 

}

(iii) 
Let $W\subset M$ be open.
 The {\rtext family} of the earliest light observation  sets with source points at $W$ is 
\beq
\firstpoint_V(W)=\{\firstpoint_V(q);\ q\in W\}\subset 2^V.\label{collection}
\eeq
\MTEXT{Note {\muutos that} $\firstpoint_V(W)$
is defined as an unindexed set, that is, for an element $\firstpoint_V(q)\in \firstpoint_V(W)$ we do not
know what is the corresponding point $q$.}
\end{definition}
Above, $2^V=\{V^\prime;\ V^\prime \subset V\}$ is the power set of $V$. 
{\muutos 
Note that when {\mmmmtext the future directed path $\mu:[-1,1]\to V$  and} the conformal type of $(V,g)$, and therefore all time-like paths in $V$ are known,
the light observation set  $\P_V(q)$ determines the earliest light observation set  $\E_V(q)$, see (\ref{eq BB a}).}

%
Below, when $\Phi:V_1\to V_2$ is a map, we 
say that the power set extension of $\Phi$ is the map
$\tilde \Phi: 2^{V_1}\to 2^{V_2}$ {\mmmmtext given by
\beq\label{power set ext}
\tilde \Phi(V^\prime) =\{\Phi(z); \, z \in V^\prime \},\quad\hbox{for $V^\prime \subset V_1$.}
\eeq
Below,} when we say that the set $V$ is given as a differentiable manifold,
we mean that we are given the set $V$ and the local coordinate charts on it for which
the corresponding transition maps are $C^\infty$-smooth.

{\mattitext 

\medskip 

\noindent
{\bf Inverse problem with passive observations:} {\it We assume that we are given the set $V$
as a differentiable manifold, the conformal class of the metric $g|_V$ on $V$,
and the {\rtext family} of the earliest light observation sets $\E_V(W)=\{\E_V(q)\subset V;\ q\in W\}$,
where $W\subset I^-(p^+)\setminus J^-(p^-)$ is a relatively compact open set.
{\gtext The inverse problem is whether} these data determine the set $W$ as a differentiable manifold and
the conformal class of the metric $g|_W$.}
\medskip

%

 A map $\Psi:(V_1,g_1)\to (V_2,g_2)$
is a conformal diffeomorphism if
$\Psi:V_1\to V_2$ is a diffeomorphism  and \refHOX{Comment 4: Pull-back used instead of push-forward.}
{\reftext $g_1(x)=e^{2f(x)}(\Psi^*g_2)(x)$ for some scalar function $f(x)$.}
{\mattitext The following theorem implies  
 that the {\rtext family}  $\firstpoint_{V}(W)$ of the earliest light observation sets determines uniquely
the conformal type of $(W,g|_W)$.} 
\canberemoved{
{\gtext In fact, it is enough to know the 
earliest light observation sets $\firstpoint_{V}(q)$ for a dense subset of points
$q$ in $W$, see Remark 2.2 in Section \ref{subsec: Obs. time}.}}

{\ggtext
\begin{theorem}\label{main thm paper part 2}
Let $(M_j,{g}_j)$, $j=1,2$ be two open, $C^\infty$-smooth, globally hyperbolic 
Lorentzian manifolds of dimension $n\ge 3$, 
 {\mmmmtext $\hat \mu_j:[-1,1]\to M_j$ be   {\rbtext smooth time-like paths}, and
$p_j^ \pm=\hat \mu_j(s_\pm)$.} Let {\itext the observation sets} $V_j\subset M_j$ be neighborhoods
of $\hat \mu_j([-1,1])$
 and
 $W_j  \subset M_j$ be relatively compact sets
 such that $\overline W_j \subset J^-(p^+_j)\setminus I^-(p^-_j).$
%
%
%
%
{\mattitext Let $\firstpoint_{V_j}(W_j)$ be
the families of the earliest light observations  sets
{\itext with source points at $W_j$}}, see (\ref{collection}). 

Assume that there is a conformal diffeomorphism
$\Phi:V_1\to V_2$ such that $\Phi(\hat \mu_1(s))=\hat \mu_2(s)$, $s\in [-1,1]$ 
and
\beq\label{eq: Es the same}
\tilde \Phi( \firstpoint_{V_1}(W_1))=  \firstpoint_{V_2}(W_2),
\eeq
{\mmmmtext where $\tilde \Phi$ is the  power set extension of $\Phi$, see  (\ref{power set ext}).}

Then there is a diffeomorphism $\Psi:W_1\to W_2$
such that the metric $\Psi^*{g}_2$ is conformal to ${ g}_1$ and
$\Psi|_{W_1\cap V_1}=\Phi|_{W_1\cap V_1}$.
\end{theorem}
}

When $M_j,\, j=1,2$, have  significant Ricci-flat parts, Theorem \ref{main thm paper part 2} can be \modified{strengthened.}

\begin{corollary}\label{coro of main thm original Einstein pre}
Assume that $(M _j,g_j)$ and $V_j$, $W_j$,  $j=1,2$ satisfy the conditions of Theorem \ref{main thm paper part 2}
with the resulting
 conformal map
$\Psi:W_1\to W_2$ as in  Theorem \ref{main thm paper part 2}. Moreover, assume that 
$W_j$ are Ricci-flat {\rtext and that $\Phi:V_1\to V_2$ is an isometry}.
\modified{ Also, assume that all {\muutos topological} components of $W_j$ intersect $V_j$, $j=1,2$. Then the map
$\Psi$ is an isometry.}
\end{corollary}

\subsubsection*{Idea of the proof of Theorem \ref{main thm paper part 2}} 
The proof is given in Section 2 where the topological structure of $W$
is reconstructed and in Section 5 where the differentiable structure of $W$ and
the conformal type of the metric $g|_W$ are reconstucted.
We will define a suitable smaller observation set $U\subset V$ and consider the {\rtext family} $\E_U(W)$.
The idea is to endow the set $\E_U(W)\subset 2^U$ with a {\mattitext Lorentzian} manifold structure that
makes it conformal to $W$. {\itext In other words, we consider the set $\E_U(W)$ as a manifold that is a ``copy''  of the manifold $W$. In the proofs we construct topological, differentialble and metric structures on  $\E_U(W)$}.
{\gtext To sketch the idea of the proof}, we
assume for simplicity that the manifold $(M,g)$
has {\reftext no} \refHOX{Comment 5: "not" changed to "no".} conjugate points or cut points and that $U\cap W=\emptyset$. Then, any light-like geodesic segment $\gamma_0$
in the light-cone $\L^+(q)$, i.e.,
$\gamma_0\subset \L^+(q)$, can be {\mattitext extended} to a  geodesic $\tilde \gamma_0\subset M$ that goes through the point $q$. Let us 
consider a  light-like geodesic segment  $\gamma_1\subset U$ and
define $\Theta(\gamma_1)$ to be the set of the elements $\E_U(q)\in \E_U(W)$ for which $\gamma_1\subset \E_U(q)$. Then, when $\gamma_1$  is continued 
to a maximal geodesic $\tilde \gamma_1\subset M$, we have
that  $\Theta(\gamma_1)$ is the image of the geodesic segment $\tilde 
\gamma_1\cap W$ in the map $q\mapsto \E_U(q)$. This means that
on the set $\E_U(W)$ 
{\rtext we can see the images of a open family of light-like geodesics of $M$ that intersect $U$.}
Using this we show that the map  $\E_U:q\mapsto \E_U(q)$
is one-to-one and defines a homeomorphism $\E_U:W\mapsto \E_U(W)$. 
{In this way we reconstruct the topological structure of $W$.
The differentiable structure can be reconstructed by using the 
earliest observation time functions $f^+_{a}(q)$ {\itext on a time-like {\muutos path} $\mu_a\subset U$.} The
function $f^+_{a}(q)$ is equal to the smallest parameter value $s$
for which $\mu_a(s)$ belongs in the light-cone $\L^+(q)$ emanating from
$q\in W$.
{\rtext We show that that, for any $q_0 \in W$, there are $a_{1}, \dots, a_{n}$ such that $f^+_{a_j}(q),$  $j=1, \dots, n,$ can be used as local coordinates near $q_0$.}
{Finally, we use the observation that {\itext for any $q\in W$ there is an open, conic set
of directions $\xi\in L_q^+M$ such that the geodesics $\gamma_{q,\xi}$ intersect
the observation set $U$.
As we can determine the images $\E_U(\gamma_{q,\xi}\cap W)$ of the geodesics $\gamma_{q,\xi}$
on the known manifold $\E_U(W)$,
we can determine the image of the light-cones, 
$d\E_U(L^+_qM)$, in the differential of the map $\E_U:q\mapsto \E_U(q)$. As this can be done 
for all $q\in W$, we see the images of the light-cones on the tangent bundle {\rtext $T(\mathcal E_U(W))$} of the manifold $\E_U(W)$.}}}
   {\rtext 
  Finally, we note that having in our possession light cones on $T(\mathcal E_U(W))$ we can determine the conformal class of the metric $g|_W$.}
}} 


\subsection{Inverse problems for active measurements}
\label{subsec: IP for active}
For active measurements, {\itext let} $(M,g)$ be a $4$-dimensional Lorentzian manifold of 
signature $(1, 3)$. 
\canberemoved{We develop a new method to solve  {\mattitext inverse problems for 
non-linear hyperbolic  equations.}
We demonstrate this method for {for a particular non-linear  wave equation,
the quadratic Klein-Gordon equation with  a time-depending metric.}
\removedEinstein{For the Einstein equations we consider for the matter fields with
 the simplest possible model,
the  scalar field equations.}
 }
 
\subsubsection{Inverse problem for the  non-linear wave equation}
{\itext Several  inverse problems encountered in applications are solved by constructing} the coefficients
of the equations using invariant {\itext techniques,} e.g.\ using travel time coordinates. 
{\itext This is why many mathematical inverse problems are  
 formulated in geometric terms}, that is, on manifolds, see e.g.\ \cite {AKKLT,BelKur,dsFKS,Eskin, GsB, GST,LeU}.
Even some linear inverse problem are not uniquely
solvable.  In fact,  counterexamples for these problems have
been based on the so-called transformation optics. This has led to  
models for  fixed frequency invisibility cloaks, see e.g.\ 
\cite{GKLU1,GLU} and references therein.
\canberemoved{These applications give one more motivation 
to study  {\mattitext the uniqueness of} inverse problems.}

Several physical models lead to non-linear differential equations.  In small perturbations, 
these equations can be approximated by linear equations,
and most of the previous results on hyperbolic inverse problems in the multi-dimensional case
with vanishing initial data concern 
linear models. {\rtext As noted above,
 the existing uniqueness results for linear hyperbolic equation with
 vanishing initial data are based on Tataru's unique 
continuation theorem \cite{Tataru1,Tataru2} that requires the coefficients
to be constant {\muutos or} real-analytic in the time variable, see \cite{Alinhac}.}
{Also, the studies  {\mattitext of} inverse problems for hyperbolic equations
 {\mattitext with} time-dependent coefficients that are based on other methods have been 
restricted to the case when only the lower order terms depend on time,
see \cite{Salazar,Stefanov1}, the monographs \cite{Petkov,Roach} and the references therein.}




Earlier studies on inverse problems  for  non-linear  equations  have  concerned  
  parabolic  equations
 \cite {Isakov1}, elliptic equations \cite {IsakovNachman,Kang,SunU},
and 1-dimensional hyperbolic equations \cite {Nakamura1}. 
The present paper differs from the earlier studies in   
 that in our approach we do not consider the non-linearity as a perturbation, whose effect
is small, but as a tool that helps us  solve the inverse problem {for
 multidimensional non-linear wave equations with vanishing initial data and time-dependent coefficients}. Indeed, 
  {\mattitext it is the non-linearity that}
 makes it possible to solve a non-linear inverse problem which linearized version is not
 yet solved. This is the key novel feature of this paper.

\subsubsection{Notations}\label{sec:notations 1}

Let $(M,g)$ be a $C^\infty$-smooth $(1+3)$-dim\-ens\-ion\-al 
globally hyperbolic Lorentzian manifold, where 
the metric signature of
$g$ is $(-,+,+,+)$.

By \cite{Bernal2}, the globally hyperbolic manifold $(\hattuM , g)$ {\mmtext is
isometric  to  {\mmmmtext  a smooth  manifold} $(\R\times N, h)$,
where $N$
is a 3-dimensional manifold and the metric $h$ has the form}
\beq\label{eq: product structure}
h=-\beta(t,y) dt ^2+\kappa(t,y).
\eeq
Here $\beta:\R\times  N\to (0,\infty)$ is a smooth function and 
$\kappa(t,\cdotp)$ is a Riemannian metric on $ N$ depending smoothly
on $t\in \R$, and the submanifolds $\{t^{\prime}\}\times N$ are $C^\infty$-smooth
Cauchy surfaces for all $t^{\prime}\in \R$. 
Let us next identify these isometric manifolds,
that is, we denote $\hattuM =\R\times N$.

{\rtext Below, we will consider wave equation on the spacetime
$
(-\infty,{\reftext T_0})\times N, 
$ \refHOX{Comment 8: Notation $t_0$ is changed to $T_0$}
where  ${\reftext T_0}>0$ is a fixed parameter, and consider solutions that vanish on 
$(-\infty,{\reftext 0})\times N$.}

\begin{center}

\arxivpreprint{
$  $\hspace{-10cm}
\includegraphics[height=4.0cm]{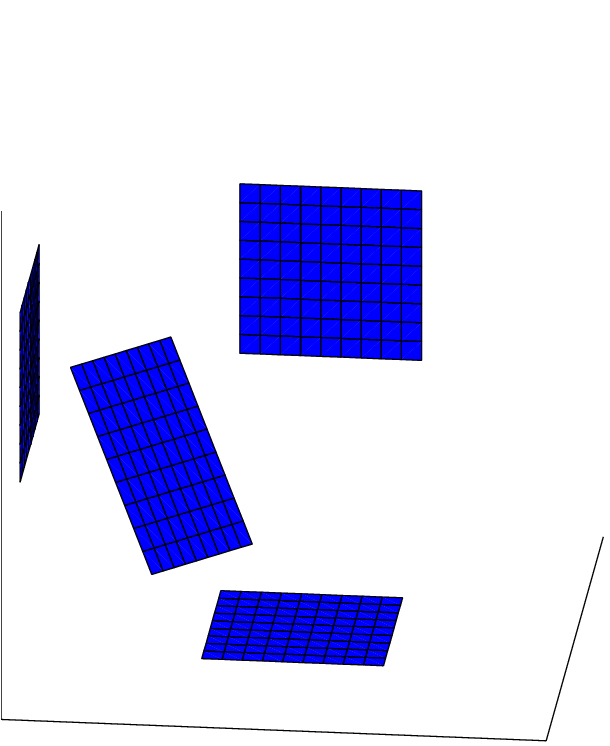}\quad 
\includegraphics[height=4.0cm]{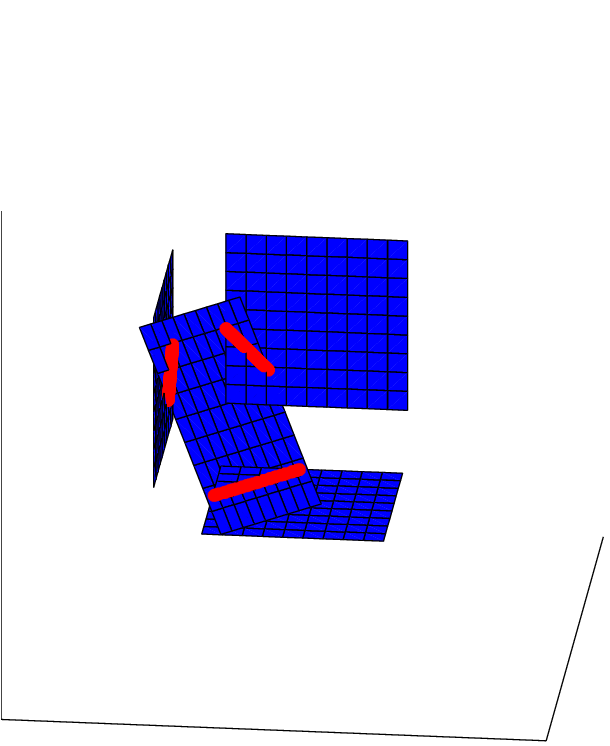}\quad
\includegraphics[height=4.0cm]{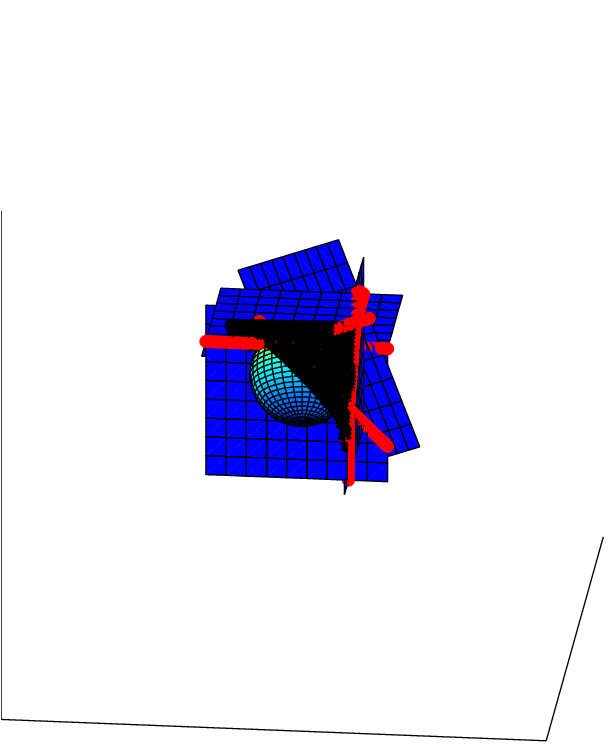}\quad
\includegraphics[height=4.0cm]{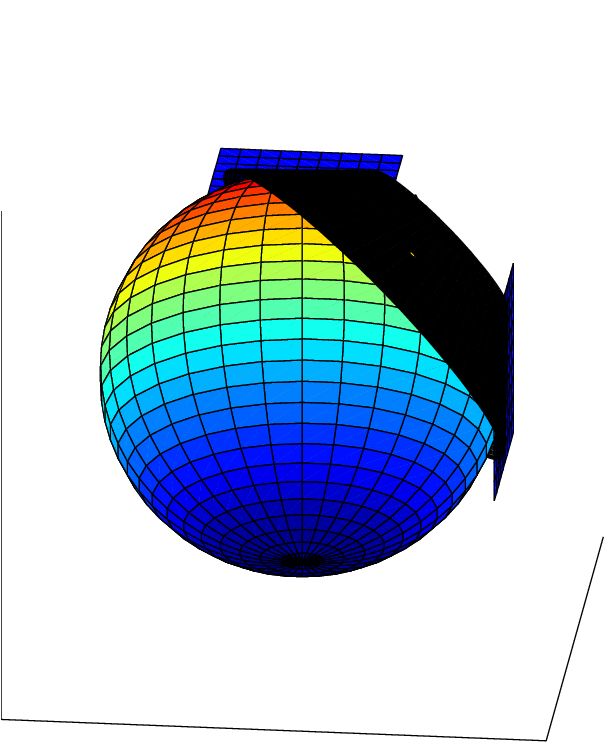}\hspace{-10cm}$  $
}
\end{center}
{\it FIGURE 3. Four plane waves propagate in space. 
When the planes intersect, the non-linearity of the hyperbolic system  produces
new waves.
The four figures show
the waves before the interaction of the waves start, when  {\mmmmtext 2}-wave interactions have  started, when all   {\mmmmtext four} waves have just interacted, and later
after the interaction. {\bf Left:} Plane waves before  interacting.
 {\bf Middle left:}
 The 2-wave interactions (red line segments) appear  but do not cause
new propagating singularities. {\bf Middle right and Right:}  All plane waves have intersected and new 
 waves have appeared. The  {\mmmmtext three} wave interactions cause 
new conic waves (black surface). Only one such wave is shown in the figure. The {\vtext interaction of four waves} causes
a point source in spacetime that sends a spherical wave in all future light-like 
directions. This spherical wave is essential in our considerations. For an animation
on these interactions, see the supplementary video \cite{Video}.} 
\medskip

\subsubsection{Main result for the inverse problem for the  non-linear wave equation}
Let $(M^{(j)},\nohat g^{(j)})$, $j=1,2$ be two globally hyperbolic  $(1+3)$  
dimension  Lorentzian 
{\mmtext manifolds that are isometric to manifolds $M^{(j)}=\R\times N^{(j)}$ having a Lorenzian metric of the form (\ref{eq: product structure})}.  {\mattitext Let $\hat \mu_j=\hat \mu_j([-1,1])
\subset M^{(j)}$ be
time-like {\rtext paths}},  $V_j\subset M^{(j)}$ be an open, relatively compact, {\muutos connected}   
{\itext neighborhood
of $\hat \mu_j([-1,1])$, and $p^+_j=\hat \mu_j(s_+),$  $p^-_j=\hat \mu_j(s_-),$ where $-1<s_-<s_+<1$.} 
{\mmtext 
We assume that 
 $V_j\subset (-\infty,{\reftext T_0})\times  
N^{(j)}$.}

{\reftext Below, we sometimes drop the index $j$ and denote $M^{(j)}$ by $M$,
 $V_j$  by $V$, $g^{(j)}$ by $g$, $a_j$ by $a$, etc.} 

{Consider the
non-linear wave equation \refHOX{Comment 6: index $j$ is dropped from the notations.}
\beq \label{eq: wave-eq general}
& &\hspace{-.5cm}  \square_{\nohat {\reftext g}}u(x)+{\reftext a}(x)\,u(x)^2={\reftext f}(x)
\quad\hbox{on }{\rbtext (-\infty,{\reftext T_0})\times N},\hspace{-.5cm}\\ \nonumber
 & &\quad 
u(x)=0,\quad \hbox{for }x\in((-\infty,{\reftext T_0})\times N)\setminus J^+_{\nohat {\reftext g}}(\supp({\reftext f})),
\eeq
where $\supp({\reftext f})\subset {\reftext V},$ $x=(x^0,x^1,x^2,x^3)=(t,x^\prime)\in (-\infty,{\reftext T_0})\times N$, and
\ba
\square_{\nohat g} u=
\sum_{p,q=0}^3(-\det (\nohat g))^{-1/2}\frac  \p{\p x^p}
\left ((-\det(\nohat g))^{1/2}
\nohat g^{pq}\frac \p{\p x^q}u(x)\right),
\ea
$\det(\nohat g)=\det((\nohat g_{pq}(x))_{p,q=0}^3)$.  {\mattitext Here,}
  ${\reftext f}\in C^6_0({\reftext V})$ is a controllable source, and ${\reftext a(x)}$ is
a  {\mattitext nowhere} vanishing $C^\infty$-smooth  function.} \refHOX{Comment 7: We have specified in what sense $f$  is small.}

{\reftext  For any given $T_0>0$, the local existence results for the non-linear hyperbolic equations imply that there is $\e=\e(T_0,N,g,a,V)>0$ {\muutos such that if }
$f\in C^6_0(V)$
satisfies $\|f\|_{C^6(\overline {\reftext V})}<\e$,} 
then 
 the equation (\ref{eq: wave-eq general}) has a unique solution $u\in C^2({\reftext M}_0)$
 see e.g. \cite{ChBook,Ringstrom}. Note that we do not consider here optimal results in terms of smoothness.

\begin{definition}
{\gtext \refHOX{Comment 9: Index $j$ is dropped from notations.}
Let {\muutos $\W=\{f\in C^6_0({\reftext V});\  \|f\|_{C^6(\overline {\reftext V})}<\e\}$,
where $\e>0$ is so small that}  
the equation (\ref{eq: wave-eq general}) has a unique solution $u\in C^2((-\infty,{\reftext T_0})\times N)$ for all $f\in {\reftext \mathcal W}$.}
The  source-to-solution map $L_{{\reftext V}}:{\reftext \mathcal W}\to C({\reftext V})$, is the non-linear
operator mapping the source ${\reftext f}$ to {\muutos the} restriction of the corresponding solution
of the wave equation $u$ {\vtext to} the observation domain ${\reftext V}$, that is,
 \beq\label{measurement operator}
L_{{\reftext V}}: {\reftext f}\mapsto u|_{{\reftext V}},\quad f\in {\reftext \mathcal W}\subset C^6_0({\reftext V}),
\eeq
where $u$ satisfies the wave equation (\ref{eq: wave-eq general}) on  
\refHOX{Comment 10: Smallness of $\mathcal W$ is now discussed above using suitable norm.}
$({\rbtext (-\infty,{\reftext T_0})\times N},g)$.
\end{definition}

{\mattitext 
\medskip 

\noindent
{\bf Inverse problem with active measurements:} {\it We assume that we are given the set $V$
as a differentiable manifold and the source-to-solution map $L_{V}: f\mapsto u|_{V}$.
The inverse problem with active measurements is whether 
 these data determine the set $I({\reftext p^-},{\reftext p^+})\subset M$ as a differentiable manifold and
the conformal class of the metric $g|_{I({\reftext p^-},{\reftext p^+})}$.}
\medskip

\refHOX{Comment 12: text on propagation of waves is discussed in detail.}
 The set $I({\reftext p^-},{\reftext p^+})$  
is the maximal set that one can reach by a causal curve that starts from $p^-$ and {\mmmmtext ends to} $p^+$.}


\refHOX{Comment 13: The details of the notation $I(p^-_j,p^+_j)$ are moved just  before Theorem 1.5.}
Below, we return to consider two manifolds  $(M^{(j)},\nohat g^{(j)})$, $j=1,2$.
%
 {\reftext   We recall that  $p^\pm_j=\hat \mu_j(s_\pm)$ and denote $I(p^-_j,p^+_j)=I^+(p^-_j)\cap I^-(p^+_j)$ on $(M^{(j)},g^{(j)})$ where
 the causality in $I(p^-_j,p^+_j)$ is defined using the metric $g^{(j)}$ {\muutos and the path $\hat \mu_j$}.}

Our main result for the the inverse problem for the non-linear wave equation is the following: 

{\ggtext
\begin{theorem}\label{main thm3 new} 
{\mmmmtext Let $(M^{(j)},\nohat g^{(j)})$, $j=1,2$ be two
smooth, globally hyperbolic    Lorentzian manifolds of  
dimension $(1+3)$ {\reftext that are represented in the form $M^{(j)}=\R\times N^{(j)}$} with a metric of the form (\ref{eq: product structure}).\refHOX{Comment 11: The time-space product structure of manifolds is clarified.}

Let 
$\hat \mu_j:[-1,1]\to  (-\infty,{\reftext T_0})\times N^{(j)}$
be smooth time-like paths,
$p^+_j=\hat \mu_j(s_+),$ $p^-_j=\hat \mu_j(s_-)$, where $-1<s_-<s_+<1$, and
  \refHOX{Comment 14: Manifold $ {N^{(j)}}_0$  is changed to $ {N^{(j)}}$.}
 $V_j\subset M^{(j)}$ be  neighborhoods
of $\hat \mu_j([-1,1])$. 

Let   $L_{V_j}$, $j=1,2$ be  the source-to-solution maps for wave equations (\ref{eq: wave-eq general}) on
manifolds $(M^{(j)},g^{(j)})$ with {\gtext nowhere vanishing} $C^\infty$-smooth functions
$a_j:M^{(j)}\to \R\setminus \{0\}$, $j=1,2$, see  (\ref{measurement operator}).

Assume that there is a %
diffeomorphism
$\Phi:V_1\to V_2$ such that $\Phi(p^-_1)=p^-_2$ and 
$\Phi(p^+_1)=p^+_2$ and the source-to-solution maps
satisfy 
\ba
((\Phi^{-1})^*\circ L_{V_1}\circ \Phi^*) f =L_{V_2}f
\ea
for $f\in \W$, where $\W$ is  a neighborhood of the zero function in $C^6_0(V_2)$.

Then there is a diffeomorphism $\Psi:I(p^-_1,p^+_1)\to I(p^-_2,p^+_2)$
and the metric $\Psi^*\nohat g^{(2)}$ is conformal to $\nohat g^{(1)}$ in  
$I(p^-_1,p^+_1)\subset M^{(1)}$,  \refHOX{Comment 15: We added that $b=1$ in $V$.}
that is, there is $b: M^{(1)}\to \R_+$ such that $\nohat g^{(1)}(x)=b(x)(\Psi^*\nohat g^{(2)})(x)$ in
$I(p^-_1,p^+_1)$. {\reftext Moreover, $b(x)=1$ for $x\in V_1$.}}
\end{theorem}
}

{\mattitext Later, in Remark 3.1 we will show that the set $V$ and the map $L_V$ determine
the metric tensor $g|_V$ and coefficient $a|_V$ in $V$.}
\medskip


\refHOX{ Comment 16: Subsection on determination of lower order terms is removed.}

\subsubsection*{Outline of the proof of Theorem \ref{main thm3 new}} 

{\revtext

The proof is given in Sections 2,3, and 4. In Sec.\ 2 we give preparatory geometrical results
to estimate the
locations of cut and conjugate points along geodesics and introduce the concepts
of first observation time functions. 
In Sections  3  and 4 we use the non-linearity to reduce the studied inverse problem to
 an inverse source problem for a linear wave equation. In particular,
we are interested in constructing ``artificial point sources'' in the spacetime. 
{\gtext This is done by {\vtext using the interaction of four distorted} plane waves, {\muutos see Fig.\ 3}. Using the
waves produced by such point sources we can determine the earliest
light observation set and use Theorem \ref{main thm paper part 2},
see Sec.\ 3-4.}

\extension{
The construction of artificial point sources is done as follows: 
We consider 4 sources $f_j$, $j=1,2,3,4$, 
supported in $V$, and their linear combination $f_{\vec\e}(x)=\sum_{j=1}^4 \e_j f_j(x)$ where
$\vec\e=(\e_1,\e_2,\e_3,\e_4)\in \R^4$ is small.
Let  $Q_{\nohat g}=\square_{\nohat g}^{-1}$ denote the inverse of the linear wave operator on $(M,\nohat g)$.
We construct sources  $f_j$ that are  non-smooth on 2-dimensional surfaces  $Y_j\subset V$ and are
such that the corresponding solutions of the linearized wave equation, $u_j=Q_{\nohat g}f_j$, are distorted plane 
waves that propagate near geodesics $\gamma_{x_j,\xi_j}$, see Fig.\ 1(Right). These linearized
waves $u_j$ are non-smooth on \HOX{Skava made a suggestion: Here the order of topics could be changed so that we first discuss the
  {\mmmmtext three} {\muutos wave interaction} and then  {\mmmmtext four} {\muutos wave interaction}, or in the opposite way.}
3-dimensional surfaces $K_j$. Due to the non-linear term in the wave  equation (\ref{eq: wave-eq general}), these waves interact. When $u_{\vec\e}(x)$ is the solution of the non-linear wave equation (\ref{eq: wave-eq general})
with the source $f_{\vec\e}(x),$
the fourth order interaction of the waves $u_j$ produce a new wave that 
has the form\HOX{Formula is changed.}
\beq\label{H term}
& &\M^{(4)}:=\p_{\e_1}\p_{\e_2}\p_{\e_3}\p_{\e_4}u_{\vec\e}(x)|_{\vec\e=0}=Q_{\nohat g}\mathcal S,\hbox{ where}\\ \nonumber
& &\mathcal S=-a\, u_4\,\cdotp Q_{\nohat g}(a\, u_3\,\cdotp Q_{\nohat g}(a\,u_2\,\cdotp u_1))\\
& &\quad\quad\quad-aQ_{\nohat g}(a\, u_4\,\cdotp u_3)\,\cdotp Q_{\nohat g}(a\,u_2\,\cdotp u_1)+\hbox{permuted terms}.\nonumber
\eeq
{For example in the Minkowski space, we can consider the case when} each linearized wave $u_j$ is non-smooth on a  {\mattitext light-like} hyperplane $K_j$. 
If the geodesics  $\gamma_{x_j,\xi_j}$ 
intersect at a single point $q\in (-\infty,{\reftext T_0})\times N$, so that $\cap_{j=1}^4 K_j=\{q\}$, see Fig.\ 3 and 4(Left),
that the 4th order interaction
creates a new source $\mathcal S$ that is similar to a point source supported at $q$.
 {\mattitext  For the visualization of this interaction and the  waves generated by it, see the supplementary video \cite{Video}.} 
Similarly, the 3rd order interaction of waves produces an artificial source $\mathcal S^{(3)}$
that is singular on a 1-dimensional variety of the spacetime. The source $\mathcal S^{(3)}$  is similar to a continuous
point source that is moving on its world line (with the speed higher that the speed of light) 
and produces a wave $\M^{(3)}:=Q\mathcal S^{(3)}$
that is  singular on a  conic set   $\mathcal Y $ of the spacetime that emanates  from the moving
point source. 
In the set $V$, the source-to-solution map determines {the wave produced
by the  {{\vtext fourth order}} interaction,}
$$
\M^{(4)}|_{V}=\p_{\e_1}\p_{\e_2}\p_{\e_3}\p_{\e_4}\,L_{V}(f_{\vec\e})|_{\vec\e=0}.
$$
However, if one would be very unlucky, it could be that the different terms in (\ref{H term})
cancel each other so that no singularities are produced by the 4th order interaction observed in $V$.
In Sec.\ 3 we show that if the intersection of these geodesics $\gamma_{x_j,\xi_j}$ happens before the
conjugate points of the geodesics  (i.e.\ the caustics of the waves), then at $x\in V\setminus \mathcal Y $,  
the amplitude of the singularities of the  {\mmmmtext four} {\muutos wave interaction} wave $\M^{(4)}(x)$ is a real-analytic function of the intersection angles
of the hypersurfaces $K_j$ at $q$ and the principal symbols  of the sources $f_j$.
Also, by analyzing the forward scattering in the  non-linear wave equation (\ref{eq: wave-eq general})
we show that this real-analytic function does not vanish in any open set
 {\mattitext of parameters}, that is, in a generic
case the 4th order interaction of intersecting plane waves sends singularities to almost all directions.
Roughly speaking, this means that by using the non-linear interaction of waves $u_j$, one
can produce an ``artificial point source'' at any point $q\in I(p^-,p^+) =I^+(p^-)\cap I^-(p^+)$.
{\mattitext When $U\subset V$ is the smaller observation domain,
we can use the singularities of the wave produced by this artificial point source at $q$
 to construct the light observation set $\P_U(q)=\L^+(q)\cap U$ and the earliest light observation set $\E_U(q)
 \subset \P_U(q).$
When $q$ varies in the set $I(p^-,p^+)$,
the {\rtext family} of these  sets  $\E_U(q)$} 
determine uniquely the topological, differentiable and the conformal structures
on  $I(p^-,p^+)$.}


As noted above, the interaction of the distorted plane waves are difficult to analyze
if the waves have caustics.
 By using global Lorentzian geometry (in Sec.\ 2) we give in Sec.\ 4 give conditions that ensure that
  no caustics affect 
the earliest observations {\gtext obtained} from the interaction of four colliding, distorted plane waves when these waves are produced
by 
appropriate sources and the collision of the waves is observed before a certain time. 
 We use this in Sec.\ 4 to give a {\gtext step-by-step} construction of the earliest light observations
 corresponding to points $q$ in the diamond set
$I(p^-,p^+)$. After this the topological, differentiable, and conformal structures in
$I(p^-,p^+)$ can be reconstructed using Theorem
\ref{main thm paper part 2}.

 In this paper we present the complete proofs of the results,
but mention for the convenience of the reader that extended versions of 
some technical computations discussed briefly in this paper and the
follow up paper \cite{E-preprint}
can be 
found in the preprint \cite{preprint}.

\medskip

\noindent
{\bf 1.3.\ Remarks and applications}
 $   $\smallskip


\noindent
{\bf Remark 1.1.} The technique developed in the proof of  Theorem \ref{main thm3 new} 
can be applied to many non-linear equations, including many semi-linear equations where
the metric $g$  depends on the solution. For example, in the follow up paper 
\cite{E-preprint},  we will show how the inverse problem 
for the coupled Einstein equations and scalar field equations can be solved the using
methods developed in this paper.

The techniques considered in this paper can be used also to study
inverse problems for non-linear hyperbolic systems encountered in applications and in problems
encountered in mathematical physics.
For instance, in medical imaging, in the
  the recently developed ultrasound elastography imaging technique
 the elastic material parameters are reconstructed by sending 
 (s-polarized) elastic waves that are imaged using (p-polarized) elastic waves,
 see  
e.g.\ \cite{Hoskins,McLaughlin1}.
This imaging method uses  interaction of waves and is based on the non-linearity of the system.

{\mmtext Passive imaging problems similar to Thm.\ \ref{main thm paper part 2} are encountered in seismic imaging based on microseismic events, where one records waves coming from
natural point sources that go off at unknown times \cite{quake2}.}
\smallskip

\noindent
{\bf Remark 1.2.}
 Theorem \ref{main thm3 new} can in some cases be improved
 so that also  the conformal factor of the metric tensor can be reconstructed.
Indeed, Theorem \ref{main thm3 new} and 
Corollary \ref{coro of main thm original Einstein pre} imply 
that if $W\subset I(p^-,p^+)$ is Vacuum, i.e., Ricci-flat,
and all points $x\in W$ 
can be connected by a  {\muutos path} $\alpha\subset W^{int}$ to points of $V$, then
under the assumptions of Theorem \ref{main thm3 new},  the
whole metric tensor $g$ in $W$ can be reconstructed.

}

\canberemoved{
  \medskip

\noindent
 {\bf Remark 1.3.}
Theorem \ref{main thm3 new} deals with an inverse problem for ``near field''  measurements. 
We remark that for inverse problems for  linear equations the measurements of the Dirichlet-to-Neumann map or the ``near field''
measurements  are
 equivalent to scattering or ``far field'' information. Analogous considerations
 for non-linear equations have not yet been done but are plausible. On
 related inverse scattering problems, see \cite{MsBV}.
 }

%
%
%

\section {Earliest observation time  functions}
\label{sec: Reconstruction of topology}

\subsection{Preliminary constructions}  \label{Observation domain}


%

{\muutos Let $(M,g)$ be  
a globally hyperbolic Lorentzian manifold of type $(1, n-1)$. 
As noted above, by \cite{Bernal2}, {\mmmmtext there is an isometry $\Phi$  from $M$  to a manifold $M=\R\times N$ having the metric of the form (\ref{eq: product structure}). }
This isometry defines a smooth time function
${\bf t}:\hattuM \to \R$ by setting
 ${\bf t}(x)=t$ if $\Phi(x)\in \{t\}\times N$.
We will use notation
\beq\label{Slava formula A}
{\muutos M_0}=(-\infty,{\reftext T_0})\times N. 
\eeq

In addition to the Lorentzian metric $g$, we introduce on $M$ a smooth Riemannian metric $g^+$, 
that obtained by changing, in local coordinates, the sign of the negative eigenvalue
of the Lorentzian metric $g$. We use the Sasaki distance induced by $g^+$  on $TM$. }

For  $W \subset M,$ let $L^+W = \bigcup_{p \in W} L^+_pM \subset TM$
{\gtext be the bundle of future pointing light-like vectors
and $L^{*,+}W = \bigcup_{p \in W} L^{*,+}_pM \subset TM$ be the 
 bundle of future pointing light-like co-vectors.} 
 {\itext Here, the covector $\eta\in T^*_xM$  is defined to be future pointing if the corresponding vector $\eta^\sharp=g^{jk}\eta_k\frac \p{\p x^k}\in T_xM$
 is  future pointing.}
 The projection from the tangent bundle $TM$ to the base point of a vector is denoted by
$\pi:TM\to M$.


Let us consider points $x, y\in M$. For $x< y$,
we define  the time separation function $\tau(x,y)\in [0,\infty)$
to be the supremum of the lengths 
$
L(\alpha)=\int_0^1 \sqrt{-g(\dot\alpha(s),\dot\alpha(s))}\,ds
$
of the piecewise smooth
causal paths $\alpha:[0,1]\to M$ from $x$ to $y$. If  the condition $x< y$ does not
hold, we define $\tau(x,y)=0$. 
We note that $\tau(x,y)$ satisfies the reverse triangle inequality
 \beq\label{eq: reverse}
 \tau(x,y)+\tau(y,z)\leq \tau(x,z),\quad\hbox{for }x\leq y\leq z.
 \eeq
As $M$ is globally hyperbolic,
the time separation function $(x,y)\mapsto \tau(x,y)$  is continuous in  $M\times M$
by \cite[{\ftext Lemma}\ 14.21]{ONeill}. By
 \cite[{\ftext Lemma}\ 14.22]{ONeill}, the sets $J^\pm(q)$ are closed.
For $q<p$ there is a causal geodesic $\gamma([0,1])$ with
 $\gamma(0)=q$ and $\gamma(1)=p$ such that $L(\gamma)=\tau(q,p)$, see
\cite[{\ftext Lemma} 14.19]{ONeill}. This geodesic, called
a {longest path} from $q$ to $p$, may not be unique.

\begin{center}
\psfrag{1}{\hspace{-2mm}$x_1$}
\psfrag{2}{\hspace{-2mm}$x_2$}
\psfrag{3}{$q$}
\psfrag{4}{$V$}
\includegraphics[height=5.5cm]{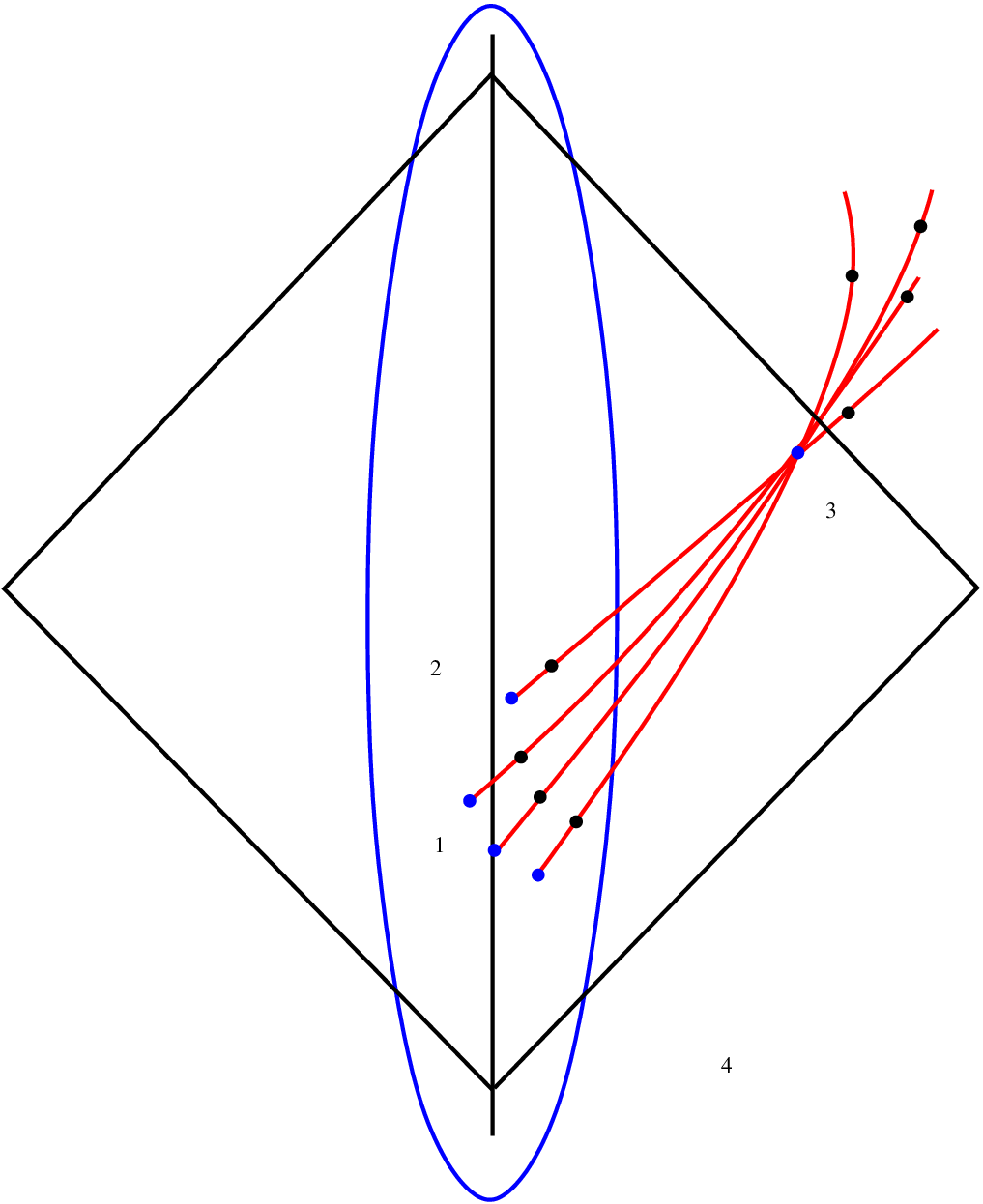}\quad\quad
\psfrag{1}{$q$}
\psfrag{2}{\hspace{-1mm}$p_j$}
\psfrag{6}{$V$}
\includegraphics[height=5.5cm]{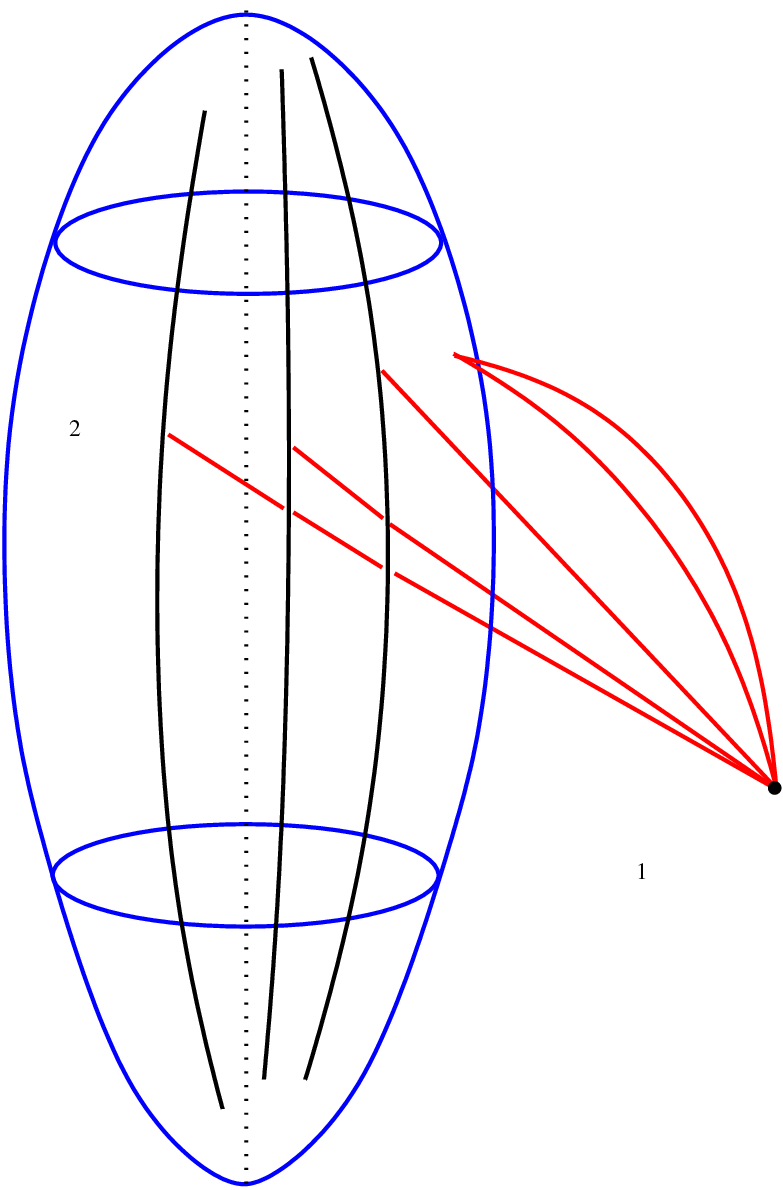}
\end{center}

{\it FIGURE 4. {\bf Left:} {We do observations in the set  $V$, marked by the blue boundary.
This set contains the set $U$, defined in (\ref{eq: Def Wg with hat}),
that is a union of time-like {\rtext paths}. In the figure,}
the four light-like geodesics $\gamma_{x_j,\xi_j}([0,\infty))$, $j=1,2,3,4$ starting at
the blue points $x_j\in U$
intersect at $q$ before the first cut points of $\gamma_{x_j,\xi_j}([0,\infty))$,
denoted by black points. The points $\gamma_{x_j,\xi_j}(t_0)$ are also shown as black points.
We use interaction of waves to produce an artificial point source at $q$.

{\bf Right:}  \refHOX{Comment 3: Figure caption is changed} The black curves are the {\reftext time-like paths $\mu_{a}\subset U,$   indexed by $a\in \A$} and
 the red curves are light-like geodesics from $q$,
{\reftext see 
Subsection \ref{subsection: Observation domain U} {\muutos for} the notation $\A$  and
Definition \ref{def: f functions} on the functions $f^+_a$}. Some
 light rays from $q$ intersect $\mu_{a}$ at {\reftext the point $p_a=\mu_{a} ( f^+_{a}(q))$,
 that is the first point of $\mu_{a}$    that is in the causal future of $q$.}
For any $q_0\in W$ we can find $a_j\in \A$, 
$j=1,2,\dots,n$ and  a neighborhood of $q_0$ where
the observation time functions $q\mapsto f^+_{a_j}(q)$
 define a smooth
coordinate system.}
\smallskip

When $(x,\xi)$ is a  non-zero vector, we define  $\T(x,\xi)\in (0,\infty]$ to
be the maximal value for which $\gamma_{x,\xi}:[0,\T(x,\xi))\to M$ is defined.

In addition to points $p^\pm=\hat \mu(s_\pm)$, we use the points $p_{\pm 2}=\hat \mu(s_{\pm 2})$
where $-1<s_{-2}<s_{-}$ and $s_+<s_{+2}<1$.

For $(x,\xi)\in L^+M$,
we define the cut locus function
\beq\label{eq: max time}
& &\rho(x,\xi)=\sup\{s\in [0,\T(x,\xi));\ \tau(x,\gamma_{x,\xi}(s))=0\},
\eeq
c.f.\ \cite[Def.\ 9.32]{Beem}.
The points $x_1=\gamma_{x,\xi}(t_1)$ and $x_2=\gamma_{x,\xi}(t_2)$,
$t_1,t_2\in [0,t_0]$, $t_1<t_2$, are cut points on $\gamma_{x,\xi}([0,t_0])$
if $t_2-t_1=\rho(x_1,\xi_1)$ where $\xi_1=\dot\gamma_{x,\xi}(t_1)$.
In particular, the point $p(x,\xi)=\gamma_{x,\xi}(s)|_{s=\rho(x,\xi)}$,
if it exists, is called the first cut point
on the geodesic  $\gamma_{x,\xi}([0,\T(x,\xi)))$.
Using \cite[Thm. 9.33]{Beem}, we see that the function $\rho(x,\xi)$ is lower semi-continuous on a globally hyperbolic 
Lorentzian manifold $(M,g)$.

%
Recall that  $\gamma_{x, \xi}(t)$ is a conjugate point on
$\gamma_{x, \xi}([0,\T(x,\xi)))$ if the differential of the map $\exp_x$
is not invertible at $t\xi$.
By \cite[Th. 9.15]{Beem}, on a globally hyperbolic manifold, $p(x,\xi)$ is either the first conjugate point along $\gamma_{x, \xi}$, or the first  point on  $\gamma_{x, \xi}$
where there is another light-like geodesic
 $\gamma_{x, \eta}$  from $x$ to $p(x,\xi)$, $\eta\not= c\xi$.

\modified{Let us return to the longest paths.}
If $q<p$ but $\tau(q, p)=0$, then there is a light-like geodesic 
$\gamma_{q,\xi}([0,t])$ 
 from $q$ to $p$ so that
 there are no cut points on $\gamma_{q,\xi}([0,t))$, see 
\cite[Thm.\ 10.51 and  Prop.\ 14.19]{ONeill}. 
Note that  if $\gamma_{q,\xi}
([0,t])$ is a light-like geodesic from $q$ to $p=\gamma_{q,\xi}(t)$
such that there are cut-points on the geodesic $\gamma_{q,\xi}
([0,t))$, \modified{(\ref{eq: reverse}) and (\ref{eq: max time}) yield}
$\tau(q, p) >0$.

 \modified{We say that 
a  {\muutos path} $\alpha([t_1,t_2])$ is a pre-geodesic if  $\alpha(t)$ is a $C^1$-smooth
 {\muutos path} such that $\dot\a(t)\not =0$ 
on $t\in [t_1,t_2]$, and $\alpha([t_1,t_2])$ can be re-parametrized so that it becomes a geodesic. A conformal diffeomorphism preserves
the light-like pre-geodesics by \cite[Th. 9.17]{Beem}.

Moreover, it follows from \cite[Prop.\ 10.46]{ONeill} that if $q$ can be connected to
$p$ with a causal path which is not a light-like pre-geodesic then $\tau(q,p)>0$. 
Let us apply this fact to a path  from $q$ to $p$ which is the union of  the future 
pointing light-like pre-geodesics
$\gamma_{q,\eta}([0,t_0])\subset M$  and $\gamma_{x_1,\theta}([0,t_1])\subset M$,
where $x_1=\gamma_{q,\eta}(t_0)$, $p=\gamma_{x_1,\theta}(t_1)$ and $t_0,t_1>0$. Let
$\xi=\dot \gamma_{q,\eta}(t_0)$. 
Then, if there is no $c>0$ such that $\xi= c\theta$, or equivalently, the union of
these geodesic is not a light-like pre-geodesics, we have  $\tau(q,p)>0$. In particular, this implies that there exists a time-like geodesic 
from $q$ to $p$. 
In the following we call this kind of argument for a union of light-like geodesics a short-cut
argument.

 \subsubsection{Smaller observation domain $U$.}
 \label{subsection: Observation domain U}
{\mattitext
 Next we define a domain $U\subset V$ that is a union of time-like paths.}

 We  assume that we {\reftext are} \refHOX{Comment 18: "are" is added.}
given 
a family of 
future pointing,  $C^\infty$-smooth, time-like paths  
$\mu_{a}:[-1,1]\to V$,
indexed by $a\in \overline \A$, where ${\rtext \A}$ is a connected metric space and
the completion of $\A$, denoted by $\overline \A$, is compact. 
We assume that there is $\hat a\in \A$  such that $\hat \mu=\mu_{\hat a}$.
Also, we assume 
that $(a,s)\mapsto \mu_{a}(s)$ defines a continuous map $\overline {\rtext \A}\times [-1,1]\to M$
and an open
map $\A\times [-1,1]\to M$.  
 Then, we define the smaller observation domain $U\subset V$ to be the  set 
 \beq\label{eq: Def Wg with hat}
U =\bigcup_{a\in \A}\mu_{a}([-1,1]).
\eeq 
Note that as 
 $(a,s)\mapsto \mu_{a}(s)$ in   is a continuous and open map
 ${\rtext \A}\times [-1,1]\to M$, the set $U$ is open.

\canberemoved{
{\gtext Let us give an example how the family $\mu_a$ of time-like paths
and the index set $\A$ can be constructed when we are
given the set $V$, as a differentiable manifold, and the conformal type of the metric $g|_V$.
Then, we could choose a smooth metric $\tilde g$ in $V$  that is conformal to 
$g|_V$ and use an
 open neighborhood $\U_{\hat z,\hat \eta}$  
of $(\hat z,\hat \eta)$ in  $TM$ as the index set $\A$.}
Here,
 $\hat z=\hat \mu(-1)$ and  $\hat \eta=\p_s \hat \mu(-1)$
 and   for  $(z,\eta)\in  \U_{\hat  z,\hat  \eta}$
the path $\mu_{z,\eta}$ is defined to be  a geodesic (i.e., a freely falling observer)
$\mu_{z,\eta}:[-1,1]\to M$ with respect to the metric $\tilde g$,
such that   $\mu_{z,\eta}(-1)=z$, and $\p_s\mu_{z,\eta}(s)|_{s=-1}=\eta$.
When $\U_{\hat z,\hat \eta}$   is a sufficiently small neighborhood of 
$(\hat z,\hat \eta)$ in  $TM$, the set $U$ given by (\ref{eq: Def Wg with hat})
  is a relatively
 open subset of $V$.
Note that then  $\hat \mu=\mu_{\hat z,\hat \eta}$.}}
%

{\mattitext
Let $s_{-2}\in (-1,s_-)$, and $s_{+2}\in (s_+,1)$.
 By replacing $\A$ in the formula (\ref{eq: Def Wg with hat}) by}  a smaller neighborhood of $\hat a$ we
may assume for all $a\in  \overline\A$ we have  
\beq \label{kaava D}
\mu_a(s_{-2}) \hspace{-1mm}\in\hspace{-1mm} I^+(\mu_{\hat a}(-1)) \cap I^-(p^-), \  \mu_a(s_{+2})\hspace{-1mm} \in \hspace{-1mm}I^-(\mu_{\hat a}(1)) \cap I^+(p^+).\hspace{-2cm}
\eeq
When $V$ is given as differentiable manifold and the 
conformal class of $g$ is given, we may define a family
of smooth time-like paths $\mu_a:[-1,1]\to V$, $a\in \overline \A$ having
the above properties and define the neighborhood
$U$ given in  (\ref{kaava D}).

\canberemoved{
{\gtext  By \cite[Prop.\ 5.7]{ONeill}, any point $x\in M$ has a neighborhood $W_x$ that is convex,
i.e., for all $x_1,x_2\in W_x$ any causal geodesic in $M$ connecting $x_1$ and $x_2$ is a subset of $W_x$. 
The compact set $\hat \mu ([s_-,s_+])$ can be covered with a finite union of such neighborhoods.
Thus when the pair $(U,g|_U)$ is given, we see using the fact that $\hat \mu$ is time-like and the continuity of 
the time-separation function $\tau(x,y)$ that one can determine 
all points $(x.y)\in W_0\times W_0$ satisfying $x\ll y$ when 
$W_0$ is some  sufficiently small open neighborhood of $\hat \mu$. Thus by making
the set $\A$ a smaller neighborhood of $\hat a$ we may assume that 
the set $U$ defined in (\ref{eq: Def Wg with hat}) is contained in
the set $W_0$. Hence
without loss of generality we may assume that the 
set
\beq\label{tau known in U}
R_U=\{(x,y)\in U\times U;\  x\ll y\}
\eeq
is known. Note that here the relation $\ll$     is defined on the whole manifold $M$.
Thus we in particular know the set $U^\prime=U\cap (I^+(p^-)\cap I^-(p^+))$.
}}
 \subsubsection{Observation time functions}
 \label{subsection: Observation time}
%
%
%
%
%

 {\mattitext 


 {\gtext Instead of the light observation sets we can consider the earliest observation
 time functions that we proceed to define.} 
 
 \begin{definition}\label{def: f functions}
 Let $a\in  \overline \A$. 
For $x\in M$
 we define $f_a^+(x),f_a^-(x)\in [-1,1]$
by setting 
\ba
f_a^+(x) &=&\inf (\{s\in (-1,1);\ \tau(x,\mu_a(s))>0\}\cup\{1\}),\\
f_a^-(x) &=&\sup (\{s\in (-1,1);\ \tau(\mu_a(s),x)>0\}\cup\{-1\}).
\ea
We call
 $f_a^+(x)$ the earliest observation time from the point $x$ on the path $\mu_a$.
 The functions
 $f_a^+:M\to \R$, $a\in  \overline\A$ are called the earliest observation time  function on the path $\mu_a$.
%

\end{definition}

}

We will show that the map 
\ba
\F:J^-(p^+)\setminus
I^-(p^-)\to C(\overline \A),\quad \F(q):a\mapsto f^+_a(q),
\ea
 that maps a point $q$ to the  earliest observation times corresponding to
 the point, is a continuous function.}
 {\muutos We will prove the following proposition  in Section \ref{subsec: Obs. time}.}

\begin{proposition}\label{propo paper part 2}
Let $(M,{g})$ be an open, $C^\infty$-smooth, globally hyperbolic 
Lorentzian manifold of dimension $n\ge 3$. Let
$\A$ be a metric space which completion $\overline \A$ is compact, 
and $\mu_{a}:[-1,1]\to M$, $a\in \overline \A$ be
 $C^\infty$-smooth, time-like paths. Let $p^-=  \mu_{\hat a}(s_-)$ 
 and $p^+=  \mu_{\hat a}(s_+)$ with
 $\hat a\in \A$ and $-1<s_-<s_+<1$. Also, assume that
$(a,s)\mapsto \mu_{a}(s)$ defines a continuous map $\overline A\times [-1,1]\to M$
and an open
map $\A\times [-1,1]\to M$.  
 
   Let $W\subset M$
be open set such that $\overline W \subset J^-(p^+)\setminus I^-(p^-)$
is compact. Then
$\F:W\to \F(W)$ is a 
 homeomorphism. {\rtext Here $F(W)$ has the metric induced by $C(\overline \A)$.}
\end{proposition}

%

%
%
%


%
%
%
%
{\mmmmtext 
In several 
 geometric inverse problems \cite{AKKLT,KKL,KLU}, in order to reconstruct an unknown manifold $W$ from a given data, one needs
to construct a copy of the manifold. The importance 
of  Proposition \ref{propo paper part 2} lies in the fact that it can be used 
to construct a homeomorphic image of  the original Lorentzian manifold $W$ embedded in the known space
$C(\overline \A)$.  After the 
homeomorphic image $\F(W)$ of the manifold is constructed, we can
construct other structures on it, e.g., the differentiable coordinates 
and a metric tensor conformal to the original metric.  
}

We need the following simple properties of the functions $ f_a^\pm(x)$.

\begin{lemma} \label{B: lemma} 
Let $a\in \overline \A$ and   $q \in  J^-(p^+)\setminus I^-(p^-)$.
Then

\smallskip

\noindent
(i) {\gtext We have that}
 $s_{-2}\leq  f^+_a(q)\leq s_{+2}$.
\smallskip 

\noindent
(ii)  We have $\mu_{a}(f^+_{a}(q))\in J^+(q)$ and
$\tau(q, \mu_{a}(f^+_{a}(q)))=0$. 
 Moreover, the function $s\mapsto \tau(q,\mu_a(s))$ is  continuous, non-decreasing on the interval $s\in [-1,1]$,
and it is  strictly increasing on $ [f^+_a(q),1]$.
\smallskip

\noindent (iii) 
Assume that $p\in U$. 
Then 
$p=\mu_{a}(f^+_{a}(q))$ with some $a\in \A$ if and only if  $p\in \P_U(q)$ and $\tau(q,p)=0$.
Furthermore,
these are equivalent to \modified{the fact} that 
there are $\xi\in L^+_qM$ and 
$t\in [0,\rho(q,\xi)]$ such that 
 $p=\gamma_{q,\xi} (t)$. 


 {\mattitext

\smallskip 

\noindent (iv)
The function ${\bf F}:(a,q) \mapsto f^+_a(q)$ is   continuous on $\overline \A\times ( J^-(p^+)\setminus I^-(p^-))$.}
%
%
%
\smallskip

{\muutos For   $q \in  J^+(p^-)\setminus I^+(p^+)$ the claims analogous with
(i)-(iv), with reversed causality, are valid for $q\mapsto f^-_a(q)$.}

\end{lemma}


\noindent  
{\bf Proof.} (i) This property follows from (\ref{kaava D}).

 (ii) Since $J^+(q)$ is closed, $\mu_{a}(f^+_{a}(q))\in J^+(q)$.
The continuity of $s\mapsto \tau(q, \mu_a(s))$ follows from the continuity of $\tau(x, y)$ on $M\times M$.

If  $\tau(q,\mu_a(f^+_a(q)))$ would be strictly positive, we would have $\mu_a(f^+_a(q))\in I^+(q)$
and there would exist $s<f^+_a(q)$ such that  $\mu_a(s)\in I^+(q)$. As this is not possible,
we have $\tau(q,\mu_a(f^+_a(q)))=0$.

Consider $s<s^{\prime}$.
Since $\mu_a$ is a time like-path, $\tau(\mu_a(s),\mu_a(s^{\prime}))>0$. 
Thus, when $s^\prime>s\geq f^+_a(q)$, 
the inequality (\ref{eq: reverse})
yields  $ \tau(q,\mu_a(s))<\tau(q,\mu_a(s^{\prime}))$. 
\modified{For $s< f^+_a(q)$} we have $\mu_a(s)\not \in J^+(q)$ and  
$\tau(q, \mu_a(s))=0$.
%

(iii) It is sufficient to prove the claim when $p\not=q$.
First, assume that  $p=\mu_{a}(f^+_{a}(q))$.  Then
$p \in J^+(q)$ and by (ii), we have $\tau(q, p)=0$.  \modified{The existence of the light-like geodesic follows from the above.}

Second, assume that $p \in J^+(q)$ and  $\tau(q, p)=0$. 
This implies by \cite[Prop. 14.19]{ONeill}  that there exists  a light-like geodesic $\gamma_{q,\xi}
([0,t])$
from $q$ to $p$. If $\gamma_{q,\xi}
([0,t))$ would have a cut-point, then $\tau(q, p) >0$ which is not possible. Thus, $t\in [0, \rho(q,\xi)]$.

Third, assume that $p=\gamma_{q,\xi}(t)$ with $\xi\in L^+_qM$ and 
$0\leq t\leq \rho(q,\xi)$. Then $\tau(q,p)=0$. Let 
$a\in  \A$ and $s_0\in [-1,1]$ be such that $p=\mu_a(s_0)$.
{\mattitext As $q\not \in  I(p^-)$, using (\ref{kaava D}) we see that 
$q\not \in I^-(\mu_a(s_{-2}))$
and hence $s_0\geq s_{-2}>-1$.}
 By (i),
$\tau(q,\mu_a(s))>0$ for $s>f^+_a(q)$
and thus  $s_0\leq f^+_a(q)\leq s_{-2}$. However, $q\leq p=\mu_a(s_0)$
and $\tau(q, \mu_a(s))>0$ for $s\in (s_0,1)$.
Thus $s_0\geq f^+_a(q)$. Hence, $s_0= f^+_a(q)$ and $p=\mu_{a}(f^+_{a}(q))$.

{\mattitext 
%

(iv) 
Assume that $x_j\to x$ in $J^-(p^+)\setminus I^-(p^-)$  and $a_j\to a$ as $j\to\infty$.
Let
$s_j=f^+_{a_j}(x_j)$ and $s=f^+_{a}(x)$. 

 Since $\tau$ is continuous,
 for any $\e>0$ we have 
$\lim_{j\to\infty} \tau(x_j,\mu_{a_j}(s+\e))=\tau(x,\mu_a(s+\e))>0$. Then for $j$ large enough 
 $\mu_{a_j}(s+\e)\in I^+(x_j)$, so that $s_j\leq s+\e.$  As $\e>0$ is arbitrary, $\limsup_{j\to \infty}s_j\leq s$
Thus ${\bf F}$ is upper-semicontinuous.

%
%
%

Next, suppose $\liminf_{j\to \infty}s_j=\tilde s< s$ and denote $\e=\tau(\mu_a(\tilde s),\mu_a(s))>0$. Then
by choosing a subsequence, we may assume that $\lim_{j\to \infty}s_j=\tilde s< s$. By continuity of  $\tau$ and (\ref{eq: reverse}),
\ba
\lim_{j\to\infty}\tau(x_j,\mu_{a_j}(s))
&\geq& 
\liminf_{j\to\infty}
\tau(x_j,\mu_{a_j}(s_j))+
\tau(\mu_{a_j}(s_j),\mu_{a_j}(s))\\
&\geq& 
0+
\tau(\mu_{a}(\tilde s),\mu_{a}(s))
=\e,\ea and
we obtain $\tau(x,\mu_a(s))=\lim_{j\to\infty}\tau(x_j,\mu_{a_j}(s))\geq \e$. This is not possible as $s=f^+_a(x)$.
Hence $\liminf_{j\to \infty}s_j\geq s$ and ${\bf F}$ is  also  lower-semicontinuous. This proves (iv).
%
%
 
  The analogous results for function $f_a^-$
follow similarly by reversing the causality. }
  \hfill \Box \medskip

Next we consider  the earliest light observation sets in the observation domains $U$ and $V$.
{\gtext We will show that without loss of generality we can take the neighborhood $V$ 
of $\hat \mu$  to
be the set $U$ defined in (\ref{eq: Def Wg with hat}).}

\begin{lemma} \label{lemma E_U} 
Let   $q \in  J^-(p^+)\setminus I^-(p^-)$.

\smallskip 

\noindent (i)
The earliest light observation set of $q$  in $U$ has the form
\beq\label{eq BB}
\E_U(q)\hspace{-1mm}=
\{\mu_{a}(f^+_{a}(q));\ a\in \A\}.
\eeq 

\smallskip  

\noindent (ii)
{\muutos Assume that we are given the sets $\E_V(q)$ and $U$ and  the paths $\mu_a,$ $a\in \A$.
These data determine the function $  f_a^+(q)$  and moreover, the 
 set $\E_U(q)$ by formula (\ref{eq BB}).}

%
%
%
%
\end{lemma}

\noindent 
{\bf Proof.} (i)
Let $x\in \E_U(q)\subset U$. Then by Lemma \ref{B: lemma} (i) there is $a\in \A$ such 
that $x=\mu_a(s)$ with $s\in [s_{-2},s_{+2}]$ and $x\in \P_U(q)$.

Assume that $\tau(q,x)>0$. Then $s>f^+_a(q)$.
By Lemma \ref{B: lemma} (iii), $y=\mu_a(f^+_a(q))\in \P_U(q)$  and
the time-like path $\mu_a([f^+_a(q),s])\subset U$ connects $
y\in \P_U(q)$ to $x$.
This is not possible by the definition of $\E_U(q)$.
This shows that $\tau(q,x)=0$. As $x\in \P_U(q)$, by Lemma \ref{B: lemma} (iii) we have $x\in
\{\mu_{a}(f^+_{a}(q));\ a\in \A\}$.

On the other hand, assume that $x=\mu_{a}(f^+_{a}(q))$ with
$a\in \A$. Then by Lemma \ref{B: lemma} (iii), $x\in \P_U(q)$ and $\tau(q,x)=0$. Then, if there would exist
$y\in \P_U(q)$ that is  connected to $x$ with a future
pointing time-like path, we would have $\tau(y,x)>0$. 
Thus (\ref{eq: reverse}) implies that $\tau(q,x)\geq \tau (q,y)+\tau (y,x)>0$.
This shows that no such $y$ can exist and
$x\in \E_U(q)$. These prove (i).

(ii)  {\muutos As  $q \in  J^-(p^+)\setminus I^-(p^-)$,
Lemma \ref{B: lemma} and inequality (\ref{eq: reverse}) yield that
the function $  f_a^+(q)$ are determined by  
\beq\label{eq: finding fa}
  f_a^+(q) =\inf  \{s\in (-1,1);\ \mu_a(s)\in \E_V(q)\}.
\eeq} 
%
%
%
%
  \hfill \Box \medskip

{\mattitext
{\muutos Due to  Lemma \ref{lemma E_U},
without loss of generality we may in Theorem \ref{main thm paper part 2}
consider the case when the set $V$
is replaced by $U$. We will do so for the remaining of this paper.}

%
%
%

}

 \subsection{Observation time representation of a Lorentzian manifold}
 \label{subsec: Obs. time}

{\mattitext In this section our main goal is to prove Proposition \ref{propo paper part 2}.}
 
 
 \subsubsection{The direction set}
 \begin{definition}\label{def: earliest element}
 Let $q \in J^-(p^+)\setminus I^-(p^-)$.
 Let 
  \beq\nonumber
\be_U(q) =\{(y,\eeta)\in L^+U&;& y=\gamma_{q,\xi}(t) \in U,\
\eeta=\dot \gamma_{q,\xi}(t),\\
\label{eq: b-sets}
& &\hbox{with some } 
\xi\in L^+_qM,\ 0\leq t\leq \rho(x,\xi)\},\\
\nonumber
\be^{reg}_U(q) =\{(y,\eeta)\in L^+U&;& y=\gamma_{q,\xi}(t) \in U,\
\eeta=\dot \gamma_{q,\xi}(t),\\
\nonumber
& &\hbox{with some } 
\xi\in L^+_qM,\ 0< t< \rho(x,\xi)\}.
\eeq 
%
%
 We say that $\be_U(q)$ is the direction set of $q$
 and $\be^{reg}_U(q)$ is the regular direction set of $q$.
 
 Then, $\firstpoint_U(q)=\pi(\be_U(q))$. We denote 
$\firstpoint^{reg}_U(q)=\pi(\be^{reg}_U(q))$ where $\pi:TU\to U$ is the canonical projection,
$\pi(y,\eeta)=y$. We say that
  $\firstpoint^{reg}_U(q)$ is {\mattitext the  regular earliest light observation set} of $q$.
\end{definition}

Note that 
 $\E_U(q) =\{ {\rbtext \mu_{a}(f^+_{a}(q))};\ a\in \A\}$ and \modified{that the lower semicontinuity
 of $\rho(x,\xi)$ implies that} 
$\firstpoint^{reg}_U(q)\subset U$ and $\be^{reg}_U(q)\subset TU$
are {\mattitext smooth $(n-1)$ and $n$ dimensional submanifolds, respectively.}

We need the following auxiliary result:

\begin{lemma}\label{lemma: properties of  E(q)}
 {\bf (i)}
Let $y\in U$,  $\eeta\in L^+_yM$, $r_1>0$, and $q\in W$  be
such that $q\not \in \gamma_{y,\eeta}([-r_1,0])$ and
$\gamma_{y,\eeta}([-r_1,0])\subset U$. 
Then  $(y,\eeta)\in \be_U(q)$
if and only if  $\gamma_{y,\eeta}([-r_1,0])\subset \E_U(q)$.
\smallskip

\smallskip

\noindent  {\bf  (ii)} Let $y\in U$,
$\eeta\in L^+_yM$, and $\hat t>0$ be the largest number such that 
the geodesic $\g_{y, \eeta}((-\hat t,0])$ is defined and has no cut points.
Then for $q\in W$ 
we have 
  $q\in \g_{y, \eeta}((-\hat t,0))$
if and only if $(y,\eeta)\in \be_{U}^{reg}(q)$.
\smallskip

\end{lemma}

\noindent 
{\bf Proof.}
(i) Suppose $(y, \eeta) \in \be_U(q)$.
Then $y \in \firstpoint_U(q)$ and $\tau(q,y)=0$.
 Since $q\not \in \gamma_{y,\eeta}([-r_1,0])\subset U$, there  is 
 $t>r_1$ such that $\gamma_{y,\eeta}(-t)=q$ and for  $\xi=\dot \gamma_{y,\eeta}(-t)$
 we have $\g_{y, \eeta}([-r_1,0])= \g_{q, \xi}([t-r_1,t])\subset \P_U(q)$. 
 If there would \modified{be} $y_1\in \g_{y, \eeta}([-r_1,0])$ such
 that $y_1\not \in  \firstpoint_U(q)$, it follows from 
  (\ref{eq BB a}) that there is $z\in  \P_U(q)$
  such that $z\ll y_1$. Then  $q\leq z\ll y_1\leq y$. These imply
  that $\tau(q,y)>0$
    and $y\in \firstpoint_U(q)$
  which is not possible by
  Lemma \ref{B: lemma} (iii) and Lemma \ref{lemma E_U}. 
  This shows
that $\g_{y, \eeta}([-r_1,0])\subset \E_U(q)$.

On the other hand, 
assume that
 $\g_{y, \eeta}([-r_1,0])\subset \E_U(q)$. 
 Then   Lemma
 \ref{B: lemma}(ii) implies that $\tau(q,y)=0$.
Denote $y_1=\g_{y, \eeta}(-r_1)$.
 Since $y_1\in \E_U(q)$ and $y_1\not =q$, there  is $\xi \in L^+_q M$
and $t_1>0$ such that $\g_{q,\xi}(t_1)=y_1$.
Then the union of the
 geodesics
 $\g_{q,\xi}([0,t_1])$ and $\g_{y,\eeta}([-r_1,0])$
 form a causal path from $q$ to $y$. Using short cut arguments,
 we see that if the union of these
 geodesics do not form one light-like pre-geodesic, we 
 have  $\tau(q,y)>0$, {\ftext that is not} possible.
 Hence $\g_{y,\eeta}([-r_1,0])$ lies in the continuation 
 of $\g_{q,\xi}([0,t_1])$, that is, there is $t>0$ such 
 that   $\g_{y, \eeta}([-r_1,0])\subset 
\g_{q,\xi}([0,t])$ and
$y= \g_{q,\xi}(t)$. Then, there is $c>0$ such that
$\eeta= c\dot \g_{q,\xi}(t)$.
Moreover, if $\g_{q,\xi}([0,t))$ would contain cut points then
 \cite[Prop. 10.46]{ONeill} implies that $\tau(q, y)>0$. This  would lead
to a contradiction with 
$y\in \firstpoint_U(q)$. Hence, $\g_{q,\xi}([0,t))$ contains no cut points. 
Therefore, we have shown that $t\leq \rho(q,\xi)$,
$y= \g_{q,\xi}(t)$, and $\eeta=c\dot \g_{q,\xi}(t)$.
These imply that   $(y, \eeta) \in \be_U(q)$.

  (ii)   Let $(y,\eeta)\in L^+U$  and $\hat t>0$ be as in the claim and $q\in W$.

First, assume that $q\in 
\gamma_{y,\eeta}(-t_1)$, $t_1\in (0,\hat t)$. Then, \modified{due to the symmetry of the cut points,} $\tau(q,y)=0$
and thus for $\xi=\dot\gamma_{y,\eeta}(-t_1)$ we have
$y=\gamma_{q,\xi}(t_1)$ and $t_1<\rho(q,\xi)$.
Thus $(y,\eeta)\in \be_U^{reg}(q)$.

Second, assume that $(y,\eeta)\in \be_U^{reg}(q)$.
Again, we see that there is $t_1>0$ such that $q=\gamma_{y,\eeta}(-t_1)$
and for $\xi=\dot\gamma_{y,\eeta}(-t_1)$ we have
$y=\gamma_{q,\xi}(t_1)$ and $t_1<R_1:=\rho(q,\xi)$.
Since $\rho$ is a lower semi-continuous, we have
that \modified{when $\e\in (0,(R_1-t_1)/2)$ is small enough,
the point $x_1=\gamma_{q,\xi}(-\e)$ and  
 $\xi_1=\dot\gamma_{q,\xi}(-\e)$ satisfy
 $\rho(x_1,\xi_1)>R_1-\e>t_1+\e$ and hence $\tau(x_1,y)=0$. 
This yields that $\hat t>t_1$. Thus $q\in 
\gamma_{y,\eeta}((-\hat t,0))$.}
\hfill \Box \smallskip

Using this result we determine the direction \modified{sets $\be_U(q)$ from $ \E_U(W)$}:

\begin{lemma}\label{lemma: from P(q) to E(q)}
Assume that we are given 
{\muutos the conformal type of $(U,g|_U)$, the paths $\mu_a:[-1,1]\to U$, $a\in \A$,
 and the set $\E_U(W)$.}
Then
\smallskip

\noindent  {\bf (i)}
For any $y\in U$, we can identify
from the set $\E_U(W)$ the element  
$\E_U(q)$ for which $q=y$, if it exists. For such elements  
$L^+_yM\subset \be_U(q)$.
\smallskip

\noindent  {\bf  (ii)} Let $q\in W$ and $(y,\eeta)\in L^+U$. Then
$(y,\eeta) \in \be_U^{reg}(q)$ if and only if there exists a  light-like pre-geodesic
$\alpha([t_1,t_2]) \subset   U$
such that $y=\alpha(t),$  $\eeta=\dot \alpha(t),$  $t_1<t<t_2$,
and $\alpha([t_1,t_2]) \subset  
 \E_U(q).$
 
\smallskip
 
 \noindent  {\bf (iii)}
When  $\E_U(q)\in \E_U(W)$ is given, one can determine
the sets $\be_U(q)$, $\be_U^{reg}(q)$,  and $\E_U^{reg}(q)$. 

\end{lemma}

\noindent 
{\bf Proof.} (i) We observe that  $q=y$ if and only if 
for $y\in \E_U(q)$ there are no $\eta\in L^+_yM$ and $t_0>0$
such that $\gamma_{y,\eta}([-t_0,0]) \subset \E_U(q)$.
Claim (i) follows from this observation.

(ii) Let $q\in W$ and $\xi\in L^+_qW$ and 
$(y,\eeta)=(\gamma_{q,\xi}(1),\dot \gamma_{q,\xi}(1)).$
 Using Definition
\ref{def: earliest element} we see that 
$(y,\eeta)\in \be_U ^{reg}(q)$ if and only if $\gamma_{q,\xi}(1)\in U$ and
$\rho(q,\xi)>1$. This is
equivalent to \modified{the fact} that
there are $t_1\in (0,1)$ and $t_2>1$ such that 
$\gamma_{q,\xi}([t_1,t_2])\subset U$
and $(\gamma_{q,\xi}(t_2),
\dot\gamma_{q,\xi}(t_2)) \in \be_U(q)$.
Also, by
 Lemma \ref{lemma: properties of  E(q)} (i)
 this is equivalent to the fact that there 
 are  $t_1\in (0,1)$ and $t_2>1$ such that 
 $\gamma_{q,\xi}([t_1,t_2])\subset \firstpoint_U(q)$. 
 This proves (ii).

(iii)  Let  $\E_U(q)$ be given.
Since the conformal class of $g|_U$ is given,
we can identify all light-like  pre-geodesics in $U$. Thus by using (ii),
we can verify for any  $(y,\eeta)\in L^+U$ whether it holds that
$(y,\eeta)\in\be^{reg}_U(q)$  or not.
Thus we can determine the set  $\be^{reg}_U(q)$.
Then the set $\be_U(q)$ can be determined as the closure
of the set $\be^{reg}_U(q)$ in $TU$. Finally, the set
$\firstpoint^{reg}_U(q)=\pi(\be^{reg}_U(q))$
can be constructed using the map
$\pi:TU\to U$.\hfill \Box \medskip

%
%
\noindent 
\subsubsection{Construction of  $W$ as a  topological  manifold.} 
%
For $q\in J^-(p^+)\setminus
I^-(p^-)$ we define the {\mattitext continuous function  $F_q:\overline \A\to \R$  by $F_q(a)= f^+_a(q)$. 
Also, we denote by $\F:J^-(p^+)\setminus
I^-(p^-)\to C(\overline \A)$ the function $\F(q)=F_q$, that maps $q$ to the 
function $F_q:\A\to \R$.} 

{\muutos By Lemma  \ref{lemma E_U} , the set
 $\E_U(q)$ determines the restriction of $F_q=\F(q)$ {\mattitext in $\A$.}}
As $F_q:\overline \A\to \R$ is continuous, this determines $F_q(a)$ for all $a\in \overline\A$.
Also,  $F_q=\F(q)$ determines $\E_U(q)$ via the formula {\muutos(\ref{eq BB})}.

{\mattitext Recall that $W$  is open and relatively compact and
$\overline W\subset  J^-(p^+)\setminus I^-(p^-)$.}
 Below,
we consider 
the sets $\F(W)=\{\F(q);\ q\in W\}\subset C(\overline \A)$ and 
$\E_U(W)=\{\E_U(q);\ q\in W\}\subset 2^U$ as two representations
for  $W$. We will construct the topological and differentiable structure of $W$ using
 $\F(W)$ and the conformal class of the metric $g|_W$ using $\E_U(W)$.
First, we consider the reconstruction of the topological type of $W$.

Now we are ready to prove Proposition \ref{propo paper part 2}.  

%
\medskip

\noindent 
{\bf Proof} (of Prop.\ \ref{propo paper part 2}).
Below, let $\overline W=\hbox{cl}\,(W)$ be the closure of $W$ in $M$,
such that $\overline W\subset  J^-(p^+)\setminus I^-(p^-)$.
{\mattitext As $\overline \A\times \overline W$ is compact and thus ${\bf F}:\A\times \overline W\to\R$ is uniformly continuous
by
Lemma \ref{B: lemma} (iv),} the 
map $ \F:\overline W\to C(\overline \A)$ is continuous. Next we show that the map 
$\F:\overline W \to \F(\overline W)$     is injective. Since $\F(q)$ determines 
the set $\E_U(q)$ uniquely,  it is enough to show that the map
$  \E_U:\,  \overline W \to \E_U(\overline W)$ is injective.  
To prove this, we assume the opposite: Assume that there are  $q_1 \neq q_2$ that satisfy  $ \E_U(q_1)= \E_U(q_2)$. By Lemma \ref{lemma: from P(q) to E(q)} (iii), this implies
\beq 
\label{e.15.1}
\be_U(q_1)=\be_U(q_2).
\eeq
Choose $a\in\A$ such that $q_i \notin \mu_a$, $i\in \{1,2\}$. Let
 $(p,\eeta) \in \be_U(q_i)$ with $p=
 \mu_{a}(f^+_{a}(q_i))$. Then there are $t_i >0$ such that
$
q_i= \g_{p,\eeta}(-t_i).
$
Since \modified{$q_1\not =q_2$}, 
we have $t_1\neq t_2$,
and let us  assume that $t_2>t_1$. Then, we see there are $\xi_i\in L^+_{q_i}M$
such that
\ba
(p, \eeta)=(\g_{q_i,\xi_i}(t_i),\, \dot\g_{q_i,\xi_i}(t_i)),\quad (q_1, \xi_1)=(\g_{q_2,\xi_2}(t_2-t_1),\, \dot\g_{q_2,\xi_2}(t_2-t_1)).
\ea
\modified{Since $\rho(q,\xi)$ is lower semicontinuous, for any $\delta_1>0$ there is $\delta_2>0$
such that $\rho(q_2,\xi_2^\prime)>\rho(q_2,\xi_2)-\delta_1$ when $\xi_2^\prime \in T_qM$ satisfies
$\|\xi_2^\prime-\xi_2\|_{g^+}<\delta_2$. 
Choosing $\delta_1$ and $\delta_2$ to be sufficiently small, we have that
 there is $\xi_2^\prime \in T_{q_2}M$ that is not parallel to $\xi_2$,
$\|\xi_2^\prime-\xi_2\|<\delta_2$,
and $t_2^\prime\in (t_2-2\delta_1,t_2-\delta_1)$
such that 
$p^\prime=\gamma_{q_2,\xi_2^\prime}(t_2^\prime)\in U$, $p^\prime\not=q_1$,
and $t_2^\prime<\rho(q_2,\xi_2^\prime).$}
%
%
%
%
%
%
%
%
Thus for  $\eeta^\prime=
\dot \gamma_{q_2,\xi_2^\prime}(t_2^\prime)$ we have $(p^\prime,\eeta^\prime)\in \be_U(q_2)$.
By (\ref{e.15.1}), $(p^\prime,\eeta^\prime)\in \be_U(q_1)$,
and hence there is $t_1^\prime>0$ such that 
 $q_1=\gamma_{p^\prime,\eeta^\prime}(-t_1^\prime)$. 
 
 Observe that $\xi_1^\prime=\dot\gamma_{p^\prime,\eeta^\prime}(-t_1^\prime)$ and $\xi_1$ are not parallel.
We have that the union of \MTEXT{the
 geodesic $\gamma_{q_2,\xi_2}([0,t_2-t_1])$ and the geodesic
$\gamma_{p^\prime,\eeta^\prime}([0,-t_1^\prime])$, oriented in the opposite direction,} form a causal  {\muutos path}
from $q_2$ to $p^\prime$ that is not a light-like pre-geodesic, and hence
$\tau(q_2,p^\prime)>0$. This is not possible as $p^\prime\in \E_U(q_2)$.
This contradiction proves that $  \E_U:\,  \overline W \to \E_U(\overline W)$ is injective.

Since $ C(\overline \A)$ is a Hausdorff space, ${\overline W}$ is a compact set, and  the map $\F:\,  \overline W \to \F(\overline W)$ is
\modified{continuous and} injective,
we have that $\F: {\overline W} \to \F({\overline W})$ is a homeomorphism. Thus $\F:W \to \F(W)$
 is a homeomorphism.
\hfill \Box 
\medskip

{\mattitext 

\canberemoved{
\noindent{\bf Remark 2.1.}
We note that 
for a dense sequence  $\{q_j\}_{j\in \Z_+}$ in  $W$ the 
{\rtext family} of light observations  $\{\E_U(q_j)\}_{j\in \Z_+}$
determines the {\rtext family} $\{\F(q_j)\}_{j\in \Z_+}\subset C(\overline\A)$.
These points are dense in  $ \F({\overline W})$. Thus the
closure of this {\rtext family} in $C(\overline \A)$
 determines the whole set   $\F({\overline W})$, or equivalently,
 the set  $\E_U({\overline W})$.
}

}

\subsubsection{Estimates for the location of the first cut point.}

 {\mattitext We finish this section by auxiliary results that are needed in the proof of  Theorem \ref{main thm3 new}.} 
{\mattitext Below,} we use  for a pair $(x,\xi)\in L^+\hattuM $  the notation 
\beq\label{eq: x(h) notation}
& &(x(t),\xi(t))=(\gamma_{x,\xi}(t),\dot \gamma_{x,\xi}(t)),\quad t\in \R_+.
\eeq

\begin{center}
\psfrag{1}{${\hat x}$}
\psfrag{2}{$x$}
\psfrag{3}{$z$}
\psfrag{6}{\hspace{-1mm}$q_1$}
\psfrag{5}{$q_0$}
\psfrag{4}{\hspace{-3mm}$p_1$}
\psfrag{7}{\hspace{-1mm}$x_1$}
\psfrag{8}{$\gamma_{{x},\xi}$}
\psfrag{9}{$\gamma_{{\hat x},\hat \xi}$}
\includegraphics[height=5.5cm]{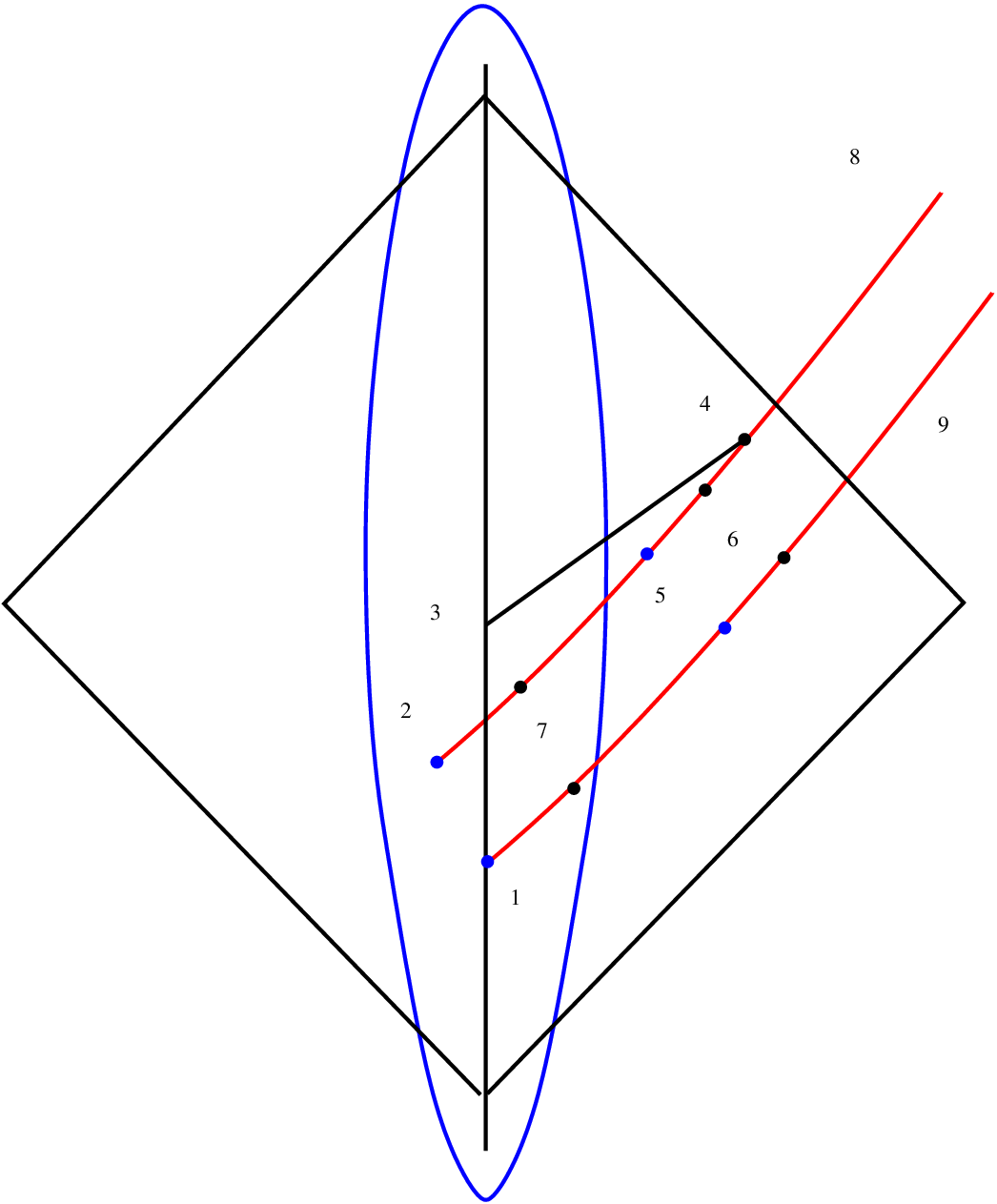}
\psfrag{1}{$x_1$}
\psfrag{2}{\hspace{-1mm}$x_2$}
\psfrag{3}{$q$}
\psfrag{4}{\hspace{-0mm}$p_1$}
\psfrag{5}{\hspace{-2mm}$p_2$}
\psfrag{6}{}
\psfrag{7}{$y_0$}
\includegraphics[height=5.5cm]{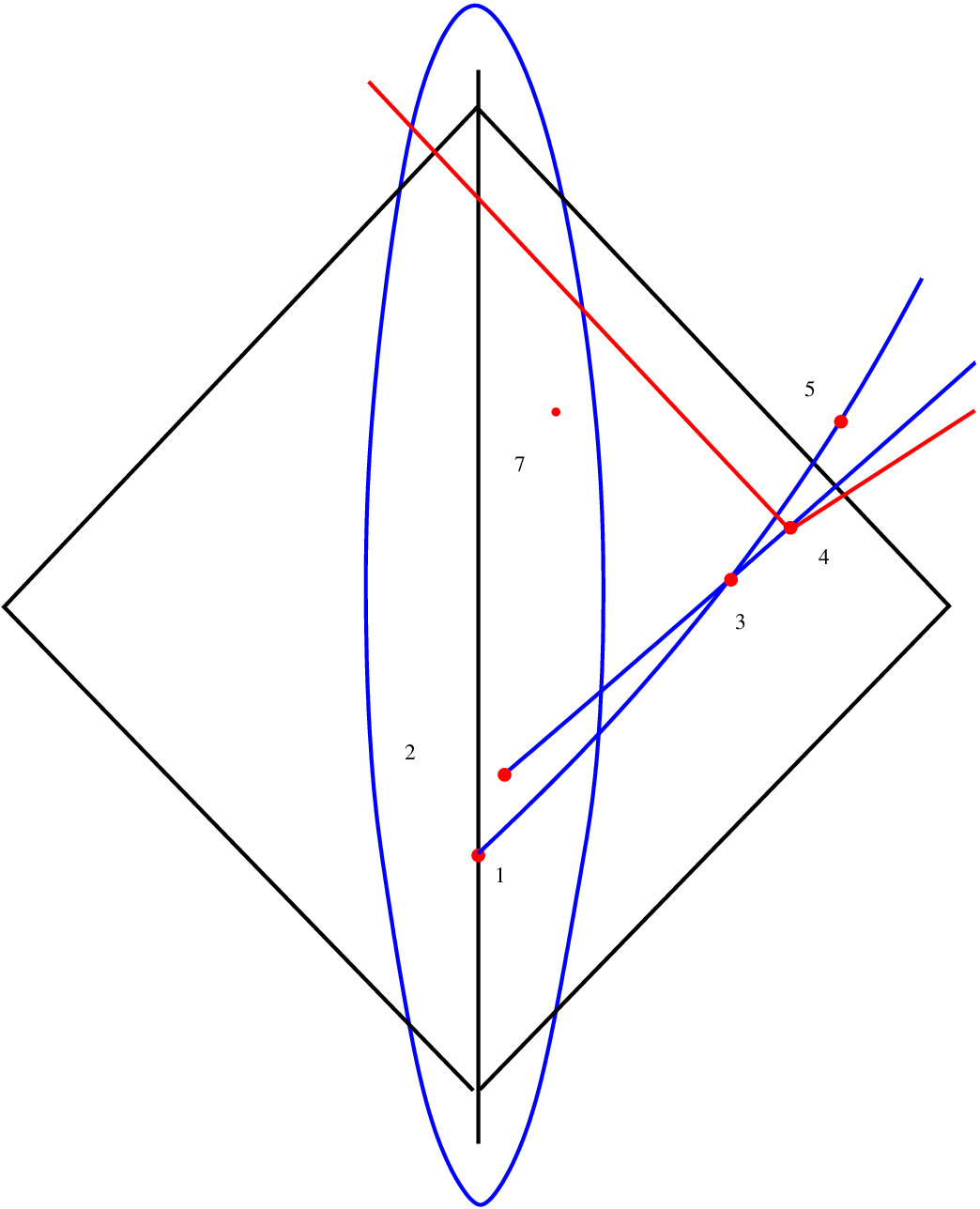}
\end{center}


\noindent
{\it FIGURE 5. {\bf Left:} The figure shows the situation in Lemma \ref{lem: detect conjugate 0}. The point ${\hat x}=\hat \mu(r_1)$ is on the time-like
{\rtext path} $ \hat \mu$ shown as a black line. 
The black diamond is the set $J^+( p^-)\cap J^-(p^+)$,
 $(x,\xi)$ is a light-like direction close to $({\hat x},\hat\xi)$,
 and $x_1=\gamma_{{x},\xi}(t_0)=x(t_0)$.
The points $q_0=\gamma_{x,\xi}(\rho(x,\xi))$ and $q_1=
\gamma_{x(t_0),\xi(t_0)}(\rho(x(t_0),\xi(t_0)))$ are the  first cut point on  
$\gamma_{{x},\xi}$ corresponding to the points $x$ and $x_1$,
respectively. The blue and black points on $\gamma_{{\hat x},\hat \zeta}$ are the corresponding
cut points on 
 $\gamma_{{\hat x},\hat \zeta}$.
Also, $p_1=\gamma_{x,\xi}(t_1)$ and $z=\hat \mu(r_1)$,  where
  $r_1=f^-_{{\hat a}}(p_1)$.
{\bf Right:}
The figure shows the configuration in formulas (\ref{eq: summary of assumptions 1}) and  (\ref{eq: summary of assumptions 2}) and in Theorem \ref{lem:analytic limits A wave}.
We send light-like geodesics $\gamma_{x_j,\xi_j}([0,\infty))$ from $x_j$,
$j=1,2,3,4.$ The boundary $\p  \V(\vec x,\vec \xi)$ is denoted by red line segments
and $y_0 \in \V(\vec x,\vec \xi)$. We assume the these  geodesics intersect
at the point $q$ before their 
 first cut points $p_j$.
}

\medskip

\motivation{Later, we will consider  waves sent from 
a point $x\in U$ that propagate
near a geodesic $\gamma_{x,\xi}([0,\infty))$. These waves may have singularities near the conjugate
points of the geodesic and due to this we analyze next how the conjugate points
move along a geodesic when its initial point is moved
from $x$ to $\gamma_{x,\xi}(t_0)$. } {\rbtext Below, let  
$T_{+}(x,\xi)=\sup \{t\geq 0\ ;\ \gamma_{x,\xi}(t)\in J^-(\hat \mu(1))\}$.}

%

 \MATTITEXT{

 \begin{lemma} \label{lem: detect conjugate 0} 
   There are 
$\vartheta_1,\kappa_1,\kappa_2>0$ such that 
for  all  $\hat x=\hat \mu(r_0)$ with $r_0\in [s_-,s_+]$,
${\hat \xi}\in L^+_{\hat x}M$, $\|{\hat \xi}\|_{{\nohat g}^+}=1$,
$t_0\in [\kappa_1,4\kappa_1]$,  
and 
$(x,\xi)\in L^+M$ satisfying $d_{{\nohat g}^+}((x,\xi),({\hat x},{\hat \xi}))\leq \vartheta_1$
the following holds:

(ii) If  $0<t\leq 5\kappa_1$, then  $f^-_{{\hat a}}(\gamma_{\hat x,{\hat \xi}}(t))=r_0$,

(ii)  If $0<t\leq 5\kappa_1$, then $\gamma_{x,\xi}(t)\in U$, 


(iii) 
Assume that there exists $t_1$ that satisfies
$t_0+\rho(\gamma_{x,\xi}(t_0),\dot \gamma_{x,\xi}(t_0))\leq t_1< T_{+}(x,\xi)$  and let
 $p_1=\gamma_{x,\xi}(t_1)$. 
Then    
 $r_1=f^-_{{\hat a}}(p_1)$ satisfies $r_1-r_0>2\kappa_2$.

\end{lemma}

Note  that above in (iii) we can choose $t_1= t_0+\rho(\gamma_{x,\xi}(t_0),\dot \gamma_{x,\xi}(t_0))$ in which case
$p_1$ is the first cut point $q_1$ of $\gamma_{x,\xi}([t_0,\infty))$, see Fig.\ 5(Left).

{\tobecheckedtext 

\noindent
{{\bf Proof.} 
  Let $B=\{({\hat x},{\hat \xi})\in L^+M;\ {\hat x}\in {\hat \mu}([s_-,s_+]),\ \|{\hat \xi}\|_{  {\nohat g}^+}=1\}$.
Since $B$ is compact, the positive and lower semi-continuous function  $\rho(x,\xi)$ obtains
its minimum on $B$.  This proves the claim (i) when $\kappa_1\in (0, \frac 15 \inf\{\rho({\hat x},{\hat \xi});\ ({\hat x},{\hat \xi})\in B\})$. 

(ii) For $ \vartheta>0$ small enough, 
 $K_\vartheta=\{(x,\xi)\in L^+M;\ d_{{\nohat g}^+}(({x},{\xi}),B)\leq \vartheta\}$
 is a compact subset of $L^+U$.
Thus yields easily  (ii) when
$ \vartheta_1$ is small enough.

(iii)  Let  $\vartheta\in (0,\vartheta_1)$ be so small
that $K_\vartheta\subset L^+U$ 
and \ba
K_\vartheta^0=\{(x,\xi)\in K_\vartheta; \ {\rho(x(\kappa_1),\xi(\kappa_1))+\kappa_1
} 
\leq T_{+}(x,\xi)\},
\quad K_\vartheta^1=K_\vartheta\setminus K_\vartheta^0.
\ea 
Using  \cite[{\ftext Lemma}\ 14.13]{ONeill}, we see that $T_{+}(x,\xi)$ is bounded in $K_\vartheta$.
 {Note that for $t_0\geq \kappa_1$ and $a>t_0$ the geodesic
$\gamma_{x,\xi}([t_0,a])$ can have a cut point only if $\gamma_{x,\xi}([\kappa_1,a])$
has a cut point and thus
$t_0+\rho(x(t_0),\xi(t_0))\geq \kappa_1+\rho(x(\kappa_1),\xi(\kappa_1))$.
}
If $K_\vartheta^0=\emptyset$,
 the claim is valid as the condition $p_1\in J^-(\hat \mu(1))$  does not hold for any $(x,\xi)\in K_\vartheta^1$.
 Thus it is enough to consider the case when  $K_\vartheta^0\not=\emptyset$.
 
Let
\ba
G_ \vartheta =\{(x,\xi,t)\in K_\vartheta \times \R_+;\ \rho(x(\kappa_1),\xi(\kappa_1))+\kappa_1\leq  t\leq T_{+}(x,\xi)\}.
\ea
As $\rho(x,\xi)$ is lower semi-continuous
and $T_{+}(x,\xi)$ is upper semi-continuous  and bounded, the sets $K_ \vartheta^0 $  and $G_ \vartheta $ are compact.

{For $(x,\xi,t)\in G_ \vartheta $,
the geodesic $\gamma_{x,\xi}([\kappa_1,t])$ has a cut point. Thus   for $y=\gamma_{x,\xi}(t)$, 
we have
$\tau(x,y)>0$. 
Hence, for $z=\hat \mu(f^-_{{\hat a}}(x))$, we have $\tau(z,y)\geq
\tau(z,x)+\tau(x,y)\geq \tau(x,y)>0$. This shows that $f^-_{{\hat a}}(y)- f^-_{{\hat a}}(x)>0$.
Since $G_ \vartheta $ is compact and $f^-_{\hat a}$
 is continuous, $\e_1:=\inf\{ f^-_{{\hat a}}(\gamma_{x,\xi}(t))- f^-_{{\hat a}}(x);\
(x,\xi,t)\in G_ \vartheta \}>0$. 

}

Then, if  $\rho(x(\kappa_1),\xi(\kappa_1))+\kappa_1\leq t_1<\mathcal T_+(x,\xi)$ 
and $p_1=\gamma_{x,\xi}(t_1)$, we have that
 $r_1=f^-_{{\hat a}}(p_1)$ and  $r_2=f^-_{{\hat a}}(x)$ satisfy 
$r_1-r_2\geq \e_1$.

As  $f^-_{{\hat a}}$ is continuous and
$\hat \mu([-1,1])$ is compact, we see that by making $\vartheta_1$ smaller if necessary,
we can assume that if ${\hat x}=\hat \mu(r_0)\in \hat \mu$ and $d_{{\nohat g}^+}(x,{\hat x})\leq \vartheta_1$ then
 $|f^-_{{\hat a}}(x)-f^-_{{\hat a}}({\hat x})|<\e_1/2$. 
%
Let $\kappa_2=\e_1/4$.
Then $r_1-r_2\geq \e_1$ and 
$r_2-r_0=|f^-_{{\hat a}}(x)-f^-_{{\hat a}}({\hat x})|<\e_1/2$ imply that $r_1-r_0>\e_1/2=2\kappa_2$.
This proves the claim.}
%
\hfill \Box \medskip
}

{\mattitext Finally, consider the case  when $(M_1,\nohat g_1)$ and $(M_2,\nohat g_2)$
are two manifolds satisfying  
(\ref{eq: Es the same}) with the sets $W_1$ and $W_2$ and
a time-orientation preserving conformal diffeomorphism $\Phi:U_1\to U_2$.
Then, if $U_1$  is defined using paths $\mu_a^{(1)}(s),$  $a\in \overline \A$, $s\in [-1,1]$,
by making $\A$ a smaller neighborhood of $\hat a$ if necessary,  we can use on 
$U_2$ the paths  $\mu_a^{(2)}(s)=\Phi(\mu_a^{(1)}(s)),$  $a\in \overline \A$.
With such paths the sets $\F(W_1)\subset C(\overline \A)$ on manifold $M_1$
and $\F(W_2)\subset C(\overline \A)$ on manifold $M_2$ coincide.}

\section{Inverse problem for active measurements}\label{Sec: Analysisforwave}

{\mattitext In this section we start the proof of Theorem \ref{main thm3 new}.} 
{\gtext Without loss of generality 
we may replace the  set $V$ where we do measurements by a smaller set $U$
of the form (\ref{eq: Def Wg with hat}). 
Also, {\muutos by redefining the path {\mmmmtext $\hat \mu_2$} in the claim of Theorem \ref{main thm3 new}, 
we can assume that {\mmmmtext$\hat \mu_2=\Phi(\hat \mu_1)$.}}
Moreover, as the proof  is constructive, and to simplify the notations, we do 
the constructions on just one Lorentzian manifold, $(M, g)$
and assume that we are given the data
\beq\label{eq: data for active}
& &\hbox{the differentiable manifold $U$ of the form  (\ref{eq: Def Wg with hat}), paths $\mu_a,$ $a\in \A$, }\\
\nonumber  & &\hbox{and the source-to-solution map $L_U$.}
\eeq 
Here, $L_U:f\mapsto u|_U$ is the source-to-observation map defined
in (\ref{measurement operator}) when the set $U\subset V$ is used as the {\muutos measurement} set.
{\itext The choice of paths $\mu_a$ are discussed in  Remark 3.1 below.}

\subsection{Asymptotic expansion for  the non-linear wave equation}\label{subset: AENWE 1}

{\revtext
Let us consider the non-linear wave equation
\beq\label{PABC eq}
& &\square_{\nohat g}u+au^2=f,\quad\hbox{in }{\muutos M_0=}(-\infty,{\reftext T_0})\times N,\\
& &u|_{(-\infty,0)\times N}=0,\nonumber
\eeq
where
$a=a(x)$ is a smooth, {nowhere} vanishing function, ${\muutos M_0=}(-\infty,{\reftext T_0})\times N\subset M=\R\times N$,
where $(M,\nohat g)$ is a globally hyperbolic Lorentzian manifold. {\vtext We denote by 
${{{\muutos \square_g^{-1}}}}$} the causal inverse  operator of $\square_{\nohat g}$.

When  {$B\subset N$  is compact and} $f$ in $C_0([0,{\reftext T_0}];H^6_0(B))\cap
C_0^1([0,{\reftext T_0}];H^5_0(B))$ 
is small enough,  
we see by using  \cite[Prop.\ 9.17]{Ringstrom},  \cite [Thm.\ III]{HKM}, or \cite[App.\ III]{ChBook}
that  the
equation (\ref{PABC eq}) has a unique solution
$u\in  C([0,{\reftext T_0}];H^{5}(N))\cap C^1([0,{\reftext T_0}];H^{4}(N))$. 
\noextension{For a detailed analysis, see Appendix B in \cite{preprint}.} \extension{(For details, see  Appendix B).} 

Let us consider the case when $f=\e f_1$ where $\e>0$
is small.
Then, we can write
\ba
u=\e w_1+\e^2 w_2+\e^3 w_3+\e^4 w_4+E_\e
\ea
where $w_j$ and the reminder term $E_\e$ satisfy {\gtext (see e.g.\ \cite[App.\ III]{ChBook})}
\beq\label{Slava formula 2}
w_1&=&{{{\muutos \square_g^{-1}}}}f_1,\\ \nonumber
w_2&=&-{{{\muutos \square_g^{-1}}}}(a\,w_1\,w_1),\\ \nonumber 
w_3
&=&2\,{{{\muutos \square_g^{-1}}}}(a\, w_1\,{{{\muutos \square_g^{-1}}}}(a\,w_1\,w_1))
,\\ \nonumber
w_4
&=&-{{{\muutos \square_g^{-1}}}}(a\, {{{\muutos \square_g^{-1}}}}(a\,w_1\,w_1)\,{{{\muutos \square_g^{-1}}}}(a\,w_1\,w_1))\\ \nonumber
& &-4\,{{{\muutos \square_g^{-1}}}}(a\, w_1\,{{{\muutos \square_g^{-1}}}}(a\, w_1\,{{{\muutos \square_g^{-1}}}}(a\,w_1\,w_1))),\\ \nonumber
& &\hspace{-1.5cm}\|E_\e\|_{C([0,{\reftext T_0}];H^{4}_0(N))\cap  
C^1([0,{\reftext T_0}];H^{3}_0(N))}\leq {\muutos C(f_1)}\,\e^5.
\eeq
{\revtext {\mattitext In particular,} we will consider sources $f_1$ for which the linearized term
 $w_1$ is a distorted plane wave.
}

\medskip

\noindent{\bf Remark 3.1.}
 {\mattitext The set $ U$, given as differentiable manifold, and the  {\slava source-to-solution} map $L_ U$ determine the linearized {\slava source-to-solution} map
 $L_ U^{lin}:{\muutos f_1}\mapsto \p_\e(L_U(\e {\muutos f_1}))|_{\e =0}$. Furthermore, this map determines all
 pairs $({\muutos f_1},w_1)$ such that $L_U^{lin}({\muutos f_1})=w_1$ and both ${\muutos f_1}$ and $w_1$ are compactly supported 
 in $U$.  Observer that then ${\muutos \square_gw_1={\muutos f_1}}$.
  In particular, for any $(x_0,\eta_0)\in T^*U$ there
  {\reftext is} \refHOX{Comment 19: "is" is added.} a pair $({\muutos f_1^\tau},w_1^\tau)$
  such that
   $ w^\tau_1(x)=e^{i\tau\phi(x)}\psi(x),
  $
where $\tau>0$, $\phi,\psi\in C^\infty_0(U)$, $d\phi(x_0)=\eta_0$, in some neighborhood of $x_0$ 
we have $\psi=1$. 
Then 
\ba
g(\eta_0,\eta_0)=-\lim_{\tau\to \infty} \frac {{\muutos \square_g w^\tau_1(x_0)}}{\tau^2}
=\lim_{\tau\to \infty} \frac {-{\muutos f_1^\tau}(x_0)}{\tau^2}.
\ea
This shows that $U$ and $L_U$ determine the metric tensor $g|_U$ in $U$.
{\itext The set of pairs $({\muutos f_1},L_U{\muutos f_1})$ that are in $C^\infty_0(U)^2$ coincide with the set of the pairs 
$\{({\muutos \square_g \phi+a\phi^2,\phi});\ \phi \in C^\infty_0(U)\}$. When $g|_U$ is known, these pairs
determine $a|_U$. Hence, $L_U$  determines also  $a|_U$.}
{\gtext Observe that by the same arguments, 
 $V$ and $L_V$ determine the metric tensor $g|_V$ in $V$, too.}}
{\muutos Also, we note that when
 $V$ and $g|_V$ are given, one can choose the time-like paths
$\mu_a:[-1,1]\to V$, $a\in \A$, appearing in (\ref{eq: Def Wg with hat}), to be perturbations of the path $\hat \mu([-1,1])$ that depend smoothly on a parameter $a$ in an open set. {\muutos Thus the data $(V,L_V)$ can be used to construct 
{\mmmmtext the paths $\mu_a$ and}
the set 
 $U\subset V$ in (\ref{eq: Def Wg with hat}).}}
%

 \subsection{Linear wave equation and distorted plane waves}
 \subsubsection{Lagrangian distributions}\label{Lagrangian distributions}

 Let us recall 
{\reftext the} \refHOX{Comment 20: "the" is added.}
definition of the {\reftext classical} \refHOX{Comment 21: "classical" is added.} 
conormal and Lagrangian distributions that we will use below, see \cite{GU1,H4,MU1}.
Let $X$ be a manifold of dimension $n$
and $\Lambda\subset T^*X\setminus \{0\}$ be
a Lagrangian submanifold.
Let  
$\phi(x,\theta)$, $(x,\theta)\in X\times \R^N$ be a non-degenerate
phase function that locally parametrizes $\Lambda$ {\muutos near a point  $(x_0,\xi_0)\in \Lambda$},
{\mattitext i.e., in a conic neighborhood $\Gamma\subset T^*X\setminus \{0\} $ {\muutos of  $(x_0,\xi_0)$,}
  the submanifold  $\Lambda$  coincides with the set $\{(x,d_x\phi(x,\theta))\in \Gamma;\ d_\theta\phi(x,\theta)=0\}$.}
We say that a distribution $u\in {\cal D}^{\prime}(X)$  
is a {\reftext classical} \refHOX{Comment 21: "classical" is added.} Lagrangian distribution associated with $\Lambda$ and denote
$u\in \I^m(X;\Lambda)$, if in local coordinates {\mmmmtext $X:W\to \R^n$,} $u$ can be represented 
 as an oscillatory integral, 
\beq\label{def: conormal}
u(x)=\int_{\R^N}e^{i\phi(x,\theta)}a(x,\theta)\,d\theta,\quad {\mmmmtext x\in W}
\eeq
where $a(x,\theta)\in S^\mu({\mmmmtext W}; \R^N)$ is a {\mattitext classical symbol of order $\mu= m+n/4-N/2$,
see  \cite{GU1,H4,MU1}.

{{\guntext 
For classical Lagrangian distributions $u\in \I^m(X;\Lambda)$ one can define a principal
symbol $\sigma_u^{(p)}(x_0,\zeta_0)$  of $u$, at $(x_0,\zeta_0)\in \Lambda$, that satisfies
$$
\sigma_u^{(p)}(x_0,\zeta_0)\in  S^{m+\frac n4}(\Lambda,\Omega^{1/2}\otimes L)/ S^{m+\frac n4-1}(\Lambda,\Omega^{1/2}\otimes L),
$$
where $L$ is the Maslov-Keller line bundle and $\Omega^{1/2}$ are  the half-densities on $X$,
on details, see \cite[Thm.\ 11.10]{GrS}. 
We note that below we {\vtext do computations }using only principal symbols of conormal distributions considered below.

}

In particular, when $S\subset X$ is a submanifold,  its conormal bundle
$N^*S=\{(x,\xi)\in T^*X\setminus \{0\};\ x\in S,\ \xi\perp T_xS\}$ is a Lagrangian submanifold.
If $u$ is a  Lagrangian distribution  associated to $\Lambda_1$
where $\Lambda_1=N^*S$, we say that $u$ is a 
{\reftext (classical)} \refHOX{Comment 21: "classical" is added.}
conormal distribution.

{\slava Let us next consider the case when $X=\R^n$,
$(x^1,x^2,\dots,x^n)$
are the Euclidean coordinates and $x^{\prime}=(x_1,\dots,x_{d_1})$,
 $S_1=\{0\}\times \R^{n-d_1}=\{x^{\prime}=0\}\subset \R^n$} and $\Lambda_1=N^*S_1$.
Then $u\in \I^m(X;\Lambda_1)$ can be represented by (\ref{def: conormal})
with $N=d_1$ and  $\phi(x,\theta)=x^{\prime}\cdotp \theta$,
that is, 
\beq\label{def: conormal e}
u(x^1,\dots,x^n)=\int_{\R^{d_1}}e^{ix^{\prime}\cdotp \theta}a(x^1,\dots,x^n,\theta)\,d\theta.
\eeq 
{For example, $\delta_{S_1}(x)\in  \I^{-n/4+d_1/2} (\R^n;N^*S_1)$, \vtext where $\delta_{S_1}(x)$ denotes
the Dirac delta distribution supported on $S_1$.}

%
{{\guntext  The principal symbol of a conormal distribution $u\in \I^m(\R^n;N^*S_1)$, represented in
the form (\ref{def: conormal e}),
can be identified with a 
 function  $c(x,\theta)$ that is $\mu$-positive homogeneous in $\theta$, such that
$a(x,\theta)-(1-\phi(\theta))c(x,\theta)\in S^{\mu-1}(X;\R^{d_1})$ 
where $\phi\in C^\infty_0(\R^{\mmmmtext d_1})$ is 1 in a neighborhood of zero. For 
 a manifold $X$  and a surface $S\subset X$, we can use this definition  to define a principal symbol of a conormal distribution $u\in \I^m(X;N^*S)$ in local coordinates. On the invariant
 nature of this definition, see 
  \cite[Sec. 18.2]{H3}.} }}

\refHOX{Comment 33: Notation of principal symbols is  clarified. }

Next we recall the definition of $\I^{p,l}(X;\Lambda_1,\Lambda_2)$, the space of
the distributions $u$ in ${\cal D}^{\prime}(X)$ associated to two 
cleanly intersecting Lagrangian manifolds $\Lambda_1,\Lambda_2\subset T^*X
\setminus \{0\}$, see \cite {CGP,GU1,MU1}. 
{\muutos We recall that $\Lambda_1$  and $\Lambda_2$ intersect cleanly if
$\Sigma=\Lambda_1\cap \Lambda_2$  is a smooth manifold and its tangent space
satisfies $T_\lambda\Sigma=T_\lambda\Lambda_1\cap T_\lambda\Lambda_2$ for all
$\lambda\in \Sigma$.}
These classes have been widely
used in the study of inverse problems, see \cite {CS2,FG}.  
Let us start with the case when $X=\R^n$.

{\slava
 Let 
$(x^1,x^2,\dots,x^n)$
be the Euclidean coordinates in $\R^n$.
Let $S_1,S_2\subset \R^n$ be the linear subspaces of codimensions
$d_1$ and $d_1+d_2$, respectively be such that
$S_2\subset S_1$.
We use in $\R^n$ the Euclidean coordinates $(x^1,x^2,\dots,x^n)=(x^{\prime},x^{\prime\prime},x^{\prime\prime\prime})$ 
where  $x^{\prime}=(x_1,\dots,x_{d_1})$,
  $x^{\prime\prime}=(x_{d_1+1},\dots,x_{d_1+d_2})$,
 $x^{\prime\prime\prime}=(x_{d_1+d_2+1},\dots,x_{n})$
and assume that}
$S_1=\{x^{\prime}=0\}$, $S_2=\{x^{\prime}=x^{\prime\prime}=0\}$.
Let us denote $\Lambda_1=N^*S_1,$  $\Lambda_2=N^*S_2$. Then
$u\in \I^{p,l}(\R^n;N^*S_1,N^*S_2)$ if and only if
\beq\label{def: conormal2}
u(x)=\int_{\R^{d_1+d_2}}e^{i (x^{\prime}\cdotp \theta^{\prime}+
x^{\prime\prime}\cdotp \theta^{\prime\prime})}a(x,\theta^{\prime},\theta^{\prime\prime})\,d\theta^{\prime}d\theta^{\prime\prime},
\eeq
where the symbol $
a(x,\theta^{\prime},\theta^{\prime\prime})$ belongs in the product type symbol class
$S^{\mu_1,\mu_2}(\R^n; (\R^{d_1}\setminus 0)\times \R^{d_2})$
that is the space of 
{\reftext functions} \refHOX{Comment 23: "function" changed to "functions".}
  $a\in C^{\infty}(\R^n\times \R^{d_1}\times \R^{d_2})$ that satisfy
\beq\label{product symbols}
|\p_x^\gamma\p_{\theta^{\prime}}^\alpha \p_{\theta^{\prime\prime}}^\beta a(x,\theta^{\prime},\theta^{\prime\prime})|
\leq C_{\alpha\beta\gamma K}(1+|\theta^{\prime}|+|\theta^{\prime\prime}|)^{\mu_1-|\alpha|}
(1+|\theta^{\prime\prime}|)^{\mu_2-|\beta|}\hspace{-1.7cm}
\eeq
for all $x\in K$, multi-indexes  $\alpha,\beta,\gamma$, and
compact sets $K \subset \R^n$. Above,
$\mu_1=p+l-d_1/2+n/4$ and
$\mu_2=-l-d_2/2$.

{{\guntext  When $X$ is a manifold of dimension $n$ and 
$\Lambda_1,\Lambda_2\subset T^*X
\setminus \{0\}$ are two 
\refHOX{Comment 24: Details of   definition of 
$\I^{p,l}(X;\Lambda_1,\Lambda_2)$  and the canonical relation of the FIO $A$ are clarified.}
cleanly intersecting Lagrangian manifolds, 
we define the class $\I^{p,l}(X;\Lambda_1,\Lambda_2)\subset \mathcal D^\prime(X)$ to consist
of locally finite sums of distributions of the form $u=Au_0$, where 
$u_0\in \I^{p,l}(\R^n;N^*S_1,N^*S_2)$ and $S_1,S_2\subset \R^n$ are the linear subspace of codimensions
$d_1$ and $d_1+d_2$, respectively, such that
$S_2\subset S_1$, 
and $A$ is a Fourier integral 
operator of order zero with a canonical relation $\Sigma$ for which
$ \Sigma\circ (N^*S_1)^\prime\subset \Lambda_1^\prime$ and 
$\Sigma\circ (N^*S_2)^\prime\subset \Lambda_2^\prime$. 
Here, for $\Lambda\subset T^*X$ we denote $ \Lambda^\prime=\{(x,-\xi)\in T^*X;\ (x,\xi)\in \Lambda\}$.
{\mtext The definition of $\I^{p,l}(X;\Lambda_1,\Lambda_2)$  is discussed in detail in  \cite{MU1},
in particular the existence of the canonical relation $\Sigma$ connecting the pair $(\Lambda_1,\Lambda_2)$ of cleanly intersecting Lagrangians
to the microlocal model $(N^*S_1,N^*S_2)$ is proven in  \cite[Prop.\ 1.3]{MU1}.}}
{When $X$  and $Y$  are manifolds and $\Sigma\subset T^*X\times  T^*Y$  we use also the notation
$\Sigma^\prime=\{(x,\xi,y,-\eta);\ (x,\xi,y,\eta)\in \Sigma\}.$}


In most cases {\reftext below,}\refHOX{Comment 25: Comma is moved.} $X=M$. We denote then
  $\I^{p}(M;\Lambda_1)=\I^{p}(\Lambda_1)$
and 
 $\I^{p,l}(M;\Lambda_1,\Lambda_2)=\I^{p,l}(\Lambda_1,\Lambda_2)$.
Also, $\I(\Lambda_1)=\cup_{p\in \R}\I^{p}(\Lambda_1)$.

{\mmtext By \cite{GU1,MU1}, if $R_1$ and $R_2$ are pseudodifferential operators of order zero on $M$  which 
are microlocally smoothing in a conic neighborhood of  $\Lambda_2$ and $\Lambda_1$,
respectively, we have} 
\beq\label{microlocally away}
& &R_1:\I^{p,l}(\Lambda_{\muutos 1},\Lambda_{\muutos 2})\to \I^{p+l}(\Lambda_{\muutos 1}),\quad 
R_2: \I^{p,l}(\Lambda_{\muutos 1},\Lambda_{\muutos 2})\to \I^{p}(\Lambda_{\muutos 2}).
\eeq
Thus the principal symbol of $u\in \I^{p,l}(\Lambda_{\muutos 1},\Lambda_{\muutos 2})$
is well defined on $\Lambda_{\muutos 1}\setminus \Lambda_{\muutos 2}$ and  $\Lambda_{\muutos 2}\setminus \Lambda_{\muutos 1}$. We denote  $\I(\Lambda_{\muutos 1},\Lambda_{\muutos 2})=\cup_{p,l\in \R}\I^{p,l}(\Lambda_{\muutos 1},\Lambda_{\muutos 2})$. 
{\muutos We recall that $(x_0,\xi_0)\in T^*M$ belong in the wave front set $\hbox{WF}\,(u)$ of a distribution
$u\in \mathcal D'(M)$ if $u(x)$ is not $C^\infty$  smooth near $x_0$  in the direction $\xi_0$ (see \cite{Duistermaat}, 
Section 1.3 for the precise definition). For Lagrangian distributions $v\in \I^p(\Lambda_1)$ and 
$u\in \I^{p,l}(\Lambda_1,\Lambda_2)$ we have
\beq\label{eq wavefronts}
\hbox{WF}\,(v)\subset \Lambda_1,\quad \hbox{WF}\,(u)\subset \Lambda_1\cup \Lambda_2.
\eeq}

Below, when $\Lambda_j=N^*S_j,$ $j=1,2$ are conormal bundles of smooth cleanly
intersecting submanifolds
$S_j\subset M$ of codimension ${\muutos d}_j$, where $\dim(M)=n$,
 we
use the traditional notations,
\beq\label{eq: traditional}\quad\quad
\I^\mu(S_1)=\I^{\mu+{\muutos d}_1/2-n/4}(N^*S_1),\quad \I^{\mu_1,\mu_2}(S_1,S_2)
=\I^{p,l}(N^*S_1,N^*S_2),\hspace{-1cm}
\eeq
where $p=\mu_1+\mu_2+{\muutos d}_1/2-n/4$ 
 and $l=-\mu_2-{\muutos d}_2/2$,
and call such distributions the conormal distributions  associated to $S_1$ or 
product type  conormal distributions  associated to
$S_1$ and $S_2$, respectively. 
 {By \cite{GU1},
$\I^\mu(X;S_1)\subset L^p_{loc}(X)$ for $\mu< -{\muutos d}_1(p-1)/p$, $1\leq p<\infty$.
Further developments
for the theory of the paired Lagrangian distributions are in \cite{GuU1,GU2}.}

\subsubsection{Inverse of the linear wave operator}
Next we will shortly discuss how {\vtext paired} Lagrangian distributions are used in
\cite{MU1,GU1} to study parametrices ({\muutos and inverses) of real-principal type operators}, in particular
the wave operator   $\square_g$ on a globally hyperbolic Lorentzian manifold $(M,g)$.
 To consider the wave operator, recall that}
the characteristic variety of 
$\square_g$ is $${\muutos \hbox{Char}(\square_g)}=\{(x,\xi)\in T^*M{\muutos \setminus 0};\ p(x,\xi)=0\},$$
where $p(x,\xi)=g^{jk}(x)\xi_j\xi_k$.}
For the
wave operator,
$\hbox{Char}\,(\square_{g})$ is the set of light-like co-vectors with respect to $g$.
Also, a bicharacteristic of $\square_{g}$ is the integral curve of the Hamiltonian vector field of ${\muutos p(x,\xi)}$ in  $T^*M$.
For $(x,\xi)\in \hbox{Char}\,(\square_{g})$, {\reftext we denote by $\Theta_{x,\xi}\subset T^*M$ the bicharacteristic of $\square_{g}$ that contains $(x,\xi)\in L^*M$.}  \refHOX{Comment 26: Definition of  bicharacteristics 
is clarified.} 
{\reftext The bicharacteristics are closely 
related to light-like {\mmmmtext geodesics}: We have 
$(y,\eta)\in \Theta_{x,\xi}$ if and only if there is $t\in\R$ such
that  
for ${\muutos v}=\eta^\sharp$ and ${\muutos w}=\xi^\sharp$ we have
 $(y,{\muutos v})=(\gamma_{x,{\muutos w}}(t),\dot \gamma_{x,{\muutos w}}(t))$
where $\gamma_{x,{\muutos w}}$ is a light-like geodesic with respect to the metric $g$
with the initial data $(x,{\muutos w})\in LM$. Here, we use {\vtext the} notations $(\xi^\sharp)^j=g^{jk} \xi_k$ and $({\muutos w}^\flat)_j=g_{jk} {\muutos w}^k$.} 

%
%
%
Let $\Lambda_1\subset T^*M$ be a Lagrangian manifold and consider the solution of ${\muutos \square_g}u_1=f_1$ with a source $f_1\in \I^m(\Lambda_1)$. 
{\muutos When the  characteristic variety ${\muutos \hbox{Char}(\square_g)}$ intersects $\Lambda_1$, this gives rise to the propagation of singularities.
%
%
 \refHOX{Comment 22: text on new singularities 
is modified.} 
{\reftext Indeed, by H\"ormander's theorem on {\ftext propagation of singularities} along bicharacteristics, \cite[Theorem 26.1.4]{H4}, see also \cite[Prop. 2.1]{GU1}, the wave front set  $\hbox{WF}(u_1)$ of $u_1$  is contained in the union of  $\Lambda_1$ and the bicharacteristics that contain points of  the intersection ${\muutos \hbox{Char}(\square_g)}\cap \Lambda_1$. When ${\muutos \hbox{Char}(\square_g)}$  and $\Lambda_1$ intersect transversally,
the union of these bicharacteristics is a {\mmmmtext Lagrangian 
manifold.} This result was extended in  \cite{MU1,GuU1} 
where it was shown that the Schwartz kernel of the inverse of the wave operator is a distribution
associated to two intersecting Lagrangian manifolds.}}
 Indeed, when
 $(M,g)$ is a globally hyperbolic manifold, the operator  
${\muutos \square_g}$ has a causal inverse operator $Q={{{\muutos \square_g^{-1}}}}$, see e.g.\  \cite[Thm.\ 3.2.11]{BGP}.
{\slava A geometric representation for its kernel is given in \cite{K2}. Below, we often use the same notation for
the operator $Q$ with its Schwartz kernel $Q(x,y)$.
By  \cite{MU1}, 
the {\slava Schwartz kernel $Q$}}
satisfies 
$Q\in \I^{p,l}(\Delta^\prime_{T^*M},\Lambda_g)$, 
$p=-\frac 32$, $l=-\frac 12$. Here,
 ${\mmmmtext \Delta_{T^*M}^\prime}=N^*(\{(x,x);\ x\in M\})$, and $\Lambda_g\subset T^*M\times T^*M$
is the Lagrangian manifold associated to the canonical relation of  the operator ${\muutos \square_g}$,
that is, 
\beq\label{eq: Lambda g}
\Lambda_g=\{(x,\xi,y,-\eta);\ (x,\xi)\in\hbox{Char}\,(\square_g),\ (y,\eta)\in \Theta_{x,\xi}\},
\eeq
where $\Theta_{x,\xi}\subset T^*M$ is the bicharacteristic of $\square_g$ containing $(x,\xi)$.  

{\muutos By \cite[Thm.\ 26.1.14]{H4}, 
${{{\muutos \square_g^{-1}}}}:H_{comp}^s(M_0)\to H_{loc}^{s+1}(M_0)$ is a bounded.}
We will repeatedly use the fact (see \cite[Prop. 2.1]{GU1}) that 
if $F\in \I^{p}(\Lambda_0)$ {\muutos is compactly supported} and $\Lambda_0$ intersects Char$({\muutos \square_g})$ transversally
so that all bicharacterestics of ${\muutos \square_g}$ intersect $\Lambda_0$ only finitely many
times, then ${\muutos \square_g^{-1}F}\in \I^{p-3/2,-1/2}(\Lambda_0,\Lambda_1)$
where $\Lambda_1^{\prime}=\Lambda_g\circ \Lambda_0^{\prime}$, that is,
\beq\label{eq: flowout}
\Lambda_1=\{(x,-\xi);\ 
(x,\xi,y,-\eta)\in \Lambda_g,\ (y,\eta)\in \Lambda_0 \}.
\eeq
{\muutos The manifold $\Lambda_1$ is called} the flowout  
{\gtext from $\Lambda_0\cap\,$Char$({\muutos \square_g})$ 
{\vtext by the Hamiltonian vector field associated to $p(x,\xi)$.}


\subsubsection{Distorted plane waves satisfying a linear wave equation}
\label{subsubsec: Distorted plane wave equation}


 Next we consider a  
 distorted plane  wave whose singular support is
concentrated near a  geodesic. \motivation{These waves, sketched in Fig.\ 1(Right),
propagate near the geodesic $\gamma_{x_0,\zeta_0}([{\rtext 0},\infty))$
and are singular on a surface $K(x_0,\zeta_0,{\muutos s_0})$, defined below in (\ref{associated submanifold}). 
\refHOX{Comment 27: We have changed $\xi_0$  to $\zeta_0$.}
The surface $K(x_0,\zeta_0,{\muutos s_0})$ is a subset
of the light cone 
{\rtext $\L^+(x_0)$ and} the parameter $s_0$ gives a ``width'' of
{\mattitext the singular support of the wave around $\gamma_{x_0,{\reftext \zeta_0}}([{\rtext 0},\infty))$}. When $s_0\to 0$, its singular support
tends to the set $\gamma_{x_0,\zeta_0}([{\rtext 0},\infty))$.
Next we will define these waves.}
%

%
%
%

\refHOX{Comment 28: Earlier $t_0$  was used in two meaning. The parameter $t_0$ used in
the introduction is now changed to $T_0$.}

{\muutos Let  $x_0\in U$, $\zeta_0\in L^{+}_{x_0}M$ and $s_0>0$ and recall that
$g^+$ is a Riemannian metric on $M$.
Also, let
\ba
 \mathcal V_{x_0,\zeta_0,s_0}=\{\eta\in T_{\rtext x_0}M:
\|\eta -{\muutos \zeta_0}\|_{\nohat g^+}<s_0,
 \  \|\eta\|_{\nohat g^+}= \|{\muutos \zeta_0} \|_{\nohat g^+}\}
\ea 
be a  neighborhood of ${\rtext \zeta_0}$ on a sphere.
%
%

We define the subset of the light cone,
$K(x_0,\zeta_0,s_0)\subset {\muutos M_0}$  associated to  the vector $(x_0,\zeta_0)$
and  $x_0\in U$ and  parameter $s_0\in \R_+ $
by 
\beq\label{associated submanifold}
& &\hspace{-.4cm}K(x_0,\zeta_0,s_0)=\{\gamma_{{\rtext x_0},\eta}(t)\in  {\muutos M_0};\ \eta\in \W_{x_0,\zeta_0,s_0},\ t\in (0,\infty)\},\hspace{-1cm}
\eeq
where $\W_{x_0,\zeta_0,s_0}=L^{+}_{x_0}M\cap \mathcal V_{x_0,\zeta_0,s_0}$, {\reftext see Figure 1.}

%


Let 
\beq
\nonumber
& &\hspace{-9mm}\Sigma({x_0,\zeta_0,s_0})= \{(x_0,r \eta^\flat)\in T^*M;\ 
\eta\in {\slava \mathcal V_{x_0,\zeta_0,s_0}},\ r\in \R\setminus\{0\}\,\},
\\ \nonumber
& &\hspace{-9mm}\Lambda(x_0,\zeta_0,s_0)=
\{(\gamma_{{\rtext x_0},\eta}(t),r \dot \gamma_{{\rtext x_0},\eta}(t)^\flat )\in T^*M;\ \eta\in \W_{x_0,\zeta_0,s_0},\\ 
& &\hspace{53mm}\ t\in (0,\infty),\ r\in \R\setminus\{0\}\,\}.\hspace{-1cm}
\label{associated submanifold 2}
\eeq
Note that $\Lambda(x_0,\zeta_0,s_0)$}
is the Lagrangian manifold that is the  flowout from  
Char$(\square_{\nohat g})\cap  \Sigma({x_0,\zeta_0,s_0})$
{\vtext by the Hamiltonian vector field of associated to $p(x,\xi)$} in the future direction,
 see (\ref{eq: flowout}). 
 {\muutos Below, we will use sources $f\in \I^{n+1}(\Sigma({x_0,\zeta_0,s_0}))$. An example of such sources are functions $A\delta_{x_0}$  where $A$ is a psudodifferential operator   \refHOX{Comment 27: Definition $K(x_0,\zeta_0;t_0,s_0)$  is motivated with few words.}
 miclocally supported near $(x_0,\zeta_0)$. For example, in local coordinates we can use 
 $A=\phi_0(x)(1-\psi_0(D))\psi_1(D)$, where  $\phi_0\in C^\infty_0(\R^n)$ is supported near $x_0$, function $\psi_0\in C_0^\infty(\R^n)$  is equal to 1 in the neighborhood of zero and $\psi_1\in C^\infty(\R^n\setminus \{0\})$ is homogeneous function 
 supported in a conic neighborhood of direction $\zeta_0$. Function  $A\delta_{x_0}$ can be considered as
 a ``directed point source'' that produces a wave $\square_g^{-1}(A\delta_{x_0})$ which singularities propagate along $\Lambda(x_0,\zeta_0,s_0)$. {\slava Outside $x_0$, such wave could be considered as  a ``piece of distorted plane wave''}, see \reftext Fig.\ 1.}\refHOX{Comment 28. Reference to figure added.}
 {\mmtext Note that $\Sigma({x_0,\zeta_0,s_0})\subset  T^*M$ is a Lagrangian submanifold that is subset of
the conormal bundle {\mmmmtext $\Sigma_{x_0}$ of the point $\{x_0\}$, considered as a 0-dimensional submanifold of $M$,}
that is,
\beq\label{Normal bundle of a point}
\Sigma_{x_0}=N^*(\{x_0\})=T^*_{x_0}M\setminus \{0\},
\eeq
and hence, $f\in \I^{n+1}(\Sigma({x_0,\zeta_0,s_0}))$ is a conormal distribution.}
 

 \modified{When $K^{reg}\subset K= K(x_0,\zeta_0,s_0)$
is the set of points $x$ that have a neighborhood $W$  such that $K\cap W$ is   a smooth 3-dimensional submanifold,
we have $N^*K^{reg}\subset \Lambda(x_0,\zeta_0,{\muutos s_0}).$ {\muutos Note that
if $(x,\xi)\in N^*K^{reg}$  then also $(x,-\xi)\in N^*K^{reg}$, and this is the reason why we used factor $r\in \R\setminus\{0\}$  in formula (\ref{associated submanifold 2}).}

%
%

}

\begin{lemma}\label{lem: Lagrangian 1 wave} 
\refHOX{Comment 32: We have removed here that $n\leq -17$.}
Let {\reftext $n$} be an integer,  ${\muutos s_0}>0$, 
$K=K(x_0,\zeta_0,s_0)$, $\Lambda_1=\Lambda(x_0,\zeta_0,s_0)$ and {\muutos $\Sigma=
\Sigma({x_0,\zeta_0,s_0})$}.
{\mmmtext Let  $({\muutos x},\xi)\in {\slava {\muutos \Sigma} \cap L^{*}M}$, $v=\xi^\sharp\in L_xM$, $r\in \R$ and
$y=\gamma_{x,v}(r)$ and $\eta=(\dot \gamma_{x,v}(r))^\flat$ 
be such that  ${\muutos x}<{\muutos y}$.}

Assume that ${f_1}\in \I^{n+1}({\muutos \Sigma})$ 
  is a {\muutos compactly supported classical {\mmmtext conormal}}
distribution.

%
%


Let us consider the restriction of $w_1=\square_{\nohat g}^{-1}{f_1}$ {\vtext to} $ {\muutos M_0}\setminus {\muutos \{x_0\}}$.
Then $w_1 |_{ {\muutos M_0}\setminus {\muutos \{x_0\}}}\in \I^{n-1/2} ( {\muutos M_0}\setminus {\muutos \{x_0\}};\Lambda_1)$.
\refHOX{Comment 33: The notations on principal symbols are now improved.}

{\guntext 
Let  ${\rtext\sigma^{(p)}_{{f_1}}({\muutos x},\xi)}$ be  the
principal symbol of ${f_1}$  at $({\muutos x},\xi)$ and ${\rtext\sigma^{(p)}_{w_1}({\muutos y},\eta)}$
 be the
principal symbol of $w_1$ at $({\muutos y},\eta)\in  \Lambda_1$.
Then 
\beq\label{eq for principal symbols}
{\rtext\sigma^{(p)}_{w_1}({\muutos y},\eta)}=
R({\muutos y},\eta,{\muutos x},\xi){\rtext\sigma^{(p)}_{{f_1}}({\muutos x},\xi)}
\eeq
 where $R=R({\muutos y},\eta,{\muutos x},\xi)$ {\mmmtext is an invertible linear operator.}
 
Moreover, when the geodesic $\gamma_{x,v}([0,r])$ has no cut points, the point
 $y$ has a neighborhood $V_0\subset M$
such that $S_1=\L^+(x)\cap V_0$ is a smooth  submanifold of codimension 1 and 
 $w_1|_{V_1}\in \I^n(S_1)$ is a conormal distribution. Then, $R$ can be considered
 as a non-zero complex number.}
 
\end{lemma}

{\mmmtext Observe that in the claim of the lemma, $(({\muutos x},\xi),({\muutos y},\eta))\in \Lambda_{\nohat g}^\prime$ and 
$({\muutos y},\eta)\in T^*M$ be  
 on the same {\mmmmtext bicharacteristic} of $\square_{\nohat g}$
as $({\muutos {\muutos x}},\xi)$.} 

We call the solution $w_1$   a distorted plane  wave associated
to the submanifold $K(x_0,\zeta_0,{\muutos s_0})$. 

\smallskip

\noindent
{\bf Proof.} 
%
{\slava Recall that the {\muutos Schwartz kernel $Q$ of the causal inverse} 
operator $Q={{{\muutos \square_g^{-1}}}}$
satisfies 
$Q\in \I^{-3/2,-1/2}(\Delta^\prime_{T^*M},\Lambda_g)$.}
%
{\muutos As {\muutos $f\in \I^{n+1}({\muutos \Sigma})$}, \cite[Prop.\ 2.1]{GU1} 
and the definition (\ref{associated submanifold 2}) of $\Lambda_1$ 
imply} that
$w_1=\square_{\nohat g}^{-1}f\in \I^{n+1-3/2,-1/2}(\Sigma,\Lambda_1)$. This yields that
 $w_1|_{ {\muutos M_0}\setminus {\muutos \{x_0\}}}\in \I^{n+1-3/2}(\Lambda_1)$.
 This implies that the restriction $w_1|_V$ is a conormal distribution in $\I^n(S_1)$. Moreover, \cite[Prop.\ 2.1]{GU1} implies} {\muutos the formula 
 (\ref{eq for principal symbols})
 for the principal symbols,}
where
$R$ is obtained by solving an ordinary differential equation along a bicharacteristic curve. Making  similar considerations
for the adjoint of the $\square_{\nohat g}^{-1}$, i.e., considering
the propagation of singularities using reversed causality,
{\muutos and by solving an ordinary differential equation along a {\mmmmtext bicharacteristic}}, we see that
 $R$ is invertible.
 
 {\mmmtext Finally, when {\vtext the} geodesic $\gamma_{x,v}([0,r])$ has no cut points, $y$ has a neighborhood $V$ where the light cone
is a smooth hypersurface.}
  \hfill \Box \medskip

 \subsection{Microlocal analysis of the non-linear interaction of waves}\label{sec: interaction}

 \MATTITEXT{ 
 Next we consider the interaction of four $C^k$-smooth 
 waves  having conormal singularities on  {\mmtext hypersurfaces, where $k\in \Z_+$ is sufficiently large. 
 Interaction of such waves produces an artificial point source at the intersection point of the hypersurfaces. We will show that such
 artificial point sources can be created to arbitrary points of $I^+(p^-)\cap  I^-(p^+)$}. {\mtext We use such
 artificial point sources on the unknown manifold $(M,g)$ to create distorted spherical waves
 that determine the earliest light observations sets.}
 
 
 {{\guntext 
First considerations on the non-linear interaction of conormal waves,  
 were done by Bony \cite{Bony}, Melrose and Ritter \cite{MR1,MR2} and Rauch and Reed,  \cite{R-R}  for  {\mtext semilinear hyperbolic equations. In particular, they analyzed
 three conormal waves and showed that 
 the interaction of three  plane waves \refHOX{Comment 34. Citation to Bony
 has been added and higher dimension is added. Citation \cite{MsBZ} to Melrose-Sa Barreto-Zworski is added.} 
in three and higher dimensional spacetimes produces ``virtual sources'' that are singular on co-dimension 3 submanifolds.
 In the three dimensional spacetime such sources correspond to point sources.}
 {\mtext In \cite{MR1,MR2}, the microlocal properties of non-linear waves are analyzed also for arbitrary many interacting waves when
 the interaction of the waves and the propagation of  {\mmmmtext singularities  take place} on a union of finitely many submanifolds
 that form so-called characteristically complete variety of finite type (the geometrical restrictions caused by this assumption
 is discussed in detail in  \cite[Section 7]{MR1}).
 {Also, the appearance of the new wavefronts due to the interaction of non-linear
 terms at caustics or due to boundary and corner diffraction have been analyzed in \cite{JsB,MsBZ,Vasy1,Vasy2,Z}.
The microlocal properties and regularity of the solutions of non-linear hyperbolic equations, that correspond to the interaction of several conormal waves, are analyzed in the monograph by Beals \cite{Beals}.
 
 
 {\ggtext As discussed in the introduction, the focus of the above papers
on the interaction of conormal singularities for non-linear hyperbolic equations is different {\vtext from our paper 
as in those it is assumed}\HOX{CHECK THIS TEXT AGAIN!}
 that the geometrical setting of the interacting singularities is a priori known.} In inverse problems, when we study waves on an unknown manifold, we do not know the geometry of the surfaces 
 on which the waves are singular. 

In this section we consider  the interaction of waves in a subset of the spacetime where we are sure that the linearized waves
 have no caustics. However, caustics may appear in the interaction of waves and these waves may
 interact with the linearized waves. {\muutos Later, in Section \ref{subsection combining} we use} global Lorentzian geometry
 to obtain  {\newtext  a procedure that marches through the diamond set $J(p^-,p^+)$ by reconstructing it
 in small pieces. {\muutos This will allow us to avoid difficulties associated with {\vtext the} appearance of caustics in the linearized waves.}
 }} }}
   
\subsubsection{Forth order interaction of  waves for the non-linear wave equation}

Next, we introduce a vector of four $\e$ variables
denoted by $\vec \e=(\e_1,\e_2,\e_3,\e_4)\in \R^4$. 
{\revtext Let $s_0>0$. {\muutos For the non-linear wave equation (\ref{eq: wave-eq general}) we denote by $u_{\vec \e}$
 its solution when the source ${f}_{\vec \e}$ is} given by 
\beq\label{eq: f vec e sources}
 {f}_{\vec \e}:=\sum_{j=1}^4\e_j {f}_{j},\quad
 {f}_{j}\in \mathcal I^{n+1} ({\muutos \Sigma}(x_j,\zeta_j,{\muutos s_0})),
 \eeq
and $(x_j,\zeta_j)$ are light-like
vectors with $x_j\in U.$}
Moreover,  we assume that  
{\rtext the sources satisfy 
 \beq\label{eq: source causality condition}
& & \supp({f}_{j})\cap J ^+(\supp({f}_{k}))=\emptyset,\quad\hbox{for all }j\not=k,\\
\nonumber
%
& &\hspace{-1cm}J^+(W)\cap J^-(W)\subset U,
\quad \hbox{\vtext where }
W=\bigcup_{j=1}^4
 \supp({f}_{j}).\hspace{-1cm}
 \eeq
  The  implies that the supports of the sources are causally independent.}
}


The sources   $ {f}_{j}$ give raise to  the solutions
of the linearized wave equations, which we denote by
\beq\label{eq: waves uj} 
u_j:=\p_{\e_j}u_{\vec \e}|_{\vec \e=0}={{{\muutos \square_g^{-1}}}} \,( {f}_{j})\in \I(M_0\setminus\{{\slava x_j}\}; \Lambda(x_j,\zeta_j,{\muutos s_0})).
\eeq
In the following we use the notations
$\p_{\vec \e}^1u_{\vec \e}|_{\vec \e=0}:=\p_{\e_1}u_{\vec \e}|_{\vec \e=0}$,
$\p_{\vec \e}^2 u_{\vec \e}|_{\vec \e=0}:=\p_{\e_1}\p_{\e_2}u_{\vec \e}|_{\vec \e=0},$
$\p_{\vec \e}^3 u_{\vec \e}|_{{\vec \e}=0}:=\p_{\e_1}\p_{\e_2}\p_{\e_3} u_{\vec \e}|_{{\vec \e}=0}$,
and
\ba
\p_{\vec \e}^4 u_{\vec \e}|_{\vec \e=0}:=\p_{\e_1}\p_{\e_2}\p_{\e_3}\p_{\e_4} u_{\vec \e}|_{{\vec \e}=0}.
\ea

{\revtext 
Below, for the non-linear wave equation, we denote the wave produced by the
fourth order interaction of waves $u_{j}$ by 
\beq \nonumber
\M^{(4)}&=&\p_{\vec \e}^4u_{\vec \e}|_{\vec\e=0}={{{\muutos \square_g^{-1}}}}\mathcal S,
\quad
{\mmmmtext \mathcal S=
\sum_{\sigma\in \Sigma(4)}\mathcal S_\sigma,}
\\
\label{4th interaction for wave eq B}
\mathcal S_\sigma&=&-\bigg(a\,  
{{{\muutos \square_g^{-1}}}}(a\,u_{\sigma(4)}\,u_{\sigma(3)})\,{{{\muutos \square_g^{-1}}}}(a\,u_{\sigma(2)}\,u_{\sigma(1)})\\  
\nonumber
& &\quad+ 4a\, u_{\sigma(4)}\,{{{\muutos \square_g^{-1}}}}(a\,  
u_{\sigma(3)}\,{{{\muutos \square_g^{-1}}}}(a\,u_{\sigma(2)}\,u_{\sigma(1)}))\bigg),
\eeq
where $\Sigma(4)$ is the set
of permutations $\sigma$ of the set $\{1,2,3,4\}$, see (\ref{Slava formula 2}).

}

\subsubsection{On the singular support of the non-linear interaction  of three waves}
\label{subsection: sing supp interactions}

We will consider the case when we send distorted plane  waves 
propagating on surfaces $K_j=K(x_j,\xi_j,{\muutos s_0})$,  ${\muutos s_0}>0$, cf.\ (\ref{associated submanifold}), and
these waves interact. 

Next we consider the {\mattitext geometry related to the} {\mmmmtext three} {\muutos wave interaction}s of the waves.
{\mmmmtext Let {\mmmtext $\X({\rtext (\vec x,\vec\xi),{\muutos s_0}})\subset L^*M$ be the set of all light-like co-vectors $(x,\xi)$}}
that are in  the normal bundles   $N^*(K_{j_1}\cap K_{j_2}\cap K_{j_3})$
{\vtext with some $1\leq j_1<j_2<j_3\leq 4$, {\rbtext that is, $\X=\X({\rtext (\vec x,\vec\xi),{\muutos s_0}})$  is
\beq\label{set X}
& &\X((\vec x,\vec\xi),s_0)=\bigcup_{1\leq j_1<j_2<j_3\leq 4} \X_{j_1j_2j_3}((\vec x,\vec\xi),s_0),
\\
\nonumber
& &\hspace{-1cm}\X_{j_1j_2j_3}((\vec x,\vec\xi),s_0)=
\bigcup_{x\in K_{j_1}\cap K_{j_2}\cap K_{j_3}} (N^*_xK_{j_1}+
N^*_xK_{j_2}+N^*_xK_{j_3})\cap L^*M_0.\hspace{-1.5cm}
\eeq
}

Moreover, we define
 \beq\label{H-set} 
& &\H((\vec x,\vec\xi),s_0)=\bigcup_{1\leq j_1<j_2<j_3\leq 4} \X_{j_1j_2j_3}((\vec x,\vec\xi),s_0),
\\
\nonumber
& &\hspace{-1cm}\H_{j_1j_2j_3}((\vec x,\vec\xi),s_0)=
\{(y,\eta)\in T^*M_0;\ \hbox{there is $(x,\zeta)\in \X((\vec x,\vec\xi),s_0)$}\hspace{-1cm}\\
&& \nonumber
\hspace{4.5cm}\hbox{such that
 $x\leq y$ and $(y,\eta)\in \Theta_{x,\zeta}$}\}
 \eeq  
and  ${\mathcal Y}({\rtext (\vec x,\vec\xi),{\muutos s_0}})=\pi(\H((\vec x,\vec\xi),s_0))$,
where $\pi:T^*M_0\to M_0$  is the projection to the base space.}
%
 Finally,
 let  
 \beq\label{Y-set} 
{\muutos {\mathcal X}(\vec x,\vec\xi)}=\bigcap_{s_0>0}{\mathcal X}({\rtext (\vec x,\vec\xi),{\muutos s_0}}),\quad  
{\muutos {\mathcal Y}(\vec x,\vec\xi)}=\bigcap_{s_0>0}{\mathcal Y}({\rtext (\vec x,\vec\xi),{\muutos s_0}}).
 \eeq  
%

The three wave interaction happens then on $\pi ( \X ({\rtext (\vec x,\vec\xi),{\muutos s_0}}))$ and, roughly speaking, this interaction sends singularities
to  $\H ({\rtext (\vec x,\vec\xi),{\muutos s_0}})$. 
\motivation{For instance, in   Minkowski space, when three plane waves (whose singular supports 
are hyperplanes) collide,
the intersections of the hyperplanes is a 1-dimensional space-like line $K_{123}=K_1\cap K_2\cap K_3$ in the 4-dimensional space-time.
This corresponds to a {\muutos virtual point source}  moving continuously in time
{\mmmmtext and creates {\vtext a ``conical'' wave} that propagates near the surface $\Y ({\rtext (\vec x,\vec\xi),{\muutos s_0}})$.
 \refHOX{Comment 35: Details of shock waves have been added.}
{\vtext To visualize this,} see the supplementary video \cite{Video} and Figure 3 that display a conic waves produced by
the interaction of three waves on $K_{123}$. The video shows also the spherical wave that is produced by the interaction of all four
waves and  that emanates from the intersection point $q\in K_1\cap K_2\cap K_3\cap K_4$.}

In this paper we do not analyze carefully the singularities produced by the  three wave interaction 
near $ \Y ({\rtext (\vec x,\vec\xi),{\muutos s_0}})$. Our goal is to  
consider the singularities produced by the four wave interaction in the domain
${\muutos M_0}\setminus \Y ({\rtext (\vec x,\vec\xi),{\muutos s_0}})$.
{\revtext 
{\vtext We consider also the limit}
when $s_0\to 0$. 
{\itext Then,
the exceptional set $\Y ({\rtext (\vec x,\vec\xi),{\muutos s_0}})$ tends to a set 
$\Y (\vec x,\vec\xi)$
whose Hausdorff dimension is at most 2.
}


%
}
 \extension{

\begin{center}

\psfrag{1}{$x^0$}
\psfrag{2}{$x^1$}
\psfrag{3}{$x^2$}
\includegraphics[width=7.5cm]{planewaves-2.eps}
\end{center}
{\it FIGURE A1: A schematic figure where the space-time is represented as the  3-dimensional set $\R^{2+1}$. 
In the figure 3 pieces of plane waves  have singularities 
on strips of hyperplanes (in fact planes) $K_1,K_2,K_3$, colored
by light blue, red, and black. These planes  have intersections, and
in the figure the sets $K_{12}=K_1\cap K_2$, $K_{23}=K_2\cap K_3$,
and $K_{13}=K_1\cap K_3$ are shown as dashed lines with
dark blue, magenta, and brown colors. These dashed lines  intersect
at a point $\{q\}=K_{123}=K_1\cap K_2\cap K_3$.}

\begin{center}

\psfrag{1}{$\pi(\mathcal U_{z_0,\eta_0})$}
\psfrag{2}{$U$}
\psfrag{3}{$J({p^-},{p^+})$}
\psfrag{4}{$W$}
\psfrag{5}{$\mu_a$}
\includegraphics[width=7.5cm]{Xpoints-before-2.eps}
\end{center}
{\it FIGURE A2: In section \ref{subsection: sing supp interactions} we consider four colliding distorted plane  waves.
In the figure we consider the Minkowski space $\R^{1+3}$ and the figure corresponds 
to a ``time-slice'' $\{T_1\}\times \R^3$.
We assume that the distorted plane  waves are sent 
so far away that the waves look like pieces of plane waves. The plane waves $u_j\in \mathcal I(K_j)$, $j=1,2,3,4$,
are conormal distributions that are  solutions of the linear wave equation and their singular supports are the sets $K_j$,
that are pieces of 3-dimensional planes in the space-time. The sets $K_j$ are not shown in the figure.
The 2-wave interaction wave $\M ^{(2)}=\p^2_{\vec\e}u_{\vec\e}|_{\vec\e=0}$ is singular on the set $\cup_{j\not =k}K_j\cap K_j$.
There are 6 intersection sets $K_j\cap K_j$ that are shown as black line segments. Note that these
lines have 4 intersection points, that is, the vertical and the horizontal black
lines do not intersect  in $\{T_1\}\times \R^3$.
The four intersection points of the black lines are the sets $(\{T_1\}\times \R^3)\cap (K_j\cap K_j\cap K_n)$.
These points correspond to  points moving in time (i.e., they are curves in the space-time) that produce singularities of the  {\mmmmtext three} {\muutos wave interaction} wave 
$\M ^{(3)}$. The points
 seem to move faster than the speed of the light (similarly, as  a shadow of a
far away  object may seem to  move faster than the speed of the light).  Such
point sources produce ``shock waves'', and due to this, $\M ^{(3)}$ is
singular on the sets $\mathcal Y((\vec x,\vec \xi),{\muutos s_0})$ defined in 
formulas (65)-(66). The set $(\{T_1\}\times \R^3)\cap \mathcal Y((\vec x,\vec \xi),{\muutos s_0})$
is the union of the four blue cones shown in the figure. 
}

\begin{center}

\psfrag{1}{$y_{-1}$}
\includegraphics[width=7.5cm]{frustum2.eps}
\end{center}
{\it FIGURE A3: A schematic figure in  3-dimensional Euclidean space  $\R^{3}$
where
we describe the geometry of the wave produced by the interaction of three
waves in the Minkowski space. The figure shows such a  {\mmmmtext three} {\muutos wave interaction} wave in a time slice  $\{T_1\}\times \R^{3}$
with two different values of the parameter $s_0$.  On the left, $s_0$ is large
and the  {\mmmmtext three} {\muutos wave interaction} wave has singular support on a cone. 
 On the right, $s_0$ is small and the   {\mmmmtext three} {\muutos wave interaction} wave is a frustum (a truncated cone, or,  a cone with its apex cut off by a plane), similar to set 
$\{(x^1,x^2,x^3)\in \R^3;\ (x^1)^2+(x^2)^2=c(x^3)^2$, $a<x^3<a+h\}$, with some
$a,c>0$ and a thickness $h=c_1s_0$.
}

\begin{center}

\psfrag{1}{$\pi(\mathcal U_{z_0,\eta_0})$}
\psfrag{2}{$U$}
\psfrag{3}{$J_{\nohat g}({p^-},{p^+})$}
\psfrag{4}{$W_{\nohat g}$}
\psfrag{5}{$\mu_{a}$}
\includegraphics[width=7.5cm]{Xpoints-after-3.eps}
\end{center}

\noindent
{\it FIGURE A4: The same situation that was described in Fig. A2 
is shown
at a later time, that is, the figure shows the time-slice $\{T_2\}\times \R^3$ with
$T_2>T_1$, when the parameter $s_0$ is quite large. The four  distorted plane waves have now collided and
produced a point source in the space-time at a point $q\in K_1\cap K_2\cap K_3\cap K_4$.
This gives  raise to the singularities of the  {\mmmmtext four} {\muutos wave interaction} wave $\M ^{(4)}$.
In the figure $T_1<t<T_2$, where $q$ has the time coordinate $t$.
The four cones in the figure, shown with solid blue and green curves
and dashed blue and yellow curves are   the intersection of  the time-slice $\{T_2\}\times \R^3$ 
and the set $\mathcal Y((\vec x,\vec \xi),{\muutos s_0})$. Inside the cones the red sphere
is the set $\mathcal L^+(q)\cap (\{T_2\}\times \R^3)$ that
corresponds to the spherical plane  wave produced by the point source at $q$.
}

\begin{center}

\psfrag{1}{$\pi(\mathcal U_{z_0,\eta_0})$}
\psfrag{2}{$U$}
\psfrag{3}{$J_{\nohat g}({p^-},{p^+})$}
\psfrag{4}{$W_{\nohat g}$}
\psfrag{5}{$\mu_{a}$}
\includegraphics[width=5.5cm]{Xpoints-aft-Alt-2.eps}
\end{center}
{\it FIGURE A5: The same situation that was described in Fig.\ A4. 
That is,  the figure shows in 
the time-slice $\{T_2\}\times \R^3$ the singularities produced by four colliding
distorted plane  waves, when the parameter $s_0$ is very small. In this case the truncated cones
degenerate to circles.
Indeed,
 the figure corresponds to the case when the parameter $s_0$ is close
to zero, that is, the pieces of the distorted plane  waves are concentrated near a single
geodesic in the space-time.
As  $s_0$ is small, the distorted plane  waves interact only during a short time.
Due to this, the sets $(K_j\cap K_n)\cap (\{T_2\}\times \R^3)$,
that were black lines in the previous figures, are now empty, and do not
appear at the figure.
In the figure the  red sphere is the set $\mathcal L^+(q)\cap (\{T_2\}\times \R^3)$ corresponding
to  the spherical  wave
 produced by a point source at $q$.
The solid blue and green circles on the sphere
and the dashed blue and yellow  circles on the sphere are
 the intersection of  the time-slice $\{T_2\}\times \R^3$ 
and the set $\mathcal Y((\vec x,\vec \xi),{\muutos s_0})$. These circles
are actually not circles but surfaces of frustums (a truncated cone, or,  a cone with its apex cut off by a plane), similar to set 
$\{(x^1,x^2,x^3)\in \R^3;\ (x^1)^2+(x^2)^2=c(x^3)^2$, $a<x^3<a+h\}$, with some
$a,c>0$ and a very small thickness $h>0$. Note that the 2-Hausdorff measure of the frustums
go to zero as $s_0\to 0$ and the 3-Hausdorff measure of the set
$\mathcal Y((\vec x,\vec \xi),{\muutos s_0})$ goes to zero as $s_0\to 0$. Thus, when all the distorted plane  waves
intersect at a point $q$, 
the singularities of $\M ^{(4)}$  are observed on a set whose 3-Hausdorff measure in
the space-time is independent of $s_0$ but the singularities of $\M ^{(3)}$ 
are observed on a set which 3-hHausdorff measure in
the space-time is smaller than $Cs_0$. 
}

}

\extension{
The source function $F_\tau$ can
be constructed in local coordinates using the asymptotic representation of the 
gaussian beam, namely, by considering first 
\ba
F^{(series)}_\tau(x)=c\tau^{-1/2}\int_\R e^{-\tau s^2}\phi(s)f_\tau(x;s)ds,
\ea
where
$f_\tau(x;s)=\square_{\nohat g}(H(s-x^0)u_\tau(x))$. Here, $\phi\in C^\infty_0(\R)$ 
has value 1 in some neighborhood of zero, and $H(s)$ is the Heaviside function,
and  we have 
\ba
e^{-i\tau p(x)}F^{(series)}_\tau(x)=e^{-i\tau p(x)}F^{(1)}(x)\tau^{-1}+e^{-i\tau p(x)}F^{(2)}(x)\tau^{-2}+O(\tau^{-3}).
\ea Then,
$F_\tau$  can be defined by $F_\tau=F^{(1)}(x)\tau^{-1}.$ 
Indeed
we see that
\ba
\square_{\nohat g}^{-1}F^{(series)}_\tau&=&
\tau^{-1}\Phi_\tau(x)\,u_\tau(x),\quad \hbox{with}\\
\Phi_\tau(x)&=&
c\tau^{1/2}\int_\R e^{-\tau s^2}\phi(s)H(s-x^0) ds\\
&=&
\left\{\begin{array}{cl}
O(\tau^{-N}),&\hbox{for }x_0<0,\\
1+O(\tau^{-N}),&\hbox{for }x_0>0.
\end{array}\right.
\ea
Moreover, by writing $\square_{\nohat g}=\nohat g^{jk}\p_j\p_k+\hat \Gamma^j\p_j$,
$\hat \Gamma^j=\hat  g^{pq}\hat \Gamma^j_{pq}$ we have  that
\ba
f_\tau(x;s)&=&
\nohat g^{00}u_\tau(x)\delta^\prime(s-x^0)\\
& &-2\sum_{j=1}^3\nohat g^{j0}(\p_ju_\tau(x))\,\delta(s-x^0)
-\hat \Gamma^0u_\tau(x)\,\delta(x^0-s).
\ea
Let us use  normal coordinates centered at the point $y$ so that
$\Gamma^j(0)=0$ and $\nohat g^{jk}(0)$ is the Minkowski metric.
Then we define $p(x)=(x^0)^2+\varphi(x)$. The leading
order term of  $F^{(series)}_\tau$ is given by
\ba
F_\tau= c\tau^{-1/2}e^{-\tau p(x)}\left(
2x^0\phi(x^0)U_0(x)+2i(\p_0\varphi(x))U_0(x)\phi(x^0)
\right),
\ea
where $\p_0\varphi|_y$ does not vanish. }}
}

\subsubsection{Wave front set of the wave produced by the interaction of four waves}

{\vtext
Next we will consider $\WF(\M^{(4)})$}
{\muutos where $\M^{(4)}$ is the} wave produced by 
the interaction} of the four linearized waves corresponding to the  sources
$ {f}_{j}\in \mathcal I^{n+1} (\Sigma({x_0,\zeta_0,s_0}))$, $j\leq 4$, {\vtext see (\ref{4th interaction for wave eq B}).}
%

\begin{definition}\label{def:  4-intersection of rays} 
We say that the geodesics corresponding to the vectors $(\vec x,\vec \xi)=((x_j,\xi_j))_{j=1}^4$
intersect  and  the intersection takes place at the point $q\in \hattuM _0$ 
if there are {\muutos $t_j>0$ such that} 
 $q=\gamma_{x_j,\xi_j}({\rtext t_j})$  
 \refHOX{Comment 41: The notation $t_0$ in the introduction is changed to $T_0$.}
 for all $j=1,2,3,4$. {\vtext We say that the intersection of geodesics is regular if
 $t_j\in (0,{\bf t}_j)$, where ${\bf t}_j=\rho({\rtext x_j,\xi_j})$ and vectors
 $\dot \gamma_{x_j,\xi_j}(t_j)\in T_qM_0$, $j=1,2,3,4$ are linearly independent.}

 \end{definition}
 
  {\reftext For $q\in {\rbtext M_0}$,} let $\Lambda_q^+$ be the Lagrangian manifold
  \refHOX{Comment 42: We added $q\in {\rbtext M}$. Also, in definition $\Lambda_q^+$ there is a change of notations.} 
  \ba
  \Lambda_q^+=\{({\muutos y},{\muutos \eta})\in T^*{\rbtext M_0}&;&\ {\muutos y}=\gamma_{q,\zeta}({\reftext 1}),\ 
  {\muutos \eta^\sharp}=r\dot\gamma_{q,\zeta}({\reftext 1}),\\\
  & & \zeta\in L^+_q{\rbtext M_0},\ {\reftext r\in \R\setminus 0}\}.
  \ea
  Note that the {\slava projection $\pi( \Lambda_q^+)$ of $\Lambda_q^+$
 on $\hattuM_0$} is  the light cone $\L^+(q)$. 

Next we consider $x_j\in U$ and $ \xi_j\in L^+_{x_j}\hattuM _0$,
  such that
  $(\vec x,\vec \xi)=((x_j,\xi_j))_{j=1}^4$
  satisfy, see Fig.\ 5(Right),
\beq\label{eq: summary of assumptions 1}
{\muutos x_j\in U\quad\hbox{and}\quad 
 x_j\not \in J^+(x_k)
  \hbox{ for  $j\not =k$.}}  \hspace{-5mm}
 \eeq
  \refHOX{Comment 43: The notation $t_0$ in the introduction is changed to $T_0$.}
We denote  
 \beq\label{eq: summary of assumptions 2} 
& & {\rtext \V(\vec x,\vec \xi)}= {\muutos M_0}\setminus \bigcup_{j=1}^4
J^+(\gamma_{{\rtext x_j,\xi_j}}({\bf t}_j)),\quad \hbox{where ${\bf t}_j:={\rtext \rho(x_j,\xi_j)}.$}
 \eeq
}
  Note that
 two geodesics $\gamma_{{\rtext x_j,\xi_j}}([0,\infty))$ can intersect at most once in ${\rtext \V(\vec x,\vec \xi)}$.
{\muutos Below, let 
\beq\label{Sets Kj}
K_j={K(x_{j},\xi_{j},{\muutos s_0})},\quad \Lambda_j={\Lambda(x_{j},\xi_{j},{\muutos s_0})},
\eeq 
{\slava cf. (\ref{associated submanifold}),  (\ref{associated submanifold 2})}, where $s_0>0$ is so small that {\mtext one of}  the following two cases are satisfied: 
\ba 
&\hbox{(A)}&
 (\cap_{j=1}^4 K_j)\cap  \V(\vec x,\vec \xi)=\emptyset,\\
&\hbox{or}& \\
&\hbox{(B)}& \hbox{ $(\cap_{j=1}^4 K_j)\cap  \V(\vec x,\vec \xi)=\{q\},$ where}\\
& &\nonumber q=\gamma_{x_j,\xi_j}(t_j)\hbox{ with $t_j>0$ for all  $j=1,2,3,4$,}
\ea
and the intersection of any  $K_i$ and $K_j$ with $i\not=j$ is  transversal. 
In the case (B) 
all geodesics $\gamma_{x_j,\xi_j}$, $j=1,2,3,4$ intersect at a point $q$ in $\V(\vec x,\vec \xi)$
and in the case (A) the geodesics do not intersect.}

Below,  $\Sigma(4)$ is the set of  permutations $\sigma:(1, 2, 3, 4)\to (1, 2, 3, 4).$

\motivation{
 
 {\reftext Observe that in the set ${\rtext \V(\vec x,\vec \xi)}$ the geodesics $\gamma_{{\rtext x_j,\xi_j}}([0,\infty))$ do not have conjugate points and thus the waves $u_j$  do not have caustics in this set.}
 \refHOX{Comment 44: The text "...before the waves $u_j$ have caustics" is  clarified.}


}

{\vtext

In the next theorem we consider the singularities of the wave $\U^{(4)}$, produced by the 
interaction of four waves $u_j$, outside a ``small'' set ${\mathcal Y}(\vec x,\vec \xi)\cup\bigcup_{j=1}^4\gamma_{x_j,\xi_j}([0,\infty))$. Essentially, we show that no such singularities can {\ftext be} detected outside the causal future $J^+(q)$ of the point
$q$  where all four plane waves interact. Also, we show that in
the case when the directions of the distorted planes at $q$ are linearly independent and we are in a generic case,  the singularities are observed in the set  $\firstpoint^{reg}_U(q)$  that is the  regular
part of the boundary $\partial J^+(q)$, see Def.\ \ref{def: earliest element}.

\medskip

\noindent {\bf Remark 3.2.} The non-linear interaction of waves may cause extraordinary singularities. For
example, M. Beals  showed in 1983 for the wave equation $\square u(t,y)+b(t,y)u(t,y)^3=0$ in $\R^4$ that there are solutions for which the singular support of the Cauchy data $(u|_{t=0},\p_tu|_{t=0})$  is the point $\{0\}$,
but the singular support of $u$  contains the entire solid cone $\{(t,y)\in \R^4;\ |y|<t\}$,
see  \cite[Thm.\ 2.10]{Beals} and  \cite{Beals-spread}. This example has similarities
 to  the above case $(B)$ when the direction vectors $\dot \gamma_{x_j,\xi_j}(t_j)$,  $j=1,2,3,4$ are not linearly independent. This happens e.g.\ when both the plane wave $u_4$  and the conic wave $w_{321}=
\square_g^{-1}(au_3\, \square_g^{-1}(au_2u_1))$, produced by the interaction
of three waves $u_1,u_2$ and $u_3$ (see Fig.\ 3), propagate along the same geodesic $\gamma_{x_4,\xi_4}\subset K_4$.
In this case it may be that the wave front set of $u_4$  contains a point $(x,\zeta)\in N^*K_4$ and 
the wave front set of $w_{321}$  contains the point $(x,-\zeta)\in N^*K_4$ with the opposite direction. In this case
it is difficult to analyze the product $u_4w_{321}$. This difficulty, as well as the possible caustics 
of $w_{321}$ (see Fig.\ 6),  are the reasons why in the claim (ii) below we restrict ourselves to a geometrically nice case.}

}

  \begin{theorem}\label{lem:analytic limits A wave}    Let $(\vec x,\vec \xi)=((x_j,\xi_j))_{j=1}^4$  be future pointing light-like vectors  
  such that (\ref{eq: summary of assumptions 1}) is satisfied. 
    Let ${\muutos y_0}\in {\rtext \V(\vec x,\vec \xi)}\cap U$ be
such that
 ${\muutos y_0}\not\in {\mathcal Y}(\vec x,\vec \xi)\cup\bigcup_{j=1}^4\gamma_{x_j,\xi_j}([0,\infty))$,
see (\ref{Y-set}) and (\ref{eq: summary of assumptions 2}).

  Assume that $s_0>0,K_j,\Lambda_j$ 
are as in (\ref{Sets Kj}) and assume that  either the above condition (A) or (B) is satisfied. In the case (B),
  we consider $t_j>0$ and $q\in M$ and co-vectors
 $b_j=(\dot\gamma_{x_j,\xi_j}(t_j))^\flat$. 
%

Let  $n\in \Z_+$ {\mtext and} 
${f}_{j}\in \mathcal I^{-n+1}({\muutos \Sigma}(x_j,\zeta_j,{\muutos s_0}))$, $j=1,2,3,4$, 
be sources satisfying (\ref{eq: source causality condition}) and
$u_j={{{\muutos \square_g^{-1}}}}f_j$
%
%
%
%
and $\M^{(4)}$ be 
{\slava the wave {\mmmmtext produced by the  4th order} interaction}
given in (\ref{4th interaction for wave eq B}). When {\mtext $n$ is large enough and} $s_0$  is small enough, the following holds:

\medskip


 \noindent
 (i) Assume that either (A) {\slava or (B) holds and that in the case (B) we have 
 {\vtext $ y_0\not \in  J^+(q)$.}}
 Then 
  ${\muutos y_0}$ has a neighborhood $W$ such that 
 $\M^{(4)}|_W$ is $C^\infty$-smooth.
\medskip

\noindent
{\mmmmtext (ii) Assume that $(B)$ holds, 
{\muutos  $b_j\in T^*_qM,$ $j=1,2,3,4$ are linearly independent} and 
{\vtext $y_0\in \firstpoint^{reg}_U(q)$, where  
  $\firstpoint^{reg}_U(q)$ is the  regular earliest light observation set} of $q$, see Def.\ \ref{def: earliest element}.}
{\vtext Also, assume that 
$ \wflat\in L^{*}_{y_0}M$ and $r\in \R$ are such that
 $\gamma_{y_0,\wflat}(r)=q$ and denote ${\mmmtext \hatzeta}=(\dot \gamma_{y_0,\wflat^\sharp}(r))^\flat\in L^*_qM$.}


%
%
%

Then 
the point ${\muutos y_0}$ has a neighborhood $W$ such that 
$\M^{(4)}$ in $W$ is a conormal distribution
associated to $S=\mathcal L^+(q)\cap W$, that is, $\M^{(4)}|_W\in  \I^m(S)$, 
with $m=-4n-4$. 
Moreover, let 
 $\zeta_j\in N^*_qK_j$  be such
that
 \beq\label{zeta 0}
 {\muutos {\mmmtext \hatzeta}}=\sum_{j=1}^4 \zeta_j. 
 \eeq %
 

{\guntext  {\slava Note that the linear independence of $b_j$ implies the uniqueness of representation (\ref{zeta 0}).}
Then {the principal symbol of $\U^{(4)}|_W\in \I^m(S)$,
  at the point  $({\muutos y_0},{\muutos \wflat})$, is}}  
 \beq\label{eq: principal symbols of U4}
 \sigma^{(p)}_{\U^{(4)}}({\muutos y_0},{\muutos \wflat})=
 R({\muutos y_0},{\muutos \wflat},q,{\muutos {\mmmtext \hatzeta}})
 a(q)^3\mathcal G _{\bf g}(\vec \zeta)\prod_{j=1}^4 \sigma^{(p)}_{u_j}(q,\zeta_j),
 \eeq 
 where 
 ${\vec \zeta}=(\zeta_j)_{j=1}^4$,
 $R({\muutos y_0},{\muutos \wflat},q,{\muutos {\mmmtext \hatzeta}} )$ is given in 
 Lemma \ref{lem: Lagrangian 1 wave} and
\beq\label{P per Qmod}
{\mathcal G} _{\bf g}(\vec \zeta)&=&\sum_{\sigma\in \Sigma(4)}\bigg(
 \frac{C_1}{G(\zeta_{\sigma(1)}+\zeta_{\sigma(2)})\,\cdotp G(\zeta_{\sigma(1)}+\zeta_{\sigma(2)}+\zeta_{\sigma(3)}) }\\
 &&\hspace{2cm}+
\frac{C_{2}}{G(\zeta_{\sigma(1)}+\zeta_{\sigma(2)})G(\zeta_{\sigma(3)}+\zeta_{\sigma(4)})}\bigg),\nonumber
\eeq
where 
$G(\xi)=g(\xi,\xi)$ and $C_1$ and $C_2$
are non-zero constants.


%

\end{theorem}

{Later, we will show that the {\muutos function $\mathcal  G_{\bf g}(\vec \zeta)$} is non-vanishing in a generic set.}
\noextension{

\medskip

\noindent
{\bf Proof.}  Due to the general geometric setting on a globally hyperbolic manifold the proof is quite long and {\vtext is}  divided to several parts. 

{\bf 1. Notations.}
\refHOX{Comment 45: "To simplify notations," is removed.}
{\muutos As  ${\muutos y_0}\not\in {\mathcal Y}(\vec x,\vec \xi)\cup\bigcup_{j=1}^4\gamma_{x_j,\xi_j}([0,\infty))$, we can assume that $s_0>0$ is {\slava so small that}
 ${\muutos y_0}\not\in {\mathcal Y}((\vec x,\vec \xi),{\muutos s_0})\cup\bigcup_{j=1}^4{K(x_{j},\xi_{j},{\muutos s_0})}$,
see (\ref{Y-set}) and (\ref{eq: summary of assumptions 2}), and that for all $i\not=j$, the surfaces $K(x_{i},\xi_{i},{\muutos s_0})$  and $K(x_{j},\xi_{j},{\muutos s_0})$ intersect transversally.}

Below, we denote $\bar N^*_yK_j=N^*_yK_j\cup \{0\}$ and $\bar L^*_yM=L^*_yM\cup \{0\}$. %
Also, {\slava let {\vtext $\V=\V(\vec x,\vex\xi)$,} $\X= \X ({   (\vec x,\vec\xi),{\muutos s_0}})$,
 {\vtext $\H= \H ({   (\vec x,\vec\xi),{\muutos s_0}})$} 
and $\Y= \Y ({   (\vec x,\vec\xi),{\muutos s_0}})$ be the sets 
given in (\ref{set X})-(\ref{Y-set}).}

We use {\vtext the} notations $K_j=K(x_{j},\xi_{j},{\muutos s_0})$,
$K_{12}=K_1\cap K_2$, $K_{123}=K_1\cap K_2\cap K_3$, 
$\Lambda_{12}=N^* K_{12}$,
etc.

{Recall that $g^+$  is the Riemannian metric obtained by changing, in {\vtext local coordinates}, the sign of the negative eigenvalue
of the Lorentzian metric $g$. On {\mmmmtext the} $\nohat g^+$-unit sphere
 bundle $S^*M$  we use the 
 Sasaki metric determined by $\nohat g^+$. 
Below we say that a conic set $\V(\e)\subset T^*M$ is a conic $\e$-neighborhood of the set $L^*M$ of the light-like co-vectors
if $\V(\e)\cap S^*M$ is the  
 $\e$-neighborhood of  $L^*M \cap S^*M$ in
 the $\nohat g^+$-unit sphere
 bundle $S^*M$. } {\muutos Note that $(x,\xi)\in \V(\e)$  if and only if $(x,-\xi)\in \V(\e)$.}

 {\mmmtext The Lorentzian volume on $(M,g)$ at point $x$ is  denoted by $dV_x$.}

 {\vtext Below, we will consider claims (i) and (ii) at the same time.
 To do that, we denote
\ba
& & \hspace{-.8cm}\V_0=\V\setminus J^+(q), \hbox{ if (B) holds and $(b_j)_{j=1}^4$ are  linearly dependent}\\
& & \hspace{-.8cm}\V_0=\V, \hbox{ if (B) holds and $(b_j)_{j=1}^4$ are linearly independent or (A) holds.}
 \ea
 We will assume below that $y_0\in U\cap \V_0$.}
 

{\muutos {\bf 2. Local coordinates.} 
Recall that the intersection of {\vtext the} surfaces $K_i$ and $K_j$ with $i\not=j$ is  transversal in $\V(\vec x,\vec \xi)$. 
To consider local coordinates, let us start with the observation that 
if three  
light-like vectors are not parallel, then  those vectors are linearly independent, see \cite[Cor.\ 1.1.5]{SW}. 
{\mmtext This implies that for any $p\in \V(\vec x,\vec \xi)$ 
and any three indexes $j_1,j_2,j_3\in \{1,2,3,4\}$
we can
choose local coordinates $X:W\to \R^4$ so that $K_{j_i}\cap W\subset \{x\in W;\ X^{j_i}(x)=0\}$ for $i=1,2,3$
and we see that $K_{j_1}\cap K_{j_2}\cap K_{j_3}$  is a smooth path in the neighborhood $W$. 
In this case we say that $X:W\to \R^4$ are adapted to {\vtext the} surface $K_{j_i}$.
Also, in the case of claim (ii), at the point $q$ we can use local coordinates $X:W\to \R^4$ such that the linearly independent co-vectors $b_j$, $j=1,2,3,4$ are the differentials of the coordinate functions at $q$ and for all $j=1,2,3,4$ we have
$K_j\cap W_0=\{x\in W_0;\ X^j(x)=0\}$. These coordinates are adapted to all $K_j$.

{\vtext As in these set $\V$  the point  $q$ is the only possible  point in $\cap_{j=1}^4K_j$,
the existence of the above coordinates imply that when  $(i_1,i_2,i_3,i_4)$ is any permutation of $\{1,2,3,4\}$,
then
in the set $\V_0$, all possible intersections $K_{i_1i_2i_3}\cap K_{i_4}$ and 
 $K_{i_1i_2}\cap K_{i_3i_4}$ and $K_{i_1i_2}\cap K_{i_3}$ are  transversal.}

{\bf 3. Testing when $({\muutos y_0,\wflat})$ is in the wave front set.}
{\mmmmtext Below, {\vtext we consider $\wflat\in L_{y_0}^{*}M$}. 
As $Q$ is the causal inverse of the wave operator, we see using
 H\"ormander's theorem {\muutos on} propagation {\vtext of} singularities along 
bicharacteristics, \cite[Theorem 26.1.4]{H4}, we see that if the point $({\muutos y_0},{\muutos \wflat})$ is in WF$(\mathcal U^{(4)})$
then
either $\wflat$  is not light-like and $(y_0,\wflat)\in \hbox{WF}(\mathcal S)$,
or,  $\wflat$  is light-like and
there is $\tsparameter\in \R$ such that
 $(\gamma_{{\muutos y_0},{\muutos \wflat^\sharp}}(\tsparameter),\dot \gamma_{{\muutos y_0},{\muutos \wflat^\sharp}}(\tsparameter)^\flat)$
 is in WF$(\mathcal S)$ and $\gamma_{{\muutos y_0},{\muutos \wflat^\sharp}}(\tsparameter)\leq y_0$.
 To apply this for the light-like singularities, we below consider a point $({\muutos x_0},\zeta_0)\in L^{*}M$ such that}}
%
%
%
%
\beq\label{y assumptions1}
({\muutos x_0},\zeta_0)=(\gamma_{{\muutos y_0},{\muutos \wflat^\sharp}}(\tsparameter),\dot \gamma_{{\muutos y_0},{\muutos \wflat^\sharp}}(\tsparameter)^\flat),
\hbox{  {\mmmmtext such that $x_0<y_0$},}\hspace{-1.4cm}
\eeq
and study {\slava if $({\muutos x_0},\zeta_0)$ belongs in the wave front $\hbox{WF}(\mathcal S)$.}
{\vtext Note that as  $y_0\in U\cap \V_0$, we have also  $x_0\in \V_0$.}

{\mmmtext 
To study the claim (ii), we see that {\vtext when $s=r$}, the point $(x_0,\zeta_0)$  coincides with $(q,{\mmmtext \hatzeta})$.}
Also, to study the claim (i), we will study several cases
when  $(x_0,\zeta_0)$ will not be in $\hbox{WF}(\mathcal S)$ and use this to show that
$({\muutos y_0},{\muutos \wflat})$ does not belong in WF$(\mathcal U^{(4)})$.

}


{\bf 4. A neighborhood of light-like directions.}
We start with {\mmmmtext some auxiliary observations.} 
First, 
note that  as ${\muutos y_0}\not \in {\mathcal Y}$, {\mmmtext formula (\ref{y assumptions1}) and definitions 
(\ref{set X}) and (\ref{Y-set}) imply that}
$({\muutos x_0},{\muutos  \zeta_0})\not \in \X$. 

{\mmmtext 
Next we will choose a small parameter $\e_1>0$ that determines
a conic neighborhood $ \V(\e_1)$ of 
 the set $L^*M$ of light-like co-vectors.
We consider separately two cases:

%
%
%
%
%
%
%

{\mmmmtext First, consider the case when 
\beq\label{case a1}
& &\hbox{the property $(B)$ is valid, so that geodesics $\gamma_{x_j,\xi_j}$ intersect at $q$,}\\ 
\nonumber
& &\hbox{$b_j\in T^*_qM,$ $j=1,2,3,4$ are linearly independent, there is} \\ \nonumber
& &\hbox{{\vtext $r\not =0$ such that $\gamma_{y_0,\wflat^\sharp}(r)=q$, and {\vtext $x_0=q$}.}}
%
\eeq

{\vtext Then, {\mmmmtext denote} $v_0=
\dot \gamma_{y_0,{\mmmtext \wflat^\sharp}}(r) \in L_{q}M$ and
$\hatzeta=v_0^\flat \in L^{*}_{q}M.$}
Let  $\zeta_j\in \bar N^*_qK_j$, $j=1,2,3,4$ be such 
that ${\mmmtext \hatzeta}=\sum_{j=1}^4 \zeta_j$.
As ${\muutos y_0}\not \in {\mathcal Y}\cup (\cup_{j=1}^4 K_j)$,
we have $(q,{\mmmtext \hatzeta})\not \in\X\cup (\cup_{j=1}^4 N^*K_j)$. This implies that 
 $\zeta_j\not =0$ for all $j=1,2,3,4$. 
Then, as ${\mmmtext \hatzeta}$ and $\zeta_j$  are light-like
we have that  ${\mmmtext \hatzeta}-\zeta_j$ is not light-like
as otherwise ${\mmmtext \hatzeta}$ and $\zeta_j\in N^*_qK_j$ would be parallel {\ftext which} is not possible.
Hence, in the case (\ref{case a1}) 
we can choose $\e_1>0$  be so small that we have
 \beq\label{eq: conditions 2 modified}
(q,{\mmmtext \hatzeta}-\zeta_j)\not \in  \V(\e_1),\quad\hbox{\vtext for $j=1,2,3,4$.}
\eeq

{\vtext Second, in the case when condition  (\ref{case a1})  does not hold, we choose $\e_1>0$ to be an arbitrary positive number.}  
}

{\bf 5. Decomposition of {\ftext the} operator $\square_{\nohat g}^{-1}$.}
Below we denote $Q=\square_{\nohat g}^{-1}$.
{\reftext  We denote also the Schwartz kernel of $Q$ by
$Q(x,y)$.} \refHOX{Comment 46: Notations on Schwartz kernels have be clarified.}
Let us next consider the map ${{} {\bQ}}:C^\infty_0({\muutos M_0})\to
 C^\infty({\muutos M_0}).$  By \cite{MU1}, {\muutos the Schwartz kernel satisfies} ${{} {\bQ}}\in I({\muutos M_0}\times {\muutos M_0};
   \Delta_{T^*\hattuM _0}^{\prime}, \Lambda_{\nohat g})$, see Sec.\ \ref{Lagrangian distributions},
 and the canonical
   relation of the operator $Q$, denoted} 
$\Lambda_{{{} {\bQ}}}^{\prime}$, {\muutos has the form} 
$\Lambda_{{{} {\bQ}}}^{\prime}= \Lambda_{\nohat g}^{\prime}\cup \Delta_{T^*\hattuM _0}$, see \cite{MU1}.
{\muutos Let $\e_1$ be as above, $\e_2\in (0,\e_1)$} and $B_{\e_1,\e_2}$ be a pseudodifferential
operator on $\hattuM _0$ which is {\muutos microlocally a smoothing
operator}  {\reftext outside} 
\refHOX{Comment 47: "outside in" changed to "outside".}
 the conic $\e_1$-neighborhood $\V(\e_1)\subset T^*\hattuM _0$ of the set of the light-like covectors
$L^*\hattuM _0$, and for which 
$(I-B_{\e_1,\e_2})$ is microlocally smoothing
operator 
 in the conic $\e_2$-neighborhood $\V(\e_2)$ of {\mattitext the bundle of the light-like co-vectors} $L^*\hattuM _0$.
 Let us decompose the operator $Q=Q_1+Q_2$
 where  $Q_1=QB_{\e_1,\e_2}$ and $Q_2=Q(I-B_{\e_1,\e_2})$.

{\muutos The Schwartz kernel $Q_1(x,y)$ of the operator $Q_1$} satisfies  $Q_1 
 \in \I({\muutos M_0}\times {\muutos M_0};
   \Delta_{T^*\hattuM _0}^{\prime}, \Lambda_{\nohat g})$,
    similarly to ${{} {\bQ}}$.
Moreover, the Schwartz kernel $Q_2(x,y)$ of the operator $Q_2$ satisfies $Q_2\in \I({\muutos M_0}\times {\muutos M_0};
   \Delta_{T^*\hattuM _0}^{\prime})$
and the operator $Q_2$  is a pseudodifferential operator that has the form
\beq\label{eq: representation of Q* p}
{\muutos ({\bQ}_2v)(y)=\int_{\mmmtext M_0\times \R^{4}}e^{i \Psi_1(y,z,\xi)}\sigma_{Q_2}(y,z,\xi)v(z)\,dzd\xi,}
\eeq
{\muutos where $\Psi_1(y,z,\xi)$ parametrises  the diagonal Lagrangian manifold $\Delta_{T^*M}^\prime$ and
$\sigma_{Q_2}(z,y,\xi)\in S^{-2}_{cl}({\muutos M_0}\times {\muutos M_0}; \R^4)$ is
a classical symbol. {\muutos When $X:W\to \R^4$ are local coordinates in an open set $W\subset M$,
the restriction  $Q_2:C^\infty_0(W)\to C^\infty(W)$, given by $v\mapsto Q_2v|_W$, can be written using 
the phase function $\Psi_1(y,z,\xi)=\sum_{j=1}^4(X^j(y)-X^j(z))\xi_j$ and} {\mmmmtext  symbol
$\sigma_{Q_2}(z,y,\xi)\in S^{-2}_{cl}({\muutos W\times W}; \R^4)$. It}  has the principal symbol}
\beq\label{AAA-formula}
\sigma_{Q_2}^{(p)}(y,z,\xi)= {\chi(z,y,\xi)}\frac 1{\nohat g^{jk}(y)\xi_j\xi_k},\quad y,z\in W,
\eeq
where $\chi(z,y,\xi)\in C^\infty$ vanishes {\muutos when $(y,\xi)\in T^*W$ is in some neighborhood of light-like co-vectors $L^*M$.}

{\bf 6. Products of $u_j$ and the singular support of   $\mathcal S$.}
{\mmtext In the computations below,  we will
represent the waves {\muutos $u_j
\in \I^{n}(K_j)=\I^{n-1/2}(N^*K_j)$} in the local coordinates $X:W\to \R^4$,  $X(x)=(X^j(x))_{j=1}^4\in \R^4$,}
that are adapted to the surface $K_j$, as
 \beq\label{U1term}
& &u_j(x)=\int_{\R}e^{i\psi_j(x,\theta)}\sigma_{u_j}(x,\theta)d\theta,\quad \sigma_{u_j}(x,\theta)\in
S^n_{cl}(W;\R),
\eeq
{\mmtext where $\psi_j(x,\theta)=\theta\,\cdotp X^j(x)$.}

{\mmtext Next, let us consider two indexes $j,k\in \{1,2,3,4\}$, $j\not =k$,
and use  local coordinates $X:W\to \R^4$
that are adapted to the surfaces $K_j$ and $K_k$.} 
Recall that  $\Lambda_j=N^*K_j$ and denote $\Lambda_{jk}=N^*(K_j\cap K_k)$. By \cite[{\ftext Lemma}\ 1.2]{GU1}, the  
pointwise product satisfies
$u_j\,\cdotp u_k\in \I(\Lambda_j,\Lambda_{jk})+ \I(\Lambda_k,\Lambda_{jk})$.
{\muutos Also, the Lagrangian manifolds $\Lambda_j$  and $\Lambda_k$ are invariant {\ftext by} the bicharacteristic flow in the future direction.
By using \cite[Prop.\ 2.2 and 2.3]{GU1}, 
we see that} ${\bQ}(a\vez u_j\,\cdotp u_k)\in \I(\Lambda_j,\Lambda_{jk})+ \I(\Lambda_k,\Lambda_{jk})$ 
can be written as  
\beq\label{Q U1U2term}
G_{jk}(x):={\bQ}(a\vez u_j\,\cdotp u_k)=\int_{\R^2}e^{i\psi_{jk}(x,\theta,\theta')}\sigma_{G_{jk}}(x,\theta,\theta')d\theta d\theta',\hspace{-1cm}
\eeq
{\mmtext where $x\in W$, $(\theta,\theta')\in \R^2$, $\psi_{jk}(x,\theta,\theta')=\theta X^j(x)+\theta' X^k(x)$} and
$\sigma_{G_{jk}}(x,\theta,\theta')$ is a sum of product type symbols, see (\ref{product symbols}).

As  $ N^*(K_j\cap K_k)
\setminus (N^*K_j\cup N^*K_k)$ consists of vectors which
are non-characteristic for $\square_{\nohat g}$,
 \MTEXT{the principal symbol 
$\sigma_{G_{jk}}^{(p)}(x,\theta,\theta')$ of $G_{jk}$ on  $N^*(K_j\cap K_k)
\setminus (N^*K_j\cup N^*K_k)$ is  given by
\beq\label{product type symbols}
& &\sigma_{G_{jk}}^{(p)}{\muutos (x,\theta,\theta')}
=
 s(x,\theta,\theta')a\vez (x)\sigma_{u_j}^{(p)}(x,\theta)\sigma_{u_k}^{(p)}(x,\theta'),
 \\ \nonumber
 & &\hspace{0cm}s(x,\theta,\theta')
   =\frac1{ g(\xi,\xi)},\quad\hbox{where $\xi=d_x\psi_{jk}(x,\theta,\theta')=\theta dX^j+\theta' dX^k$} 
   \hspace{-1cm}
\eeq
 and $g(\xi,\xi)=g^{jk}(x)\xi_j\xi_k$.

{\mmmmtext Let us next consider {\ftext the} singular supports of {\ftext the} functions $S_{\sigma}$ given
in (\ref{4th interaction for wave eq B}). 
   Let us start with {\ftext the} case when the permutation $\sigma$  is the identity map.   
   
%
   
   As we showed above, for $i\not =j$  we have
 ${\bQ}(a\vey u_i\,\cdotp u_j)\in \I(\Lambda_i,\Lambda_{ij})+\I(\Lambda_i,\Lambda_{ij})$,
 so that by (\ref{eq wavefronts}), the wave front set of this function is
a {\mmmmtext subset of $N^*K_i\cup N^*K_j\cup N^*K_{ij}$.
Thus,
\beq \label{WF of Qu1u2}
\hbox{singsupp}(\bQ(a\vey u_i\,\cdotp u_j))\subset K_i\cup K_j.
\eeq
Moreover, as $\hbox{WF}(u_3)\subset N^*K_3$ {\vtext and the intersection of $K_{12}$  and $K_3$ is transversal,}
{\ftext the theorem} for the wave front set of the pointwise product of distributions, \cite[Thm.\ 1.3.6]{Duistermaat}, yield that $F_{321}=au_3\,\cdotp {\bQ}(a\vey u_2\,\cdotp u_1)$ satisfies 
\beq\label{WF of F321}\\ \nonumber
\hbox{WF}(F_{321})\cap T^*_{x}M\subset \P^{(123)}_x= (\bigcup_j N^*_{x}K_j)\cup
(\bigcup_{j,k} N^*_{x}K_{jk})\cup (\bigcup_{j,k,l} N^*_{x}K_{jkl}),\hspace{-14mm}
\eeq
where $x\in M_0$ and  {\mmmtext $j,k,l\in \{1,2,3\}$ and} we interpret $N^*_{x}K_j$  to be an empty set if $x\not\in K_j$ etc. 
Thus, $\hbox{WF}(F_{321})\cap L^*M \subset \X\cup (\cup_{j=1}^4 N^*K_j)$. As $Q$  is the causal inverse of the 
wave operator, 
we see  by using
 H\"ormander's theorem on {\ftext propagation of singularities} along bicharacteristics, \cite[Theorem 26.1.4]{H4}, that
%
%
{\vtext
 \beq\label{WF QF321}
 \hbox{WF}(QF_{321})
 \subset \hbox{WF}(F_{321})\cup  {\H}_{123}\cup (\bigcup_{j=1}^3 N^*K_j)\hspace{-1cm}
 \eeq
 where $ {\H}_{123}=  {\H}_{123}((\vec x,\vec\xi),s_0)$, see (\ref {H-set}).
  In particular, this implies that
 $\hbox{singsupp}(QF_{321})\subset  \Y\cup (\cup_{j=1}^3 K_j).$ 

 Formulas (\ref{WF of Qu1u2}) and (\ref{WF QF321}) {\ftext give} that   
  for $\sigma=Id$ we have
 \beq\label{singsup S123}
 \hbox{singsupp}\,(S_{\sigma})\subset  \Y\cup (\bigcup_{j=1}^4 K_j).
   \eeq
 The same arguments yield that  (\ref{singsup S123})   holds for all permutations $\sigma$.}

%

 }

{\bf 7. Decomposition of the source term  $\mathcal S$.}
{\muutos Below we will analyze  the wave front set of
the source $\mathcal S$  that {\mmmtext is produced by} the 
{\mmmmtext fourth order} {\muutos interaction}. 
To this end, we use
the decomposition $Q=Q_1+Q_2$, and  write the source $\mathcal S$  in the form}}
\beq\label{wave eq: m indicator A1} 
\mathcal S={   {\mathcal S}^{(1)}}+{   {\mathcal S}^{(2)}}+{   {\mathcal S}^{(3)}},
\quad {   {\mathcal S}^{(p)}}
=
\sum_{\sigma\in \Sigma(4)} 
{   {\mathcal S}^{(p)}_{\sigma\vel}},\ p\in \{1,2,3\}\hspace{-1cm}
\eeq
where $\Sigma(4)$ is the set of permutations of the set $\{1,2,3,4\}$
and
\beq\label{eq- S sources}
{\mathcal S}^{(1)}_{\sigma\vel}
&=&-4
 a\vex  u_{\sigma(4)} \, \cdotp   Q_1(
a\vey  u_{\sigma(3)}\,\cdotp Q( a\vez  u_{\sigma(2)}\,\cdotp u_{\sigma(1)})),\\
\nonumber
{\mathcal S}^{(2)}_{\sigma\vel}
&=&-4
 a\vex  u_{\sigma(4)} \, \cdotp   Q_2(
a\vey  u_{\sigma(3)}\,\cdotp Q( a\vez  u_{\sigma(2)}\,\cdotp u_{\sigma(1)})),\\
\nonumber
{   {\mathcal S}^{(3)}_{\sigma\vel}}
&=&- a\vex Q(a\vey  u_{\sigma(4)}\,\cdotp  \,u_{\sigma(3)})\,\cdotp
  Q(a\vez  u_{\sigma(2)}\,\cdotp \, u_{\sigma(1)}).
\eeq
}

%
%

Later, we consider the terms (\ref{eq- S sources}) 
{\mmtext with {\ftext the} permutation $\sigma=Id$. 
 Note that the terms corresponding to the other permutations $\sigma$  can be analyzed similarly by renumbering the indexes.}}
 
{ 
 

%



  \begin{center}

\psfrag{1}{$(x_1,\xi_1)$}
\psfrag{2}{\hspace{-5mm}$(x_2,\xi_2)$}
\psfrag{3}{\hspace{-5mm}$(x_3,\xi_3)$}
\psfrag{4}{\hspace{-5mm}$(x_4,\xi_4)$}
\psfrag{5}{\hspace{6mm}$y$}
\psfrag{6}{$z$}
\psfrag{7}{$x_0$}
\psfrag{8}{\hspace{-3mm}$\omega^\sharp$}
\psfrag{9}{\hspace{-2mm}$w_4^\sharp$}
\psfrag{0}{}
\includegraphics[height=5cm]{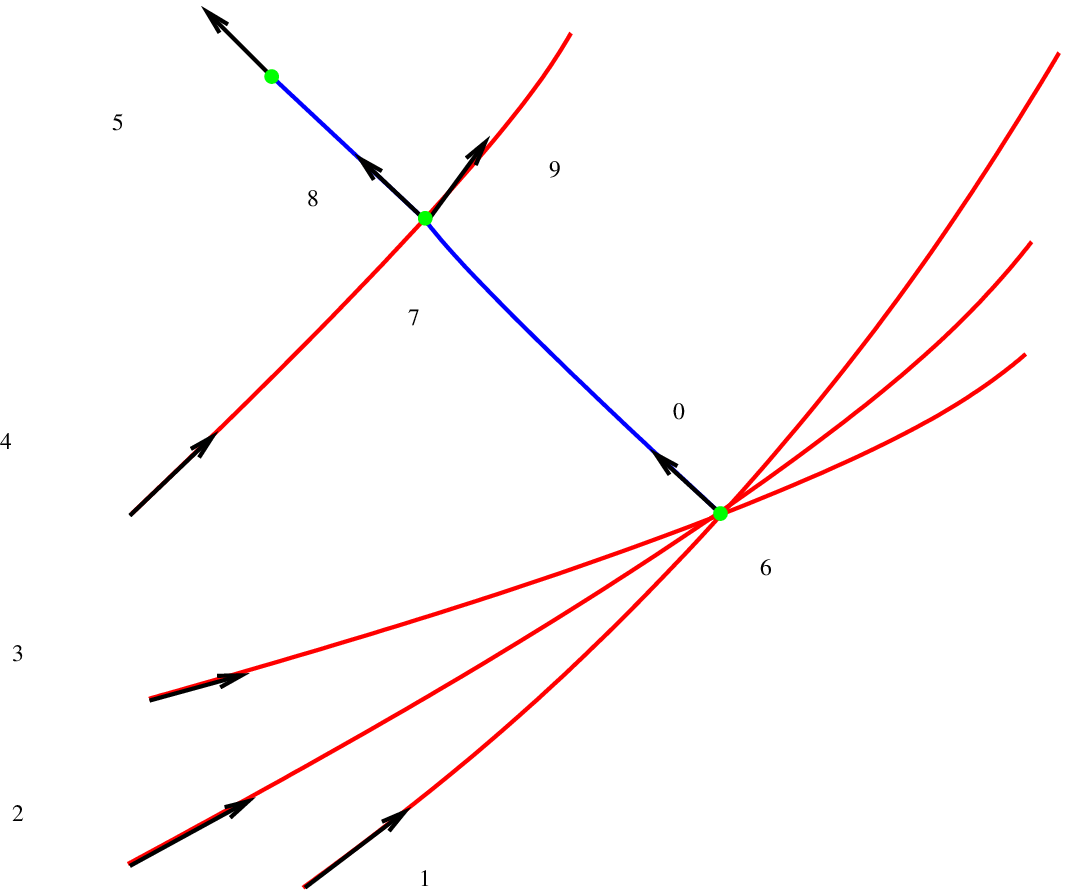}
\psfrag{1}{$q$}  
\psfrag{2}{$y$}  
\quad
\includegraphics[height=5cm]{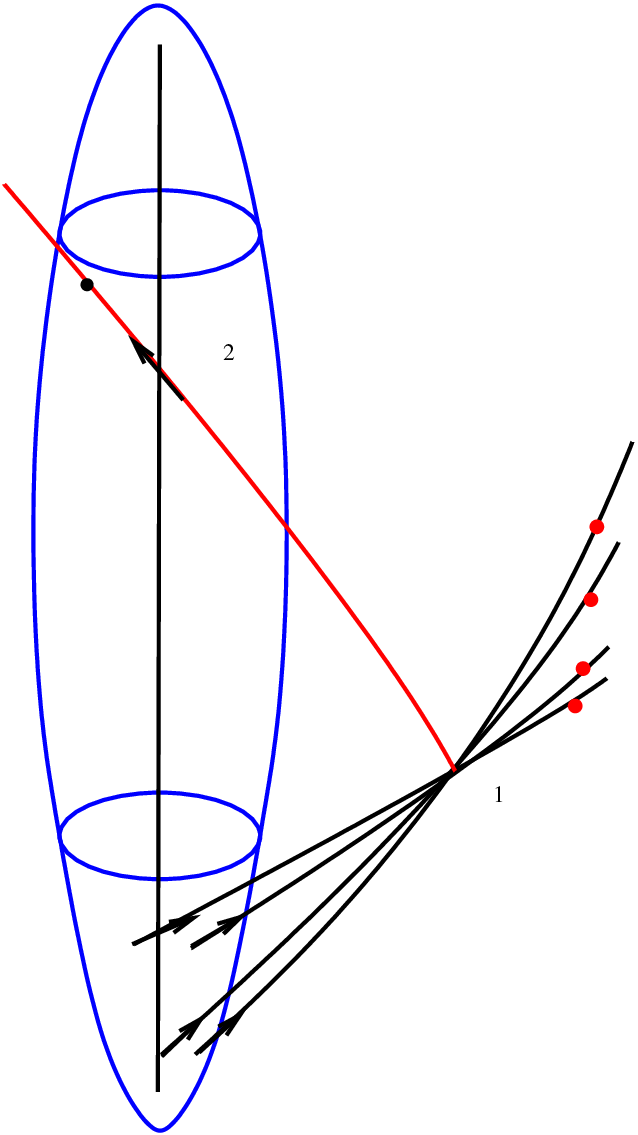}  

\end{center}

\noindent
{\it FIGURE 6. {\bf Left:} {\vtext The figure shows the case  when three geodesics intersect at $z$ and 
the waves propagating near these geodesics interact and create a wave that hits
to the fourth geodesic at the point $x_0$.  The produced singularities propagate to  {\muutos the point ${\muutos y}\in \Y$}. Note that $z$ and $x_0$ may be conjugate points on the geodesic connecting
them {\mmtext or the waves propagating from $z$ to $x_0$ may have caustics}.}
{\bf  Right:} {\muutos Geodesics corresponding to directions $(x_j,\xi_j)$, $j=1,2,3,4$ intersect at the point $q$  {\mmmtext and $b_j$, $j=1,2,3,4$ are linearly independent}} and the Condition I is valid for the point $y$ with vectors $(x_j,\xi_j)$ and with the parameter $q$. The red points are the conjugate points of $\gamma_{x_j,\xi_j}([{\rtext 0},\infty))$ and $\gamma_{q,w}([0,\infty))$.
}
\medskip

    {\bf 8. Analysis the  wave front set of source $\mathcal S^{(1)}_{\sigma\vel}$.}   
{\mmmmtext       In this step we {\ftext consider the case} when $(x_0,\zeta_0)$ is in the wave front 
   set of  source functions $\mathcal S^{(1)}_{\sigma\vel}$, see Steps 2 and 5.

{\vtext 
 
 We start in the case when $\sigma=Id$. Then, $\mathcal S^{(1)}_{Id}=au_4\,\cdotp Q_1F_{321}.$
 To analyze the product $u_4\,\cdotp Q_1F_{321}$ in the set $\V_0$,
 we will first show that the  wave front sets satisfy for $z\in \V_0$
\beq\label{WF sums}
\hbox{If $ (z,\hatomega)\in
\hbox{WF}(Q_1F_{321})$ and $ (z,w_4)\in \hbox{WF}(u_4)$ then $\hatomega+w_4\not =0$.}\hspace{-15mm}  \eeq
  To show this, we assume the opposite, that there are  $z_0\in \V_0$ and
\ba
(z_0,\hatomega)\in \hbox{WF}(Q_1F_{321})\quad\hbox{and}\quad
 (z_0,w_4)\in \hbox{WF}(u_4)\subset N^*K_4,
 \ea
such that
  $\hatomega+w_4=0$.
  We consider different cases for  $(z_0,\hatomega)$ that are given by
  equations (\ref{WF of F321}) and (\ref{WF QF321}).

First, we consider  the case when $(z_0,\hatomega) \in N^*K_{123}$.
Since  $z_0\in K_4$, this yields $z_0\in \cap_{j=1}^4K_j$.
Thus we are case $(B)$ and $z_0=q$. However, as  $z_0=q\in \V_0$,
the vectors $b_j$, $j=1,2,3,4$ are linearly independent,
and it is not possible that $\hatomega=-w_4 \in N^*_qK_{123}\cap N^*_qK_4$.
Thus we see that  $(z_0,\hatomega) \not\in N^*K_{123}$.

Second,  we consider the case when $(z_0,\hatomega) \in N^*K_{j_1j_2}\setminus (N^*K_{j_1}\cup N^*K_{j_2})$
where $j_1,j_2\in \{1,2,3\}$, $j_1\not= j_2$.
Then, $\hatomega$  is not light-like and so it is not possible that $\hatomega=-w_4\in N^*_{z_0}K_4$.
Thus we {\ftext conclude} that $(z_0,\hatomega) \not\in N^*K_{j_1j_2}\setminus (N^*K_{j_1}\cup N^*K_{j_2})$.

Third, we see that 
 it is not possible that $\hatomega=-w_4\in (\bigcup_{j=1}^3 N^*_{z_0}K_j)\cap N^*K_4$ as
 the surfaces $K_i$ and $K_j$, $i\not =j$ intersect transversally.

 Fourth, we consider the remaining case when $(z_0,\hatomega)\in {\H}_{123}\subset L^*M$.
Then there is 
  $(x,\zeta)\in  \X_{123}\subset N^*K_{123}$, $x\in K_{123}$ such that 
  $x\le z_0$ and the bicharacteristic $\Theta_{x,\zeta}$  passes through $(z_0,\hatomega)$
  and $\Theta_{x,\zeta}=\Theta_{z_0,\hatomega}$.
   Also, $(z_0,\hatomega)=(z_0,-w_4)\in L^*K_4=\Lambda_4$.
  As $x\in  J^+(x_1)$, we see using (\ref{eq: summary of assumptions 1}) that $x_4\not \in J^+(x)$.
Thus,  as  $(x,\zeta)\in \Theta_{z_0,\hatomega}=\Theta_{z_0,-w_4}$, we have $(x,\zeta)\in \Lambda_4=N^*K_4$ and $x\in K_4$.
  These imply that $x\in \cap_{j=1}^4 K_j$. Hence, the case $(B)$  is valid and $x=q$. Moreover, as $(x,\zeta)$  is in the intersection
  of $N^*_qK_{123}$ and $N^*_qK_4$, we {\ftext get} that
{\ftext the} vectors $b_j$, $j=1,2,3,4$  are not linearly independent. As $x=q$ and $x\leq z_0$,
this  implies $z_0\not \in \V_0$, that not possible by our assumptions. 

As none of the above four case is possible, we {\ftext obtain} that it is not possible that 
  $\hatomega+w_4=0$. Hence, (\ref{WF sums}) is true.

 Next we consider the question, can 
 $({\muutos x_0},\zeta_0)$, given in (\ref{y assumptions1}), be in $\hbox{WF}(\S^{(1)}_{id})$.
Recall that  $\zeta_0$  is light-like.

 Due to  (\ref{WF sums}) we can use the formula for the wave front set of the pointwise product of distributions, \cite[Thm.\ 1.3.6]{Duistermaat}.
 It implies that  if $({{x_0}},\zeta_0)\in \hbox{WF}(\S^{(1)}_{id})$,
 where $\S^{(1)}_{id}= au_4\,Q_1(F_{321})$, then 
 there are
 \ba
 ({{x_0}},\hatomega)\in\hbox{WF}( Q_1F_{321})
\cup (M_0\times \{0\})\hbox{ and  }({{x_0}},w_4)\in \hbox{WF}(u_4)
 \cup (M_0\times \{0\})
 \ea
such that
\ba
\zeta_0=\hatomega+w_4.
\ea

Let us use the fact that $w_4$ is light-like or zero and $\zeta_0=\hatomega+w_4$ is light-like. If $\hatomega$ is light-like or zero, then $\zeta_0$ 
has to be parallel either to $\hatomega$ or   $w_4$. Then
 $(x_0,\zeta_0)\in \H\cup (\cup_{j=1}^4 N^*K_j)$, see (\ref{H-set}) and (\ref{WF QF321}).
 However, this is not possible {\ftext since} $y_0\not \in \Y\cup (\cup_{j=1}^4 K_j)$.
 Hence, $\hatomega$ is not light-like or zero.
 
Since $(x_0,\hatomega)\in\hbox{WF}( Q_1F_{321})$ is not zero and 
we see by using the definition of $Q_1$
that 
\beq\label{Ne1 formula}
(x_0,\hatomega)\in \V(\e_1).
\eeq
 Also, as $\hatomega$ is not light-like and zero and $\zeta_0=\hatomega+w_4$ is light-like
we have  $w_4\not=0$. Then $x_0\in K_4$.

As the above implies that  $\hatomega\in \P^{(123)}_{x_0}\setminus\{0\}$, we can consider separately
the different cases given by definition of $P^{(123)}_{x_0}$
in (\ref{WF of F321}).

First, as $\hatomega$  is not light-like, we have $\hatomega\not \in N^*_{x_0}K_{j}$ for $j=1,2,3$.

Second, if $\hatomega\in N^*_{x_0}K_{j_1j_2}$,  with $j_1,j_2\in \{1,2,3\}$, we have $x_0\in K_{j_1j_2}\cap K_4$ and $(x_0,\zeta_0)=(x_0,\hatomega+w_4)\in N^*(K_{j_1j_2}\cap K_4)$. As $\zeta_0$  is light-like, we have
 $(x_0,\zeta_0)\in \X$. This implies $y_0\not \in \Y$ {\ftext which is not} possible by our assumptions.
 Thus $\hatomega\not \in N^*_{x_0}K_{j_1j_2}$.
 
Third, consider the case when $\hatomega\in N^*_{x_0}K_{123}$. Then $x_0\in \cap_{j=1}^4 K_j$. 
Hence, $(B)$ is valid and $x_0=q$.
Again, as $x_0=q\in \V_0$, we see that
 the vectors $b_j$, $j=1,2,3,4$ have to be linearly independent. 
 Then,  the condition (\ref{case a1}) is valid.  Note that as $x_0=q$, we have $\zeta_0=\hatzeta$.
%
We recall that in the case (\ref{case a1})
we chose $\e_1>0$  to be so small that (\ref{eq: conditions 2 modified}) is valid.

 
As  $b_j$, $j=1,2,3,4$ are linearly independent,  $  \hatzeta$ has a unique
representation $\hatzeta=\sum_{j=1}^4 \zeta_j$ where $\zeta_j\in N^*_qK_j$. Then,
$\hatomega=\sum_{j=1}^3\zeta_j$ and $w_4=\zeta_4$. 
Then,
 by (\ref{eq: conditions 2 modified}), we have that $(x_0,\hatomega)=(x_0,\hatzeta-\zeta_4)\not \in  \V(\e_1).$
 This is not possible {\ftext since} $({{x_0}},\hatomega)\in \V(\e_1)$ by (\ref{Ne1 formula}).
Hence, we have  $(x_0,\zeta_0)\not \in \hbox{WF}(\S^{(1)}_{\sigma})$.}

 {{} Above, we have analyzed the case when the permutation $\sigma$ is the identity.
 For other permutations $\sigma\in \Sigma(4)$, the same computations are valid with a renumbering of indexes, and we 
 {\ftext conclude} that
 \beq\label{eq S1 result}
({\muutos x_0},\zeta_0)\not \in \hbox{WF}( \mathcal S^{(1)}).
 \eeq
 }

 {\bf 9. Analysis of the term $ \mathcal S^{(2)}_\sigma$.}
Let start by analyzing the case when $\sigma=Id$.  
{\mmmmtext Recall that the wave front set of $F_{321}=au_3\,\cdotp {\bQ}(a\vey u_2\,\cdotp u_1)$ satisfies
(\ref{WF of F321}).}
Also, as $Q_2$ is a pseudodifferential operator,
$$
\hbox{WF}(Q_2F_{321})\subset \hbox{WF}(F_{321}).
$$
{\vtext Recall, in the set $\V_0$  the intersections $K_4\cap K_{123}$ and 
$K_{j_1j_2}\cap K_{j_1}$ are transversal,
 where $j_1,j_2,j_3\in \{1,2,3\}$. 
Thus as  $\hbox{WF}(u_4)\subset N^*K_4$ and $\S^{(2)}_{Id}=au_4\,\cdotp Q_2F_{321}$, 
we can apply}
the formula for the wave front set of the pointwise product of distributions, 
\cite[Thm.\ 1.3.6]{Duistermaat}, 
{\mmmtext and see {\vtext for $x\in \V_0$ that}
\beq\label{Z-space argument}
\hbox{WF}(\S^{(2)}_{Id})\cap T^*_{x}M\subset Z_{x}
\eeq
where 
\beq\label{Z-space}
Z_{x}=
(\bigcup_j N^*_{x}K_j)\cup
(\bigcup_{j,k} N^*_{x}K_{jk})\cup (\bigcup_{j,k,l} N^*_{x}K_{jkl})\cup (\oplus_{j=1}^4 N^*_{x}K_j)\hspace{-14mm}
\eeq
{\ftext and}} {\mmmtext $j,k,l\in \{1,2,3,4\}$.} \ftext We} interpret $N^*_{x}K_j$  to be an empty set if $x\not\in K_j$ etc. 

{\mmmtext {\vtext Consider light-like co-vector $({\muutos x_0},\zeta_0)$ given in (\ref{y assumptions1}).
First, we consider the cases when (A) holds or (B) holds and $x_0\not =q$. 
 Then,} 
\beq \label{formula F} 
{\mmmtext Z_{x_0}}
\subset (
{\mmmtext \bigcup_j N^*_{x_0}K_j)\cup
(\bigcup_{j,k} N^*_{x_0}K_{jk})\cup (\bigcup_{j,k,l} N^*_{x_0}K_{jkl}}).
\eeq 
Then, if
$(x_0,\zeta_0)\in Z_{x_0}$, formula
(\ref{formula F}) yields  $(x_0,\zeta_0)\in\X\cup 
(\cup_{j=1}^4 N^*_{x_0}K_j)$. By  (\ref{y assumptions1}), this yields
that  ${\muutos y_0}\not\in {\mathcal Y}((\vec x,\vec \xi),{\muutos s_0})\cup(\cup_{j=1}^4K_j)$
{\ftext {\ftext {\ftext which is not}}} possible by our assumptions. Thus}   in all the above mentioned cases
 $(x_0,\zeta_0)\not \in \hbox{WF}(\S^{(2)}_{Id})$.}
%
%
%
%
%
  In particular, this holds under the assumption of claim (i).
The similar analysis holds for a  general permutation $\sigma$.

{\mmmtext 
{\vtext We have above considered the cases when (A) holds or (B) holds and $x_0\not =q$.
It remains to consider the case when (B) holds and $x_0=q$. As $x_0\in \V_0$, then the condition (\ref{case a1})
has to be valid. {\vtext Moreover, as $x_0=q$,} we have $\zeta_0={\mmmtext \hatzeta}$, and}
\beq\label{y is the intersection point B2}
\hbox{there are $\zeta_j\in \bar N^*_qK_j$, $j=1,2,3,4,$ so that ${\mmmtext \hatzeta}=\sum_{j=1}^4 \zeta_j\in T^*_qM$,}\hspace{-1.4cm}
\eeq
and  such $\zeta_j\in \bar N^*_qK_j$ are uniquely determined by ${\mmmtext \hatzeta}$.
Recall that $({\mmmtext q},{\mmmtext \hatzeta})\not \in \mathcal X\cup\bigcup_{j=1}^4N^*K_j$. 
This implies $\zeta_j\not =0$ for all $j=1,2,3,4$, that is,  $\zeta_j\in N^*_qK_j$.
Also, note that as ${\mmmtext \hatzeta}$ and $\zeta_{\sigma(4)}$ are light-like vectors that are not parallel,
we have that $\zeta_{\sigma(1)}+\zeta_{\sigma(2)}+\zeta_{\sigma(3)}={\mmmtext \hatzeta}-\zeta_{\sigma(4)}$ is not light-like or zero.

{\mmtext {\vtext Next, in a neighborhood of} the point $q$ we use local coordinates $X:W_0\to \R^4$ that are adapted to all surfaces $K_j$, see Step 2.}

Let $R_{j,k}=R_{j,k}(x,D)$ be a pseudodifferential operator which symbol is one
in a conic neighborhood $\Sigma_{j,k}\subset T^*M_0$ of $(q,\zeta_j+\zeta_k)$ and vanishes in a conic
 neighborhood of sets $N^*K_{i}$, $i=1,2,3,4$. Also, let
 $R_{i,j,k}=R_{i,j,k}(x,D)$ be a pseudodifferential operator which symbol is one
in a conic neighborhood 
in $\Sigma_{i,j,k}\subset T^*M_0$ of  $(q,\zeta_i+\zeta_j+\zeta_k)$ and vanishes 
in a conic
 neighborhood of all sets $N^*K_{i_1}$ and $N^*K_{i_1i_2}$.
 Moreover, let $R_{4321}=R_{4321}(x,D)$ be a pseudodifferential operator 
 which symbol is one
in a conic neighborhood of $(q,{\mmmtext \hatzeta})$ and 
vanishes in a conic
 neighborhood of all sets $N^*K_{i_1}$, $N^*K_{i_1i_2}$ and
$N^*K_{i_1i_2i_3}$. We 
 denote $R^\sigma_{jk}=R_{\sigma(j),\sigma(k)}$ and $R^\sigma_{ijk}=R_{\sigma(i),\sigma(j),\sigma(k)}$.
{\mmtext Also, assume that the Schwartz kernels of these operators {\mmmmtext and that of $R_{4321}$} are supported in $W_0\times W_0$.}
Then we define
\ba
{\mathcal S}^{(2),1}_{\sigma\vel}
&=&-4
 a\vex  u_{\sigma(4)} \, \cdotp   Q_2(
a\vey  u_{\sigma(3)}\,\cdotp R^\sigma_{21} Q( a\vez  u_{\sigma(2)}\,\cdotp u_{\sigma(1)})),\\
{\mathcal S}^{(2),2}_{\sigma\vel}
&=&-4
 a\vex  u_{\sigma(4)} \, \cdotp   R^\sigma_{321}Q_2(
a\vey  u_{\sigma(3)}\,\cdotp R^\sigma_{21} Q( a\vez  u_{\sigma(2)}\,\cdotp u_{\sigma(1)})),\\
{\mathcal S}^{(2),3}_{\sigma\vel}
&=&-4R_{4321}(
 a\vex  u_{\sigma(4)} \, \cdotp   R^\sigma_{321}Q_2(
a\vey  u_{\sigma(3)}\,\cdotp R^\sigma_{21} Q( a\vez  u_{\sigma(2)}\,\cdotp u_{\sigma(1)}))).
\ea
Using (\ref{eq wavefronts}) and \cite[Thm.\ 1.3.6]{Duistermaat} we observe that then
\ba
{\mathcal S}^{(2)}_{\sigma\vel}-{\mathcal S}^{(2),1}_{\sigma\vel}
&=&-4
 a\vex  u_{\sigma(4)} \, \cdotp   Q_2(
a\vey  u_{\sigma(3)}\,\cdotp (I-R^\sigma_{21}) Q( a\vez  u_{\sigma(2)}\,\cdotp u_{\sigma(1)})),
\ea
satisfies 
\ba
& &\hbox{WF}( \mathcal S_{\sigma\vel}^{(2)}-\mathcal S_{\sigma\vel}^{(2),1})\cap T^*_{q}M
\subset\\
& &\quad \bar N^*_{q}K_{\sigma(4)}+\bar N^*_{q}K_{\sigma(3)}+
((\bar N^*_{q}K_{\sigma(2)}+\bar N^*_{q}K_{\sigma(1)})\setminus \Sigma_{\sigma(2),\sigma(1)}).
\ea
As ${\mmmtext \hatzeta}$ is written in  (\ref{y is the intersection point B2}) in {\ftext a} unique way as a sum of  vectors in $\bar N^*_qK_j$, $j=1,2,3,4$, we see that
$
(q,{\mmmtext \hatzeta})\not \in \hbox{WF}( \mathcal S_{\sigma\vel}^{(2)}-\mathcal S_{\sigma\vel}^{(2),1}).
$
Similarly,
\ba
& &(q,{\mmmtext \hatzeta})\not \in \hbox{WF}( \mathcal S_{\sigma\vel}^{(2),1}-\mathcal S_{\sigma\vel}^{(2),2}),\quad (q,{\mmmtext \hatzeta})\not \in \hbox{WF}( \mathcal S_{\sigma\vel}^{(2),2}-\mathcal S_{\sigma\vel}^{(2),3}).
\ea
Hence,
\beq\label{final2}
(q,{\mmmtext \hatzeta})\not \in \hbox{WF}( \mathcal S_{\sigma\vel}^{(2)}-\mathcal S_{\sigma\vel}^{(2),3}).
\eeq

{\mmtext 
To simplify notations, we next consider   {\ftext the case when}  $\sigma=Id$. Using formulas (\ref{microlocally away}) and (\ref{eq: traditional}), we see that 
 $R^\sigma_{21} Q( a\vez  u_2\,\cdotp u_1)\in \I(K_{12})$. Then 
 the results for the products of conormal distributions, \cite[Lemma 1.1]{GU1},  imply
$u_3\,\cdotp R^\sigma_{21} Q( a\vez  u_2\,\cdotp u_1)\in \I(K_{12},K_{123})+ \I(K_{3},K_{123})$. By (\ref{microlocally away}), we have $R^\sigma_{321}(u_3R^\sigma_{21} Q( a\vez  u_2\,\cdotp u_1))\in \I(K_{123})$.
Repeating the arguments, and computing above the orders of  {\ftext the} distributions,  we see that   $\mathcal S_{\sigma\vel}^{(2),3}\in 
\I^{-4n-3}(\Sigma_q)$, where $\Sigma_q=N^*\{q\}=T^*_qM\setminus \{0\}$, see (\ref{Normal bundle of a point}). 
The other permutations can be analyzed in the same way.}

 {\bf 9. Analysis of the term $ \mathcal S^{(3)}_\sigma$.}
We start by considering the case when $\sigma=Id$.  
{\muutos 
The term $\S^{(3)}_{Id}$ is the product of 
  the terms ${\bQ}(a\vey u_2\,\cdotp u_1)$ 
  and ${\bQ}(a\vez u_4\,\cdotp u_3)$. As above, we see that 
${\bQ}(a\vey u_2\,\cdotp u_1)\in \I(\Lambda_1,\Lambda_{12})+\I(\Lambda_2,\Lambda_{12})$. Similarly, 
${\bQ}(a\vez u_4\,\cdotp u_3)\in \I(\Lambda_3,\Lambda_{34})+\I(\Lambda_4,\Lambda_{34}) $. Then (\ref{eq wavefronts}) and  \cite[Thm. 1.3.6]{Duistermaat} {\mmmtext yield that {\vtext for  $x\in \V_0$,}
$$
{\vtext \hbox{WF}(\S^{(3)}_{Id})\cap T^*_{x}M\subset Z_{x}}.
$$
{\vtext 
Consider next the cases when  (A) holds or  (B) holds and $x_0\not =q$.}
Then, the same arguments that were used {\mmmmtext above} to analyze the formula (\ref{Z-space argument})
yield that 
%
 $(x_0,\zeta_0)\not \in \hbox{WF}(\S^{(3)}_{Id})$.}} 
%

{\mmmtext {\vtext Again, as we have considered 
the cases when  (A) holds or  (B) holds and $x_0\not =q$,
it remains to consider the case (\ref{case a1}). Then, as $x_0=q$, we have} $\zeta_0={\mmmtext \hatzeta}$.}
We use the same notations as in Step 8.}

Let ${\mathcal S}^{(3),0}_{\sigma\vel}={\mathcal S}^{(3)}_{\sigma\vel}$ and
\ba
\nonumber
{\mathcal S}^{(3),1}_{\sigma\vel}
&=&-  a\vex  Q(a\vey  u_{\sigma(4)}\,\cdotp  \,u_{\sigma(3)})\,\cdotp
  R^\sigma_{21}Q(a\vez  u_{\sigma(2)}\,\cdotp \, u_{\sigma(1)}),\\
\nonumber
{\mathcal S}^{(3),2}_{\sigma\vel}
&=&-  a\vex  R^\sigma_{43}Q(a\vey  u_{\sigma(4)}\,\cdotp  \,u_{\sigma(3)})\,\cdotp
  R^\sigma_{21}Q(a\vez  u_{\sigma(2)}\,\cdotp \, u_{\sigma(1)}),\\
\nonumber
{\mathcal S}^{(3),3}_{\sigma\vel}
&=&-  R_{4321}(a\vex  R^\sigma_{43}Q(a\vey  u_{\sigma(4)}\,\cdotp  \,u_{\sigma(3)})\,\cdotp
  R^\sigma_{21}Q(a\vez  u_{\sigma(2)}\,\cdotp \, u_{\sigma(1)})).
\ea
Again, using (\ref{eq wavefronts}) and \cite[Thm.\ 1.3.6]{Duistermaat} we see that for $i=0,1,2$
\ba
(q,{\mmmtext \hatzeta})\not \in \hbox{WF}( \mathcal S_{\sigma\vel}^{(3),i}-\mathcal S_{\sigma\vel}^{(3),i+1}).
\ea
so that
\beq\label{final3}
(q,{\mmmtext \hatzeta})\not \in \hbox{WF}( \mathcal S_{\sigma\vel}^{(3)}-\mathcal S_{\sigma\vel}^{(3),3}).
\eeq
Using formula (\ref{microlocally away}) and \cite[Lemma 1.1]{GU1}, we see that $\mathcal S_{\sigma\vel}^{(3),3}\in 
\I^{-4n-3}(\Sigma_q)$.

{\bf 10. Proof of claim (i).}
Above we have shown  in the case of claim (i) that $({\muutos x_0},\zeta_0)\not \in \hbox{WF}( \mathcal S_{\sigma\vel}^{(p)})$ for all $p=1,2,3$,
and hence $(y_0,\wflat)$ can not be in $ \hbox{WF}( \mathcal U^{(4)})$.
{\mmmmtext As $y_0\not \in \cup_{j=1}^4 K_j$  and $\wflat\in L^{*}_{y_0}M$  can be arbitrary,  
 (\ref{singsup S123}) and H\"ormander's theorem on {\ftext propagation of singularities} along bicharacteristics, \cite[Theorem 26.1.4]{H4}, prove the claim (i).}

  {\bf 11.  Proof of claim (ii).} {\mmmmtext Assume that condition (\ref{case a1}) is valid.} 
   Let 
 \ba
{\mathcal S}^{mod}= \sum_{\sigma\in\Sigma(4)}   {\mathcal S}^{(2),mod}_{\sigma}+  {\mathcal S}^{(3),mod}_{\sigma},
\quad  {\mathcal S}^{(2),mod}_{\sigma}= {\mathcal S}^{(2),3}_{\sigma},\ 
{\mathcal S}^{(3),mod}_{\sigma}= {\mathcal S}^{(3),3}_{\sigma}. 
 \ea
 {\mmmmtext Since $\mathcal S^{mod}\in 
\I^{-4n-3}(\Sigma_q)$, we have  $\hbox{WF}(\mathcal S^{mod})\subset T^*_qM$.  In steps 7,8 and 9 we have shown that in the case 
when $(x_0,\zeta_0)\in \Theta_{y_0,\wflat}$ is not equal to $(q,{\mmmtext \hatzeta})$, we have
$(x_0,\zeta_0)\not \in 
 \hbox{WF}(
  {\mathcal S}^{(p)}_{\sigma})$ for $p=1,2,3$. Also, 
  by  (\ref{eq S1 result}), (\ref{final2}) and (\ref{final3})  we have $(q,{\mmmtext \hatzeta})\not \in \hbox{WF}(
  {\mathcal S}_{\sigma}- {\mathcal S}_{\sigma}^{mod})$. These show that the bicharasteristic $\Theta_{y_0,\wflat}$
  does not intersect $\hbox{WF}(
  {\mathcal S}- {\mathcal S}^{mod})$.}}
  Thus, using 
 H\"ormander's theorem on {\ftext propagation of singularities} along bicharacteristics, \cite[Theorem 26.1.4]{H4},
   we see 
 {\muutos that  $({\muutos y_0},{\muutos \wflat})$ is not in the wave front set of the function}
$\M^{(4)}-{\bQ}\S^{mod}={\bQ}(\S-\S^{mod})$.
{\vtext By replacing $\wflat$ by $-\wflat$, the arguments  {\ftext above}  show also that 
$({\muutos y_0},-{\muutos \wflat})$ is not in the wave front set of ${\bQ}(\S-\S^{mod})$.}
{\vtext Since $y_0\in 
  \firstpoint^{reg}_U(q)$, see Def.\ \ref{def: earliest element},
%
%
the light-like geodesic from $q$  to $y_0$  has no cut points and there is only
one  light-like geodesics connecting $q$ to ${\muutos y_0}$.}
{\muutos Moreover, as $y_0\not \in \Y\cup \bigcup _{1\leq j\leq 4} K_j,$ we  {\ftext obtain}  from (\ref{singsup S123})  that
 {\ftext the}  function
${\mathcal S}$ {\muutos is smooth in a neighborhood of ${\muutos y_0}$.
As $\mathcal S^{mod}\in 
\I^{-4n-3}(\Sigma_q)$, also ${\mathcal S}^{mod}$ is smooth in a neighborhood of ${\muutos y_0}$.
{\mmmtext Thus the above and \cite[Theorem 26.1.4]{H4}
 show that  $(q,w)\not \in \hbox{WF}(Q({\mathcal S}-{\mathcal S}^{mod}))$ for all $w\in T^*_{y_0}M$.}
Hence, $y_0$  has a neighborhood $W$ such that ${\bQ}(\S-\S^{mod})$ is $C^\infty$-smooth in $W$.}  

%
%
%
%
%
%
%

By \cite[Prop.\ 2.1]{GU1}, ${\bQ}:\I^{-4n-3}(\Sigma_q)\to \I^{-4n-3-3/2,-1/2}(\Sigma_q,\Lambda_q^+)$. This and (\ref{microlocally away})
imply {\muutos that} ${\muutos y_0}$ as a neighborhood $W$ such that ${\bQ}\S^{mod}|_W$ and thus $\M^{(4)}|_W$ are in $\I^{-4n-3-3/2}(\Lambda_q^+)= \I^{m}(S)$,
where $S=\L^+(q)\cap W$. Here, $S$  is a smooth surface as the light-like geodesic from $q$  to $y_0$   {\ftext does not have}  cut points.

{\mmmmtext We consider  {\ftext now}  the case when $x_0 =q$,
  so that $\zeta_0={\mmmtext \hatzeta}$, and compute
the principal symbol of $\mathcal S^{mod}\in \I^{-4n-3}(\Sigma_q)$}.
The principal symbol of 
$F^{mod}=R_{3,2,1}(au_3\,\cdotp R_{2,1}{\bQ}(a\vey u_2\,\cdotp u_1))\in \I(K_{123})$
at $(q,\xi)\in N^*_qK_{123}$, where $\xi=\zeta_1+\zeta_2+\zeta_3$, is by (\ref{product type symbols})
$$
\sigma^{(p)}_{F^{mod}}({\muutos q},\xi)=
\frac{Ca(q)^2}{ g(
\zeta_1+\zeta_2,\zeta_1+\zeta_2)}
 \prod_{j=1}^3 \sigma_{u_j}^{(p)}(q,\zeta_j),
$$
and the principal symbol of $Q_2$  at $(q,\xi)$ is given by (\ref{AAA-formula}).
  Using these, we  {\ftext obtain}  that the principal symbols of the sources ${\mathcal S}^{(2),mod}_{\sigma}$ and ${\mathcal S}^{(3),mod}_{\sigma}$, with $\sigma=Id$, at $({\muutos q},{\mmmtext \hatzeta})$, are given by
\beq\label{T-example2 wave A}
& &\sigma^{(p)}_{{\mathcal S}^{(2),mod}_{Id}}({\muutos q},{\mmmtext \hatzeta})=
\frac{Ca(q)^3}{\nohat g(\zeta_1+\zeta_2,\zeta_1+\zeta_2)\, \nohat g(
\sum_{j=1}^3
\zeta_j , \sum_{k=1}^3 \zeta_k  )}
 \prod_{j=1}^4 \sigma_{u_j}^{(p)}(q,\zeta_j),
\hspace{-0.6cm}
\eeq
and
\beq\label{T-example2 wave B}
& &\sigma^{(p)}_{{\mathcal S}^{(3),mod}_{Id}}({\muutos q},{\mmmtext \hatzeta})=
\frac{Ca(q)^3}{g(\zeta_1+\zeta_2,\zeta_1+\zeta_2)\, 
g(\zeta_3+\zeta_4,\zeta_3+\zeta_4)}
\prod_{j=1}^4 \sigma_{u_j}^{(p)}(q,\zeta_j),\hspace{-0.6cm}
\eeq
where we recall that ${\mmmtext \hatzeta}=\sum_{j=1}^4\zeta_j$, where $\zeta_j\in N^*_{{\muutos q}}K_j$.
%


 Lemma \ref{lem: Lagrangian 1 wave} implies
that the principal symbol of $\M^{(4)}$  at $({\muutos y_0},{\muutos \wflat})$ is the product of a non-zero function $R({\muutos y_0},{\muutos \wflat},{\muutos q},{\mmmtext \hatzeta})$ times the principal symbol of $\S^{mod}$ at $({\muutos q},{\mmmtext \hatzeta})$.
%
{\muutos    This and formulas (\ref{T-example2 wave A}) and (\ref{T-example2 wave B}),
written for general permutation $\sigma$, yield formula  (\ref{eq: principal symbols of U4}). This proves the claim (ii).}}
\hfill \Box \medskip


\extension{

\noindent
{\bf Proof.}  
Below, to simplify notations,
we denote $K_j={   K(x_{j},\xi_{j},{\muutos s_0})}$ and 
$K_{123}=K_1\cap K_2\cap K_3$ and $K_{124}=K_1\cap K_2\cap K_4$, etc.
We will denote $ \Lambda_j={   \Lambda(x_j,\xi_j,{\muutos s_0})}$
to consider also the singularities of $K_j$ related to conjugate points.

By our assumption, we can assume that following linear independency condition is satisfied:
\medskip

(LI) Assume if that  $J\subset \{1,2,3,4\}$ and $y\in J^-(x_6)$
are such that for all $j\in J$  we have $\gamma_{x_j,\xi_j}(t^\prime_j)=y$
with some $t^\prime_j\geq 0$, then  the vectors
$\dot\gamma_{x_j,\xi_j}(t_j^\prime )$, $j\in J$ are linearly independent.
\medskip


By the definition of ${\bf t}_j$,  if
the intersection of geodesics, $\gamma_{{\muutos y_0},{\muutos \wflat^\sharp}}(\R_-) \cap(\cap_{j=1}^4\gamma_{   x_j,\xi_j}((0,{\bf t}_j)))
$ is non-empty, it can
contain only one point. In the case that such a
point exists, we denote it by $q$. When $q$ exists,
the intersection of $K_j$ at this point are transversal and
we see that when $s_0$ is small enough, the set
$\cap_{j=1}^4 K_j$ consists only of the point $q$.
We assume below that $s_0$ is such that  this is true.

We define    two local coordinates $Z:W_0\to \R^4$
and  $Y:W_1\to \R^4$ such that $W_0,W_1\subset {   \V(\vec x,\vec \xi)}$,
 see definition after (\ref{eq: summary of assumptions 2}).
 We assume that these local coordinates are such that
$K_j\cap W_0=\{x\in W_0;\ Z^j(x)=0\}$ and $K_j\cap W_1=\{x\in W_1;\ Y^j(x)=0\}$ for $j=1,2,3,4$.
In the Fig.\ 6, $W_0$ is a neighborhood of $z$ and
$W_1$ is a neighborhood of $y$. 
We note that the origin $0=(0,0,0,0)\in \R^4$ does not
necessarily belong to the set $Z(W_0)$ or the set $Y(W_1)$,
for instance in the case when the four geodesics $\gamma_{x_j,\xi_j}$ do not intersect.
However, this is the case when the  four geodesics intersect at the point $q$, we need
to consider the case when $W_0$ and $W_1$ are neighborhoods of $0$. 
 %
Note that
$W_0$ and $W_1$ do not contain any cut points of the geodesic
 $\gamma_{x_j,\xi_j}([{   0},\infty)$.

We will denote $z^j=Z^j(x)$.
We assume that  $Y:W_1\to \R^4$ are similar coordinates
and denote  $y^j=Y^j(x)$. We also denote
 below $dy^j=dY^j$ and  $dz^j=dZ^j$.
 
 Let us next considering the map ${   {\bQ}}:C^\infty_0(W_1)\to
 C^\infty(W_0).$  By \cite{MU1}, ${   {\bQ}}\in I(W_1\times W_0;
   \Delta_{T\hattuM _0}^{\prime}, \Lambda_{\nohat g})$ is an  operator 
   with a classical symbol and its canonical
   relation 
$\Lambda_{{   {\bQ}}}^{\prime}$ is associated to a union of two  intersecting Lagrangian manifolds,
$\Lambda_{{   {\bQ}}}^{\prime}= \Lambda_{\nohat g}^{\prime}\cup \Delta_{T\hattuM _0}$,
intersecting cleanly \cite{MU1}.
Let $\e_2>\e_1>0$ and $B_{\e_1,\e_2}$ be a pseudodifferential
operator on $\hattuM _0$ which is  a microlocally smoothing
operator (i.e., the full  symbol vanishes  in local coordinates)
outside of the $\e_2$-neighborhood $\V_2\subset T^*\hattuM _0$ of the set of the light like covectors
$L^*\hattuM _0$ and for which 
$(I-B_{\e_1,\e_2})$ is microlocally smoothing
operator 
 in the $\e_1$-neighborhood $\V_2$ of $L^*\hattuM _0$.
 The neighborhoods here are defined 
 with respect to the Sasaki metric of $(T^*\hattuM _0,\nohat g^+)$ and
$\e_2,\e_1$ are chosen later in the proof. 
 Let us decompose the operator ${   {\bQ}}={   Q}_2+{   Q}_1$
 where  ${   Q}_2={   {\bQ}}(I-B_{\e_1,\e_2})$ and   ${   Q}_1={   {\bQ}}B_{\e_1,\e_2}$.
 As  $\Lambda_{{   {\bQ}}}= \Lambda_{\nohat g}\cup \Delta_{T\hattuM _0}^{\prime}$,
we see that then there is 
a neighborhood $ \W_2=\W_2(\e_2)$ of $L^*\hattuM _0\times L^*\hattuM _0\subset (T^*\hattuM _0)^2$
such that the Schwartz
kernel ${   Q}_1(r,y)$ of the operator ${   Q}_1$ satisfies
\beq\label{eq: W_2 neighborhood}
\hbox{WF}({   Q}_1)\subset \W_2.
\eeq
Moreover,  $\Lambda_{{   Q}_2}\subset \Delta_{T\hattuM _0}^{\prime}$ implying
that ${   Q}_2$ is a pseudodifferential operator with a classical symbol,
 ${   Q}_2\in I(W_1\times W_0;
   \Delta_{T\hattuM _0}^{\prime})$,
 and
 ${   Q}_1 
 \in I(W_1\times W_0;
   \Delta_{T\hattuM _0}^{\prime}, \Lambda_{\nohat g})$
   is a Fourier integral operator (FIO) associated
   to two cleanly intersecting Lagrangian manifolds, similarly to ${   {\bQ}}$. 
In the case when  $p=1$ we can write ${   {\bQ}}_p$ as 
\beq\label{eq: representation of Q* p}
({   Q}_2v)(z)=\int_{\R^{4+4}}e^{i(y-z)\cdotp \xi}\sigma_{Q_2}^{(p)}(z,y,\xi)v(y)\,dyd\xi,
\eeq
where 
{
a classical symbol $\sigma_{Q_2}^{(p)} (z,y,\xi)\in S^{-2}(W_1\times W_0; \R^4)$ with a real
valued principal symbol
\ba
\hat \sigma_{Q_2}^{(p)} (z,y,\xi)=\frac {\chi(z,\xi)}{\nohat g^{jk}(z)\xi_j\xi_k}
\ea
where $\chi(z,\xi)\in \C^\infty$ is a cut-off function vanishing
in a neighborhood of  the set where $g(\xi,\xi)=0$. Note that
then $Q_2-Q_2^*\in \Psi^{-3}(W_1\times W_0)$.}

Furthermore, let us decompose 
${   Q}_2={   {\bQ}}_{1,1}+{   {\bQ}}_{1,2}$ corresponding to the decomposition
$q_1(z,y,\xi)=q_{1,1}(z,y,\xi)+q_{1,2}(z,y,\xi)$
 of the symbol,
where
 \beq\label{eq: q_{1,1} and q_{1,2}} \hspace{-.3cm}
 & &q_{1,1}(z,y,\xi)=q_{1}(z,y,\xi)\psi_R(\xi),\\
 & & \nonumber q_{1,2}(z,y,\xi)=q_{1}(z,y,\xi)(1-\psi_R(\xi)),
 \hspace{-2cm}
 \eeq
and $\psi_R\in C^\infty_0(\R^4)$
is a cut-off function that is equal to one in  a ball $B(R)$ of radius $R$ specified below.

Next we start to consider the terms $T_{\tau}^{(4)}$ and $\tilde T_{\tau}^{(4)}$ of the type
 (\ref{T-type source wave}) 
and
(\ref{tilde T-type source wave}). In these terms,
we can represent the gaussian beam $u_{\tau}(z)$ in $W_1$ in the form 
\beq\label{eq:representation of the wave}
u_\tau(y)=e^{i\tau \varphi (y)}{{\Phi}}(y,\tau)
\eeq
where the function $ \varphi$ is a complex phase function having non-negative imaginary
part such that $\im   \varphi$, defined on $W_1$, vanishes exactly on the geodesic
 $\gamma_{x_5,\xi_5}\cap W_1$. Note that 
  $\gamma_{x_5,\xi_5}\cap W_1$ may be empty. Moreover,
 ${{\Phi}}\in
S^{0}_{clas}(W_1;\R)$ is a classical symbol.

 Also, for $y=\gamma_{x_5,\xi_5}(t)\in W_1$
we have that $d\varphi(y)=c\dot \gamma_{x_5,\xi_5}(t)^\flat$,  where $c\in \R\setminus \{0\}$, is light-like.



We consider first the asymptotics of 
 terms $T_\tau^{(4)}$ and $\tilde T^{(4)}_\tau$ of the type
 (\ref{T-type source wave}) 
and
(\ref{tilde T-type source wave})
   where  the symbols $\sigma_{u_j}(z,\theta_j)$, $j\leq 3$ and  
    $\sigma_{u_j}(y,\theta_j)$, $j\in \{4,5\}$ are 
   symbols written in the $Z$ and $Y$ coordinates.

Let us consider functions
$U_j\in \I^{p_j}(K_j)$, $j=1,2,3,$ supported in $W_0$ and 
$U_4\in \I(K_4)$,  supported in $W_1$,  having classical symbols. They have the form
 \beq\label{U1term}
& &U_j(x)=\int_{\R}e^{i\theta_jx^j}\sigma_{u_j}(x,\theta_j)d\theta_j,\quad a_j\in
S^{p_j}_{clas}(W_{\kappa(j)};\R), 
\eeq
for all $j=1,2,3,4$ (Note that here the phase function is  $\theta_jx^j=\theta_1x^1$ for $j=1$ etc,
that is, there is no summing over index $j$). We may assume that $\sigma_{u_j}(x,\theta_j)$ vanish near $\theta_j=0$.

Since  $x_5\in  {   \V(\vec x,\vec \xi)}\cap U_{\nohat g}$ and 
$W_0,W_1\subset  {   \V(\vec x,\vec \xi)}$, we see that
an example of functions (\ref{U1term}) are $U_j(z)=\Phi_{\kappa(j)}(z)u_j(z)$, $j=1,2,3,4$,
where $u_j(z)$ are the distorted plane  waves. Here and below, 
 $\kappa(j)=0$ for $j=1,2,3$  and $\kappa(4)=1$ and we
also denote $\kappa(5)=1$.
Note that $p_j=n$ correspond to the case when $U_j\in \I^n(K_j)=\I^{n-1/2}(N^*K_j)$.

Denote $\Lambda_j=N^*K_j$ and $\Lambda_{jk}=N^*(K_j\cap K_k)$. By \cite[{\ftext Lemma}\ 1.2 and 1.3]{GU1}, the  
pointwise product 
$U_2\,\cdotp U_1\in \I(\Lambda_1,\Lambda_{12})+ \I(\Lambda_2,\Lambda_{12})$  and thus
by \cite[Prop.\ 2.2 and 2.3]{GU1},
${\bQ}(U_2\,\cdotp U_1)\in \I(\Lambda_1,\Lambda_{12})+ \I(\Lambda_2,\Lambda_{12})$ 
and  it
can be written as 
\beq\label{Q U1U2term}
{\bQ}(U_2\,\cdotp U_1)=\int_{\R^2}e^{i(\theta_1z^1+\theta_2z^2)}c_1(z,\theta_1,\theta_2)d\theta_1d\theta_2.
\eeq
Note that here $c_1(z,\theta_1,\theta_2)$ is sum of product type symbols, see (\ref{product symbols}).
As  $ N^*(K_1\cap K_2)
\setminus (N^*K_1\cup N^*K_2)$ consists of vectors which
are non-characteristic for 
the wave operator, that is, the wave operator is elliptic in
a neighborhood of this subset of the cotangent bundle,   the principal symbol 
$\hat c_1$ of ${\bQ}(U_2\,\cdotp U_1)$ on  $N^*(K_1\cap K_2)
\setminus (N^*K_1\cup N^*K_2)$ is  given by
\beq\label{product type symbols}
& &\hat  c_1(z,\theta_1,\theta_2)\sim
 s(z,\theta_1,\theta_2)\sigma_{u_1}(z,\theta_1)\sigma_{u_2}(z,\theta_2),
 \\ \nonumber
 & &\hspace{-1cm}s(z,\theta_1,\theta_2)
=1/{ \nohat g(\theta_1b^{(1)}+\theta_2b^{(2)},\theta_1b^{(1)}+\theta_2b^{(2)})}
 =1/({2 \nohat g(\theta_1b^{(1)},\theta_2b^{(2)})}). \hspace{-1cm}
 \eeq
Note that $s(z,\theta_1,\theta_2)$ is a smooth function  
on $N^*(K_1\cap K_2)
\setminus (N^*K_1\cup N^*K_2)$ and homogeneous of order $(-2)$ in $\theta=(\theta_1,\theta_2)$.
Here, we use $\sim$ to denote that the symbols have the same principal symbol.
%
%
%
%
Let us next
make computations in the case when 
 $\sigma_{u_j}(z,\theta_j)\in C^\infty(\R^4\times \R)$ is positively homogeneous for $|\theta_j|>1$, that is, we have 
 $ \sigma_{u_j}(z,s)=a^{\prime}_j(z)s^{p_j}$, where $p_j\in \N$  and $|s|>1$.
 We consider $T_\tau^{(4)}=\sum_{p=1}^2T^{(1)}_{\tau,p}^{(4)}$
 where $T^{(1)}_{\tau,p}^{(4)}$ is defined as 
 $T_{\tau}^{(4)}$ by replacing the term ${   {\bQ}}(U_4\,\cdotp u_{\tau})$
  by ${   {\bQ}}_p(U_4\,\cdotp u_{\tau})$.

 Let us now consider the case $p=2$ and choose the parameters $\e_1$ and $\e_2$
that determine the
decomposition ${   {\bQ}}={   Q}_1+{   Q}_2$. 
First, we observe that for $p=2$ we can write using $Z$ and $Y$ coordinates
 \beq\label{eq; T tau asympt, p=2}
T^{(1)}_{\tau,2}^{(4)}
&=&
\tau^{4}
\int_{\R^{12}}e^{i\tau \Psi_2(z,y,\theta)}
c_1(z,\tau \theta_1,\tau \theta_2)\cdotp\\ \nonumber
& &\hspace{-1.3cm}
\cdotp  \sigma_{u_3}(z,\tau \theta_3)
{   Q}_1(z,y)
\sigma_{u_4}(y,\tau \theta_4) {{\Phi}}(y,\tau )\,d\theta_1d\theta_2d\theta_3d\theta_4 dydz,\\
\nonumber \Psi_2(z,y,\theta)
&=&\theta_1z^1+\theta_2z^2+\theta_3z^3+\theta_4y^4+\varphi(y).
\eeq

Denote $\theta=(\theta_1,\theta_2,\theta_3,\theta_4)\in \R^{4}$.
Consider the case when $(z,y,\theta)$ is a critical point of $\Psi_2$
satisfying $\im \varphi(y)=0$. Then we have
 $\theta^{\prime}=(\theta_1,\theta_2,\theta_3)=0$
and $z^{\prime}=(z^1,z^2,z^3)=0$, 
$y^4=0$,
  $d_y\varphi(y)=(0,0,0,-\theta_4)$,
implying that $y\in K_4$ and 
$(y,d_y\varphi(y))\in N^*K_4$.
As  $\im \varphi(y)=0$, we have that 
  $y=\gamma_{x_5,\xi_5}(\tsparameter)$ with some $\tsparameter\in\R_-$.
 As we have $\dot \gamma_{x_5,\xi_5}(\tsparameter)^\flat=d_y\varphi(y)\in N^*_yK_4$,
we obtain $\gamma_{x_5,\xi_5}([\tsparameter,0])\subset K_4$. However, 
this 
is not possible by our assumption 
$x_5\not \in \cup_{j=1}^4 K(x_{j},\xi_{j},s_0)$ when $s_0$ is small enough.
Thus  the phase function 
$ \Psi_2(z,y,\theta)$ has no critical points satisfying $\im \varphi(y)=0$.

%

When the orders $p_j$ of the symbols $a_j$ are small enough, the integrals
in the $\theta$ variable in  (\ref{eq; T tau asympt, p=2}) are convergent in the classical sense.
We use  now properties of the wave front set to 
compute the asymptotics of  oscillatory integrals
and to this end we introduce the function
\beq
\tilde {   Q}_1(z,y,\theta)={   Q}_1(z,y),
\eeq
that is, consider ${   Q}_1(z,y)$ as a constant function
in $\theta$. Below, denote $\psi_4(y,\theta_4)=\theta_4y^4$ and $r=d\varphi(y)$. 
Note that then $d_{\theta_4}\psi_4=y^4$ and $d_y\psi_4=(0,0,0,\theta_4)$.
Then  in $W_1\times W_0\times \R^4$
\ba
d_{z,y,\theta}\Psi_2&=&(\theta_1,\theta_2,\theta_3,0;
r+d_y\psi_4(y,\theta_4),
z^1,z^2,z^3,d_{ \theta_4}\psi_4(y,\theta_4))
\\
&=&(\theta_1,\theta_2,\theta_3,0;
d\varphi(y)+(0,0,0,\theta_4),
z^1,z^2,z^3,y^4)
\ea
and we see that if $((z,y,\theta),d_{z,y,\theta}\Psi_2)\in \hbox{WF}(\tilde {\bQ}_2^*)$
and  $\im \varphi(y)=0$, we have $(z^1,z^2,z^3)=0$, $y^4=d_{\theta_4}\psi_4(y,\theta_4)=0$
and $y\in  \gamma_{x_5,\xi_5}$. Thus 
$z\in K_{123}$ and $y\in  \gamma_{x_5,\xi_5}\cap K_4$.

Let  us use the following notations
\beq\label{new notations}
& &z\in K_{123},\quad \omega_\theta:=(\theta_1,\theta_2,
 \theta_2,0)=\sum_{j=1}^3\theta_j dz^j\in T^*_z\hattuM _0,\\
 \nonumber
& &y\in K_4\cap \gamma_{x_5,\xi_5},\quad (y,w):=(y,d_y\psi_4(y,\theta_4))\in N^*K_4,
\\
 \nonumber
& & r=d\varphi(y)=r_jdy^j\in T^*_y\hattuM _0,\quad \kappa:=r+w.
 \eeq
Then, $y$ and $ \theta_4$ satisfy $y^4=d_{ \theta_4}\psi_4(y,\theta_4)=0$ and
$w=(0,0,0,\theta_4)$.
 
 Note that  by definition of the $Y$  coordinates $w$ is a light-like covector.
 By definition of the  $Z$ coordinates,
 $\omega_\theta\in 
 N^*K_1+N^*K_2+N^*K_3=N^*K_{123}$.  
 
 Let us first consider what happens if $\kappa=r+w=d\varphi(y)+(0,0,0,\theta_4)$ is light-like.
In this case, all vectors $\kappa$, $w$, and $r$ are light-like
and satisfy $\kappa =r+w$. This is possible
only if  $r\parallel w$, i.e., $r$ and $w$ are parallel, see \cite[Cor.\ 1.1.5]{SW}.  Thus  $r+w$ is light-like if and only
if $r$ and $w$ are parallel.

Consider next the case when  $(x,y,\theta)\in W_1\times W_0\times \R^4$ 
is such that $((x,y,\theta),d_{z,y,\theta}\Psi_2)\in \hbox{WF}(\tilde {\bQ}_2^*)$
and   $\im \varphi(y)=0$.
Using the above notations (\ref{new notations}), we obtain
$
d_{z,y,\theta}\Psi_2=(\omega_\theta,
r+w;(0,0,0,d_{ \theta_4}\psi_4(y, \theta_4)))=(\omega_\theta,d\varphi(y)+(0,0,0,\theta_4);(0,0,0,y^4))$,
 where $y^4=d_{ \theta_4}\psi_4(y, \theta_4)=0$,
and thus we have
\ba
((z,\omega_\theta),(y,
r+w))\in \hbox{WF}( {   {\bQ}}_w)=\Lambda_{{   Q}_1}.
\ea

\begin{center}

\psfrag{1}{$(x_1,\xi_1)$}
\psfrag{2}{\hspace{-5mm}$(x_2,\xi_2)$}
\psfrag{3}{\hspace{-5mm}$(x_3,\xi_3)$}
\psfrag{4}{\hspace{-5mm}$(x_4,\xi_4)$}
\psfrag{5}{\hspace{-5mm}$(x_5,\xi_5)$}
\psfrag{6}{$z$}
\psfrag{7}{$y$}
\psfrag{8}{$r$}
\psfrag{9}{$w$}
\psfrag{0}{$\omega_\theta$}
\includegraphics[height=5.5cm]{three_and_one_int5.eps}
\psfrag{1}{${\hat x}^\prime$}  
\psfrag{2}{$p_2$}  
\psfrag{3}{$J({p^-},{p^+})$}  
\psfrag{4}{${\hat x}$}  
\psfrag{5}{}  
\psfrag{6}{${\hat x}_1$}  
\psfrag{6}{${\hat x}_1$}  
\psfrag{7}{$\tilde p$}  
\psfrag{8}{$$}  
\includegraphics[width=5.5cm]{iteration-2_mod.eps}  

\end{center}

\noindent
{\it FIGURE 6. {\bf Left:} In the figure we consider the case (A1) where three geodesics intersect at $z$ and 
the waves propagating near these geodesics interact and create a wave that hits
 the fourth geodesic at the point $y$, the produced singularities are detected by the gaussian beam
source at the point $x_5$. Note that $z$ and $y$ can be conjugate points on the geodesic connecting
them. In the case (A2) the points $y$ and $z$ are the same.
{\bf Right:} Condition I  is valid.
}
\medskip

Since $\Lambda_{{   Q}_2}\subset \Lambda_{\nohat g}\cup \Delta_{T\hattuM _0}^{\prime}$,
this implies that one of the following conditions are valid:
\ba
& &(A1)\ ((z,\omega_\theta),(y,
r+w))\in \Lambda_{\nohat g}, \hspace{5cm} \
\\ & &\hspace{-0.5cm}\hbox{or}\\ 
& &(A2)\ ((z,\omega_\theta),(y,
r+w))\in \Delta_{T\hattuM _0}^{\prime}.
\ea
Let
$\gamma_0$ 
be the geodesic with
$
\gamma_0(0)=z,\ \dot \gamma_0(0) =\omega_\theta^\sharp.
$
Then (A1) and (A2) are equivalent to the following conditions:
\ba
& &(A1)\ \hbox{There is $\tsparameter\in \R$ such
that $(\gamma_0(\tsparameter), \dot \gamma_0(\tsparameter)^\flat) =(y,\kappa)$ and},\\
& &\quad\quad \ \ \hbox{
 the vector $\kappa$ is  light-like,}
\\ & &\hspace{-0.5cm}\hbox{or}\\ 
& &(A2)\ z=y\hbox{ and }
\kappa=-\omega_\theta.
\ea

Consider next the case when (A1)  is valid. 
As  $\kappa$ is light-like,  $r$ and $w$ are parallel. 
Then, since $(\gamma_{x_5,\xi_5}(t_1),\dot \gamma_{x_5,\xi_5}(t_1)^\flat)=(y,r)$ we see that $\gamma_0$ is a continuation of the geodesic $ \gamma_{x_5,\xi_5}$,
that is, for some $t_2$ we
have $(\gamma_{x_5,\xi_5}(t_2),\dot \gamma_{x_5,\xi_5}(t_2))=(z,\omega_\theta)
\in N^*K_{123}$. This implies that $x_5\in {\mathcal Y}$  {\ftext {\ftext {\ftext which is not}}} possible
by our assumptions. Hence (A1) is not possible.

Consider next the case when (A2)  is valid.Then we would also have 
that $r\parallel w$  then $r$ is parallel to $\kappa=-\omega_\theta\in N^*K_{123}$,
 and since $(\gamma_{x_5,\xi_5}(t_1),\dot \gamma_{x_5,\xi_5}(t_1)^\flat)=(y,r)$ we
 would have  that $x_5\in {\mathcal Y}$. As this is not possible by our assumptions,
 we see that $r$ and $w$ are not parallel. This implies that
 $\omega_\theta=-\kappa$ is not light-like.


For any given  $(\vec x,\vec \xi)$  and $(x_5,\xi_5)$ there exists
 $\e_2>0$ so that $(\{(y, \omega_\theta)\}\times T^*M)\cap \W_2(\e_2)=\emptyset$,
 see (\ref{eq: W_2 neighborhood}), and thus
 \beq\label{eq: y and omega}
 (\{(y, \omega_\theta)\}\times T^*M)\cap \hbox{WF}({\bQ}_2^*)=\emptyset.
 \eeq
 Next we assume that $\e_2>0$ and also $\e_1\in (0,\e_2)$ are chosen
 so that  (\ref{eq: y and omega}) is valid.
 Then there are no $(z,y,\theta)$ such that 
  $((z,y,\theta),d\Psi_2(z,y,\theta))\in  \hbox{WF}(\tilde {\bQ}_2^*)$ 
 and  $\im \varphi(y)=0$. Thus by 
  Corollary 1.4 in \cite{Delort} or
  \cite[{\ftext Lemma}\ 4.1]{Ralston}
 yields  $T^{(1)}_{\tau,2}^{(4)}=O(\tau^{-N})$ for all $N>0$.
Alternatively, one can use 
 the complex version of \cite[Prop. 1.3.2]{Duistermaat}, 
 obtained using combining the proof of  \cite[Prop. 1.3.2]{Duistermaat} and 
 the method of stationary phase with a complex phase, see \cite[Thm.\ 7.7.17]{H1}.
 
 Thus to analyze the asymptotics of  $T_{\tau}^{(4)}$ 
 we need to consider only $T^{(1)}_{\tau,{\muutos \vec \ell,\beta}}^{(4)}$.
Next, we analyze the case when $U_4$ is a conormal
distribution and has the form (\ref{U1term}).

%
%

 Let us thus consider the case $p=1$. Now
 \beq
U_4(y)\,\cdotp u_{\tau}(y)&=&
\int_{\R^{1}}e^{i\theta_4y^4+i\tau \varphi(y)}
\sigma_{u_4}(y,\theta_4){{\Phi}}(y,\tau )
\,d\theta_4. \label{QU4Ut term PRE}
\eeq
 We obtain
by (\ref{eq: representation of Q* p}) 
\beq\nonumber
({   Q}_2(U_4\,\cdotp u_{\tau}))(z)&=&
\int_{\R^{9}}e^{i((y-z)\cdotp \xi+\theta_4y^4+\tau \varphi(y))}q_1(z,y,\xi)\cdotp\\
& &\quad
\cdotp \sigma_{u_4}(y,\theta_4){{\Phi}}(y,\tau )
\,d\theta_4\,dy d\xi. \label{QU4Ut term}
\eeq
Then
 $T^{(1)}_{\tau,1}^{(4)}=T^{(1)}_{\tau,1,1}^{(4)}+T^{(1)}_{\tau,1,2}^{(4)}$, 
 cf.\ (\ref{eq: q_{1,1} and q_{1,2}}),
 where 
  \beq\label{eq; T tau asympt pre}
& &\hspace{-.4cm}T^{(1)}_{\tau,1,k}^{(4)}
=
\int_{\R^{16}}e^{i(\theta_1z^1+\theta_2z^2+\theta_3z^3+(y-z)\cdotp \xi+\theta_4y^4+\tau \varphi(y))}
c_1(z,\theta_1,\theta_2)\cdotp
\hspace{-1cm}
\\ \nonumber
& &\cdotp  \sigma_{u_3}(z,\theta_3)
q_{1,k}(z,y,\xi)\sigma_{u_4}(y,\theta_4) {{\Phi}}(y,\tau )\,d\theta_1d\theta_2d\theta_3d\theta_4dydz d\xi,
\eeq
or
 \beq\label{eq; T tau asympt}
& &\hspace{-.1cm}T^{(1)}_{\tau,1,k}^{(4)}=\tau^{8}
\int_{\R^{16}}e^{i\tau (\theta_1z^1+\theta_2z^2+\theta_3z^3+(y-z)\cdotp \xi+\theta_4y^4+\varphi(y))}
c_1(z,\tau \theta_1,\tau \theta_2)\cdotp\hspace{-1cm}\\ \nonumber
& &\hspace{-0.3cm}
\cdotp  \sigma_{u_3}(z,\tau \theta_3)
q_{1,k}(z,y,\tau \xi)\sigma_{u_4}(y,\tau \theta_4) {{\Phi}}(y,\tau )\,d\theta_1d\theta_2d\theta_3d\theta_4dydz d\xi.\hspace{-1cm}
\eeq

%
Let $(\cirt z,\cirt \theta,\cirt y,\cirt \xi)$ be a critical point of
 the phase function
  \beq\label{Psi3 function}
 {\mmtext \Psi_2}(z,\theta,y,\xi)=\theta_1z^1+\theta_2z^2+\theta_3z^3+(y-z)\,\cdotp \xi+\theta_4y^4+\varphi(y).\hspace{-1cm}
 \eeq
 Denote $\cirt  w=(0,0,0,\cirt \theta_4)$ and 
$ \cirt r=d\varphi(\cirt y)=\cirt r_jdy^j $. 
Then
 \beq\label{eq: critical points} 
  & &\p_{\theta_j}{\mmtext \Psi_2}=0,\ j=1,2,3\quad\hbox{yield}\quad \cirt z\in K_{123},\\
  \nonumber
& &\p_{\theta_4}{\mmtext \Psi_2}=0\quad\hbox{yields}\quad \cirt  y\in K_4,\\
  \nonumber
 & &\p_z{\mmtext \Psi_2}=0\quad\hbox{yields}\quad\cirt  \xi=\omega_{\cirt  \theta},\\
   \nonumber
& &\p_{\xi}{\mmtext \Psi_2}=0\quad\hbox{yields}\quad\cirt  y=\cirt z,\\
  \nonumber
 & &\p_{y}{\mmtext \Psi_2}=0\quad\hbox{yields}\quad \cirt \xi=
 -d\varphi(\cirt  y)-\cirt  w.
  \eeq
  The critical points we need to consider for the asymptotics satisfy also
\beq\label{eq: imag.condition} \\ \nonumber
\im \varphi(\cirt y)=0,\quad\hbox{so that }\cirt y\in \gamma_{x_5,\xi_5},\ \im d\varphi(\cirt y)=0,\
\re d\varphi(\cirt y)\in L^{*,+}_{\cirt y}\hattuM _0.\hspace{-1cm}
\eeq


Next\hiddenfootnote{Alternatively, to analyze  the term $T^{(1)}_{\tau,1,2}^{(4)}$
we could use the fact that that  $e^{{\mmtext \Psi_2}}=|\xi|^{-2} (\nabla_z-\omega_\theta)^2 e^{{\mmtext \Psi_2}}$
where $\omega_\theta=(\theta_1,\theta_1,\theta_1,0)$. This can be used in the 
integration by parts in the $z$ variable. 
Possibly, we could use this to analyze directly
$T^{(1)}_{\tau,1}^{(4)}$ and omit the decomposition ${   Q}_2={   {\bQ}}_{1,1}+{   {\bQ}}_{1,2}$.  
} we analyze      the terms $T^{(1)}_{\tau,1,k}^{(4)}$ starting with $k=2$.
 Observe that 
 the 3rd and 5th equations in (\ref{eq: critical points})  imply that at the critical points
 $\xi=( \p_{y_1}\varphi(y), \p_{y_2}\varphi(y),$ $ \p_{y_3}\varphi(y),0)$.
 Thus the critical points are bounded in the $\xi$ variable. 
 Let us now fix the parameter $R$ determining $\psi_R(\xi)$ in
 (\ref{eq: q_{1,1} and q_{1,2}}) so that 
 $\xi$-components of the critical
 points are in a ball $B(R)\subset \R^4$.
  Using the identity  $e^{{\mmtext \Psi_2}}=|\xi|^{-2} (\nabla_z-\omega_\theta)^2 e^{{\mmtext \Psi_2}}$
where $\omega_\theta=(\theta_1,\theta_1,\theta_1,0)$
we can include the operator $|\xi|^{-2}(\nabla_z-\omega_\theta)^2$
in
the integral (\ref{eq; T tau asympt}) with $k=2$
and integrate by
parts. Doing this  two times we can show
that this oscillatory  integral (\ref{eq; T tau asympt}) with $k=2$ becomes an integral of a Lebesgue-integrable
 function. 
 Then, by using method of stationary phase
 and  
%
  the fact that $\psi_R(\xi)$ vanishes at all critical points of ${\mmtext \Psi_2}$
 where $\im {\mmtext \Psi_2}$ vanishes, we see that   $T^{(1)}_{\tau,1,2}^{(4)}=O(\tau^{-n})$
 for all $n>0$.  
%
%


Above, we have shown that  the term $T^{(1)}_{\tau,1,1}^{(4)}$  has the same
asymptotics as $T_\tau^{(4)}$. Next we analyze  this term.
Let $(\cirt z,\cirt \theta,\cirt y,\cirt \xi)$ be a critical point
of $ \Psi_3(z,\theta,y,\xi)$  such that $y$ satisfies (\ref{eq: imag.condition}).
Let us next use the same notations (\ref{new notations}) which
we used above.
Then (\ref{eq: critical points}) and (\ref{eq: imag.condition}) imply
\beq\label{eq: psi1 consequences}
\cirt z=\cirt y\in \gamma_{x_5,\xi_5}\cap \bigcap_{j=1}^4 K_j,\ \cirt \xi=\omega_{\cirt \theta}=-\cirt r-\cirt w.
\eeq
Note that in this case all the four geodesics $\gamma_{x_j,\xi_j}$ intersect
at the point $q$ and by our assumptions, $\cirt r=d\varphi(\cirt y)$ is such a co-vector
that in the $Y$-coordinates $\cirt r=(\cirt r_j)_{j=1}^4$ with $\cirt r_j\not =0$ for all $j=1,2,3,4$.
In particular, this shows that
 the existence of the critical point  of $ \Psi_3(z,\theta,y,\xi)$
implies that there exists 
an intersection point of $\gamma_{x_5,\xi_5}$ and $ \bigcap_{j=1}^4 K_j$.
Equations (\ref{eq: psi1 consequences}) imply also that 
\ba
\cirt r=\sum_{j=1}^4\cirt r_jdy^j=-\omega_{\cirt \theta}-\cirt w=-\sum_{j=1}^3\cirt \theta_jdz^j-
\cirt \theta_4 dy^4.
\ea

%


To consider the case when $\cirt y=\cirt z$, let us assume for a while that that $W_0=W_1$
and that the $Y$-coordinates and $Z$-coordinates coincide, that is, $Y(x)=Z(x)$.
Then the covectors $dz^j=dZ^j$ and $dy^j=dY^j$ coincide for $j=1,2,3,4$.
Then we have
\beq\label{eq: r is -theta}
\cirt r_j=-\cirt \theta_j,\quad\hbox{i.e.,\ } \cirt \theta:=\cirt \theta_jdz^j=-\cirt r=\cirt r_jdy^j\in T^*_{\cirt y}\hattuM _0.
\eeq


Let us apply
the method of stationary phase to $T^{(4)}_{\tau,1,1}$ as $\tau\to \infty$. 
Note that
as $c_1(z,\theta_1,\theta_2)$ is a product type symbol, we need to use the fact $\theta_1\not=0$ and $\theta_2\not=0$
for the critical points as we have by (\ref{eq: r is -theta}) 
and the fact that $\cirt r=d\varphi(y)|_{\cirt y}\not \in N^*K_{234}\cup N^*K_{134}$
as $x_5\not \in  {\mathcal Y}$ and assuming that $s_0$ is small enough.

{
In local $Y$ and $Z$ coordinates where $\cirt z=\cirt y=(0,0,0,0)$ we
can use the method of  stationary phase,
 similarly to the proofs of \cite[Thm.\ 1.11 and 3.4]{GrS},
 to compute 
the asymptotics of  (\ref{eq; T tau asympt}) with $k=1$.
Let us explain this  computation in detail.
To this end, let us start with some preliminary considerations.

Let $\phi_1(z,y,\theta,\xi)$ be a smooth bounded
 function in $C^\infty(W_0\times W_1\times \R^4\times \R^4)$ that is 
homogeneous of degree zero in the $(\theta,\xi)$ variables in the set  $\{|(\theta,\xi)|>R_0\}$ with some $R_0>0$. 
Assume that $\phi_1(z,y,\theta,\xi)$ is equal to one  in 
a conic neighborhood, with respect to $(\theta,\xi)$, of the points where some of the 
$\theta_j$ or $\xi_k$ variable is zero. Note that in this set 
the positively homogeneous  functions $\sigma_{u_j}(z,\theta_j/|\theta_j|)$ may be non-smooth.
Let $\Sigma\subset W_0\times W_1\times \R^4\times \R^4$ 
be a conic 
 neighborhood 
of the critical points of the phase function $\Psi_3(y,z,\theta,\xi)$.
Also, assume that the function  $\phi_1(y,z,\theta,\xi)$ vanishes in 
the intersection of $\Sigma$ 
and the set  $\{|(\theta,\xi)|>R_0\}$.
Let  
\ba
\Psi_{(\tau)}(z,y,\theta,\xi)= 
\theta_1z^1+\theta_2z^2+\theta_3z^3+(y-z)\cdotp \xi+\theta_4y^4+\tau \varphi(y)\ea 
 be the phase function appearing in (\ref{eq; T tau asympt pre}).
 Note that $\Sigma$ contains  the critical points of $\Psi_{(\tau)}$ for all $\tau>0$.
 Let
$$L_\tau=\frac {\phi_1(z,y,\tau^{-1}\theta,\tau^{-1} \xi)}{ |d_{z,y,\theta,\xi}
\Psi_{(\tau)}(z,y,\theta,\xi)|^{2}}\, (d_{z,y,\theta,\xi}
\overline{\Psi_{(\tau)}(z,y,\theta,\xi)})\cdotp d_{z,y,\theta,\xi},$$
so that 
\beq\label{kaav A}
L_\tau \exp(\Psi_{(\tau)})= \phi_1(z,y,\tau^{-1} \theta,\tau^{-1}\xi)\,\exp(\Psi_{(\tau)}).
\eeq
Note that if $\im \varphi(y)=0$, then $y\in \gamma_{x_5,\xi_5}$
and hence $d\varphi(y)$ does not vanish and that when $\tau$ is large enough,  the function $\Psi_{(\tau)}$ has 
no critical points in 
the support of $\phi_1$. Using
these we see that $ \gamma_{x_5,\xi_5}\cap W_1$
has a neighborhood $V_1\subset W_1$ where $|d\varphi(y)|>C_0>0$ and
there are $C_1,C_2,C_3>0$ so that if
  $\tau>C_1$, $y\in V_1$, and $(z,y,\theta,\xi)\in \supp(\phi_1)$
then
\beq\label{str est1}
\,|d_{z,y,\theta,\xi}
\Psi_{(\tau)}(z,y,\theta,\xi)|^{-1}\leq \frac {C_2}{\tau-C_3}.
\eeq
After these preparatory steps, we are ready to  compute the asymptotics of
 $T^{(4)}_{\tau,1,1}$. To this end, we first transform the integrals in  (\ref{eq; T tau asympt})
 to an integral of a Lebesgue integrable function by using integration by parts of $|\xi|^{-2}(\nabla_z-\omega_\theta)^2$
 as explained above. Then
 we decompose  $T^{(4)}_{\tau,1,1}$ into three
 terms $T^{(4)}_{\tau,1,1}=I_1+I_2+I_3$. To obtain the first term $I_1$ we include the factor
 $(1-\phi_1(z,y,\theta,\xi))$ in the  integral  
 (\ref{eq; T tau asympt}) with $k=1$. 
  The integral $I_1$ can then be computed
 using the method of stationary phase similarly to the proof of \cite[Thm.\ 3.4]{GrS}. 
 Let $\chi_1\in C^\infty(W_1)$ be a function that 
 is supported in $V_1$ and vanishes on $\gamma_{x_5,\xi_5}$. 
 The terms $I_2$ and $I_3$
 are obtained by 
 including the factor
 $\phi_1(z,y,\tau^{-1}\theta,\tau^{-1}\xi)\chi_1(y)$ and  $\phi_1(z,y,\tau^{-1}\theta,\tau^{-1}\xi)(1-\chi_1(y))$
  in the  integral 
 (\ref{eq; T tau asympt pre}) with $k=1$, respectively.
(Equivalently, 
 the terms $I_2$ and $I_3$
 are obtained by 
 including the factor
 $\phi_1(z,y,\theta,\xi)\chi_1(y)$ and  $\phi_1(z,y,\theta,\xi)(1-\chi_1(y))$
  in the  integral 
 (\ref{eq; T tau asympt}) with $k=1$, respectively.)
 Using integration by parts in the integral (\ref{eq; T tau asympt pre}) and inequalities  (\ref{kaav A}) and (\ref{str est1}), we see that that $I_2=O(\tau^{-N})$
 for all $N>0$. Moreover, the fact that $\im \varphi(y)>c_1>0$ in  
$W_1\cap V_1$ implies that $I_3=O(\tau^{-N})$
 for all $N>0$.}


 Combining the above we obtain
the asymptotics  
 \beq\label{definition of xi parameter1}
& &\quad \quad T^{(4)}_\tau
\sim \tau^{4+4-16/2-2+\rho-2}\sum_{k=0}^\infty c_k
\tau^{-k} =\tau^{-4+\rho}\sum_{k=0}^\infty c_k\tau^{-k}
, \\ \nonumber
& &\hspace{-8mm}c_0=h(\cirt z)c_1(0,-\cirt r_1,-\cirt r_2)\sigma_{u_3}(0,-\cirt r_3)
\hat q_1(0,0,-(\cirt r_1,\cirt r_2,\cirt r_3,0))\sigma_{u_4}(0,-\cirt r_4){{\Phi}}(0,1),\hspace{-1cm}
\eeq
where $\rho=\sum_{j=1}^5 p_j$
and $a_j$ is the principal
symbol of $U_j$ etc.   The factor $h(\cirt z)$ is non-vanishing and is
determined by the determinant of the Hessian of the phase function
$\varphi$ at $q$.
A direct computation shows that $\det(\hbox{Hess}_{z,y,\theta,\xi} 
\Psi_{3}(\cirt z,\cirt y,\cirt \theta,\cirt \xi))=1$. Above,
$(\cirt z,\cirt y,\cirt \theta,\cirt \xi)$ is the critical point satisfying
(\ref{eq: critical points}) and (\ref{eq: imag.condition}), 
where in the local coordinates $(\cirt z,\cirt y)=(0,0)$ and
$h(\cirt z)$ is constant times powers of values of the cut-off functions $\Phi_0$
and $\Phi_1$ at zero. 
 The term $c_k$ depends on the 
derivatives of the symbols $a_j$ and $q_1$ of order less
or equal to $2k$ at the critical point.
If  $ \Psi_3(z,\theta,y,\xi)$  has no critical points, that is,
 $q$ is not an intersection point of all five geodesics $\gamma_{x_j,\xi_j}$, $j=1,2,3,4,5$ we
obtain the asymptotics $T^{(4)}_{\tau,1,1}=O(\tau^{-N})$ for all $N>0$.

 For  future reference we note that if we use the method of  stationary
 phase in the last integral of (\ref{eq; T tau asympt}) only in the integrals with respect
 to $z$ and $\xi$, yielding that at the critical point we have  $y=z$ and $\xi=\omega_\beta(\theta)=(\theta_1,\theta_2,\theta_3,0)$, we 
 see that 
 $T^{(1)}_{\tau,1,1}^{(4)}$ can
 be written as
 \beq \label{eq; T tau asympt modified}
& &T^{(1)}_{\tau,1,1}^{(4)}
=c\tau^{4}
\int_{\R^{8}}e^{i(\theta_1y^1+\theta_2y^2+\theta_3y^3+\theta_4y^4)+i\tau \varphi(y)}
c_1(y, \theta_1, \theta_2)\cdotp\\ \nonumber
& &
\cdotp \sigma_{u_3}(y, \theta_3)
q_{1,1}(y,y,\omega_\beta(\theta))\sigma_{u_4}(y,\theta_4) {{\Phi}}(y,\tau )\,d\theta_1d\theta_2d\theta_3d\theta_4dy\\ \nonumber
& &=c\tau^{8}
\int_{\R^{8}}e^{i\tau (\theta_1y^1+\theta_2y^2+\theta_3y^3+\theta_4y^4+\varphi(y))}
c_1(y, \tau \theta_1, \tau \theta_2)\cdotp\\ \nonumber
& &
\cdotp \sigma_{u_3}(y,\tau  \theta_3)
q_{1,1}(y,y,\tau  \omega_\beta(\theta))\sigma_{u_4}(y,\tau \theta_4) {{\Phi}}(y,\tau )\,d\theta_1d\theta_2d\theta_3d\theta_4dy.
\eeq

Next, we consider the terms $\tilde T^{(4)}_\tau$ of the type
(\ref{tilde T-type source}). Such term
is an integral of the product of $u_\tau$ and two other factors 
${\bQ}(U_2\,\cdotp U_1)$ and ${\bQ}(U_4\,\cdotp U_3)$.
As  the last two factors can be written in the form (\ref{Q U1U2term}),
one can see using the method of stationary phase that  $\tilde T^{(4)}_\tau$ 
has similar asymptotics to  $T^{(4)}_\tau$ as $\tau\to \infty$,
with the leading order coefficient $\tilde c_0=\tilde h(\cirt z)\hat c_1(0,-\cirt r_1,-\cirt r_2)
\hat c_2(0,-\cirt r_3,-\cirt r_4)$ where $\hat c_2$ is given as in (\ref{product type symbols})
with symbols $a_3$ and $a_4$, and moreover,
$\tilde h(\cirt z)$ is a constant times powers of values of the cut-off functions $\Phi_0$
and $\Phi_1$ at zero.

This proves the claim in the special case where
$u_j$ are conormal distributions supported in the coordinate neighborhoods $W_{\kappa(j)}$,
$a_j$ are positively homogeneous scalar valued symbols.

By using a suitable partition of unity and 
summing the results of the  above computations, 
similar results to the above follows when  
$a_j$ are general classical symbols.

{\revtext
Let us also comment how the compuations go later for the Einstein equations.
There, the symbols \HOX{13.10.2014: This is modified, please check.}
 are $\mathcal B$-valued and 
the waves $u_j$ are supported on
$J^+_{\nohat g}(\supp({f}_j))$.
Also, $Q$ can be replaced by operators $S_j$ of type (\ref{extra notation}) and
the multiplication by function $a(x)$ can replaced by differential operators $\B_j$ without other essential changes
expect that the highest order power of $\tau$ changes.
Then, in the asymptotics of terms $T^{(4)}_\tau$  the function
$h(\cirt z)$ in (\ref{definition of xi parameter1}) is a section in  dual bundle $(\B_L)^4$. The coefficients
of $h(\cirt z)$ in local coordinates are polynomials of $\nohat g^{jk}$, $\nohat g_{jk}$, 
$\hat \phi_{\ell}$, and 
their derivatives at $\cirt z$. Similar representation is obtained for the asymptotics
 of the terms $\tilde T^{(4)}_\tau$.   

Let us next consider the scalar wave equation and the Einstein equations together.}
Above we integrated by parts two times the operator
$(\nabla_z-\omega_\theta)^2$. {\newcorrection For the scalar wave equation the multiplication by $a(x)$
are zeroth order operators and for the Einstein equation the total order of 
of $\B_j$ is less or equal to 6. Thus we see for both equations that it is
enough to assume above that the symbols $\sigma_{u_j}(z,\theta_j)$
are of  order $(-12)$ or less.
The leading order asymptotics come from the term
where and $p_j=n$ for $j=1,2,3,4$,
$p_5=0$, and for the Einstein equations,  the sum of orders of $\B_j$ is 6. This
is why we use  $m=-4-\rho+6=4n+2$ {\revtext for Einstein equations.
For the wave equation, we have $m=4n-4$.}
We also see that the terms containing permutation $\sigma$ of the indexes
of the distorted plane  waves can be analyzed analogously. This proves (\ref{indicator}).}

Note that above $\cirt r$  depends only on $\nohat g$ and ${\bf b}$ and we write
$\cirt r=\cirt r({\bf b})$, omiting the $\nohat g$ depdency.
Making the above computations explicitly, we obtain an explicit formula for
the leading order coefficient $s_{m}$ in (\ref{indicator})  in terms of
$\bsequence$ and ${\bf v}$.
This show that $s_{m}$
coincides with some real-analytic function $G(\bsequence,{\bf v})$.
%
%
This proves the claims (i) and (ii). 


(iii)  Let us fix $(x_j,\xi_j)$, $j\leq 4$, $s_0$, and the
 waves $u_j\in \I(K_j)$, $j\leq 4$.

 First we observe that if the condition (LI) is not valid, we see similarly to the above that
$\M^{(4)}$ is $C^\infty$ smooth in $\V\setminus (\Y\cup \bigcup_{j=1}^4 K_j)$.
Thus to prove the claim of the proposition we can assume that (LI) is valid.S

Let us decompose  $\S$, {\newcorrection given by (\ref{4th interaction for wave eq B}) for the wave
equation, and by
 (\ref{w1-3 solutions}) and (\ref{M4 terms})-(\ref{tilde M4 terms}) for the Einstein equation,}
 as  $\S=\S_{1}^{(4)}+\S_{2}^{(4)}$
 where
   $\S_{p}$ is defined similarly to
 $\S$ in (\ref{4th interaction for wave eq B}), or in (\ref{w1-3 solutions}) and 
 (\ref{M4 terms})-(\ref{tilde M4 terms})
by modifying these formulas as follows: \HOX{Here is a modification, please check.}
{\revtext For the wave equation ${\bQ}$ is replaced by ${\bQ}_p$.
 For the Einstein equation,}
 the operator $S_1^\beta$ is replaced by 
$S_{1,p}^\beta$, where $S_{1,p}^\beta={\bQ}_p$, when
 $S_{1}^\beta={\bQ}$, and $S_{1,p}^\beta=(2-p)I$, when $S_{1}^\beta=I$.
 Here, the operators ${\bQ}_p$  are defined as above using the parameters
  $\e_2$ and $\e_1$ defined below.
 
%

 Using formulas  (\ref{U1term}), (\ref{Q U1U2term}), and (\ref{QU4Ut term})
 we see that near $q$ in the $Y$ coordinates $\M_1^{(4)}={\bQ} \S_1^{(4)}$ 
 can be calculated using that
 \beq\label{f4-formula}
 \S_1^{(4)}(y)=\int_{\R^4} e^{iy^j\theta_j}b(y,\theta)\,d\theta,
 \eeq
 where $K_j$ in local coordinates is given by $\{y^j=0\}$ and
  $b(y,\theta)$  is a finite sum of terms  that are products of some of the following
 terms: at most one  product type symbol $c_l(y,\theta_j,\theta_k)\in S(W_0;\R\times (\R\setminus \{0\}))$ 
 and one ore more term which is either
 the symbols $\sigma_{u_j}(y,\theta_j)\in S^n(W_0;\R)$,
or the functions $q_1(y,y,\omega_\beta(\theta))$,
  cf.\  (\ref{eq; T tau asympt modified}),
 where $\omega_\beta(\theta)$ is equal to some of the vectors
$(\theta_1,\theta_2,\theta_3,0)$,
 $(\theta_1,\theta_2,0,\theta_4)$, $(\theta_1,0,\theta_3,\theta_4)$, or
 $(0,\theta_1,\theta_3,\theta_4)$, depending on the permutation $\sigma$.

Let us consider next the source $F_\tau$ is determined by the functions $(p,h)$ in (\ref{Ftau source}). 
 Then using the  method of  stationary phase gives  the asymptotics, c.f. (\ref{eq; T tau asympt modified}),
 \ba
 \bra u_\tau,\S_1^{(4)}\cet \hspace{-1mm}\sim \hspace{-1mm}\tau^8\hspace{-1mm}\int_{\R^8} e^{i\tau (\varphi(y)+y^j\theta_j)}({{\Phi}}(y,\tau),b(y,\tau \theta))_{\nohat g}d\theta dy
\sim \sum_{k=m}^\infty  
 s_k(p,h) \tau^{-k}
 \ea
where $\nohat g$ {\revtext for the wave equation product in $\R$ and for the Einstein equation it} is a Riemannian metric
on the fiber of $\B^L$ at $y$, that is isomorphic to $\R^{10+L}$, and
the critical point of the phase function is $y=0$ and $\theta=-d\varphi(0)$.
  As we saw above, we have that
 when $\e_2>0$ is small enough then for $p=2$ we have
$\bra u_\tau,\S_p^{(4)}\cet=O(\tau^{-N})$ for all $N>0.$
%

Let us choose sufficiently
small $\e_3>0$ and choose a function $\chi(\theta)\in C^\infty(\R^4)$ that 
vanishes in a $\e_3$-neighborhood (in the $\nohat g^+$ metric) of
$\mathcal A_q$, 
\beq\label{set Yq}
\mathcal A_q&:=&
N_q^*K_{123}
\cup N_q^*K_{134}\cup N_q^*K_{124}\cup N_q^*K_{234}
\eeq
and is equal to 1 outside the $(2\e_3)$-neighborhood of this set.

Let $\phi\in C^\infty_0(W_1)$ be a function that is one near $q$.
Also, let  \ba\hbox{$b_0(y,\theta)=\phi(y)\chi(\theta)b(y,\theta)$}\ea be a classical symbol,
$p=\sum_{j=1}^4 p_j$,
and let  $
\S^{(4),0}(y)
\in \I^{p-4}(q) 
$ be the
 conormal distribution 
 that is given by the formula (\ref{f4-formula}) with
$b(y,\theta)$ being replaced by 
 $b_0(y,\theta)$.

When $\e_3$ is small enough (depending on the point $x_5$), we have that 
$F_\tau$ is determined by functions $(p,h)$ 
and the corresponding 
gaussian beams $u_\tau$
propagating on the geodesic $\gamma_{x_5,\xi_5}(\R)$ such
that the geodesic
passes through 
$x_5\in  U$, we have 
\ba
\bra u_\tau,\S^{(4),0}\cet\sim \sum_{k=m}^\infty s_k(p,h) \tau^{-k},
\ea
that is, we have $\bra u_\tau,\S^{(4),0}\cet-\bra u_\tau,\S\cet=O(\tau^{-N})$
for all $N$.
When $\gamma_{x_5,\xi_5}(\R)$  does not pass through $q$,
we have that  $\bra u_\tau,\S\cet$ and $\bra u_\tau,\S^{(4),0}\cet$ 
are both of order $O(\tau^{-N})$ for all $N>0$.

Let $W\subset \ U((\vec x,\vec \xi),t_0)\setminus \bigcup_{j=1}^4\gamma_{x_j,\xi_j}([0,\infty))$, see   (\ref{eq: summary of assumptions 2}), be an open set.
By varying  the source $F_\tau$,  defined in  (\ref{Ftau source}),  
 we see, by multiplying the solution  with a smooth
 cut of function and  using  Corollary 1.4 in \cite{Delort} in local coordinates,
 or \cite{MelinS},
%
 we have
 that the function
$\M^4-{\bQ}\S^{(4),0}$ has no wave front set in $T^*(W)$ and it is
thus $C^\infty$-smooth function in $W$.

Since by \cite{GU1}, ${\bQ}:\I^{p-4}(\{q\})\to \I^{p-4-3/2,-1/2}(N^*(\{q\}),\Lambda_q^+)$, the above implies that
\beq\label{eq: conormal}
\M^4|_{W\setminus {\mathcal Y}}\in \I^{p-4-3/2}(W\setminus  {\mathcal Y};\Lambda_q^+),
\eeq
where $  {\mathcal Y}=  {\mathcal Y}((\vec x,\vec \xi),t_0,s_0)$. When $x_5$ if fixed, choosing $s_0$ to
be small enough, we obtain the claim (iii).
\hfill \Box \medskip
}
 
 \refHOX{Earlier Remark 3.2 on first order perturbation is removed.} 
 
Next we will show that $\mathcal G_{\bf g}(\vec\zeta)$ in  (\ref{eq: principal symbols of U4}) is not  vanishing
identically. 
{\gtext This implies that at the point of interaction of the four waves a {\mmmmtext spherical wave is produced in a generic case.}}

%

\generalizations
{{\bf Remark 3.6.A.} 
Similar considerations to the 
proof of  Prop.\ \ref{lem:analytic limits C}, can be used to show that $\M^3|_{(U\setminus \cup_j K_j)\cap I^-(x_5)}$
is a Lagrangian distribution. As an example of this, let us analyze the term
\ba
m={\bQ}(\B_3 u_3\,\cdotp v_{12}),\quad \hbox{where }v_{12}={\bQ}(\B_2 u_1\,\cdotp \B_2 u_3),
\ea
where $u_j\in \I^{p_j}(K_j)$ are solutions
of the linear wave equation.
Let $S=(\bigcup_j N^*K_j)\cup(\bigcup _{j,k}N^*K_{jk})$.
Next, we do computations in local coordinates $X:W\to \R^n$,  $x^j=X^j(x)$, such that $W\subset  I^-(x_6)$ 
and $K_j\cap W =\{x\in W;\ x^j=0\}$. 
As  $v_{12}=v_{12}^1+v_{12}^2\in \I(K_1,K_{12})+\I(K_2,K_{12})$, we can write the function
$f=\B_3 u_3\,\cdotp v_{12}$  in
$W$ as an oscillatory integral
\beq\label{def; f}
f(x)=\int_{\R^3}e^{i(\theta_1x^1+\theta_2x^2+\theta_3x^3)}b(x,\theta_1,\theta_2,\theta_3)d\theta_1 d\theta_2d\theta_3.
\eeq
Here, $b(x,\theta_1,\theta_2,\theta_3)\in S_{cl}^{p_1+p_2+p_3-2}(W;\R^3)$ away from
the set where some  some of the variables $\theta_j$ vanish and $\xi=\sum_{j=1}^3 \theta_jdx^j$
is light-like.
Let $ {\mmmtext \tilde \theta} (x,\theta_1,\theta_2,\theta_3)=\phi(x,\sum_{j=1}^3 \theta_jdx^j),$
where $\phi$ is  a cut-off function
that
vanishes  in a conic neighborhood $S$ 
and is one outside a suitable conic neighborhood of $S$.
Then 
 $c(x,\theta_1,\theta_2,\theta_3)= {\mmmtext \tilde \theta} (x,\theta_1,\theta_2,\theta_3)b(x,\theta_1,\theta_2,\theta_3)$ is a classical symbol. Writing 
$f=f_1+f_2$ where $f_1$ and $f_2$ have the representation 
(\ref{def; f}) with the symbol $b$ is replaced by $c$ and $b-c$, correspondingly.
Then $f_1\in \I^{p_1+p_2+p_3-2}(K_{123})= \I^{p_1+p_2+p_3-5/2}(N^*K_{123})$ and 
we see that ${\bQ}f_1\in
\I^{p_1+p_2+p_3-4}(Y;\tilde \Lambda),$  is a Lagrangian distribution on  $Y=I^-(x_5)\setminus ( \cup_j K_j)$,
where $\tilde \Lambda$ is the flow-out of $N^*K_{123}$,
and ${\bQ}f_2$ is a $C^\infty$-smooth in  $Y$.
Hence, 
\beq\label{3-singularity order}
\M^3|_{Y}\in \I^{\tilde p-10}(Y;\tilde \Lambda),\quad \tilde p=p_{1}+p_{2}+p_{3}.
\eeq
}

\subsubsection{Non-vanishing of the function $\mathcal G_{\bf g}(\vec\zeta)$ in a generic set}

\begin{proposition}\label{singularities in Minkowski space} 
Let ${\muutos {\mmmtext \hatzeta}}\in L^{+,*}_qM$  
and 
\ba
\mathcal Z_0({\muutos {\mmmtext \hatzeta}})&=&\{(\zeta_j)_{j=1}^4 \in (L^{*}_qM)^4;\ (\zeta_j)_{j=1}^4\hbox{ are linearly independent},\\
& &\sum_{j=1}^4\zeta_j={\muutos {\mmmtext \hatzeta}},\hbox{ and } {\muutos {\mmmtext \hatzeta}\not =\zeta_k}
\hbox{ for all } k=1,2,3,4\}.
 \ea
{\muutos Then $\mathcal Z_0({\muutos {\mmmtext \hatzeta}})$ is a real analytic manifold having several topological components and} $\mathcal G_g(\vec\zeta)$, given in (\ref{P per Qmod}), is non-vanishing for $\vec\zeta =(\zeta_j)_{j=1}^4$ in an open and dense subset of 
$\mathcal Z_0({\muutos {\mmmtext \hatzeta}})$.
\end{proposition}
%


\noextension{ \noindent{\bf Proof.} 
{\muutos By its definition,
 $\mathcal Z_0({\muutos {\mmmtext \hatzeta}})$ is a real analytic manifold having several topological components.} 
In the proof, we use in $T^*_qM$ a basis where the metric tensor $g$ is the standard Minkowski metric
$\diag(-1,1,1,1)$. {\mmmmtext
 Also, without loss of generality we can assume that  ${\mmmtext \hatzeta}=(1,1,0,0)$.} {\ftext Moreover,} we 
identify the space $T^*_qM$ with $\R^4$.
{\mtext Also, note that for $(\zeta_j)_{j=1}^4\in \mathcal Z_0({\muutos {\mmmtext \hatzeta}})$
all $\zeta_j$  are non-zero and hence
${\mmmtext \hatzeta}$ does not belong in the span of any three co-vectors $\zeta_j$, $j=1,2,3,4$.}

%
%
{\rbtext 
Denote
$\mathcal B=(L^{+,*}_qM)^4$.
Let ${\muutos {\mmmtext \hatzeta}}\in L^{+,*}_qM$ and 
\ba
\mathcal B_0({\muutos {\mmmtext \hatzeta}})&=&\{(b_j)_{j=1}^4 \in L^{+,*}_qM;\ (b_j)_{j=1}^4\hbox{ are linearly independent},\\
& &\ \hbox{there are $a_j\in \R\setminus \{0\}$ such that }\sum_{j=1}^4a_jb_j={\muutos {\mmmtext \hatzeta}}\hbox{ and}\\
& &
\ {\muutos {\mmmtext \hatzeta}\not =a_jb_j}, 
\hbox{ for all } j=1,2,3,4\}.
 \ea
 We {\muutos observe from this that  $\mathcal B_0({\muutos {\mmmtext \hatzeta}})$ is a real analytic manifold that} has several topological components but it is contained in {\mtext the} connected {\muutos real-analytic}
 manifold
 $\mathcal B$.

When $(b_j)_{j=1}^4\in \mathcal B_0({\muutos {\mmmtext \hatzeta}})$, let ${\alpha}_j =\a_j(\vec b,{\muutos {\mmmtext \hatzeta}})$, $j=1,2,3,4$ be such that
$$
{\muutos {\mmmtext \hatzeta}}=\sum_{j=1}^4{\alpha}_j b_j.
$$ 
{\mtext Since $b_j$ are linearly independent, $\a_j(\vec b,{\muutos {\mmmtext \hatzeta}})$  are uniquely determined.}
Considering $b_j$ as elements of $\R^4$ and using Cramer's rule, we obtain
  $$
 {\alpha_1}(\vec b,{\muutos {\mmmtext \hatzeta}})=\frac{\det({\muutos {\mmmtext \hatzeta}},b_{2},b_{3},b_{4})}{\det(b_{1},b_{2},b_{3},b_{4})},
 \quad  {\alpha_2}(\vec b,{\muutos {\mmmtext \hatzeta}})=\frac{\det(b_{1},{\muutos {\mmmtext \hatzeta}},b_{3},b_{4})}{\det(b_{1},b_{2},b_{3},b_{4})}.$$
 {\ftext Similar}  formulas hold for $ {\alpha_j}(\vec b,{\muutos {\mmmtext \hatzeta}})$ with $j=3,4$.
Using these formulas, we define 
$$
F(\vec b,{\muutos {\mmmtext \hatzeta}})=(\alpha_1(\vec b,{\muutos {\mmmtext \hatzeta}})b_1,\alpha_2(\vec b,{\muutos {\mmmtext \hatzeta}})b_2,\alpha_3(\vec b,{\muutos {\mmmtext \hatzeta}})b_3,
\alpha_4(\vec b,{\muutos {\mmmtext \hatzeta}})b_4).
$$
 Note that $\zeta_j={\alpha}_j b_j$ 
and ${\muutos {\mmmtext \hatzeta}}=\sum_{j=1}^4 \zeta_j$ so that for all $\sigma\in \Sigma(4)$,
\beq\label{3-trick}
\zeta_{\sigma(1)}+\zeta_{\sigma(2)}+\zeta_{\sigma(3)}={\muutos {\mmmtext \hatzeta}}-\zeta_{\sigma(4)}.
\eeq
{\muutos {\mmmmtext Then, if $\zeta_{\sigma(1)}+\zeta_{\sigma(2)}+\zeta_{\sigma(3)}$ would be light-like {\mmmmtext or zero}
then ${\muutos {\mmmtext \hatzeta}}-\zeta_{\sigma(4)}$  would be light-like  {\mmmmtext or zero} that is possible only if  
${\muutos {\mmmtext \hatzeta}}=\zeta_{\sigma(4)}$ {\mmmmtext and this can not happen  {\ftext since}  $(b_j)_{j=1}^4\in \mathcal B_0({\muutos {\mmmtext \hatzeta}})$.}}
We recall also that the inner product of two light-like vectors is zero if and only if
the vectors are parallel.
Then, we consider the function ${\mathcal G} _{\bf g}(\vec \zeta)$ given in
(\ref{P per Qmod}) with $\vec\zeta=F(\vec b,{\muutos {\mmmtext \hatzeta}})$. {\mmmmtext We} denote 
$\tilde {\mathcal G}(\vec b,{\muutos {\mmmtext \hatzeta}})={\mathcal G} _{\bf g}(F(\vec b,{\muutos {\mmmtext \hatzeta}}))$
 and
observe that 
$$
\tilde {\mathcal G}(\vec b,{\muutos {\mmmtext \hatzeta}})={P^{}{}(\vec b,{\muutos {\mmmtext \hatzeta}})}/{Q(\vec b,{\muutos {\mmmtext \hatzeta}})}
$$
where ${\mmmmtext \vec b\mapsto} P(\vec b,{\muutos {\mmmtext \hatzeta}})$ and
${\mmmmtext \vec b\mapsto}Q(\vec b,{\muutos {\mmmtext \hatzeta}})$ are
real analytic functions defined on the whole set $\mathcal B$ and
$Q(\vec b,{\muutos {\mmmtext \hatzeta}})\not =0$  for $\vec b\in  \mathcal B_0({\muutos {\mmmtext \hatzeta}})$.

%
%
%
Let us next show that $\tilde {\mathcal G}(\vec b,{\muutos {\mmmtext \hatzeta}})$ obtains a non-zero finite value at
some $(\vec b,{\muutos {\mmmtext \hatzeta}}).$
%
 Let
\beq\label{b distances}
  & &b_1=(1+\rhoepsilon_1^2,1-\rhoepsilon_1^2,2\rhoepsilon_1,0),\quad
 b_{2}=(1+\rhoepsilon_2^2,1-\rhoepsilon_2^2,0,2\rhoepsilon_2),\\ \nonumber
    & &b_{3}=(1+\rhoepsilon_3^2,1-\rhoepsilon_3^2,2\rhoepsilon_3,0),\quad
 b_{4}=(1+\rhoepsilon_4^2,1-\rhoepsilon_4^2,0,2\rhoepsilon_4),\\ \nonumber
& &{\muutos {\mmmtext \hatzeta}}=(1,1,0,0),
\eeq
where  $\rhoepsilon_j\in (0,1)$ are small parameters.
{\mtext Below in this proof,} we use the parameters $\rhoepsilon_j$ {\mmmmtext given by}
\beq\label{eq: ordering of epsilons wave equation}
\rhoepsilon_4=\rhoepsilon_3^{100},\  
\rhoepsilon_3=\rhoepsilon_2^{100},\hbox{ and  
}\rhoepsilon_2=\rhoepsilon_1^{100}.
\eeq
We denote $\vec\rhoepsilon\to 0$ when
$\rhoepsilon_1\to 0$ and}
$\rhoepsilon_2,$ $\rhoepsilon_3,$ and $\rhoepsilon_4$ are defined using
 (\ref{eq: ordering of epsilons wave equation}).
{Note that $\rhoepsilon_4<\rhoepsilon_3<\rhoepsilon_2<\rhoepsilon_1$. 
 
 The vectors $b_k$ 
are  light-like. {\mmmmtext  For small} $\rho_1$ we have $\vec b=(b_1,b_2,b_3,b_4)\in  \mathcal B_0({\muutos {\mmmtext \hatzeta}})$
 and 
\beq\label{BBB eq}
& &{\nohat g({\muutos {\mmmtext \hatzeta}},b_j)=-2 \rhoepsilon_j^2},\quad \nohat g( b_k,b_j)=
-2 (\rhoepsilon_k^2+ \rhoepsilon_j^2+O(\rhoepsilon_k\rhoepsilon_j)),\hspace{-1cm}
\eeq
for $j,k=1,2,3,4$.


{Below, we denote ${\alpha}_j =\a_j(\vec b,{\muutos {\mmmtext \hatzeta}})$ and see that}
\beq\label{formula to be corrected}
& & {\alpha_1}=\frac {\rho_2\rho_4}{\rho_1^2}(1+O(\rhoepsilon_1)),\quad {\alpha_2}=-\frac {\rho_4}{\rho_2}(1+O(\rhoepsilon_1)),\\ \nonumber
& & {\alpha_3}=-\frac {\rho_2\rho_4}{\rho_1\rho_3}(1+O(\rhoepsilon_1)),\quad
{\alpha_4}=1+O(\rhoepsilon_1).
\eeq

Then
$\tilde {\mathcal G}(\vec b,{\muutos {\mmmtext \hatzeta}})=\sum_{\sigma \in \Sigma(4)}H^{-1}(M_\sigma^1+M_\sigma^2)$, where 
%
%
$$
M_\sigma^1=\frac{C_1\a_{\sigma(3)}}{g({\muutos {\mmmtext \hatzeta}},b_{\sigma(4)})g(b_{\sigma(1)},b_{\sigma(2)})} \quad\hbox{and}\quad  M_\sigma^2=\frac{C_2}
{g(b_{\sigma(3)},b_{\sigma(4)})g(b_{\sigma(1)},b_{\sigma(2)})},
$$
and  $H=\prod_{j=1}^4 {\a}_j$.

When  {\mtext $\sigma(4) =4$, see that} the leading order asymptotics of $M^1_{\sigma}$ is given by
\beq\label{M1 vika}
& &M^1_{Id}\sim C_0\frac 1{\rho_4^2}\frac 1{\rho_1^2}\frac {\rho_2\rho_4}{\rho_1\rho_3}=C_0\rho_4^{-1}\rho_3^{-1}\rho_2^{+1}\rho_1^{-3},\\ \nonumber
& &M^1_{(1,3,2,4)}\sim C\frac 1{\rho_4^2}\frac 1{\rho_1^2}\frac {\rho_4}{\rho_2}=C\rho_4^{-1}\rho_3^{0}\rho_2^{-1}\rho_1^{-2},\\  \nonumber
& &M^1_{(3,2,1,4)}\sim C\frac 1{\rho_4^2}\frac 1{\rho_2^2}\frac {\rho_2\rho_4}{\rho_1}=C\rho_4^{-1}\rho_3^{0}\rho_2^{-1}\rho_1^{-2},\\ \nonumber
& &M^1_{(3,1,2,4)}\sim C\frac 1{\rho_4^2}\frac 1{\rho_1^2}\frac {\rho_4}{\rho_2}=C\rho_4^{-1}\rho_3^{0}\rho_2^{-1}\rho_1^{-2},\\ \nonumber
 &&M^1_{(2,3,1,4)}\sim C\rho_4^{-1}\rho_3^{0}\rho_2^{-1}\rho_1^{-2},\\  \nonumber
& &M^1_{\sigma_1}=M^1_{id},
\eeq
where $\sigma_1={(2,1,3,4)}$  and $C_0\not =0$.

{\mtext Using  the formulas for $\a_j$} given by  Cramer's rule, 
we  {\ftext have}  {\mtext that}  {\ftext the} terms $M^1_{\sigma}$ with $\sigma(4)\not =4$  {\vtext do} not contain  {\ftext the} factor $\rho_4^{j}$
with $j\leq -1$. Also $M^2_{\sigma}$  {\vtext do} not contain  {\ftext the} factor $\rho_4^{j}$
with $j\leq -1$. Hence,  $M^1_{\sigma_1}=M^1_{id}$ has the strongest {\mmmmtext asymptotics} 
when   $\vec\rho\to 0$.
More precisely, for all $\sigma\not \in \{Id,\sigma_1\}$, we have 
that $M^1_{\sigma}/M^1_{Id}\to 0$ as ${\mtext \vec\rho \to 0}.$
Also, for all $\sigma$ we have 
that $M^2_{\sigma}/M^1_{Id}\to 0$ as ${\mtext \vec\rho \to 0}.$
As $C_0\not=0$, this implies that
 \beq\label{last limit}
\hbox{$\tilde {\mathcal G}(\vec b,{\muutos {\mmmtext \hatzeta}})/M^1_{Id}\to 2$ as ${\mtext \vec\rho \to 0}.$}
\eeq

Recall that  $\tilde {\mathcal G}(\vec b,{\muutos {\mmmtext \hatzeta}})$ is a quotient of two real analytic functions
${P^{}{}(\vec b,{\muutos {\mmmtext \hatzeta}})}$ and ${Q(\vec b,{\muutos {\mmmtext \hatzeta}})}$.
Since for  {\ftext the} vectors given in 
(\ref{b distances}) with small $\rho_1$ we have
$\vec b\in  \mathcal B_0({\muutos {\mmmtext \hatzeta}})$,
 {\vtext we see using (\ref{M1 vika}) and (\ref{last limit})} that there is a point $\vec b\in  \mathcal B_0({\muutos {\mmmtext \hatzeta}}){\mmmmtext \subset \mathcal B}$
for which  $Q(\vec b,{\muutos {\mmmtext \hatzeta}})\not =0$  and   $P(\vec b,{\muutos {\mmmtext \hatzeta}})\not =0$,
that is, $Q(\vec b,{\muutos {\mmmtext \hatzeta}})$  and   $P(\vec b,{\muutos {\mmmtext \hatzeta}})$
are not
 identically vanishing. As  $\mathcal B$ is a connected, real-analytic manifold, 
we  {\ftext get}  that the functions  ${P^{}{}(\vec b,{\muutos {\mmmtext \hatzeta}})}$ and ${Q(\vec b,{\muutos {\mmmtext \hatzeta}})}$
 are not identically vanishing in 
any open subset of $\mathcal B$. As $\mathcal B_0({\muutos {\mmmtext \hatzeta}})\subset \mathcal B$ is open and dense, this implies that 
{\mtext ${P^{}{}(\vec b,{\muutos {\mmmtext \hatzeta}})}$ and ${Q(\vec b,{\muutos {\mmmtext \hatzeta}})}$, and thus {\mmmmtext also}
$\tilde {\mathcal G}(\vec b,{\muutos {\mmmtext \hatzeta}})$, are} non-vanishing
in an open and dense subset of $\mathcal B_0({\muutos {\mmmtext \hatzeta}})$, too. Since $F_{\muutos {\mmmtext \hatzeta}}:\mathcal B_0({\muutos {\mmmtext \hatzeta}})\to \mathcal Z_0({\muutos {\mmmtext \hatzeta}})$, $F_{\muutos {\mmmtext \hatzeta}}(\vec b)=F(\vec b,{\muutos {\mmmtext \hatzeta}})$,
is {\muutos an open, continuous and surjective map,} 
we  {\ftext conclude}  that
${\mathcal G}_{\bf g}(\vec \zeta)$ is non-vanishing
in an open and dense subset of $\mathcal Z_0({\muutos {\mmmtext \hatzeta}})$.
\hfill \Box \medskip
}



}

\subsection{Detection of singularities}\label{subsubsec: detection of singularities wave}

{\revtext We use now the  {\mtext above} results to  {\mtext  detect in the set $U$ the  singularities that are produced by the interaction of four waves}.}

{\vtext First we show that for all  $ q\in I^+(p^-)\cap I^-(p^+)$ there are
$(\vec x,\vec \xi)$ such that $q$ is
the intersection point of the geodesics corresponding to $(\vec x,\vec \xi)$.}

\begin{lemma}\label{lemma: existence of x,xi}
 {\mtext Let $ q\in I^+(p^-)\cap I^-(p^+)$. Then

(i) There are $({\mtext z},\zeta)\in L^+M$ and $0<r {\mmmmtext <}\rho({\mtext z},\zeta)$ such that $z\in 
I^+(p^-)\cap U$
and $q=\gamma_{{\mtext z},\zeta}(r)$. 

(ii) 
%
In any} neighborhood of $({\mtext z},\zeta)$ there are $(x_j,\xi_j)\in L^+U$, $j=1,2,3,4$   such that
{\vtext the geodesics corresponding to 
$(\vec x,\vec \xi)$ intersect regularly at $q$, see Def.\ \ref{def:  4-intersection of rays}. Moreover,
 we have that} $q\in \V(\vec x,\vec \xi)$, where  
 $\V(\vec x,\vec \xi)$ is defined in  (\ref{eq: summary of assumptions 2}),
and the points $x_j$  satisfy the
 condition (\ref{eq: summary of assumptions 1}).
%
\end{lemma}

\noindent
{\bf 
Proof.} {\mtext 
(i)  In the case when $q\in \hat \mu$,
let $\zeta_0\in L^+_q(M)$, and let $r_0>0$ be so small that
the geodesic $\gamma_{q,\zeta_0}([-r_0,0])\subset U\cap I^+(p^-)$ has no cut points.
Then, we define  $z=\gamma_{q,\zeta_0}(-r_0)$ and $\zeta=\gamma_{q,\zeta_0}(-r_0)$.

In the case when  $q\not \in \hat \mu$,
let $z_1=\hat \mu(f^-_{\hat a}(q))\in U\cap I^+(p^-)$.
By Lemma \ref{B: lemma} (iii) there is
$\zeta_1\in L^+_{z_1} M$
such that $\gamma_{z_1,\zeta_1}$ is one of the longest light-like {\mmmmtext geodesics} connecting $y$  and $q$,
$q=\gamma_{z_1,\zeta_1}(r_1)$  where  $0<r_1\leq \rho({\mtext z}_1,\zeta_1)$.
Then, let $r_{\mmmmtext 2}>0$ be so small that
$\gamma_{z_1,\zeta_1}([0,r_{\mmmmtext 2}])\subset U$,
and let   $z=\gamma_{z_1,\zeta_1}(r_{\mmmmtext 2})$ and $\zeta=\gamma_{z_1,\zeta_1}(r_{\mmmmtext 2})$.
 
 In both cases $(z,\zeta)\in L^+U$ has the properties required in (i).}

(ii) Note that
$q=\gamma_{{\mtext z},\zeta}(r)\in {\mmmmtext I}^-({p^+})$ with ${\muutos 0}<r\leq \rho({\mtext z},\zeta)$. {\mtext Let
$ {\rtext \rho_0}\in (0,r)$ be such that 
$\gamma_{{\mtext z},\zeta}([0, {\rtext \rho_0}])\subset U$.} 
Let $W\subset TM$ be  a 
neighborhood of $({\mtext z},\zeta)$ {\muutos such that $\pi(W)\subset U$}.

Let $\theta=-\dot\gamma_{{\mtext z},\zeta}(r)\in T_qM$. 
Then the geodesic
$\gamma_{{\mtext z},\zeta}([ {\rtext \rho_0},r])=\gamma_{q,\theta}([0,r- {\rtext \rho_0}])$ has no
cut points. 
 {\vtext Thus,}  consider four  geodesics that emanate  from $q$ to the past, 
in the light-like direction $\eta_1=\theta$ and in the light-like directions
 $\eta_j\in  T_qM$,  $j=2,3,4$ that are close to the
direction $\theta$, {\mtext such that {\mmmmtext $(\eta_j)_{j=1}^4$}  are linearly independent.}
Let  
$\gamma_{q,\eta_j}(r_j)$ be the intersection points  of $\gamma_{q,\eta_j}$ with
the surface $\{x\in M;\ {\bf t}(x)=c_0\}$, on which the time function  ${\bf t}(x)$ 
has the constant value 
$c_0:={\bf t}({\muutos {\mtext z}})$, {\mmmmtext see subsection \ref{Observation domain}}.
Let $x_j=\gamma_{q,\eta_j}({\muutos r_j})$ and
$\xi_j=-\dot\gamma_{q,\eta_j}({\muutos r_j})$.
 When $\eta_j$ are sufficiently close to $\eta_1$,
such $r_j=r_j(\eta_j)$ exist by the inverse function theorem, 
we have that $\rho{\rtext (x_j,\xi_j})\vtext{>}r_j$,  $(x_j,\xi_j)$ are in $W$.
Then the obtained points $((x_j,\xi_j))_{j=1}^4$  {\mtext satisfy
condition (\ref{eq: summary of assumptions 1}) and  geodesics corresponding to 
$((x_j,\xi_j))_{j=1}^4$ intersect regularly at $q\in \V(\vec x,\vec \xi)$}. 
\hfill \Box \medskip

 We say
that the {\it interaction condition} (I) is satisfied for $y\in U$ with light-like
vectors $(\vec x,\vec \xi)=((x_j,\xi_j))_{j=1}^4$ 
{\mtext and} parameters  $(q,{\muutos w},t)$,  
 if
 \medskip


%
%

\noindent
{\bf (I)} {\mtext {\it There exist
 $q\in \bigcap_{j=1}^4\gamma_{ {\mtext x_j,\xi_j}}((0,{\bf t}_j))$,  where ${\bf t}_j=\rho {\mtext (x_j,\xi_j)}$, and
 ${\muutos w}\in L^+_q\hattuM _0$ and $t\geq 0$ 
 such that $y=\gamma_{q,{\muutos w}}(t)$.}}


 \medskip


Below, in  $T\hattuM _0$ 
we use the Sasaki metric corresponding to the Riemannian metric $\nohat g^+$.
\rtext Moreover, let  
 $B_j\subset  U$ be open sets
such that, 
cf.\  
(\ref{eq: source causality condition}),  
  we have 
  \beq\label{eq: admissible Bj}
& &\hbox{ 
$\overline B_j\subset  U$ and $B_j\cap J ^+(B_k)=\emptyset$ for all 
$j\not =k$.}
\eeq
}



{\revtext Next we formulate 
 a condition ($D_{}$) that 
{\gtext is valid  when we can detect
singularities at {\mtext a point} $y\in U$.}
 {We say that a function $v(x)$ is $C^\infty$-smooth at $y$ if there is an open neighborhood $W$ of $y$  such that
 $v|_W\in C^\infty(W)$.} 
 
  We define that point $y\in U$
   satisfies the singularity {\it detection condition} ($D_{}$)  
with light-like directions $(\vec x,\vec \xi)$
   and $ \hat s>0$
  if 
  \medskip  
  
  \noindent
(${\bf D}_{}$) {\it For any $s_0,{\mtext s_1}\in (0,\hat s)$ and   $j=1,2,3,4$ {\mtext and sufficiently large $n$}
{\ftext there exist $(x_j^{\prime},\xi_j^{\prime})$
 in the
{$s_1$-neighborhood} of $(x_j,\xi_j)$, {open sets $B_j\subset B_{\nohat g^+}({\rtext x_j},{\mtext s_1})$ satisfying (\ref{eq: admissible Bj}),
and source functions ${f}_{j}\in {\mathcal  
I}^{n+1}({\muutos \Sigma(x_j^{\prime},\xi_j^{\prime},s_0)})$, $\supp(f_j)\subset B_j$,} such  
that the following holds:  
When  $u_{\vec \e}$ is the solution  of the non-linear wave equation (\ref{PABC eq})
with the source ${f}_{\vec \e}=\sum_{j=1}^4
\e_j{f}_{j}$, then
the function
$\p_{\vec \e}^4u_{\vec \e}|_{\vec \e=0}$ is not
$C^\infty$-smooth at $y$.}}
  \medskip


\begin{lemma}\label{lem: sing detection in in normal coordinates 2 wave}
Let $(\vec x,\vec \xi)$,   and ${\bf t}_j$ with $j=1,2,3,4$ 
satisfy (\ref{eq: summary of assumptions 1})-(\ref{eq: summary of assumptions 2}).
%
{\mtext  
Assume  that   $y\in   {\rtext \V(\vec x,\vec \xi)}\cap U$ satisfies $y\not \in 
 {\mathcal Y}(\vec x,\vec \xi)\cup\bigcup_{j=1}^4
 \gamma_{x_j,\xi_j}(\R)$.} 
Then


{\vtext (i) If the geodesics corresponding to $(\vec x,\vec \xi)$ either do not intersect in $\V(\vec x,\vec \xi)$,
or the geodesics intersect at a point $q\in \V(\vec x,\vec \xi)$ and $y\not \in J^+(q)$,
then $y$ does not satisfy
condition ($D_{}$) with  $(\vec x,\vec \xi)$ and  any $\hat s>0$.}

(ii) Assume  
$y\in U$ satisfies condition (I) with  $(\vec x,\vec \xi)$ and {\mtext the} parameters $q,{\muutos w}$, and  $0<t<\rho(q,{\muutos w})$.
Then $y$ satisfies condition ($D_{}$) with  $(\vec x,\vec \xi)$  {\gtext for} any sufficiently
small $\hat s>0$. 

(iii) Using {\mmmmtext the source-to-solution operator} $L_{V}$ one can determine
whether the condition ($D_{}$) is valid  for {\ftext the given $y\in U$, $(\vec x,\vec \xi)$ and $\hat s$.}

   \end{lemma}

{\revtext
\noindent {\bf Proof.}
(i) {\vtext Assume that $y\in U$ satisfies the conditions stated in (i) and}
let $s_0,$ $s_1$, $(x_j^\prime,\xi_j^\prime)$, $B_j$, 
and $f_{j}\in  
\I^{n+1}({\muutos \Sigma(x_j^{\prime},\xi_j^{\prime},s_0)})$ be as in condition
($D_{}$).  {\muutos When {\mtext $s_0$ and $s_1$ are small enough, 
we see that 
{\vtext 
$y\in \V(((x_j^\prime,\xi_j^\prime))_{j=1}^4)$ and if the geodesics corresponding to
$((x_j^\prime,\xi_j^\prime))_{j=1}^4$  intersect at a point $q'\in \V(((x_j^\prime,\xi_j^\prime))_{j=1}^4)$,
then $y\not \in J^+(q')$. 
Hence,} the point $y$, the vectors
$(x_j^\prime,\xi_j^\prime)$ and the sources $f_j$ satisfy the assumptions of the claim (i) in 
Thm.\ \ref{lem:analytic limits A wave}.}

Let
 $u_{\vec \e}$ be the solution  of (\ref{PABC eq})
with  ${f}_{\vec \e}=\sum_{j=1}^4
\e_j{f}_{j}$ and
$\M^{(4)}=\p_{\vec \e}^4u_{\vec \e}|_{\vec \e=0}$.
 Then Theorem \ref{lem:analytic limits A wave} (i) implies that
$\M^{(4)}$ is  $C^\infty$-smooth at   
$y$. {\muutos Thus $y$ does not satisfy condition ($D_{}$).}

(ii) Let $y\in U$ satisfy condition (I) with  $(\vec x,\vec \xi)$ and {\mtext the} parameters $q,{\muutos w}$, 
and  $0<t<\rho(q,{\muutos w})$, {\muutos so that $\gamma_{q,{\muutos w}}(t)=y$.} {\vtext Let $s_0,s_1>0$.} 
{\vtext Note that then $y\in \firstpoint^{reg}_U(q)$.
Let ${\mmmmtext \hatzeta}={\muutos w}^\flat\in L^{*,+}_qM_0$.} 
For $j=1,2,3,4$, let 
 $t_j>0$ be such that 
$\gamma_{x_j,\xi_j}(t_j)=q$.
{\muutos Let $b_j=\dot \gamma_{x_j,\xi_j}(t_j)^\flat$. We can make {\mtext an arbitrarily small perturbation to the co-vectors $b_j$
to obtain co-vectors 
$\hat b_j\in L^{*,+}_q{M}$ such} that
$(\hat b_j)_{j=1}^4$ are linearly independent and
${\mmmmtext \hatzeta}$ is not in the space spanned by any three {\mmmmtext of the} vectors $b_j$, $j=1,2,3,4$.
Then, there are $\alpha_j\in \R\setminus \{0\}$ such that $\hat \zeta_j=\alpha_j\hat b_j$  satisfy
$\sum_{j=1}^4 \hat \zeta_j={\mmmmtext \hatzeta}$.
Then,
$(\hat \zeta_j)_{j=1}^4 \in \mathcal Z_0({\muutos {\mmmmtext \hatzeta}})$.
%
%
Furthermore, by using
Prop.\ \ref{singularities in Minkowski space}, we see that  {\mmmmtext
there are arbitrarily small perturbations
$\zeta^{\prime}_j\in L^*_q{M}$ of the co-vectors  ${\hat \zeta_j}$,} such that 
$(\zeta^{\prime}_j)_{j=1}^4 \in \mathcal Z_0({\muutos {\mmmmtext \hatzeta}})$
and
$\mathcal G_{\bf g}((\zeta^{\prime}_j)_{j=1}^4)\not =0$. Let {\mmmmtext ${\mtext b^{\prime}_j}= \frac 1{\alpha_j}\zeta^{\prime}_j$,  $j=1,2,3,4$  and
 $x_j^\prime=\gamma_{q,(b^{\prime}_j)^\sharp}(-t_j)$ and  $\xi_j^\prime=\dot \gamma_{q,(b^{\prime}_j)^\sharp}(-t_j)$.}}
{\itext 
When the perturbations above are small enough, we have that $x_j^\prime \in B_{\nohat g^+}({x_j},s_1/2)$ and {\ftext the} points $x_j^\prime$ satisfy
conditions (\ref{eq: summary of assumptions 1}) {\mtext and $\rho(x_j^\prime,\xi_j^\prime)>t_j$.}
Then,
using Lemma \ref{lem: Lagrangian 1 wave} we  {\ftext obtain}  that {\mtext when $n\in \Z_+$ is sufficiently large,}
  there are $f_{j}\in  
\I^{n+1}({\muutos \Sigma(x_j^{\prime},\xi_j^{\prime},s_0)})$
for which the principal symbol of ${{{\muutos \square_g^{-1}}}}f_j$ at 
{\mtext $(q,{\mmmmtext \zeta'_j})$} are non-vanishing
and $f_j$ are supported in sets  $B_j\subset B_{\nohat g^+}({\rtext x_j},{\mtext s_1})$ satisfying (\ref{eq: admissible Bj}).}}
By {\rbtext Thm.\ \ref{lem:analytic limits A wave} {\mmmmtext (ii)}, 
  the function  $\M^{(4)}$ is not 
  $C^\infty$-smooth {\mtext at $y$.  Hence, {\vtext as $s_0,s_1>0$ are above arbitrary,} we  {\ftext conclude}  that condition ($D_{}$) is} valid.
 
 (iii) 
 The {\mtext non-linear source-to-solution} map $L_U$ determines the functions
 $\M^{(4)}=\p_{\vec \e}^4u_{\vec \e}|_{\vec \e=0}$ in $ U$. This yields (iii).
 \hfill \Box \medskip
}

\canberemoved{
{\gtext 

\noindent
{\bf Remark 3.3.} 
Let us consider the case when the manifold $({\rbtext M},g)$ has no cut points.
Let $U^\prime= U\cap (I^+(p^-)\cap I^-(p^+))$.
Using Lemma \ref{lem: sing detection in in normal coordinates 2 wave} we see that 
the geodesics corresponding 4-tuples $(\vec x,\vec \xi)\in (L^+U^\prime)^4$ 
 intersect at some point $q\in I^+(p^-)\cap I^-(p^+)$ if and only if
 condition (D) is valid for some point $y\in \hat \mu((s_-,s_+))$.
Thus, by Lemma \ref{lem: sing detection in in normal coordinates 2 wave}
the pair $(U,g|_U)$ and the source-to-solution map $L_U$ determine 
all 4-tuples $(\vec x,\vec \xi) \in (L^+U^\prime)^4$ for which the corresponding geodesics
 intersect at some point $q\in I^+(p^-)\cap I^-(p^+)$.
Moreover, for such 4-tuples $(\vec x,\vec \xi)$ these data determine all points $y\in U$
that satisfy condition (I)  with  $(\vec x,\vec \xi)$ and  some $ {\rtext \rho_0}$, that is,
all points $y$ in  the earliest light observation set $\E_U(q)$. Thus by 
varying the 4-tuples $(\vec x,\vec \xi)$ we can determine  sets $\E_U(q)$ for all 
$q\in I^+(p^-)\cap I^-(p^+)$.  Hence, in
 the case when the manifold $({\rbtext M},g)$ has no cut points
 we have already reduced the proof of Theorem  \ref{main thm3 new} to proving Theorem \ref{main thm paper part 2}.}}

\section{Determination of the earliest light observation sets }\label{subsection combining}

  
  \removedEinstein{
We will do considerations for the non-linear wave equation and for the Einstein-scalar field
 equations. Due to this, we assume that we are given the data set ${\mathcal D}$
 that is  the set ${\mathcal D}_{wave}$ for the wave equation, see (\ref{eq: wave data set}),
 or ${\cal D}(g,\hat \phi,\e)$ for the Einstein equations, see (\ref{eq: main data}).
 Also, we denote conditions $(D_{wave})$ and $(D_E)$ by $(D)$.
 }

{\mattitext In this section we 
reduce {\gtext the proof of} Theorem  \ref{main thm3 new} to proving Theorem \ref{main thm paper part 2} that is
{\mtext proven} later in Section
\ref{sec: Proof for passive observations}.
\canberemoved{As explained in  Remark 3.3, we have already done this reduction in Section 3 in the case when
the manifold $(M,g)$ has no cut points. Next we consider the general case which
leads to some technical considerations involving the cut points and light-like 
geodesics that may intersect several times.}
Below, we assume that we are given 
 $(U,g|_U)$ and the {\slava source-to-solution} map $L_U$.
%
%
%

 \subsection{Surfaces of the earliest singularities} 
  
  Next we consider the determination of the
  the earliest light observation sets $\E_U(q)$, see Def.\ \ref{def. O_U}. To this end we need the following
  notation:
  
  
 \begin{definition}\label{def: first observations}
%
For a closed set $S\subset U$, we define the earliest points of set $S$ on 
the  {\muutos path} $\mu_{a}=\mu_{a}([-1,1])$ to be 
\beq\label{Earliest element sets}
& &\pointear_{a}(S)\hspace{-1mm}=\hspace{-1mm}\{\mu_{a}(\inf\{s\in [-1,1];\ \mu_{a}(s)\in S\})\},
\hbox{ if }\mu_{a}\cap S\not =\emptyset, \hspace{-10mm}\\
\nonumber 
& &\pointear_{a}(S)\hspace{-1mm}=\hspace{-1mm}\emptyset,
\quad\hbox{if }\mu_{a}\cap S =\emptyset.
\eeq

\end{definition}

  Note that by  the above definition, $\mu_{a}(f^+_{a}(q))=\pointear_{a}(\P_U(q))$  {\mtext for $q\in J(p^-,p^+).$}


Our next  
aim is to consider the global problem of constructing the set of the earliest light observations in $U$ of all points $q\in J(p^-,p^+)$. To this end, we need to handle the technical problem that  
{\vtext in the set $ {\mathcal Y}(\vec x,\vec\xi)$, see {(\ref{Y-set}), we  
have not analyzed if we observe singularities. Moreover, we have not analyzed the wave $\U^{(4)}$ in the set $J^+(q)$, if
{\vtext the geodesics intersect at $q$  and} the velocity vectors of the geodesics
at the point $q$ are not linearly independent, or in the set $M_0\setminus \V(\vec x,\vec \xi)$ that 
is the causal future of first conjugate points of the geodesics $\gamma_{x_j,\xi_j}$. As discussed in Remark 3.2,
the waves created by the non-linear interaction can be very complicated in these sets.}
To avoid these  difficulties, 
 we make the following definition} }  
    

  \begin{definition}\label{def: Sreg}  
Let $(\vec x,\vec \xi)=((x_j,\xi_j))_{j=1}^4$ be a collection of light-like  
vectors with $x_j\in U$.
 We define 
 \ba
 \S(\vec x,\vec \xi)=\{ y\in U&;& \hbox{\ftext there is $\hat s>0$ such that the condition ({D}) }
 \\
& & \hbox{\ftext is valid for $y$, $(\vec x,\vec \xi)$  and $\hat s$}\}.
 \ea  
 Moreover, {\vtext let   
$\S_{H} (\vec x,\vec \xi)$ be the set  
of such points  $y_0\in U$ that for every neighborhood $W\subset U$
of $y_0$  
{\vtext the Hausdorff dimension of the intersection $W\cap \S(\vec x,\vec \xi)$  
is at least  3.}  Note that $\S_{H} (\vec x,\vec \xi)$ is closed in the relative topology of $U$.
We denote   
(see (\ref{Earliest element sets}) and Def.\ \ref{def. O_U}) 
\beq  \label{EE1}
& &{\Scle}(\vec x,\vec \xi)\hspace{-0.5mm}=
\bigcup
 _{a\in \A}\pointear_{a}(\S_{H} (\vec x,\vec \xi)).
\eeq
We call ${\Scle}(\vec x,\vec \xi)$ the surface of the earliest stable singularities
produced by the interaction of  four waves.}
\end{definition}

%

   \begin{lemma}\label{lem: aux construction}
The {\vtext path $\hat \mu:[-1,1]\to U$, the} manifold $(U,g|_U)$ and $L_U$  determine the set
${\Scle}(\vec x,\vec \xi)$ for all  $(\vec x,\vec \xi)\in (L^+ U)^4 $.
{\rtext Moreover, these data determines the sets $\E_U({\mtext p})$ for all $p\in U$,
 the causality relation $R_U^<=\{(p_1,p_2)\in  {\mtext U\times U}
:\ p_1<p_2\}$, where $<$ is the causality
relation of $(M_0,g)$, {\mmmtext and the set}}
\beq\label{The set G}
G_{\hat\mu}=\{(x,\xi)\in L^+M_0;\ x\in U,\gamma_{x,\xi}(\R_+)\cap \hat \mu\not=\emptyset\}.
\eeq
   \end{lemma}
   
 \noindent
 {\bf Proof.} {\rtext 
 The first claim follows from
 {\mattitext  Lemma \ref {lem: sing detection in in normal coordinates 2 wave}
 and Def.\ \ref{def: Sreg}. 
 
 The derivative of the map $L_U$ at zero is the map $DL_U|_0:f\mapsto Qf|_U$,
defined for the distributions $f$ supported in $U$.
 {\mtext This map coincides with the source-to solution map for the linearized wave equation.}
 {\mmmtext By \cite[Thm.\ 8.1.4]{H1}, for any
 $(x,\xi)\in L^*M$  there is a distribution $f_1$  such that WF$(f_1)$ is
 the {\mmmmtext half-line} $\{(x,s\xi)\in T^*M;\ s>0\}.$ Then 
 by using H\"ormander's theorem on {\ftext propagation of singularities} along bicharacteristics, \cite[Theorem 26.1.4]{H4}, 
 we see that the singular support of $DL_U|_0(f_1)$ is $\gamma_{x,\xi}([0,\infty))\cap U$.
 Thus, we can determine the set $G_{\hat\mu}$.
%
%
%
%
{\mmmmtext Also,} for $p\in U$ we can find} the {\mmmtext sets $\mathcal L^+({\mtext p})$ and} $\pointear_{a}(\mathcal L^+({\mtext p}))$ and 
{\mmmmtext the} values
$f^+_a(p)$ for $a\in \A$. These determine the sets $\E_U({\mtext p})$ and $J^+(p)\cap U$ for all $p\in U$.
The latter determines all pairs $(p_1,p_2)\in U^2$ such that $p_1<p_2$.}}
\hfill \Box \medskip  
   

\medskip

\medskip


{\mattitext {\gtext Next we show} that
if the geodesics emanating from $(\vec x,\vec \xi)$ intersect before
their first cut points at $q_0$ then
 ${\Scle}(\vec x,\vec \xi)$ coincides with the set $\E_U(q_0)$, see
Def.\ \ref{def. O_U}.} 

  \extension{
  
Our goal is to show that ${\Sclo}(\vec x,\vec \xi)$ coincides with the intersection   
of the light cone $\L^+_{g}(q)$ and $U$ where $q$ is the intersection  
point of the geodesics corresponding to $(\vec x,\vec \xi)$, see Fig.\ 4(Left).

Let us next motivate the analysis we do below:   
We will consider how to create an artificial point source using interaction  
of distorted plane  waves propagating along light-like geodesics $\gamma_{x_j,\xxi_j}(\R_+)$  
where $(\vec x,\vec \xxi)=((x_j,\xxi_j))_{j=1}^4$ are perturbations  
of a light-like $(y,\zeta)$, $y\in \hat \mu$. 
We will use the fact that for all $q\in J(p^-,p^+)\setminus \hat \mu$ there  
is a light-like geodesic $\gamma_{y,\zeta}([0,t])$ from $y=\hat \mu(f^-_{\hat a}(q))$  
to $q$ with $t\leq \rho(y,\zeta)$.  
We will next  
show that   
when we choose $(x_j,\xxi_j)$ to be suitable   
perturbations of $\p_t \gamma_{y,\zeta}(t_0)$, $t_0>0$,  
it is possible that all geodesics $\gamma_{x_j,\xxi_j}(\R_+)$ intersect at $q$ before   
their first cut points,  
that is, $\gamma_{x_j,\xxi_j}(r_j)=q$, $r_j<\rho(x_j,\xxi_j)$.   
We note that we cannot analyze the interaction of the waves  
if the geodesics intersect after the cut points as then  
 the distorted plane   
waves can have caustics. Such interactions of wave caustics can, in principle,  
cause propagating singularities. Thus the sets ${\Sclo}((\vec x,\vec \xi),t_0)$  
contain singularities propagating along the light cone $\L^+_{g}(q)$  
and in addition that they may contain singularities produced by caustics that we do not know how to analyze (that could  
be called ``messy waves''). Fortunately, near an open and dense set of geodesics  
$\mu_{a}$ the nice singularities propagating along the light cone $\L^+_{g}(q)$  
arrive before the ``messy waves''. This is the reason why we consider below the  
first observed singularities on geodesics  
$\mu_{a}$.  
 Let us now return to the rigorous analysis.  
 
Below in this section we fix $t_0$ to have the value  
 ${t_0}=4\kappa_1$, cf.\ Lemma \ref{lem: detect conjugate 0}. 
 Recall the notation that  
\ba  
& &(x({t_0}),\xxi({t_0}))=(\gamma_{x,\xxi}({t_0}),\dot \gamma_{x,\xxi}({t_0})),\\  
& &(\vec x({t_0}),\vec \xxi({t_0}))=((x_j({t_0}),\xxi_j({t_0})))_{j=1}^4.  
\ea  
\medskip

\noindent
{\bf Lemma 5.1.A}  
{\it
Let $\vartheta>0$ be arbitrary, $q\in J(p^-,p^+)\setminus {\hat \mu}$ and   
let $y=\mu_{\hat a}(f_{\hat a}^-(q))$  
and $\zzeta\in L^+_y\hattuM _0$, $\|\zeta\|_{g^+}=1$ be such that $\gamma_{y,\zzeta}([0,r_1])$, $r_1>t_0=4\kappa_1$ is   
a longest causal (in fact, light-like) geodesic connecting $y$ to $q$.   
Then there exists  a set $\mathcal G$ of 4-tuples   of light-like vectors $(\vec x,\vec \xxi)=((x_j,\xxi_j))_{j=1}^4$   
 such that the points  $x_j$  
 and the directions $\xxi_j$ and the points $x_j(t_0)=\gamma_{x_j,\xxi_j}(t_0)$
  have the following properties:  
\begin{itemize}  
\item[(i)] $x_j(t_0)\in U$, $x_j(t_0)\not \in J^+(x_k(t_0))$ for $j\not =k$,  
  
\item[(ii)] $d_{g^+}  
((x_l,\xxi_l),(y,\zeta))<\vartheta$ for $l \leq 4$,  
  
\item[(iii)] $q=\gamma_{x_j,\xxi_j}(r_j)$ and $\rho(x_j(t_0),\xxi_j(t_0))+t_0>r_j$,  
  
\item[(iv)] when $(\vec x,\vec \xxi)$ run through the set  $\mathcal G$,  
the directions  
 $ (\dot \gamma_{x_j,\xxi_j}(r_j))_{j=1}^4$ form an open  
 set in $(L^+_qM)^4$.  
 
 \end{itemize}  
   
In addition, $\mathcal G$ contains elements $(\vec x,\vec \xxi)$ for which $(x_1,\xxi_1)=  
(y,\zzeta)$.  
}
 \medskip
 
 \noindent 
{\bf Proof.}   
Let  $\eta=\dot \gamma_{y,\zzeta}(r_0)\in L^+_q\hattuM _0$.  
Let us choose light-like directions  $\eta_j\in T_q\hattuM _0$, $j=1,2,3,4,$   
close to $\eta$ so that $\eta_j$ and $\eta_k$ are not parallel for $j\not  =k$.  
In particular, it is possible (but not necessary) that   
$\eta_1=\eta$. Let ${\bf t}:M\to \R$ be a time-function on $M$ that can be  
used to identify $M$ and $\R\times N$. Moreover, let  us choose  
 $T\in ((\gamma_{y,\zeta}(r_0-t_0)),{\bf t}(\gamma_{y,\zeta}(r_0-2t_0)))$ and  for $j=1,2,3,4$, let $s_j>0$   
be such that   
${\bf t}(\gamma_{q,\eta_j}(-s_j))=T$. Choosing  first $T$ to be sufficiently close ${\bf t}(\gamma_{y,\zeta}(r_0-t_0))$  
and then all $\eta_j$, $j=1,2,3,4$ to be sufficiently  
close to $\eta$     
 and defining  
$x_j=\gamma_{q,\eta_j}(-s_j-t_0)$ and $\xxi_j=\dot\gamma_{q,\eta_j}(-s_j-t_0)$  
we obtain the pairs $(x_j,\xxi_j)$ satisfying the properties stated in the claim.  
Indeed, this follows from the fact that $\rho(x,\xi)$ is lower semicontinuous function,
$\rho(q,\eta)\geq r_0$
and the geodesics $\gamma_{q,\eta_j}$ can not intersect before their first
cut points.
As vectors $\eta_j$ can be varied in sufficiently small open sets   
so that the properties stated in the claim stay valid, we obtain   
the  claim  concerning the open set of light-like directions.  
  
The last claim follows from the fact that $\eta_1$ may be equal to $\eta$ and $T={\bf t}(y)$.  
 \hfill \Box \medskip

Next we  analyze the   
set ${\Sclo}((\vec x,\vec \xi),t_0)$.

  }

 If the set 
 $\cap_{j=1}^4 \gamma_{x_j,\xi_j}([{\rtext 0},\infty))$ is non-empty we denote  its earliest point by  
$Q(\vec x,\vec \xi)$.  
If such intersection point does not exists, we define $Q(\vec x,\vec \xi)$  
to be the empty set.  \MTEXT{Next we consider the relation of ${\Scle}(\vec x,\vec \xi)$ and $\E_U(q)$ {\mtext for}
 $q=Q(\vec x,\vec \xi)$.}

\begin{lemma}\label{lem: sing detection in in normal coordinates 3}  
Let {\mtext $(\vec x,\vec \xi)=((x_j,\xi_j))_{j=1}^4\in (L^+U)^4$ be such that the points $x_j$   
satisfy} (\ref{eq: summary of assumptions 1}).  
Let $\V=  {\rtext \V(\vec x,\vec \xi)}$ be the set defined in  (\ref{eq: summary of assumptions 2}).
Then 
 \smallskip  
  
{(i)  
Assume that $y\in \V\cap U$ satisfies the condition (I) with $(\vec x,\vec \xi)$  and parameters   
$q,{\muutos w}$, and $t$ such that $0\leq t\leq \rho(q,{\muutos w})$. 
 Then   
$y\in  {\vtext \S_{H}}(\vec x,\vec \xi)$. 
 \smallskip   
  

 {\vtext (ii) {\vtext Assume that  the} geodesics corresponding to $(\vec x,\vec \xi)$ either do not
 intersect in $\V$  or those intersect at $q$ and $y\not \in J^+(q)$.  
Then $y\not \in {\vtext \S_{H}}(\vec x,\vec \xi)$.}

\smallskip   
  
(iii)  {\mtext The sets ${\Scle}(\vec x,\vec \xi)$ satisfy 
\ba
 {\Scle}(\vec x,\vec \xi)=\E_U(q)\subset \V,&&\hbox{if  $Q(\vec x,\vec \xi)\not =\emptyset$ and $q=Q(\vec x,\vec \xi)\in\V$,}\\
 {\Scle}(\vec x,\vec \xi)\subset M_0\setminus \V,\hspace{6mm}&& \hbox{if  $Q(\vec x,\vec \xi)\cap \V =\emptyset.$}
\ea
}}
%
   \end{lemma}  
      
   \noindent  
   {\bf Proof.}    {\vtext (i) 
    Consider sets ${\mathcal X}(\vec x,\vec\xi)$ and ${\mathcal Y}(\vec x,\vec\xi)$
 defined in (\ref{Y-set}).
 When $p\in \V$ and $(p,\zeta)\in
 {\mathcal X}(\vec x,\vec\xi)$,
 we see that there are $1\leq i_1<i_2<i_3\leq 4$
 such that $p=\gamma_{x_{i_k},\xi_{i_k}}(t_{i_k})$  with some $t_{i_k}>0$,
 that is, $p$  is an intersection point of some three geodesics. When
 $v_k=(\dot \gamma_{x_{i_k},\xi_{i_k}}(t_{i_k}))^\flat$, $k=1,2,3$,
 we see that ${\mathcal X}(\vec x,\vec\xi)\cap T^*_pM=\{v=\sum_{k=1}^3a_kv_k\in T^*_pM\setminus \{0\};\
 g(v,v)=0\}$. Then, we see that 
 ${\mathcal X}(\vec x,\vec\xi)\cap T^*_{\vtext p}M$
 is a union of two 2-dimensional cones.
This implies that the Hausdorff dimension of the set ${\mathcal Y}(\vec x,\vec \xi)\cap \V$ is at most 2. 
   
   Assume first that the point  $y$ satisfies conditions in  (i), $y\not\in {\mathcal Y}(\vec x,\vec\xi)$
   and $t<\rho(q,{\muutos w})$. Then {\mtext $y\in {\mmmmtext \firstpoint^{reg}_U}(q)$, see (\ref{eq: b-sets}),
   and $y$  has a neighborhood $W\subset U$  such that ${\mmmmtext \firstpoint^{reg}_U}(q)\cap W$
   is a smooth 3-dimensional submanifold. Then 
    the 
   assumptions  of Lemma
   \ref{lem: sing detection in in normal coordinates 2 wave} (ii)
   are valid for all points $y^\prime\in {\mmmmtext \firstpoint^{reg}_U}(q)\cap W$, and Lemma
   \ref{lem: sing detection in in normal coordinates 2 wave} (ii)
   implies that $y\in {\vtext \S_{H}} (\vec x,\vec \xi)$.}
   
   Consider next a general point $y$ satisfying  the
   assumptions in (i) and let $q=Q(\vec x,\vec \xi).$ 
   Then $y\in \E_U(q)$.
    {\gtext Recall that  $\rho(x,\xi)$ is
   lower semi-continuous. Then 
   the set ${\mmmmtext \firstpoint^{reg}_U}(q)\setminus {\mathcal Y}(\vec x,\vec\xi)$ is  dense in $\E_U(q)$ and
   we have that $y$ is a limit point of points {\mtext $y_n \in {\mmmmtext \firstpoint^{reg}_U}(q)\setminus {\mathcal Y}(\vec x,\vec\xi)$}. As ${\vtext \S_{H}}(\vec x,\vec \xi)$ is closed in the relative topology of $U$,   this yields that $y\in  {\vtext \S_{H}} (\vec x,\vec \xi)$.

 (ii) In the case when the geodesics corresponding to $(\vec x,\vec \xi)$ do
 intersect at $q\in \V$, denote $\V_1=\V\setminus J^+(q)$. Also, in the case when the geodesics corresponding to $(\vec x,\vec \xi)$ do not
 intersect in $\V$, denote $\V_1=\V$. 
 Then by Lemma
   \ref{lem: sing detection in in normal coordinates 2 wave} (i), condition (D) is not
   valid for any point in the set $\V_1\setminus ({\mathcal Y}(\vec x,\vec \xi)\cup\bigcup_{j=1}^4
 \gamma_{x_j,\xi_j}(\R))$. 
 As the Hausdorff dimension of the set $({\mathcal Y}(\vec x,\vec \xi)\cup\bigcup_{j=1}^4 \gamma_{x_j,\xi_j}(\R))\cap \V$ is at most 2, 
$ \V_1$ does not intersect 
 ${\vtext \S_{H}}(\vec x,\vec \xi)$. This 
 yields the claim (ii).
   

  (iii)   
Suppose
    $q=Q(\vec x,\vec \xi)\in  \V$ and
     $y\in \E_U(q)\setminus
  \bigcup_{j=1}^4\gamma_{x_j,\xi_j}([{\rtext 0},\infty))$. Let $\gamma_{q,\eta}([0,l])$ be a light-like
  geodesic that is {\mtext one of the longest causal geodesics from $q$ to $y$. Then $l\leq \rho(q,\eta)$. Let 
  $p_j=\gamma_{x_j,\xi_j}({\rtext {\bf t}_j})$, ${\bf t}_j=\rho( {\mtext x_j,\xi_j})$, be the first cut point on the geodesic $\gamma_{x_j,\xi_j}([{\rtext 0},\infty))$.
 To show that 
   $y$ is in $\V$, we assume the opposite, $y\not \in \V$. Then  for some $j$ there is
  a causal geodesic $\gamma_{p_j,\theta_j}([0,l_j])$ from $p_j$ to $y$.  Now we can use
  a short-cut argument: Let $q=\gamma_{x_j,\xi_j}(t^\prime)$. 
  As  $q\in \V$, we have $t^\prime<{\rtext {\bf t}_j}$. 
  Moreover, as 
  $y\not \in \gamma_{x_j,\xi_j}([{\rtext 0},\infty))$, the union of the geodesic
  $\gamma_{x_j,\xi_j}([t^\prime,{\rtext {\bf t}_j]})$ from $q$ to $p_j$ and 
  $\gamma_{p_j,\theta_j}([0,l_j])$ from $p_j$ to $y$ does not form a light-like geodesic
  and thus $\tau(q,y)>0$. As     $y\in \E_U(q)$, this is not possible. Hence
  $y\in \V$.   Thus by (i), $y\in {\vtext \S_{H}}(\vec x,\vec \xi)$
  and hence $\E_U(q)\setminus ( \bigcup_{j=1}^4\gamma_{x_j,\xi_j}([{\rtext 0},\infty)))\subset{\vtext \S_{H}}(\vec x,\vec \xi)$. 
  Since  the set $\E_U(q)\setminus (
  \bigcup_{j=1}^4\gamma_{x_j,\xi_j}([{\rtext 0},\infty)))$ is dense in the closed set
  $\E_U(q)\subset U$, the above shows that  $\E_U(q)\subset {\vtext \S_{H}}(\vec x,\vec \xi)$.
  Also
  by (ii),  $ {\vtext \S_{H}}(\vec x,\vec \xi)\subset {\vtext J^+(q)}.$ 
 Using {\vtext Lemma \ref{lemma E_U}, Def.\ \ref{def: first observations},
 and  (\ref{EE1}), we {\vtext conclude} 
 that $ {\Scle}(\vec x,\vec \xi)= \E_U(q).$}

  On the other hand, if $Q(\vec x,\vec \xi)\cap \V =\emptyset$, we can apply
  (ii) for all $y\in \V\cap U$ and {\vtext obtain}  that  ${\vtext \S_{H}}(\vec x,\vec \xi)\cap \V=\emptyset$.
  This and (\ref{EE1}) prove (iii).}}}
\hfill \Box \medskip

  {\rbtext

 {\mtext Next we show  the sets ${\Scle}(\vec x,\vec\xi)$, where $ (\vec x,\vec\xi)\in (L^+U)^4$,
 determine the  family of the earliest light observation sets. As we believe that this type of result} is useful for inverse problems for a wide range of non-linear {\mmmmtext partial differential equations}, we formulate this in more general terms using two geometric properties,
 denoted below by (P1) and (P2), {\vtext and the first cut points $\gamma_{x_j,\xi_j}(\rho(x_j,\xi_j))$ of the geodesics
 $\gamma_{x_j,\xi_j}([0,\infty))$.} 
 We note that Theorem \ref{prop. A} below is also valid for Lorentzian manifolds of dimension $n\geq 3$.

  
\begin{theorem}\label{prop. A}
(i) Assume that for  all $(\vec x,\vec\xi)=((x_j,\xi_j))_{j=1}^4\in (L^+U)^4$,
such that $(x_j)_{j=1}^4$ satisfy  (\ref{eq: summary of assumptions 1}),
{\mmmmtext there are} sets
  $S_0(\vec x,\vec\xi)\subset U$ {\mmmmtext that} have the following properties:
 
   \medskip

\noindent (P1) {\mtext If there is  $q\in J^-(p^+)$ such that $q=\gamma_{x_j,\xi_j}(t_j)$  with $t_j\in (0,\rho(x_j,\xi_j))$, for all $j=1,2,3,4$,  
 then   
 $S_0(\vec x,\vec\xi)= \E_U(q)$,}
    \medskip

\noindent 
 (P2) {\mtext If there are no such  $q\in J^-(p^+)$,} then 
 {\muutos $S_0(\vec x,\vec\xi)\subset M\setminus \V,$ where 
 $M\setminus \V=\cup_{j=1}^4 J^+(\gamma_{x_j,\xi_j}(  \rho({\rtext x_j,\xi_j})))$.}
  \medskip
  
Assume that we are given $(U,g|_U)$, the causality relation $R_U^<$  in $U$,
{{\vtext the set $G_{\hat\mu}$, see (\ref{The set G}), and}
the family $\{S_0(\vec x,\vec\xi);\ (\vec x,\vec\xi)\in (L^+U)^4\}$.}

Then these data
determine uniquely the {\rtext family} $\{\E_U(q);\ q\in I^+(p^-)\cap I^-(p^+)\}$ of
the earliest light observation
sets.

(ii) 
The properties (P1) and (P2) are valid
for all set  ${\Scle}(\vec x,\vec\xi)$ such that $(\vec x,\vec\xi)\in (L^+U)^4$
and $(x_j)_{j=1}^4$ satisfy  (\ref{eq: summary of assumptions 1}).
\end{theorem}
}
{\mtext We call the sets  $S_0(\vec x,\vec\xi)$   the generalized observations, and emphasise 
that for such a set we do not a priori know it is of the type considered in (P1) or (P2).
{\ftext In particular, we do not know a prior if a given set $S_0(\vec x,\vec\xi)$ corresponds to the interaction of waves
that has started to happen before or after the conjugate points of the geodesics $\gamma_{x_j,\xi_j}([0,\infty))$.}
By claim (ii),  the observations related to wave equation are an example of generalized observations

The proof of the claim (ii) of
Theorem \ref{prop. A} 
is obtained immediately from Lemma \ref{lem: sing detection in in normal coordinates 3}  (iii)}.

We prove the claim of Theorem \ref{prop. A} (i) in the next subsection. {\mtext To give the idea, 
before proving}   Theorem \ref{prop. A} (i) for general globally hyperbolic manifolds,  we consider a simpler case where the proof of  Theorem \ref{prop. A} is {\ftext easier}.}\medskip

  \noindent
{\bf Proof of Theorem \ref{prop. A} (i) in \ftext a special case.} 
{\rtext Let us consider a special case when the light-like geodesics
  in {\mtext $I(p^-,p^+)=I^+(p^-)\cap I^-(p^+)$ do not contain cut points.
  Then, we consider all   $(\vec x,\vec \xi)\in (L^+U)^4$ such that $(x_j)_{j=1}^4$ satisfy conditions (\ref{eq: summary of assumptions 1}) and $x_j\in  I^+(p^-)\cap U$. When there are no cut points, 
 all intersection points of geodesics in $I(p^-,p^+)$ are automatically in the set $\V(\vec x,\vec \xi)$, see} (\ref{eq: summary of assumptions 2}). {\mtext Then (P1) and (P2) imply that
  $S_0(\vec x,\vec \xi)=\E_U(q)$
  if the geodesics corresponding to $(\vec x,\vec \xi)$ intersect at some point $q\in I(p^-,p^+)$. {\mtext Moreover,} $\mathcal S_e(\vec x,\vec \xi)$  is an empty set if no such intersection point exists in $I^+(p^-)$. 
  
  {\mtext  Consider a point $q\in{\mtext I^+(p^-,p^+)}$. By Lemma \ref{lemma: existence of x,xi}
  there are $(\vec x,\vec \xi)$ satisfying (\ref{eq: summary of assumptions 1}) such that $x_j\in I^+(p^-)\cap U$  and that the corresponding geodesics intersect at $q$. Also, we have that $S_0(\vec x,\vec \xi)$ intersects the set $U\cap I^-(p^+)$. 
{\ftext Then, let us consider the family
\ba
\{S_0(\vec x,\vec \xi)&;&   \hbox{$(\vec x,\vec \xi)\in (L^+U)^4$ are such that $x_j\in I^+(p^-)\cap U$}\\
& & \hbox{satisfy conditions (\ref{eq: summary of assumptions 1})
and $S_0(\vec x,\vec \xi)\cap (U\cap I^-(p^+))\not=\emptyset$}\}.
\ea
The above yields that} this family coincides with the family $\{\E_U(q);\ q\in I^+(p^-,p^+)\}$ of
the earliest light observation
sets.}  This completes the proof of Theorem \ref{prop. A} in the special case 
when the light-like geodesics
  in $ {\mtext I(p^-,p^+)}$ do not contain cut points.} \hfill \Box \medskip  

{\mtext Theorem \ref{prop. A} reduces the active inverse problem considered in Theorem \ref{main thm3 new}
to the passive inverse problem.}
Later, in Section \ref{sec: Proof for passive observations} we finish
the proof on the uniqueness for the inverse problem with passive observations.
In the rest of this section we will prove  Theorem \ref{prop. A} for general globally hyperbolic spacetimes and remark that a reader who is interested {\mmmmtext just} in spacetimes having no cut points can
move to Section \ref{sec: Proof for passive observations}. }

 \subsection{Determination of the earliest light observation set}

 {\rtext In this subsection we consider the proof of   Theorem \ref{prop. A} (i) in the general case. The proof will be quite technical due to the reason that
 above we have analyzed interaction of waves that propagate near light-like geodesics
 and intersect before the first conjugate to cut points of the geodesics. 
 If the geodesics intersect after or at the conjugate points, the waves may have caustics and
 the waves produced by the non-linear interaction may be very complicated.
 However, we do {\ftext not} know the manifold $(M,g)$ and thus we do not a priori 
 know when the geodesics have conjugate points. Therefore, we have to determine from our observations
 when we are sure that the interaction of the waves has taken place before the the conjugate points. {\mmmmtext Then we can} remove 
 from our data all observations that may be caused by caustics.}

\medskip

%
%

} 
%
%
%
%
%
%

  \noindent
{\bf Proof of Theorem \ref{prop. A} (i).}
   Below, we {\rtext use the numbers} $\vartheta_1,\kappa_1,\kappa_2>0$  appearing in Lemma \ref{lem: detect conjugate 0}
   and denote ${t_0}=4\kappa_1$.
%
%

First, consider the set 
 \beq\label{set Kt0}
 \K_{t_0}=\{x \in {\mtext M}&;&
x= \gamma_{{\hat x},\xi}(r),\
{\hat x}=\hat \mu(s), \ s\in [s^-,s^+],\\
\nonumber
& & 
 \xi\in L^+_{\hat x}M,\ \|\xi\|_{g^+}=1,\ r\in [0,2{t_0})\}{\mtext \subset U}.
\eeq that is the closure of 
a small neighborhood of $\hat \mu$. As we will consider waves propagating
near geodesics
$\gamma_{\hat x,\hat \zeta}([t_0,\infty))$ where $\hat x\in \hat \mu$,
we have to consider the earliest light observation sets
corresponding to  the points $\gamma_{\hat x,\hat \zeta}([0,t_0))$ separately.
This is why the set  $ \K_{t_0}$ is introduced.

\begin{center}  
  
 \psfrag{1}{\hspace{-2mm}$J^+(B_1)$}  
\psfrag{2}{$J^+(B_2)$}  
\psfrag{3}{$B_1$}  
\psfrag{4}{$B_2$}  
\psfrag{5}{$p^-$}

\psfrag{1}{$\hspace{-2mm}{\hat x}_2$}  
\psfrag{3}{\hspace{-1.5mm}$p_0$}  
\psfrag{2}{$p_0^\prime$}  
\psfrag{4}{\hspace{-2mm}${\hat x}$}  
\psfrag{5}{y}  
\psfrag{6}{${\hat x}_1$}  
\psfrag{7}{$z$}  
\psfrag{8}{$y^\prime$}  
\psfrag{9}{$\hspace{-2mm}p^-$}  
\psfrag{0}{$p^+$}  
\psfrag{P}{$p$}  
\includegraphics[width=5.5cm]{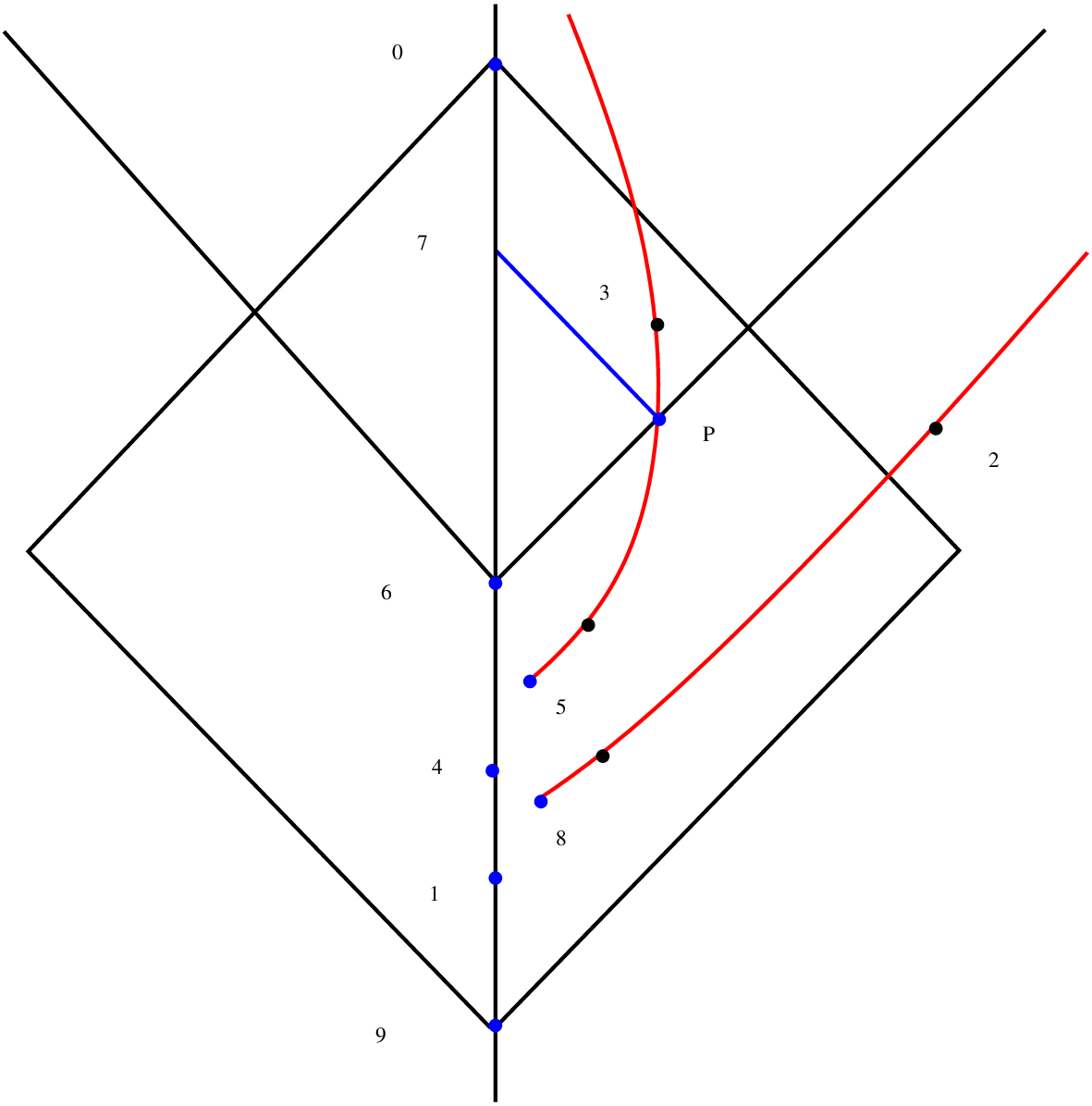}  
\end{center}  
{\it FIGURE 7. 
%
The blue points on $\hat \mu$ are ${\hat x}_1=\hat \mu(s_1)$, ${\hat x}_2=\hat \mu(s_2)$, and ${\hat x}=\hat \mu(s)$.
The blue points $y$ and $y^\prime$ are close to $\hat x$.  
{\gtext The boundary of $J^+({\hat x}_1)$ is marked by black.}
 We consider the 
 geodesics $\gamma_{y,\zeta}([0,\infty))$ and $\gamma_{y^\prime,\zeta^\prime}([0,\infty))$.   
 These geodesics  corresponding to the cases when
the geodesic  $\gamma_{y,\zeta}([0,\infty))$ enters in $J^-({p^+})\cap J^+({\hat x}_1)$,
 and the case
when the geodesic  $\gamma_{y^\prime ,\zeta^\prime}([0,\infty))$ does not
enter this set. The point $p_0$ is the cut point of  $\gamma_{y,\zeta}([t_0,\infty))$
and  $p_0^\prime$ is the cut point of  $\gamma_{y^\prime ,\zeta^\prime}([t_0,\infty))$.
At the point $z=\hat \mu(\mathbb{S} (y,\zeta,s_1))$  we observe 
for the first time  on the geodesic $\hat \mu$  that  
the geodesic $\gamma_{y,\zeta}([0,\infty))$ has entered   
  $J^+({\hat x}_1)$.  {\gtext The entering in the set  $J^+({\hat x}_1)$ happens at the point $p$.}
}  


   \MTEXT{
   
\smallskip

\extension{ 
Recall that $\U_{z_0,\eta_0}=\U_{z_0,\eta_0}(\hat h)$  
 was defined using the parameter $\hat h$. We see 
 using compactness of $\hat \mu ([s_-,s_+])$, the continuity of $\tau(x,y)$, and 
the existence of convex neighborhoods \cite[Prop.\ 5.7]{ONeill},
  c.f.\ Lemma 2.A.1,
 that if $t_0^\prime >0$  
 and  $\hat h^\prime\in (0, \hat h)$ are small enough and $a\in \A(\hat h^\prime)$,  
 then the  
 longest geodesic from $x\in \overline \K_{t_0^\prime}$ to the point  
 $\E_{a}(x)$ is contained in $U$
 and hence we can determine the point $\E_{a}(x)$ for such $x$ and $(z,\eta)$.  
 Let us replace the parameters $\hat h$ and $t_0$ by $\hat h^\prime$  
 and $t_0^\prime$, correspondingly in our considerations below. Then we  
 may  assume  
 that in addition to the data given in the original formulation of the problem,  
 we are given also the set $\be_U(\K_{t_0})$. Next we do this.
}
\noextension{
{\mattitext 
\canberemoved{Recall that by (\ref{tau known in U}) the 
set $R_U=\{(x,y)\in U\times U;\  x\ll y\}$
is assumed to be known.}

{\itext {\rtext As
we are given $(U,g|_U)$ and the causality relation {\vtext $R_U^<$} in $U$, we can determine the subset
$\K_{t_0}\cap J^+(\hat \mu(s))$ for
all $s\in [s_-,s_+]$. As this set is a subset of $U$, {\mtext by Lemma \ref{lem: aux construction}} we can determine the earliest light observation sets {\mtext in}}
 $\E_U(\K_{t_0}\cap J^+(\hat \mu(s)))$ for
all $s\in [s_-,s_+]$. Here, recall that for a set $W\subset M$, we denote  $\E_U(W)=\{\E_U(q);\ q\in W\}\subset 2^{U}$.}
%
Also, we may assume below that $\vartheta_1$ is so \MTEXT{small that
$\gamma_{y ,\zeta}([0,t_0]) \cap J(p^-,p^+)\subset \K_{t_0}$ when 
$y\in J(p^-,p^+)$, $d_{g^+}(y,\hat \mu)<\vartheta_1$ and
$\zeta\in L^+_yM$, $\|\zeta\|_{g^+}\leq 1+\vartheta_1$.}

{Let  $s_0\in [s_-,s_+]$   
be so close to $s_+$  
 that $J^+(\hat \mu(s_0))\cap J^-({p^+})\subset \K_{t_0}$. Then the {\mtext data given in the claim} 
  determine $\E_{U}(J^+(\hat \mu(s_0))\cap J^-({p^+}))$.


\subsubsection{Determination of the time when a geodesic is observed to enter in to the 
already reconstructed set.}

{\mattitext 
Let us describe the rough idea of the construction that we do next:
We will consider the point  $p=\gamma_{y,\zeta}(r)$ 
where a geodesic $\gamma_{y,\zeta}$ {\ftext enters for the first time} (see Fig.\ 7) in the
``already reconstructed'' set $J^+({\hat x}_1)\cap J^-({p^+})$ or exits
the set $J^-({p^+})$. In particular we consider the time $s=
\mathbb{S} (y ,\zeta ,s_1)$ where the light coming from the point $p$
is observed on $\hat \mu$ {\ftext for the first time}. 
The time $\mathbb{S} (y ,\zeta ,s_1)$ will be essential for us as we {\mtext are} sure
that {\mtext before this time} we do not {\mtext observe} on $\hat \mu$ any strange signals 
that caustics may have produced.
The idea {\mtext to find $\mathbb{S} (y ,\zeta ,s_1)$} is that when
$r^\prime>r$ is close to $r$, the observations from an artificial point sources produced
at the point $\gamma_{y,\zeta}(r^\prime)$ coincide with some of the observations
that we have made earlier.  
Next we present details of this construction.}

We recall the notation that $\hat \mu=\mu_{\hat a}$ {\mtext and we use $s_+<s_{+2}<1$ such 
that $p^+=\hat \mu(s_+)$.}

{\rbtext 
Next we consider $\vartheta\in (0,\vartheta_1)$, the
points $x_j\in U$ and the  directions $ \xi_j\in L^+_{x_j}\hattuM _0$,
denoted by
  $(\vec x,\vec \xi)=((x_j,\xi_j))_{j=1}^4$, that
  satisfy, see Fig.\ 5(Right), 
\beq\label{eq: summary of assumptions 1 version B}& &\hspace{-1.0cm}(i)\ \gamma_{x_j,\xi_j}([0,{\rtext t_0}])\subset U,\ 
 {\rtext \gamma_{x_j,\xi_j}({\rtext t_0})}\not \in J^+(\gamma_{x_k,\xi_k}({\rtext t_0})),
  \hbox{ for  $j\not =k$,}  \hspace{-5mm}\\ & &\nonumber
 \hspace{-1.0cm}
 {\rtext (ii)\ 
\hbox{For all $j,k\leq 4$,
 $d_{\nohat g^+}((x_j,\xi_j),(x_k,\xi_k))<\vartheta$,}} \hspace{-5mm}\\
 & &\nonumber
 \hspace{-1.0cm}{\rtext (iii)\ 
\hbox{There is $\hat x\in  \hat \mu$ such that 
for all $j\leq 4$,
 $d_{\nohat g^+}(\hat x,x_j)<\vartheta$.}} \hspace{-5mm}
 \eeq
  \refHOX{Comment 43: The notation $t_0$ in the introduction is changed to $T_0$.}
}

\extension{

\vfill

\goodbreak 
\Egin{center}  
  
\psfrag{1}{$\tilde y^\prime$}  
\psfrag{2}{$p_2$}  
\psfrag{3}{$J_{g}(\hat p^-,\hat p^+)$}  
\psfrag{4}{$\tilde y$}  
\psfrag{5}{}  
\psfrag{6}{$y_1$}  
\psfrag{7}{}  
\psfrag{8}{}
\psfrag{A}{\hspace{-3mm}$\mu_{a}$}  
\psfrag{B}{$\mu_{\hat a}$}  
\includegraphics[width=7cm]{strategy_3.eps}  
\end{center}  
  {\it FIGURE A8: A schematic figure where the space-time is represented as the  2-dimensional set $\R^{1+1}$.   
We consider  
the geodesics emanating from a point $x(r)=\gamma_{\tilde y,\tilde \zeta}(r)$,  
When $r$ is smallest value for which $x(r)\in J^+(y_1)$, a light-like geodesic  (black line segment)  
emanating from $x(r)$ is observed at the point $\tilde p\in \hat \mu $.  
Then $\tilde p= \hat\mu(\mathbb{S} (\tilde y,\tilde \zeta,s_1))$.  
When $r$ is small enough so that $x(r)\not \in J^+(y_1)$,  
the light-geodesics (red line segments) can be observed on $\hat \mu$ in the set  $J^-(\tilde p)$.  
Moreover, when  $r$  is such that $x(r)\in J^+(y_1)$, the light-geodesics can be observed at $\hat \mu$ in the set  $J^+(\tilde p)$. The golden point is the cut point  on $\gamma_{\tilde y,\tilde \zeta}([t_0,\infty))$ and the singularities on the light-like geodesics   
starting before this point (brown line segments) can be   
analyzed, but after the cut point the singularities on the light-like geodesics (magenta line segments)  
are not analyzed in this paper.  
}

\vfill

\goodbreak 

\psfrag{1}{$x_1$}  
\psfrag{2}{$x_2$}  
\psfrag{3}{$q$}  
\psfrag{4}{$p$}  
\psfrag{5}{}  
\psfrag{6}{$y_1$}  
\psfrag{7}{$x$}  
\psfrag{8}{$z$}  
\begin{center}  
\includegraphics[width=6.5cm]{garbage_2D_2.eps}  
\end{center}  
{\it FIGURE A9:  
A schematic figure where the space-time is represented as the  2-dimensional set $\R^{1+1}$.   
In section 5 we consider geodesics $\gamma_{x_j,\xi_j}([0,\infty))$, $j=1,2,3,4,$ that all intersect for the first time at a point  
$q$. 
When we consider geodesics $\gamma_{x_j,\xi_j}([t_0,\infty))$, $j=1,2$ with $t_0>0$, they may  
have cut points at $x^{cut}_j=\gamma_{x_j,\xi_j}({\bf t}_j)$.  
In the figure $\gamma_{x_j,\xi_j}([0,{\bf t}_j))$ are colored by black and the geodesics   
 $\gamma_{x_j,\xi_j}([{\bf t}_j,\infty))$ are colored by green. In the figure  
 the geodesics intersect at the point $q$ before the cut point and for  
the second time at $p$  
 after the cut point. We can analyze the singularities caused by the distorted plane  waves   
 that interact at $q$ but not the interaction of waves after the cut points of the geodesics.  
 It may be that e.g. the intersection at $p$ causes new singularities to appear and we observe  
 those in $U$. As we cannot analyze these singularities, we consider these singularities  
 as ``messy waves''. However, the ``nice'' singularities caused by the interaction at $q$ propagate along  
 the future light cone of the point $q$ and in
$U\setminus \cup_{j=1}^4\gamma_{x_j,\xi_j} $ these ``nice'' singularities are observed before the ``messy waves''. Due to this the first singularities we observe near $\hat \mu$ 
come from the point $q$.  
}  
}

\begin{definition}\label{def: hat s 1}  
Let $s_-\leq s_2\leq s<s_1\leq s_+$ satisfy  $s_1<s_2+\kappa_2$, 
 ${\hat x}_j={\hat \mu}(s_j)$, $j=1,2$, and  ${\hat x}={\hat \mu}(s)$,
 ${\hat \zeta} \in L^+_{{\hat x}}U$, 
    $\|{\hat \zeta}\|_{g^+}=1$.   
   \MTEXT{Let $(y,\zeta) \in L^+U$ be in $\vartheta_1$-neighborhood of
    $({\hat x},{\hat \zeta})$ such that
     $y\in J^+({\hat x}_2)$ and 
  the geodesic $\gamma_{y,\zeta}(\R_+)$ does not intersect $\hat \mu$.  
Define
\ba 
& &\rkaksi(y,\zeta,s_1)=\inf\{r>0;\ \gamma_{y,\zeta}(r)\in J^+(\hat\mu (s_1))\},\\
& &\ryksi(y,\zeta)=\inf\{r>0;\ \gamma_{y,\zeta}(r)\in M\setminus I^-(\hat\mu(s_{+2}))\},\\
& &r_0(y,\zeta,s_1)=\min(\ryksi(y,\zeta),\rkaksi(y,\zeta,s_1)).
\ea}
 
When  
 $\gamma_{y ,\zeta }(\R_+)$ intersects $J^+(\hat \mu(s_1))\cap J^-({p^+})$  
we define  
\beq\label{special S formula}  
\mathbb{S} (y ,\zeta ,s_1)=f^+_{\hat a}(q_0),
\eeq  
where $ q_0=\gamma_{y ,\zeta }(r_0)$ and   
$r_0=r_0(y,\zeta,s_1).$
In the case when  
$\gamma_{y ,\zeta }(\R_+)$ does not intersect $J^+(\hat \mu(s_1))\cap  {J}^-({p^+})$,  
we define $\mathbb{S} (y ,\zeta ,s_1)={\mtext s_+}$.

\end{definition}  
  
{We note that above 
  $\ryksi(y,\zeta)$ is finite by \cite[{\ftext Lemma}\ 14.13]{ONeill}.}
{\rtext Note that by {\ftext the assumptions of the claim, we can check  if $\gamma_{x,\xi}\cap\hat \mu=\emptyset$
for a given $(x,\xi)$}}.
%
%
%

\begin{definition}\label{def: hat s}  
Let $0<\vartheta<\vartheta_1$ and 
$\qP_{\vartheta} (y ,\zeta )$  be the set of $(\vec x,\vec \xxi)=((x_j,\xxi_j))_{j=1}^4$ 
that satisfy {\ftext the} {\mtext conditions} 
in the formula {\rtext (\ref{eq: summary of assumptions 1 version B})} {and $(x_1,\xi_1)=(y,\zeta)$.}
  We say that the set $S\subset U$ is a 
  {\mtext repeated} observation
 associated to the geodesic $\gamma_{y,\zeta}$ if  there is $\hat \vartheta>0$
 such that for all $\vartheta\in (0,\hat \vartheta)$
  there are  $(\vec x,\vec \xxi)\in \qP_{\vartheta} (y ,\zeta )$ such that
  $S={{\Smod}}(\vec x,\vec \xi)$.
\end{definition}


%

Next we use \MTEXT{a step-by-step construction}: We    
consider  
 $s_1\in (s_-,s_+)$ and assume that  we are given  $\E_{U}(J^+({\hat x}_1)\cap J^-({p^+}))$ with ${\hat x}_1={\hat \mu}(s_1)$.} Then, let  $s_2\in (s_1-\kappa_2,s_1)$.  
Our next aim is to find the earliest light observation sets $\E_{U}(J^+({\hat x}_2)\cap J^-({p^+}))$ with ${\hat x}_2={\hat \mu}(s_2)$.

{\rtext Let us consider {\ftext now} four light-like future pointing directions $(x_j,\xi_j)$, $j=1,2,3,4$, and
use below the notation $(\vec x(t),\vec\xi(t))$ defined in (\ref{eq: x(h) notation}).}

{\begin{lemma}\label{conjugatepoints are nice}   
\MTEXT{{\gtext Assume that} $\max(s_-,s_1-\kappa_2)\leq s<s_1<s^+$
and let
${\hat x}=\hat \mu(s)$, ${\hat x}_1=\hat \mu(s_1)$, and ${\hat \zeta}\in L^+_{{\hat x}}M$, $\|{\hat \zeta}\|_{g^+}=1$.
Moreover, let
  $(y,\zeta)$ be in  a $\vartheta_1$-neighborhood of $({\hat x},{\hat \zeta})$. 
  {\gtext Assume also} that the geodesic $\gamma_{y,\zeta}(\R_+)$ does not intersect $\hat \mu$.}
  {\mtext Then,}

 \noindent (A) {\mtext The} cut point  $p_0=\gamma_{y(t_0),\zeta(t_0)}({\bf t}_*)$, ${\bf t}_*= \rho(y(t_0),\zeta(t_0))$ of   
 the geodesics $ \gamma_{y(t_0),\zeta(t_0)}([0,\infty))$, if it exists, satisfies either  

  \smallskip  
   
(i)  $p_0\not \in J^-(\hat \mu(s_{+2}))$,  
 \smallskip  
  
\noindent

or  
 \smallskip

(ii)  $r_0=r_0(y,\zeta,s_1)<\ryksi(y,\zeta)$ 
and $p_0\in I^+({\hat x}_1)$.
 \smallskip  
 
 \noindent
(B)
 There is $\vartheta_2(y,\zeta,s_1)\in (0,\vartheta_1)$   
   such that if $0<\vartheta<\vartheta_2(y,\zeta,s_1)$, 
 $(\vec x,\vec \xi)\in \qP_\vartheta(y,\zeta)$,
  \MTEXT{ and the geodesics $ \gamma_{ {\rtext  x_j, \xi_j}}([0,\infty))$,
   $j\in \{1,2,3,4\}$,
 has a cut point  
 $p_j=\gamma_{  {\mtext  x_j(t_0), \xi_j(t_0)}}({\mtext \rho(x_j(t_0), \xi_j(t_0))})$,
 then the following holds: 
   
If either the point $p_0$ does not exist or it exists and (i) holds
then 
  $p_j\not \in J^-({p^+})$.
On the other hand, if $p_0$ exists and
 (ii) holds, then $f^+_{\hat a}(p_j)
>f^+_{\hat a}(q_0)$, where $q_0=\gamma_{y,\zeta}(r_0(y,\zeta,s_1))$.}

Note that
$f^+_{\hat a}(q_0)=\mathbb{S} (y ,\zeta ,s_1)$.  

  
  \end{lemma}

\noindent  {\bf Proof.}  
(A) Assume that (i) does not hold, that is, $p_0
=\gamma_{y,\zeta}(t_0+{\bf t}_*) \in J^-(\hat \mu(s_{+2}))$. 
By  Lemma \ref{lem: detect conjugate 0}  (iii) we have
$f^-_{\hat a}(p_0)>s+2\kappa_2\geq s_1$
that yields $p_0
\in I^+({\hat x}_1)$.
Thus, the geodesic
$\gamma_{y(t_0),\zeta(t_0)}([0,\rho(y(t_0),\zeta(t_0)))$ intersects
 $J^+({\hat x}_1)\cap J^-(\hat \mu(s_{+2}))$. Hence the alternative (ii)
 holds {with $0<r_0
<\ryksi(y,\zeta)$ and moreover, $r_0< t_0+\rho(y(t_0),\zeta(t_0))$.}

(B) If (i) holds, the claim follows since the function $(x,\xi)\mapsto \rho(x,\xi)$ is lower semi-continuous and $(x,\xi,t)\mapsto \gamma_{x,\xi}(t)$ is continuous.  

In the case (ii), we saw above that $r_0<t_0+\rho(y(t_0),\zeta(t_0))$.
Let  $q_0=\gamma_{y,\zeta}(r_0)$. Then by using a short cut argument and the 
fact that $\gamma_{y,\zeta}(\R_+)$ does not intersect $\hat \mu$
we see similarly to the above that $f^+_{\hat a}(p_0)>f^+_{\hat a}(q_0)=\mathbb{S} (y ,\zeta ,s_1)$.
Since the function $(x,\xi,t)\mapsto f^+_{\hat a}(\gamma_{x,\xi}(t))$ is continuous
and $t\mapsto f^+_{\hat a}(\gamma_{x,\xi}(t))$ is non-decreasing, and  the function $(x,\xi)\mapsto \rho(x,\xi)$ is lower semi-continuous, 
 we have that the function $(x,\xi)\mapsto f^+_{\hat a}(\gamma_{x(t_0),\xi(t_0)}(\rho(x(t_0),\xi(t_0))))$
 is lower semi-continuous, and
 the claim follows.
  \hfill \Box \medskip  
}

  {\mattitext We recall that below we assume that the set $ \E_{U}(J^+({\hat x}_1)\cap J^-({p^+}))$ is already constructed.}

   \begin{definition}  \label{def 56}
\MTEXT{Let $s_-\leq s_2\leq s<s_1\leq s_+$ satisfy  $s_1<s_2+\kappa_2$, 
 ${\hat x}_j={\hat \mu}(s_j)$, $j=1,2$, and  ${\hat x}={\hat \mu}(s)$,
 ${\hat \zeta} \in L^+_{{\hat x}}U$, 
    $\|{\hat \zeta}\|_{g^+}=1$.   
 Also, let $(y,\zeta) \in L^+U$ be in $\vartheta_1$-neighborhood of
    $({\hat x},{\hat \zeta})$ and $\mathcal R(y ,\zeta ,s_1)$  be the set
    of the
     {\mtext repeated} observations  $S\subset U$
associated to the geodesic $\gamma_{y,\zeta}$ 
such that $S\in \E_{U}(J^+({\hat x}_1)\cap J^-({p^+}))$.
Moreover, define
  $ \mathbb{S}^{obs} (y ,\zeta ,s_1)$ to be the infimum of $s^\prime\in [-1, s_+]$ 
such that $ \hat \mu(s^\prime)\in S\cap \hat \mu$ with
some  $S\in \mathcal R(y ,\zeta ,s_1)$.
If no such $s^\prime$ exists,
we define $\mathbb{S}^{obs} (y ,\zeta ,s_1)=s^+$.}
  \end{definition}  

Let us next  consider  $(\vec x,\vec\xi)\in \qP_\vartheta(y,\zeta)$   
where
$0<\vartheta<\vartheta_2(y,\zeta,s_1)$. Here, 
$\vartheta_2(y,\zeta,s_1)$ is defined in  
Lemma \ref{conjugatepoints are nice}.    
  Assume that for some $j=1,2,3,4$ we have that $\rho(x_j({t_0}),\xxi_j({t_0}))<\mathcal T(x_j({t_0}),\xxi_j({t_0}))$,
  {\mtext see Sec.\ \ref{Observation domain},}  and consider the cut point  
$p_j=\gamma_{x_j({t_0}),\xxi_j({t_0})}(\rho(x_j({t_0}),\xxi_j({t_0})))$.  
Then either the case (i) or (ii) of Lemma \ref{conjugatepoints are nice} holds.
If (i) holds, {\mattitext by Lemma \ref{conjugatepoints are nice} (B),} $p_j$  
 satisfies  
$p_j\not\in J^-({p^+})$ and thus $f^+_{\hat a}(p_j)>s_+\geq \mathbb{S} (y ,\zeta,s_1)$.
 If (ii) holds, there exists $r_0=r_0(y,\zeta,s_1)<\rkaksi(y,\zeta)$ such that  
  $q_0=\gamma_{y,\zeta}(r_0)\in J^+({\hat x}_1)$ and 
$f^+_{\hat a}(p_j)
  >f^+_{\hat a}(q_0)=\mathbb{S} (y ,\zeta ,s_1).$ Thus {\ftext in both cases} (i) and (ii) we have 
 \beq\label{eq: cut points observed after S}
  f^+_{\hat a}(p_j)>  \mathbb{S} (y ,\zeta ,s_1).
 \eeq

 {\rtext 
 {\gtext We consider next}  a point 
$q=\gamma_{y,\zeta}(r)\in J^-({p^+}),$
where $t_0<r\leq r_0=r_0(y,\zeta,s_1)$.
By Lemma \ref{lem: detect conjugate 0} (iii), the geodesic
$\gamma_{y,\zeta}([t_0,r])$ has no
cut points.


By Lemma \ref{lemma: existence of x,xi},
we see that
when  $\vartheta_3\in (0,\vartheta_2(y,\zeta,s_1))$ is small enough,
 for all $\vartheta\in (0,\vartheta_3)$, there is
$(\vec x,\vec\xi)=((x_j,\xi_j))_{j=1}^4\in \qP_\vartheta(y,\zeta)$   
such
that the geodesics corresponding to $(\vec x,\vec\xi)$ intersect
at $q$. As  the set $(U,g)$  is known, 
that for sufficiently small $\vartheta$
one can verify {\mtext using $(U,g|_U)$ and the causality relation in $U$}, if {\ftext the} given vectors $(\vec x,\vec\xi)$ satisfy 
$(\vec x,\vec\xi)=((x_j,\xi_j))_{j=1}^4\in \qP_\vartheta(y,\zeta)$. 
Also, note that as then $\vartheta<\vartheta_2(y,\zeta,s_1)$, 
the inequality (\ref{eq: cut points observed after S}) yields 
 $\tilde x{\mtext =\hat \mu(\mathbb{S} (y ,\zeta,s_1))}\in  {\rtext \V(\vec x,\vec \xi)}$ and thus $q\in 
 J^-(\tilde x)\subset  {\rtext \V(\vec x,\vec \xi)}$.
 Then  {\rtext property (P2) implies that} 
 ${{\Smod}}(\vec x,\vec \xi)=\E_{U}(q)$. }


 As   $\vartheta\in (0,\vartheta_3)$
 above can  be arbitrarily small, 
  {\mtext we see} for 
  any $q=\gamma_{y,\zeta}(r)\in J^-({p^+})$
where $t_0<r\leq r_0=r_0(y,\zeta,s_1)$, we {\mtext have}
\beq\label{aaaa}
& &\hspace{-1cm}
\hbox{
$S=\E_{U}(q)$ is a {\mtext repeated} 
observation associated to  $\gamma_{y,\zeta}$ and}
\hspace{-1cm}
 \\
\nonumber
& &\hspace{-1cm}
S\cap 
\hat \mu=\{\hat \mu(\hat s)\},
\ \hat s:=f^+_{\hat a}(q)\leq 
\mathbb{S} (y ,\zeta,s_1).\hspace{-1cm}
\eeq

  {
\begin{lemma}\label{lem: Determination of hat s}  
Assume that $\gamma_{y,\zeta}(\R_+)$ does not intersect $\hat \mu$. Then
we have  $\mathbb{S}^{obs} (y ,\zeta ,s_1)=\mathbb{S} (y ,\zeta,s_1)$. 
  
\end{lemma}

\noindent{\bf Proof.}   
  Let us first prove that $\mathbb{S}^{obs} (y ,\zeta ,s_1)\geq \mathbb{S} (y ,\zeta,s_1)$.
To this end, let $s^\prime<\mathbb{S} (y ,\zeta,s_1)$
and  $x^\prime =\hat \mu(s^\prime)$.  
Assume  $S\in  \E_{U}(J^+({\hat x}_1)\cap J^-({p^+}))$ is a {\mtext repeated} 
observation associated to the geodesic $\gamma_{y,\zeta}$ 
and $S\cap \hat \mu=\{\hat \mu(s^\prime)\}$.
 Let   $q \in J^+({\hat x}_1)\cap J^-({p^+})$ be such that
$S=\E_U(q)$.

Then for arbitrarily small  $0<\vartheta<\vartheta_2(y,\zeta,s_1)$ 
there is
 $(\vec x,\vec\xi)\in \qP_\vartheta(y ,\zeta)$
 satisfying   
 ${{\Smod}}(\vec x,\vec \xi)=S$.  Let $\V= {\rtext \V(\vec x,\vec \xi)}$.
 Then by (\ref{eq: cut points observed after S}), we have 
 $x^\prime \in \V$.

 If the geodesics  corresponding to $(\vec x,\vec\xi)$
   intersect  
at some point $q^\prime\in J^-(x^\prime)$, then  by  {\mtext (P1)} 
 we have
$S=\E_{U}(q^\prime)$. Recall that $\hat \mu=\mu_{\hat a}$.
Then 
$\hat \mu(s^\prime)=\hat \mu(f^+_{\hat a}(q^\prime))$
implying   {\ftext that} $ f^+_{\hat a}(q^\prime)=s^\prime$. Moreover, we have
 then that $\E_{U}(q)\
=\E_{U}(q^\prime)$ and 
Proposition \ref{propo paper part 2}  
 yields 
$q^\prime=q$.
Since $(\vec x,\vec\xi)\in \qP_\vartheta(y,\zeta)$ 
implies $(x_1,\xi_1)=(y,\zeta)$, 
we see that $q^\prime\in \gamma_{y,\zeta}([t_0,\infty))$.
As  $q\in J^+({\hat x}_1)$,
we {\gtext have} that $q=q^\prime\in 
\gamma_{y,\zeta}([t_0,\infty)) \cap J^+({\hat x}_1)=
\gamma_{y,\zeta}([r_0(y ,\zeta,s_1),\infty))$.
However, then $f^+_{\hat a}(q^\prime)\geq
\mathbb{S} (y ,\zeta ,s_1)  >s^\prime$ and thus $S\cap \hat \mu=\E_{U}(q)\cap \hat \mu$
can not be equal to $\{\hat \mu(s^\prime)\}$.

On the other hand, if the geodesics  corresponding to $(\vec x,\vec\xi)$
do not   intersect  at any point in 
 $ J^-(x^\prime)\subset \V$, \MTEXT{then either they intersect in some
 $q^\prime_1\in (M\setminus J^-(x^\prime))\cap \V$,    {\ftext they} do not intersect at all, or    {\ftext they} intersect at 
 $q^\prime_2\in M\setminus \V$. In the first case, $S=\E_U(q^\prime_1)$ do not satisfy
 $S\cap \hat \mu\in \hat\mu((-1,s^\prime))$. In the other cases,}
 {\rtext property (P2)} 
 yields 
 ${{\Smod}}(\vec x,\vec \xi)\cap \V=\emptyset$.
As  $x^\prime=\hat \mu(s^\prime)\in \V$, we see that     $S\cap \hat \mu$
can not be equal to $\{\hat \mu(s^\prime)\}$. Since above 
$s^\prime< \mathbb{S} (y ,\zeta,s_1)$ is arbitrary, this shows that
$\mathbb{S}^{obs} (y ,\zeta ,s_1)\geq \mathbb{S} (y ,\zeta,s_1)$.

Let us next show that $\mathbb{S}^{obs} (y ,\zeta ,s_1)\leq \mathbb{S} (y ,\zeta,s_1)$.
\MTEXT{Assume the opposite. Then, if $\mathbb{S} (y ,\zeta,s_1)=s_+$,
 we see by Def.\ \ref{def 56} that 
$\mathbb{S}^{obs} (y ,\zeta ,s_1)=\mathbb{S} (y ,\zeta,s_1)$ which leads to a contradiction. 
However, if 
$\mathbb{S} (y ,\zeta,s_1)<s_+$,  by Def.\ \ref{def: hat s 1},
we have
(\ref{eq: cut points observed after S}).  
This implies the existence of  $q_0=\gamma_{y,\zeta}(r_0)$, 
$r_0=r_0(y,\zeta,s_1)$ such that 
 $q_0\in  J^+({\hat x}_1)
\cap J^-({p^+})$ and
by (\ref{aaaa}), 
$S=\E_{U}(q_0)$ is a {\mtext repeated} 
observation associated to the geodesic $\gamma_{y,\zeta}$.
By Lemma \ref{conjugatepoints are nice} (ii), 
$\mathbb{S} (y ,\zeta,s_1)=f^-_{\hat a}(q_0)$ which implies,
by Def.\ \ref{def 56}, that $\mathbb{S}^{obs} (y ,\zeta ,s_1)\leq \mathbb{S} (y ,\zeta,s_1)$.
%
%
Thus, $\mathbb{S}^{obs} (y ,\zeta ,s_1)=\mathbb{S} (y ,\zeta,s_1)$.}
 \hfill \Box \medskip  
  }

The above means that {\gtext the} function $\mathbb{S} (y ,\zeta,s_1)=\mathbb{S}^{obs} (y ,\zeta ,s_1)$
can be reconstructed from the data {\mtext given in the claim}.

\subsubsection{Construction of the {\rtext family} of the earliest light observation sets}
Next we will {\mtext collect  together 
all $\E_{U}(q)$ where $q$ is in an appropriate} geodesic segment.

\begin{lemma}\label{lem: Determination of be-sets} 
Let $s_-\leq s_2\leq s<s_1\leq s_+$ with  $s_1<s_2+\kappa_2$, let  
 ${\hat x}_j={\hat \mu}(s_j)$, $j=1,2$, and  ${\hat x}={\hat \mu}(s)$,
 ${\hat \zeta} \in L^+_{{\hat x}}U$, 
    $\|{\hat \zeta}\|_{g^+}=1$.   
    Let $(y,\zeta) \in L^+U$ be in the $\vartheta_1$-neighborhood of
    $({\hat x},{\hat \zeta})$ such that $y\in J^+({\hat x}_2)$.  Assume that $\gamma_{y,\zeta}(\R_+)$ does not intersect
    $\hat \mu$.   
{\rtext Then, under    {\ftext the} assumptions of Theorem \ref{prop. A} we can determine}  the {\rtext family} $\{\E_{U}(q) 
;\ q\in G_0(y,\zeta,s_1)\}$, where $G_0(y,\zeta,s_1)=\{q\in \gamma_{y ,\zzeta }([{t_0},\infty))\cap  (I^-({p^+}) \setminus J^+({\hat x}_1))\}$.
\end{lemma}  
  
 \noindent
{\bf Proof.} Let $s^\prime=\mathbb{S} (y ,\zzeta ,s_1)$, $x^\prime=\hat\mu(s^\prime)$, and
$\Sigma$ be the set of all {\mtext repeated} observations $S$ associated
to the geodesic $\gamma_{y,\zeta}$ such that $S$ intersects
$\hat\mu([-1,s^\prime))$.
%

Let $q=\gamma_{y,\zeta}(r)\in G_0(y,\zeta,s_1)$. Since
$\gamma_{y,\zeta}(\R_+)$ does not intersect
    $\hat \mu$,  
using a short cut argument for the geodesics from  $q$ to $q_0=\gamma_{y,\zeta}(r_0(y,\zeta,s_1))$
and from $q_0$ to $x^\prime$,
 we see that $f^-_{\hat a}(q)<s^\prime$.
 Then, 
$q\in I^-({p^+}) \setminus J^+({\hat x}_1)$ and $r<r_0(y,\zeta,s_1)$,
and we {\mtext see} using  (\ref{aaaa}) that
$S=\E_{U}(q)$ is a {\mtext repeated} 
observation associated to the geodesic $\gamma_{y,\zeta}$
and $S\cap \hat \mu=\{\hat \mu(f^-_{\hat a}(q))\}$ with $
f^-_{\hat a}(q)<s^\prime$. Thus
$\E_{U}(q)\in \Sigma$ and we conclude that $\E_{U}(G_0(y,\zeta,s_1))\subset \Sigma.$

%

Next, suppose $S\in \Sigma$. Then
there is $\hat \vartheta\in (0,\vartheta_2(y,\zeta,s_1))$ 
 such that for all $\vartheta\in (0,\hat \vartheta)$ 
 there is $(\vec x,\vec\xi)\in \qP_\vartheta(y,\zeta)$   so that
${{\Smod}}(\vec x,\vec \xi)=S$. 
Observe that by (\ref{eq: cut points observed after S}) we have  $\hat\mu([-1,s^\prime))\subset J^-(x^\prime)\subset 
  {\rtext \V(\vec x,\vec \xi)}$. 

First, consider the case when the geodesics corresponding to $(\vec x,\vec\xi)$ do not intersect
 at any point in $I^-(x^\prime)$. Then {\mtext properties (P1) and (P2) yield}
 that ${{\Smod}}(\vec x,\vec \xi)$ is either empty or does not intersect $I^-(x^\prime)$. 
 Thus $S\cap I^-(x^\prime)={{\Smod}}(\vec x,\vec \xi)\cap I^-(x^\prime)$ is empty 
 and $S$ does not intersect
  $\hat\mu([-1,s^\prime))$. Hence $S$ cannot be in $\Sigma$.

Second, consider the case when the geodesics corresponding to $(\vec x,\vec\xi)$
intersect at some point $q\in I^-(x^\prime){\mtext\subset \V(\vec x,\vec\xi)}$. Then,  {\mtext property (P1)}
yields $S=\E_{U}(q)$. 
 {Since $(\vec x,\vec\xi)\in \qP_\vartheta(y,\zeta)$ 
implies $(x_1,\xi_1)=(y,\zeta)$,
 the intersection point $q$ has
a representation $q=\gamma_{x_1,\xi_1}(r)$. As  $q\in I^-(x^\prime)$, this yields
$q\in G_0(y,\zeta,s_1)$ and $S\in \E_{U}(G_0(y,\zeta,s_1))$.
 Hence $\Sigma\subset \E_{U}(G_0(y,\zeta,s_1))$.
 }

%
%
%
%

Combining the above arguments, we    {\ftext conclude}  that $\Sigma=\E_{U}(G_0(y,\zeta,s_1))$.
As  $\Sigma$ is determined by the given data, the claim follows.
\hfill \Box \medskip



{\rtext Now we can complete the proof of Theorem \ref{prop. A}.}
 Let $B(s_2,s_1)$ be the set of all $(y,\zeta,t)$ such that  there are
 $ {\hat x}={\hat \mu}(s)$, $s\in [s_2,s_1)$, $\hat x_j=\hat \mu(s_j)$, $j=1,2,$ and 
 ${\hat \zeta} \in L^+_{{\hat x}}U$,      $\|{\hat \zeta}\|_{g^+}=1$ so that
 $(y,\zeta) \in L^+U$ in $\vartheta_1$-neighborhood of    $({\hat x},{\hat \zeta})$, $y\in J^+({\hat x}_2)$,  and $t\in  [t_0,r_0(y,\zeta,s_1)]$.
   Moreover, let   $B_0(s_2,s_1)$ be the set of all $(y,\zeta,t)\in B(s_2,s_1)$
      \MTEXT{such that  $t<r_0(y,\zeta,s_1)$}
     and 
    $\gamma_{y,\zeta}(\R_+)\cap \hat \mu=\emptyset$.
\MTEXT{Using Lemma \ref{lem: Determination of be-sets} {\rtext we then can determine 
the collection}
$\Sigma_0(s_2,s_1):=\{\E_{U}(q);  
\ q=\gamma_{y ,\zeta }(t),\ (y,\zeta,t)\in B_0(s_2,s_1)\}$. 
We denote also $\Sigma(s_2,s_1):=\{\E_{U}(q);  
\ q=\gamma_{y ,\zeta }(t),\ (y,\zeta,t)\in B(s_2,s_1)\}$.

{\mtext Recall} that the sets $\E_{U}(q)\subset {U}$, where $q\in J:=J^-(p^+)\cap J^+(p^-)$, can be identified
with the {\mattitext continuous} function, $F_q:\overline \A\to \R$,
$F_q(a)=f^+_{a}(q)$, c.f.\ (\ref{kaava D}). 
When we  endow the set $C(\overline \A)$ of {\mattitext  continuous} maps $\overline \A\to \R$
with the topology of {\mattitext uniform} convergence, 
Lemma \ref{B: lemma} yields that $F:q\mapsto F_q$ is continuous
map $F:J\to C(\overline \A)$. 
By
Proposition \ref{propo paper part 2},  
$F:J\to F(J)$ is homeomorphism.
Next, we identify $\E_{U}(q)$ and $F_q$. {\mtext Also, on
the space {\mmmmtext $\E_{U}(J)$} we will use the topology
that makes the map $\E_{U}\circ F^{-1}:F(J)\to \E_{U}(J)$ a homeomorphism}.
%

{Using standard results of differential topology, we {\gtext have} that any neighborhood
of $(y,\zeta)\in L^+U$ contains  $(y^\prime,\zeta^\prime)\in L^+U$
such that the geodesic $\gamma_{y^\prime,\zeta^\prime}([0,\infty))$
does not intersect $\hat \mu$. 
Since $(y,\zeta)\mapsto r_0(y,\zeta,s_1)$ is lower semicontinuous, this implies that $\Sigma_0(s_2,s_1)$ is dense in $\Sigma(s_2,s_1)$. 
Hence we obtain
the closure $\overline \Sigma(s_2,s_1)$ of $\Sigma(s_2,s_1)$
as  
the limits points of $ \Sigma_0(s_2,s_1)$.

\MTEXT{Then, we obtain the set $\E_U( J^+({\hat \mu}(s_2 ))\cap J^-({p^+}))$
as
the union  $\overline\Sigma(s_2,s_1)\cup
\E_U(J^+({\hat \mu}(s_1))\cap J^-({p^+}))\cup \E_U(\K_{t_0}\cap J^+({\hat \mu}(s_2 )))$,  
see (\ref{set Kt0}).

Let  
$s_0,\dots,s_K\in [s_-,s_+]$ be such that  
 $s_j>s_{j+1}>s_j-\kappa_2$ and $s_K=s_-$.
Then, by iterating the above construction so that the values of the parameters $s_1$ and $s_2$ are replaced by
$s_j$ and $s_{j+1}$, respectively,
we can construct the set $\E_{U}( J^+(\hat \mu (s_-))\cap J^-(\hat \mu (s_+)))$.}}

Moreover, similarly to the above construction, we can find the   
sets $\E_{U}( J^+(\hat \mu (s^\prime))\cap J^-(\hat \mu (s^{\prime\prime}))$  
for all $s_-<s^\prime<s^{\prime\prime}<s_+$,
and taking their union, we {\gtext construct} the set 
$\E_{U}( I(\hat \mu (s_-),\hat \mu (s_+))$. This proves Theorem \ref{prop. A} (i).
{\mtext As the claim (ii) was already proven earlier, this finishes the proof of Theorem \ref{prop. A}.}
   \hfill \Box \medskip

{\mtext
   {\ftext After}  Theorem \ref{main thm paper part 2} is proven,
Theorem  \ref{main thm3 new} will follow from  Lemma \ref{lem: aux construction} and Theorem \ref{prop. A}.
Thus we consider next the proof of Theorem \ref{main thm paper part 2}
and return later to the proof of Theorem  \ref{main thm3 new}.}}


\section{Solution of the inverse problem for passive observations} 
\label{sec: Proof for passive observations}

In this section we prove Theorem \ref{main thm paper part 2}. 
 {\mattitext Earlier we have shown that the conformal type of $(U,g|_U)$
and {\gtext the} {\rtext family} $\E_U(W)$  determine the {\rtext family} of the 
earliest observation times $\F(W)\subset C(\overline\A)$.
 In the proof of  Proposition \ref{propo paper part 2}
we have shown that the topology of $C(\overline\A)$ induces 
on the set $\F(W)$ a topology that makes it homeomorphic
to $W$. In this section we 
construct on $\F(W)$ smooth coordinates 
   {{\gtext and show that then $\F(W)$ is}
diffeomorphic to $W$. After this we construct a metric on $\F(W)$ 
that makes it conformal to $(W,g|_W)$.}  

The proof  is constructive and to simplify the notations, we do    {\ftext the} 
constructions on just one Lorentzian manifold, $(M, g)$
and assume that we are given the data 
\beq\label{eq: data}
& &\hbox{the differentiable manifold $U$, the conformal class of $g|_U$, 
}\\ \nonumber
& &\hbox{the paths $\mu_a:[-1,1]\to U$, $a\in \A$,
 and the set $\E_U(W)$,}
\eeq 
where $W$ is a relatively compact open set such that
 $\overline  W \subset I^-(p^+) \setminus J^-(p^-)$.


In this section we do not use {\gtext the functions} $f^-_a$ and denote 
\ba
f_a(x)=f^+_a(x).
\ea
{\mattitext Also we use the notations defined in Section \ref{sec: Reconstruction of topology}.}

\subsection{Construction of the differentiable structure}
Let us next consider the set ${\mathcal Z}=\{(q,p)\in W\times U;\ p\in \E_U^{reg}(q)\}$. 
For every $(q,p)\in {\mathcal Z}$ there is a unique $\xi\in L^+_{q}M$ such that $\gamma_{q,\xi}(1)=p$
and $\rho(q,\xi)>1$.  We will denote 
$\Theta(q,p)=(q,\xi)$ that defines a map
by $ \Theta:{\mathcal Z}\to L^+W$. 
Below, let $\W_\e(q_0,\xi_0)\subset TM$ be an $\e$-neighborhood
of $(q_0,\xi_0)$ with respect to the Sasaki-metric induced by $g^+$ on $TM$.

\begin{lemma}
Let $(q_0,p_0)\in {\mathcal Z}$ and $(q_0,\xi_0)=\Theta(q_0,p_0)$.
When $\e>0$ is small enough, the map 
\beq \label{e.16.0}
& &X: \W_\e(q_0,\xi_0) \to  M\times M,\quad
X(q, \xi)= (q,\exp_q(\xi))
\eeq
is open and defines a diffeomorphism $X: \W_\e(q_0,\xi_0) \to \U_\e(q_0,p_0):=X(\W_\e(q_0,\xi_0))$.  When $\e$ is small enough, $\Theta$ coincides  in $\modified{\mathcal Z\cap} \U_\e(q_0,p_0)$ with
the inverse map of $X$.  Moreover, \modified{$\mathcal Z$ is a $(2n-1)$-dimensional manifold and} the map  $\Theta:{\mathcal Z}\to L^+M$  is $C^\infty$-smooth.
\end{lemma}

\noindent{\bf Proof.}
{\mattitext We start the proof with  {\gtext some} technical considerations.}

Let $p_{+2}=\hat \mu(s_{+2})$, $s_+<s_{+2}<1$.
For $(x,\xi)\in L^+M$, $x\in J^-(p^+)$, the value $T_{+2}(x,\xi)=\sup \{t\geq 0\ ;\ \gamma_{x,\xi}(t)\in J^-(p_{+2})\}$,
 is finite by \cite[{\ftext Lemma}\ 14.13]{ONeill}.
Since $J_2=J^-(p_{+2})$ is closed and $\gamma_{x,\xi}$, $(x,\xi)\in L^+M$, are future-pointing  {\muutos path}s, we have that $T_{+2}:L^+J_2 \to \R$ is upper semicontinuous. As $W\subset J^-(p_{+2})$ is relatively compact, the 
set  
\beq\label{eq: set K}
K=\{(x,\xi)\in L^+M;\ x\in \hbox{cl}\,(W),\ \|\xi\|_{g^+}=1\}
\eeq is compact and there is 
$c_0\in \R_+$ such that $T_{+2}(x,\xi)\leq c_0$ for all $(x,\xi)\in K$.

{\mattitext Let us now start the proof of the claim.}
Since the geodesic $\gamma_{q_0,\xi_0}([0,1])$ does not contain cut points and \modified{thus
 conjugate points}, we see that when $\e>0$ is small enough, the set
$\U_\e(q_0,p_0)=X(\W_\e(q_0,\xi_0))\subset M\times M$ is open and
the map
$X:  \W_\e(q_0,\xi_0)\to \U_\e(q_0,p_0)$
has a $C^\infty$-smooth inverse map $X^{-1}:\U_\e(q_0,p_0)\to  \W_\e(q_0,\xi_0)$.
Thus $\modified{X:} \W_\e(q_0,\xi_0) \to  \U_\e(q_0,p_0)$  is a diffeomorphism. 
\modified{Note that
$X^{-1}(q_0, p_0) =(q_0, \xi_0)$.}

{{\mattitext First, we prove that   $ \Theta:{\mathcal Z}\to L^+W$ is continuous.}
If  $ \Theta:{\mathcal Z}\to L^+W$ would not be continuous at $(q_0,p_0)\in{\mathcal Z}$,
there would exists a sequence $(q_k,p_k)\in{\mathcal Z}$ converging to 
$(q_0,p_0)$ as $k\to \infty$, such that $\Theta(q_k,p_k)\in L^+M$
does not converge to $(q_0,\xi_0)=\Theta(q_0,p_0)$.

Since $p_k\in J^-(p_{+2})$ and the function 
$T_{+2}$ is bounded by $c_0\in \R_+$ in the set $K$ given in (\ref {eq: set K}),
the sequence $\|\Theta(q_k,p_k)\|_{g^+}$ is uniformly bounded. By
considering a subsequence we may assume that
$\Theta(q_k,p_k)\to (q_0,\eta)\in L^+M$ as $k\to \infty$ and  $\eta\not =\xi_0$.
In this case the geodesics $\gamma_{q_0,\xi_0}([0,1])$  and 
$\gamma_{q_0,\eta}([0,1])$ would be two light-like geodesics
connecting $q_0$ to $p_0$ so that $\rho(q_0,\xi_0)\leq 1$.
This would be in contradiction with
the assumption that $p_0\in \E_U^{reg}(q_0)$. This shows that 
 $ \Theta:{\mathcal Z}\to L^+W$ is continuous at $(q_0,p_0)$.
 
 Let $\e_1\in (0,\e)$  and ${\mathcal Z}\cap  \U_{\e_1}(q_0,p_0)$ be a neighborhood of $(q_0,p_0)$ in the 
relative topology of ${\mathcal Z}\subset W\times U$.
When $\e_1$  is small enough, we have 
$\Theta({\mathcal Z}\cap  \U_{\e_1}(q_0,p_0))\subset \W_\e(q_0,\xi_0)$.
Then for $(q,p)\in {\mathcal Z}\cap  \U_{\e_1}(q_0,p_0)$ and $(q,\xi)=\Theta(q,p)\in \W_\e(q_0,\xi_0)$
  we have $\exp_q(\xi)=p$, and hence
 $X(\Theta(q,p))=(q,p)$. 
 Since $\Theta(q,p)\in \W_\e(q_0,\xi_0)$, we have
 $\Theta(q,p)=X^{-1}(q,p)$.
 {\gtext Therefore} for $(q,p)\in {\mathcal Z}\cap  \U_{\e_1}(q_0,p_0)$ 
the function $\Theta:{\mathcal Z}\cap  \U_{\e_1}(q_0,p_0)\to TM$ coincides with the smooth
function  $X^{-1}:{\mathcal Z}\cap  \U_{\e_1}(q_0,p_0)\to TM$. 
{\gtext Given that} $(q_0,p_0)\in {\mathcal Z}$ is arbitrary, this shows that
\modified{$\mathcal Z$ is a $(2n-1)$-dimensional manifold and
 $\Theta:{\mathcal Z}\to L^+M$} is $C^\infty$-smooth.
\hfill \Box \medskip

%
%
%
%
}

%
%
%

\begin{proposition} \label{lemma coordinates} Let $q_0\in  I^-(p^+)\setminus J^-(p^-) $
and  $(q_0,p_j)\in {\mathcal Z}$, $j=1,2,\dots,n$ and $\xi_j \in L^+_{q_0}M$
be such that $\gamma_{q_0,\xi_j}(1)=p_j$. Assume that $\xi_j$,
$j=1,2,\dots,n$ are linearly independent. Then, if $a_j\in \A$ and
$\vec a=(a_j)_{j=1}^n$  are such that $p_j\in \mu_{a_j}$,
there is a neighborhood $V_1\subset M$ of $q_0$ such that
the corresponding observation time functions
\ba
\bbf_{\vec a}(q)=(f_{a_j}(q))_{j=1}^n 
\ea
define  $C^\infty$-smooth coordinates  in $V_1$. 
Moreover,  $\nabla f_{a_j}|_{q_0}$, the gradient of $f_{a_j}$ \modified{with respect to $q$} at $q_0$
satisfies $\nabla f_{a_j}|_{q_0}=c_j \xi_j $
{\gtext for} some $c_j\not =0$.

\end{proposition}
\noindent
{\bf Proof.} 
Let $(q_0,p_0)\in {\mathcal Z}$ and $\xi_0\in L^+_{q_0}M$ such that $\gamma_{q_0,\xi_0}(1)=p_0$.
Moreover, let 
$\e>0$ be so small that the map
$X:  \W_\e(q_0,\xi_0)\to \U_\e(q_0,p_0)$
 has a $C^\infty$-smooth inverse, see (\ref{e.16.0}).
  We denote this inverse map by 
  \beq\label{xi map}
  X^{-1}(q, p) =(q, \xi(q, p)).
  \eeq {\mattitext Recall that
  we denote $\W=\W_\e(q_0,\xi_0)$ and $\U=\U_\e(q_0,p_0)$.}

We associate with any $(q, p) \in \W$ the energy
$E(q, p)= E(\g_{q, \xi(q,p)}([0, 1]))$ of the geodesic segment $ \g_{q, \xi(q,p)}([0, 1])$ from $p$ to $q$.
Here, the energy of a piecewise smooth  {\muutos path} $\alpha:[0,l]\to M$
is defined by
\ba
E(\a)=\frac 12\int_0^l g(\dot \alpha(t),\dot \alpha(t))\,dt.
\ea
Observe that the sign of $E(q, p)$ depends on the causal nature of $\g_{q, p}$. In particular,
$E(q,p)=0$ if and only if  $\xi(q,p)$ is light-like. Moreover, since $X^{-1}$ is $C^\infty$-smooth on $\U$, also $E(q,p)$ is $C^\infty$-smooth in $\U$.

Let us return to consider $(q_0,p_0) \in{\mathcal Z}$ and let ${a_0}\in \A$ be such that $p_0\in \mu_{a_0}$.
 Then $p_0=\mu_{a_0}(s_0)$ with $s_0=f_{a_0}(q)$.

 \modified{Let $V_0\subset W$ be an open neighborhood of $q_0$ and $t_1,t_2\in (s_{-2},s_{+2})$,  $  t_1<s_0 <t_2$
 be such that}
$V_0\times \mu_{a_0}([t_1, t_2])\subset \U$. Then for $q\in V_0$ and $s\in (t_1,t_2)$  the function
${\bf E}_{a_0}(q, s):=E(q, \mu_{a_0}(s))$ is well defined and smooth.
 Using the first variation formula for ${\bf E}_{a_0}(q,s)$, see e.g.\ \cite[Prop. 10.39]{ONeill},  we obtain
 \beq \label{e.16.3}
 & & \modified{\frac{\p {\bf E}_{a_0}(q_0, s)}{\p s}\bigg|_{s=s_0}=  g\left(\eeta,\, \dot \mu_{a_0}( f_{a_0}(q_0)) \right),
 \quad
\nabla {\bf E}_{a_0}(q, s_0)\bigg|_{q=q_0}=- \xi_0} ,
 \eeq
 where $\xi_0=\xi(q_0,p_0)$ and $\eeta=  \dot \g_{q_0, \xi_0}(1)$, see (\ref{xi map}). 
Since $\dot \mu_{a_0}(s)$ is time-like and future-pointing and
  $\eeta$ is light-like and future-pointing, 
 $ \frac{\p {\bf E}_{a_0}}{\p s}(q_0, s_0) <0$.
 
It follows from the implicit function
 theorem that  there is  an open neighborhood $V_{a_0}\subset V_0$ 
 of $q_0$ 
and a smooth function $q\mapsto s(q,{a_0})$ defined for $q \in V_{a_0} $ such that 
  $s(q_0,{a_0})= f_{a_0}(q_0)$  and
 ${\bf E}_{a_0}(q, s(q,{a_0}))=0$.  \MTEXT{Then
 $q\mapsto s(q,{a_0})$ and $q\mapsto f_{a_0}(q)$ coincide in $V_{a_0}$,  and
 it follows from (\ref{e.16.3}) that 
 \beq
 \label{e.16.4} 
   \modified{\nabla f_{a_0}(q)\bigg|_{q=q_0} = \frac{1}{c(q_0,{a_0})} \xi_0,\,\, c(q_0,{a_0})= \frac{\p {\bf E}_{a_0}}{\p s}(q_0, s))\bigg|_{s=f_{a_0}(q_0)}}.\hspace{-1cm}
 \eeq

Next we choose
$p_1,p_2,\dots,p_n\in \E_U^{reg}(q_0)$ and 
let 
 $\xi_1, \dots, \xi_n \in L^+_{q_0}(M)$ 
be such that $p_i=\gamma_{q_0,\xi_i}(1)$. We assume that 
 $\xi_1, \dots, \xi_n \in L^+_{q_0}(M)$ are linearly independent.
 Moreover, let 
 $a_j\in \A$ be such
 that $p_j\in \mu_{a_j}$ and  $\vec a=(a_j)_{j=1}^n$. Finally,
 we denote by $q\mapsto s(q,a_j)$ the above constructed smooth
 functions that are defined in some neighborhoods $V_{a_j}\subset W$ of $q_0$

 Let  $V_{\vec a}=\bigcap_{j=1}^n V_{a_j} $ and consider the map
\ba
\bbf_{\vec a}: V_{\vec a}\to \R^n,\quad \bbf_{\vec a}(q)=(f_{a_1}(q), \dots, f_{a_n}(q)).
\ea
 It follows from (\ref{e.16.4}) that the map $\bbf_{\vec a}$ has {\gtext an invertible} differential at $q_0$ and, therefore, the function 
 $\bbf_{\vec a}: V_{\vec a}\to \R^n$
 defines a $C^\infty$-smooth coordinate
system in some neighborhood of $q_0$.
\hfill \Box \medskip

\subsubsection{Properties of the $C^0$ and $C^\infty$ smooth coordinates}

\begin{definition}\label{def:-observation  coordinates}
Let $\vec a=(a_j)_{j=1}^n\in \A^n$,  
${\mathcal O}\subset \F(W)$ be an open set, $s_{a_j}=f_{a_j}\circ \F^{-1}$,
and  $\bs_{\vec a}=\bbf_{\vec a}\circ \F^{-1}$.
We say that $({\mathcal O},\bs_{\vec a})$
are $C^0$-observation  coordinates on $\F(W)$  if the map
$\bs_{\vec a}:{\mathcal O}\to \R^n$ is an open and injective map.
Also,   we say that  $({\mathcal O},\bs_{\vec a})$
are $C^\infty$-observation  coordinates on $\F(W)$ if
$ \bs_{\vec a}\circ \F:\F^{-1}({\mathcal O})\to \R^n$ are
$C^\infty$-smooth local coordinates on $W\subset M$,
 see Fig.\ 1(Right).
\end{definition}

\modified{Note that, by the invariance of domain theorem, the above continuous map $\bs_{\vec a}:{\mathcal O}\to \R^n$ is open if it is injective.} Even though for a given $\vec a\in \A^n$
there are several sets ${\mathcal O}$  for which 
$({\mathcal O},\bs_{\vec a})$ form $C^0$-observation  coordinates, to clarify
the notations, we sometimes 
denote the coordinates $({\mathcal O},\bs_{\vec a})$  by $({\mathcal O}_{\vec a},\bs_{\vec a})$.

Since $\F:W\to \F(W)$ is a homeomorphism, we can
determine all $C^0$-observation  coordinates on $\F(W)$ using
data (\ref{eq: data}).  
 Next we will consider $\F(W)$ as a topological manifold 
 endowed with the $C^0$-observation  coordinates and denote $\F(W)=\bW$.
 We  denote the points of this manifold by $\q=\F(q)$.
Next we construct  a differentiable structure on $\bW$ \modified{that is compatible with that of $W$}.

\subsubsection{Construction of the $C^\infty$ smooth coordinates}

\begin{lemma}\label{lemma: testing coordinates}
Assume that we are given data (\ref{eq: data}). Then for  any $C^0$-observation  coordinates
 $({\mathcal O}_{\vec a},\bs_{\vec a})$ with 
  $\vec a\in \A^n$ we can determine if  $({\mathcal O}_{\vec a},\bs_{\vec a})$ 
are $C^\infty$-observation  coordinates on $\bW$. Moreover,
for any $\q\in \bW$ there exists $C^\infty$-observation  coordinates
$({\mathcal O}_{\vec a},\bs_{\vec a})$ such that $\q\in {\mathcal O}_{\vec a}$.
%
\end{lemma}

%

 \noindent
{\bf Proof.}  
Let $q\in W$.  \MTEXT{We say that $p\in \E_U(q)$ and $a\in \A$ are associated
if $p\in \mu_a$. Next, consider  $p\in \E_U^{reg}(q)$ and $a\in \A$ that  are associated.
Note that then $q\not\in \mu_a$.}
By (\ref{e.16.4}), 
the function $f_a(q)$ satisfies 
\ba
\modified{\nabla f_a(q)=c(q,a)}\, \xi(q,y),\quad c(q,a) \not =0,
\ea
where $y=\mu_{a}(f^+_{a}(q))$, see (\ref{xi map}).
Let 
\ba
K(q)=\{(\xi_j)_{j=1}^n;\  \xi_j\in L^+_qM,\ 
\modified{ \rho( {q,\xi_j})>1},\ \gamma_{q,\xi_j}(1)\in U\}
 \ea 
 and $H:K(q)\to U^n$
 be the map $H((\xi_j)_{j=1}^n)=(p_j)_{j=1}^n,$ where $p_j= \gamma_{q,\xi_j}(1)$.
 Then $p_j\in \E_U^{reg}(q)$ and $\xi_j=\Theta(q,p_j)$. 
 {\gtext Given that} $\rho$ is lower semi-continuous,
 we {\gtext have} that $K(q)\subset (L^+_qM)^n$ is open. 
 Clearly, $H$ is  $C^\infty$-smooth. {\gtext Since} $ \Theta:{\mathcal Z}\to L^+W$ is continuous and injective, we see that
 $H:
K(q)\to H(K(q))=(\E_U^{reg}(q))^n$ is a  homeomorphism.
We denote below $Y(q)=(\E_U^{reg}(q))^n$. Note that
for all $\q\in \bW$
 the  data (\ref{eq: data})
determine the set $Y(q)\subset U^n$, where $q=\F^{-1}(\q)$.

  Let us consider the set 
 \ba
 K_0(q)=\{(\xi_j)_{j=1}^n\in K(q);\hbox{
 $\xi_j$, $j=1,\dots, n$ are linearly independent}\}.
 \ea
 Clearly, the set  $K_0(q)$ is dense and open in $K(q)$,
 and hence $Y_0(q):=H(K_0(q))$ is open and dense in
 $Y(q)$.

%

Let   $({\mathcal O}_{\vec a},\bs_{\vec a})$, $\vec a\in \A^n$
be $C^0$-observation  coordinates on $\bW$,  $\q\in {\mathcal O}_{\vec a}$,
and $q=\F^{-1}(\q)$.
\MTEXT{Also, let  $(p_j)_{j=1}^n\in Y(q)$  be
 such that {\gtext the $p_j$'s} are associated with $a_j$. Similarly, let  $({\mathcal O}_{\vec b},\bs_{\vec b})$, $\vec b\in \A^n$
be another $C^0$-observation  coordinates on $\bW$ such that
$\q\in {\mathcal O}_{\vec b}$, and let 
 $(z_j)_{j=1}^n\in Y(q)$ be such that $z_j$ are associated with $b_j$.
 Note that then 
 $p_j=\firstpoint^{reg}_{a_j}(q)$ and  $z_j=\firstpoint^{reg}_{b_j}(q)$.}


In the case when $z_j)_{j=1}^n\in Y_0(q)$,
$q$ has a neighborhood $V_1\subset W$ in which the function  $\bbf_{\vec b}:V_1\to \R^n$
give $C^\infty$-smooth local coordinates. 
Thus, if $(z_j)_{j=1}^n\in Y_0(q)$, then it holds that $(p_j)_{j=1}^n\in Y_0(q)$ if and only if
 \ba
  & &\hspace{-5mm}(i)\hbox{The functions $s_{a_j}\circ \bs_{\vec b}^{-1}$, $j=1,2,\dots,n$ are
$C^\infty$-smooth at $\bs_{\vec b}(\q)$  and}
\\
 & &\hbox{ the Jacobian determinant $\det(D(\bs_{\vec a}\circ \bs_{\vec b}^{-1}))$ 
 at  $\bs_{\vec b}(\q)$ is non-zero.}
\ea

%
%
%
%
%
\MTEXT{
Denote $\vec p=(p_j)_{j=1}^n\in Y(q)$, and define
$\X_{\vec p}\subset Y(q)$ to be the set of those $(z_j)_{j=1}^n\in Y(q)$,
 for which there are   $\vec b\in \A^n$ and $C^0$-observation  coordinates 
$({\mathcal O}_{\vec b},\bs_{\vec b})$ such that $\q=\F(q)\in  {\mathcal O}_{\vec b}$,
$z_j$ are associated with $b_j$ for $j=1,2,\dots,n$,
and the condition (i) is satisfied.
If  ${\vec p}$ is in $Y_0(q)$, we see that  $Y_0(q)\subset \X_{\vec p}$.
On the other hand,  if ${\vec p}$ is  not in $Y_0(q)$, we have  $Y_0(q)\cap \X_{\vec p}=\emptyset$.
Since the set $Y_0(q)$ is open and dense in $Y(q)$,
we {\gtext observe} that  ${\vec p}\in Y_0(q)$ if and only if the interior of set $\X_{\vec p}$ is a dense
subset of $Y(q)$.}
%
%
%
This in particular implies that using the data (\ref{eq: data}) 
 we  can determine whether  $(p_j)_{j=1}^n$ is in $Y_0(q)$ or not. {\mattitext As we can do the 
 above considerations {\gtext for} all $\q\in  {\mathcal O}_{\vec a}$, we {\gtext get} that}
the $C^0$-observation  coordinates $({\mathcal O}_{\vec a},\bs_{\vec a})$, $\vec a\in \A^n$ 
are $C^\infty$-observation  coordinates on $\bW$ if and only if
for all   $\q\in  {\mathcal O}_{\vec a}$, $q=\F^{-1}(\hat q)$, and  $p_j=
\mu_{a_j}(f^+_{a_j}(q))$, $j=1,2,\dots,n$
 we have $(p_j)_{j=1}^n\in Y_0(q)$.
Thus we can determine all $C^0$-observation  coordinates
 $({\mathcal O}_{\vec a},\bs_{\vec a})$
 on $\bW$ that  are $C^\infty$-observation  coordinates. Moreover,
 since for all $q\in W$ the set $Y_0(q)$ is non-empty, we {\gtext obtain} that
 any $\q=\F(q)\in \bW$ belongs in the domain some $C^\infty$-observation  coordinates.
 \hfill \Box \medskip

%

We endow $\bW=\F(W)$ with the differentiable structure provided by all $C^\infty$-observation  coordinates
on $\bW$.  
By Lemma \ref{lemma: testing coordinates} and \cite[{\ftext Lemma}\ 1.42]{ONeill}
the $C^\infty$-observation  coordinates
 make $\bW$ a differentiable manifold and
its the differentiable structure is uniquely  determined.
Since the differentiable structure of $W$ is determined by
the functions $f_{\vec a}$ that are $C^\infty$-smooth local coordinates, 
we {\gtext have} using Def.\ \ref{def:-observation  coordinates} that 
the map
\beq \label{e.18.1}
\F: W \to \bW=\F(W)
\eeq
is a diffeomorphism.

\noindent

\subsection{Construction of the conformal type of the metric}
Let us  denote by ${ \underhat  g}=\F_*g$ the metric on $\underhat  W= \F(W)$ 
that makes $\F$  an isometry. Next we show
that the set $\F(W)$, the paths $\mu_a$ and the conformal
class of the metric $g$ on $U$ determine the conformal class of 
${ \underhat  g}$ on ${ \underhat  W}$.

\begin{lemma}
The data (\ref{eq: data}) determines 
a metric $G$ on $\bW=\F(W)$ that is conformal to ${ \underhat  g}$
and the time orientation on $\bW$ that makes $\F:W\to \bW$
a causality preserving map.
\end{lemma}

 \noindent
{\bf Proof.} Let $({\mathcal O}_{\vec a},\bs_{\vec a})$
be $C^\infty$-observation  coordinates on $\bW$.
Then by (\ref{e.16.4}) the co-vectors
$-ds_{a_1}|_\q$ and $-ds_{a_2}|_\q$ are non-parallel future-pointing light-like \modified{co-vectors}.
Thus their sum determines a  future-pointing time-like co-vector field on ${\mathcal O}_{\vec a}$. Using a suitable
partition of unity we can construct {\gtext a}  future-pointing time-like co-vector field $X$ on $\bW$.

Let  
$({\mathcal O}_{\vec a},\bs_{\vec a})$ be $C^\infty$-observation  coordinates on $
\underhat W$.
Let  $\q\in {\mathcal O}_{\vec a}$ and 
$q\in W$ be such that $\q=\F(q)$.
Using the data (\ref{eq: data}), the function $F_q=\F(q):\A\to \R$,  and the formula (\ref{eq BB},
we can determine the set $\E_U(q)\subset U$.
By Prop.\ \ref{lemma: from P(q) to E(q)} (iii),
this further determines  the set
$\be_U^{reg}(q)$.
%
%

%

{Then, let  us fix a point $\q=\F(\qpoint) \in {\mathcal O}_{\vec a}$. 
Let $(y,\eeta)\in \be_U^{reg}(\qpoint)$ and  let $\hat t>0$
be the largest number such that the geodesic $\gamma_{y,\eeta}((-\hat t,0])\subset M$ is defined and has no cut points.
%
%
For $q\in W$, Proposition \ref{lemma: properties of  E(q)} (ii) yields that
 $q\in \gamma_{y,\eeta}((-\hat t,0))$ if and only if
$(y,\eeta)\in \be_U^{reg}(q)$.
}
%
%
Hence $(y,\eeta)$ and  the data (\ref{eq: data})  determine the set 
\ba
\beta=\{\q\in {\mathcal O}_{\vec a};\ \q=\F(q),\  \be_U^{reg}(q)\ni (y,\eeta)\}=\modified{\F(\gamma_{y,\eeta}((-\hat t,0)))\cap
{\mathcal O}_{\vec a}.}
\ea
This implies that on ${\mathcal O}_{\vec a}\subset \bW$ we can find the
 image, in the map $\F$, of 
the light-like geodesic segment 
$\gamma_{y,\eeta}((-\hat t,0))\cap \F^{-1}({\mathcal O}_{\vec a})$ that contains  $\qpoint=\gamma_{y,\eeta}(-t_1)$. 
Let $\a(s)$, $s\in (-s_0,s_0)$ be a smooth path on ${\mathcal O}_{\vec a}$ such that $\p_s\a(s)$
does not vanish, $\alpha((-s_0,s_0))\subset \beta$, and $\alpha(0)=\hat \qpoint$.
Such smooth path  $\a(s)$ can \modified{be} obtained  e.g.\ by parametrizing $\beta$ by
 arc-length with respect to
some auxiliary smooth Riemannian metric on ${\mathcal O}_{\vec a}$.
Then $\underhat \xi=\p_s\a(s)|_{s=0}\in T_{\hat \qpoint}\bW$  has the form
$\underhat \xi=c\F_*(\dot\gamma_{y,\eeta}(t_1))$ where $c\not=0$.
Since we can do the above construction for all points  $(y,\eeta)\in \be_U^{reg}(\qpoint)$, 
we determine   in the tangent space $T_{\hat \qpoint}\bW$ 
the set \modified{$\Gamma=\F_*(\{c\xi\in L_{\qpoint}M;\ \exp_{\qpoint}(\xi)\in \E_U^{reg}(\qpoint),\
c\in \R,\ c\not=0\}),$}
 that is an open, non-empty subset
of the light cone at $\hat \qpoint$ associated to the metric $\underhat g$.
Let us now consider the set $\Gamma$ in the coordinates of
$T_{\hat \qpoint}\bW$ associated to $\bs_{\vec a}$.
Since   the light cone is determined by a quadratic equation in the
tangent space, having an open set $\Gamma$ of  the light cone we can  uniquely determine the whole light cone.
Using this construction with all points $\hat \qpoint\in {\mathcal O}_{\vec a}$, 
we can determine all light-like vectors
in the tangent space  $T_{\hat \qpoint}{\mathcal O}_{\vec a}$ for all  $\hat \qpoint\in {\mathcal O}_{\vec a}$.
The collections of  light-like vectors at tangent
spaces of $\underhat W$ determine uniquely the conformal
class of the tensor ${ \underhat  g}=\F_*g$ in the manifold $ \underhat W$,
see  \cite[Thm.\ 2.3]{Beem} (or  \cite[{\ftext Lemma}\ 2.1]{Beem} \modified{for} a constructive procedure).

The above shows that 
the data (\ref{eq: data})  determines the conformal class of the  metric tensor 
 ${ \underhat  g}$. In particular, we can construct a metric  $G$ on $ \underhat W$ that is conformal to 
 ${ \underhat  g}$ and satisfies $G(X,X)=-1$.
 \hfill \Box \medskip

We have shown  that the data (\ref{eq: data})
determine the topological and the differentiable structures on 
$\underhat W=\F(W)$ and a metric $G$ on it
that makes the map $\F:(W,g|_W) \to (\bW,G)$ a diffeomorphism
 and a conformal map. Moreover, we determine the time-orientation on
 $\underhat W$ that makes $\F$ a causality preserving map.

 Finally, by Prop.\ \ref{lemma: from P(q) to E(q)}  {\it (i)}, for any $y\in U$
 we can {\mmmmtext verify} \modified{if $y=q\in W$ and find the corresponding element $\F(q)\in \F(W)$.}
Thus we can find  the set $\F(W\cap U)$
and the map $\F^{-1}:
\F(W\cap U)\to W\cap U$.  This yields {\mattitext the claim (ii) of} Thm.\ \ref{main thm paper part 2}.
Thus  Theorem \ref{main thm paper part 2} is proven.
\hfill \Box \medskip

\subsubsection{Construction of the conformal factor in the vacuum spacetime}
$  $
\smallskip

\noindent
{\bf Proof of Corollary \ref{coro of main thm original Einstein pre}.}
By Theorem \ref{main thm paper part 2}, there is a  conformal
diffeomorphism $\Psi:(W_1,g^{(1)})\to (W_1,g^{(2)})$. By our assumptions,   $\Phi:(V_1,g^{(1)})\to (V_1,g^{(2)})$
is an isometry, the
Ricci curvature of $g^{(j)}$ is zero 
in $W_j$, and 
  any point $x_1\in W_1$ is connected to some point $y_1\in V_1\cap
W_1$ with a 
\modified{piecewise smooth path $\mu_{y_1,x_1}([0,1])\subset W_1$,
$\mu_{y_1,x_1}(0)=y_1$.
Note that then $\Psi(\mu_{x_1,y_1}([0,1]))\subset W_2$
connects $x_2=\Psi(x_1)$ to  $y_2=\Psi(y_1)$.}

To simplify notations we  denote $\hat g= g^{(1)}$ and 
$g= \Psi^*g^{(2)}$. 
Since $\Psi$ is conformal, there is $f:W_1\to \R$ such that \modified{$\hat g=e^{2f} g$} on $W_1$, and as $\Phi:V_1\to V_2$ is an isometry,  \modified{$f=0$ in ${V_1}$}.
 By \cite[formula (2.73)]{Rendall},
 the
 Ricci tensors $ \hbox{Ric}_{jk}(g)$ of $g$ and ${\hbox{Ric}}_{jk}(\hat g)$ of $ \hat g$ satisfy
 on  $W_1$ 
 \ba
0={\hbox{Ric}}_{jk}(\hat g)&=&{\hbox{Ric}}_{jk}(g)-2\nabla^{}_{j}\nabla_{k}^{}f+2(\nabla_{j}^{}f)(\nabla_{k}^{}f)\\
\hspace{2cm}& &- (g^{pq}\nabla^{}_{p}\nabla_{q}^{}f+2g^{pq}(\nabla_{p}^{}f)(\nabla_{q}^{}f))g_{jk}
\ea
where $\nabla=\nabla^g$. For the scalar curvature this yields
 \ba
0=e^{2f}\hat g^{pq}{\hbox{Ric}}_{pq}(\hat g)=
g^{pq}{\hbox{Ric}}_{pq}(g)-3g^{pq}\nabla^{}_{p}\nabla_{q}^{}f.
\ea
Combining the above \modified{with the fact that $\hbox{Ric}_{jk}(g)=0$, we obtain
\ba
\nabla^{}_{j}\nabla_{k}^{}f\hspace{-1mm}-\hspace{-1mm}(\nabla_{j}^{}f)(\nabla_{k}^{}f)
\hspace{-1mm}+\hspace{-1mm}g^{pq}(\nabla_{p}^{}f)(\nabla_{q}^{}f)g_{jk}\hspace{-1mm}=\hspace{-1mm}0.
\ea
This equation gives a system of  first order ordinary differential equations for the vector field $Y=\nabla f$
along $\mu_{y_1,x_1}([0,1])$ with initial value $Y(y_1)=\nabla f(y_1)=0$,
that has the unique solution $Y=0$. \MTEXT{As $f(y_1)=0$, we obtain} $f(\mu_{y_1,x_1}(t))=0$ for $t\in [0,1]$.
Since all points $x\in W_1$ are connected in $W_1$  to the set $V_1$ by piecewise smooth paths,
this shows that $f=0$}.
\hfill \Box \medskip

Finally, we are ready to complete the proof of the main theorem for
   active measurements.\medskip

   \noindent
 {\bf Proof.} (of Theorem \ref{main thm3 new})
  Theorem  \ref{main thm3 new} follows from 
  Theorem \ref{prop. A}  and Theorem \ref{main thm paper part 2}.
   \hfill \Box \medskip

\generalizations{
 \section{Generalizations and outlook}\label{sec: Generalization}

\HOX{T23. This section needs to be checked and made less messy. - Matti }
{
In this section we present a sketch of the analysis how the above results can
be applied for the single measurement inverse problem, that is, how we can  obtain an approximation
of the space-time with a single measurement. Recall that  
 Theorem \ref{main thm Einstein} concerns the case when all possible
measurements near the freely falling observer $\mu$ are known. In principle this implies that the perfect
spacetime cloaking with a smooth globally hyperbolic metric is not possible, that is, no two
space-times having different conformal structure appear the same in
all measurements. However, one can ask if 
one can make an approximative image of the space-time knowing
 only one measurement.  In general, in many
inverse problems several measurements can be packed together to 
one measurement.
For instance, for the wave equation with a  time-independent simple metric
this is done in \cite{HLO}. Similarly, Theorem \ref{main thm Einstein} and its proof
make it possible to do approximate reconstructions as discuss below.
\HOX{Add a remark that we can produce any number $N\in \Z_+$ of point
sources in the space time.}

\subsection{Pre-compact collection of Lorentzian manifolds}
For $m\in \Z_+$ and $\Lambda_{m}\in \R_+$, let $\M_{m,\Lambda}$
be the set of pointed  globally hyperbolic Lorentzian manifolds $(M,g,x_0)$ such that
the corresponding Riemannian manifold $(M,g^+)$
has the property that in the closure of the ball $B_{g^+}(x_0,m)$ with center $x_0$
and radius $r=m$
the covariant derivatives of the curvature tensor $R=R(g^+)$ of the Riemannian
metric $g^+$ satisfy $\|\nabla^j R\|_{g^+}\leq \Lambda$
for all $j\leq m$, the  injectivity radius $i_0(x)$ of $(M,g^+)$ satisfies
$i_0(x)>\Lambda^{-1}$ for all $x\in \hbox{cl}(B_{g^+}(x_0,m))$. 
Let us choose some numbers $\Lambda_m$ for all $m\in \Z_+$,
and let 
\ba
\M=\bigcap_{m\in \Z_+}\M_{m,\Lambda_m}
\ea
be the set endowed with the initial topology for which all
identity maps $\iota_m:\M\to \M_{m,\Lambda_m}$ are continuous.
Using the Cheeger-Gromov theory and the Cantor diagonalization
procedure, we see that the set $\B_R$ of the balls  $B_{g^+}(x_0,R)$ of radius $R$ of the pointed manifolds  
 $(M,g,x_0)\in \M$ form a precompact set in the  topology given by the Gromov-Hausdorff metric
 and the closure of $\B_R$, denoted $\hbox{cl}_{GH}(\B_R)$, consists of smooth manifolds,
 see \cite{And1,Cheeger}, see  also the applications for inverse problems in  \cite{AKKLT}.
 Indeed, we see that the balls of radius $R$ of the pointed manifolds in $\M_{m,\Lambda_m}$ are
 precompact for all $m$ and thus  using a Cantor diagonalization argument
 we see claimed is precompactness.
  
Let us fix $R>0$ and define  $\mathcal N$ to be the set 
$((M,g),(x_0,\xi_0))$  such that  $(M,g,x_0)$ is a pointed  Lorentzian manifold for which  $B_{g^+}(x_0,R)\in \hbox{cl}_{GH}(\B_R)$,
and a time-ike vector $(x_0,\xi_0)\in TM$ is such that the  
 freely falling observer $\hat \mu$ in $(M,g,x_0)$, such that $\hat \mu(-1)=x_0$ and 
 $\p_s\hat \mu(-1)=\xi_0$ satisfies
$U\cup J_{g}(\hat \mu(-1-\e_0),\hat \mu(1+\e_0))\subset B_{g^+}(x_0,R)$, with fixed parameters 
$\e_0>0$ and  $\hat h>0$ (that determines $U$), 
and finally that $\hat \mu([-1,1])$ is such that no light-like geodesic has cut points
in $J_{g}(\hat \mu(-1),\hat \mu(1))$.

\subsubsection{Geometric preparations for single source measurement}
Let $S^{g^+}(y,r)$ denote the sphere of $T_yM$ of radius $r$ with respect to the metric
$g^+$ and  $S^{g^+}(y)= S^{g^+}(y,1)$.
For  $(y,\xi)\in L^+M$ we define spherical surfaces
\ba
& &\L((y,\xi);\rho)=\{\exp_y(t\theta);\ \theta\in L_y^+M,\ 
 \|\theta-\frac 1{\|\xi\|_{g^+}}\xi\|_{g^+}<\rho,\ t>0\},\\
 & &\L^c((y,\xi);\rho)=\{\exp_y(t\theta);\ \theta\in L_y^+M,\  
 \|\theta-\frac 1{\|\xi\|_{g^+}}\xi\|_{g^+}\geq \rho,\ t>0\}.
  \ea
  We say that $y$ is the center point of these surfaces.

In particular,
$((M,g),(x_0,\xi_0))\in \mathcal N$ implies that
no light-like geodesic starting from $\mu_{a}$ can intersect $\mu_{a}$ again.

 Using compactness of $\mathcal N$, \HOX{We need to check (\ref{consequence of compactness})
 and existence of $\rho>0$  and add parameter $r_0$ below. Idea below is first
 to choose $\e$, i.e. the accuracy on which we want to find $(M,g)$, then
 suitable $L$, then perturbation so that light-like geodesics do not intersect
 ${\hat x}_a$, then $r_0$ and then $\rho$.}
for any $r_0>0$ we can choose $\rho=\rho(r_0)>0$ such that 
for all  $((M,g),(x_0,\xi_0))\in \mathcal N$  and any ${\hat x}_j\in J_{g}(\hat \mu(s_-),\hat \mu(s_+)) $, $j=1,2,3,4$ and $\xi_j\in L^+_{{\hat x}_j}M$ such that 
\beq\label{far away}
d_{g^+}({\hat x}_j,
\gamma_{{\hat x}_k,\xi_k}(\R_+))\geq r_0,\quad \hbox{ when $j\not   =k$},
\eeq
the  set
\beq\label{consequence of compactness}
\bigcap_{j=1}^4 \L(({\hat x}_j,\xi_j);\rho)\cap J_{g}(\hat \mu(s_-),\hat \mu(s_+))\hbox{ contains at most one point.}\hspace{-2.1cm}
 \eeq
 Indeed, if there are no such $\rho$ we find a sequence of manifolds $(M_k,g_k)$
 and $(\vec {\hat x}_k,\vec \xi_k)\in (TM_k)^4$ such that the sets
 (\ref{consequence of compactness}) contain at least two points with the parameter $\rho_k$
 and $\rho_k\to 0$ as $k\to \infty$. Considering suitable subsequences of $(M_k,g_k,x_k)$
 that converge to $(M,g,x)$ 
 and $(\vec {\hat x}_k,\vec \xi_k)$ that converge to $(\vec y,\vec \xi)$ we 
 see that some  geodesics in $(M,g,x)$ have cut points in $J_{g}(\hat \mu(s_-),\hat \mu(s_+))$
 {\ftext which is not} possible due to the definition of $\mathcal N$. Hence required parameter $\rho>0$ exists.

Next, for $((M,g),(x_0,\xi_0))\in  \hbox{cl}( \mathcal N)$ we denote 
\ba
J((M,g),(x_0,\xi_0))=J_{g}(\hat \mu(s_-),\hat \mu(s_+))
\ea
and 
\ba
\mathcal J=\{(J((M,g),(x_0,\xi_0)),[g],x_0,\xi_0);\ ((M,g),(x_0,\xi_0))\in \hbox{cl}( \mathcal N)\},
\ea
where $[g]$ denotes the conformal class of $g$. Equivalence classes $[g]$
are closed sets and we use Hausdorff distance for to define
the distance for the equivalence classes $[g_1]$ and $[g_2]$ using representatives
of the equivalence classes for which the determinant of the metric tensor in the Fermi
coordinates are equal to the determinant of the metric tensor of the Euclidean
metric in the Fermi coordinates associated to a line. We call such representative of 
the equivalence class $[g]$ a Fermi-normalized metric and denote it by $g^{(n)}$.
The corresponding Riemannian metric is denoted by $g^{(n),+}$

If $p\in J_{g}(\hat \mu(s_-),\hat \mu(s_+))$, then for given $(z,\eta)$ there
is unique point $y\in \mu_{a}(s)$
for which there exists a light like geodesic
 $\gamma_{y,\xi}([0,t])$ that connects $y$ to $p$. Moreover, then
 $s=f^-_{a}(p)$. Let $(z_j,\eta_j)\in \U_{z_0,\eta_0}$, $j=1,2,\dots$,
 be dense set 
 $(z_j,\eta_j)\not= (z_0,\eta_0)$. We consider these points to be defined so that they have fixed 
 representation in the Fermi coordinates. \HOX{Add details on the points in the Fermi coordinates}
 Using compactness of $J_{ g}(\hat \mu(s_-),\hat \mu(s_+))$,
  we see that  there is $m_0=m_0(M,g,x_0,\xi_0)$ such that
 for all $p\in J_{ g}(\hat \mu(s_-),\hat \mu(s_+))$ there are four pairs $(z_{j_k},\eta_{j_k})$
 with $j_k\leq m_0$ such
 that $q\mapsto ( f^-_{z_{j_k},\eta_{j_k}}(q))_{k=1}^4$ forms smooth coordinates
 in a neighborhood of $p$.

Let choose dense set 
 of numbers $s_j^k\in (-1,1)$, $j,k\in \Z_+$ and let
 ${\hat x}_{j,k}=\gamma_{z_j,\eta_j}(s_j^k)$. For each
 ${\hat x}_{j^\prime,k^\prime}$ we choose a dense
 set of directions, $\xi_{j^\prime,k^\prime,n}$, $n=1,2,\dots,$ of $S^{g^+}({\hat x}_{j^\prime,k^\prime})\cap L^+_{{\hat x}_{j^\prime,k^\prime}}M$. 
 After this, let us use the pairs $({\hat x}_{j^\prime,k^\prime},\xi_{j^\prime,k^\prime,n})$ to define
the pairs $({\hat x}_{j^\prime,k^\prime}(t_{j^\prime,k^\prime}),\xi_{j^\prime,k^\prime,n}(t_{j^\prime,k^\prime}))$ with some 
 $t_{j^\prime,k^\prime}>0$, and 
 re-enumerate the obtained pairs 
 $({\hat x}_{j^\prime,k^\prime}(t_{j^\prime,k^\prime}),\xi_{j^\prime,k^\prime,n}(t_{j^\prime,k^\prime}))\in L^+M$,
 $j^\prime,k^\prime,n\in \Z_+$ as 
$({\hat x}_a,\xi_a)$, $a\in \Z_+$. 
 
 Let $\A$ be the set of 4-tuples $\vec a=(a_\ell)_{\ell=1}^4\in \Z_+^4$ that contain
 four pairwise non-equal positive integers.  
 Let
 \ba
 q(\vec a)=\bigcap_{\ell=1}^4 \L(({\hat x}_{a_\ell},\xi_{a_\ell});\rho)\cap J_{ g}(\hat \mu(s_-),\hat \mu(s_+)).
 \ea
 Note that $q(\vec a)$ contains then only one point or is empty
 if (\ref{far away}) is satisfied with some $r_0$ and $\rho<\rho(r_0)$. Also,
 let 
 \ba
 & & \Psi^{(1)}_k(\vec a)=\prod_{\ell,\kappa\in \{1,2,3,4\},\ \ell\not=\kappa } \phi(\frac 1k \dist_{g^{(n),+}}({\hat x}_{a_\kappa},
\gamma_{{\hat x}_{a_\ell},\xi_{a_\ell}}(\R_+))),\\
& & \Psi^{(2)}_k(\vec a)=\prod_{\ell=1}^4 \phi(\frac 1k \dist_{g^{(n),+}}(q(\vec a), 
 \L^c(({\hat x}_{a_\ell},\xi_{a_\ell});\rho)),\\
& & \Psi_k(\vec a)=\Psi^{(1)}_k(\vec a)\,\cdotp \Psi^{(2)}_k(\vec a), 
  \ea
where $\phi\in C(\R\cup\{\infty\})$ is a function with $\phi(t)=0$ for $t\leq 1$ and $\phi(t)=1$
for $t>2$ and $t=\infty$ and $\dist_{g^{(n),+}}$ is the distance
with respect to the metric $g^{(n),+}$ that is conformal to the metric $g^+$
and which determinant in the Fermi coordinates is normalized as above.

Using the existence
 of $m_0(M,g,x_0,\xi_0)$ and local coordinates associated to four pairs $(z_{j_k},\eta_{j_k})$,
we see that $\{q(\vec a);\ \vec a\in \A,\ q(\vec a)\not=\emptyset\}$ is a dense set of $J_{ g}(\hat \mu(s_-),\hat \mu(s_+))$.

 Let us then consider the continuous map $G:\mathcal J\to \R^{\Z_+\times \Z_+}$,
 \ba
G: (J((M,g),(x_0,\xi_0)),[g],x_0,\xi_0)\mapsto (\Psi_k(\vec a)\,f^+_{z_j,\eta_j}(q(\vec a)))_{\vec a\in \A,\ j\in \Z_+,\ k\in \Z_+}.
 \ea
  \HOX{Check that $G$ is continuous}
In the case when the set $q(\vec a)$ is empty, we define 
 $f^+_{z_j,\eta_j}(q(\vec j(n),\vec k(n))$ to be equal to $\sup\{s;\ \mu_{z_j,\eta_j}(s)\in 
 J_{ g}(\hat \mu(s_-),\hat \mu(s_+))\}$. In fact, we can define this value
 to be arbitrary because of the "cut-off" function $\Psi_k(\vec a)$.
 
  Note that the sets $\P_V(q)$ are  closed. We see using the arguments of the main text that the map $G$ is injective in $\mathcal J$. 
Since 
 $\mathcal J$ is compact and $G$ is continuous, we see for $G_N:\mathcal J\to \R^{\Z_+\times \Z_+}$,
 \ba
G_N: (J((M,g),(x_0,\xi_0)),[g],x_0,\xi_0)\mapsto (\Psi_k(\vec a)\,f^+_{z_j,\eta_j}(q(\vec a)))_{|\vec a|\leq N,\ j\leq N,\ k\leq N}
 \ea
that for any $\e>0$ there is  $N$ 
 such that $G_N(g)$ determines the conformal type $(J((M,g),(x_0,\xi_0)),[g])$ up to $\e$-error in the  metric  in $\mathcal J$. 
 Here, for on $ \mathcal J $  we use
 the distance function defined using 
 $((B_1,[g_1]),(B_2,[g_2]))\mapsto
 d_{GH}((B_1,g_1^{(n),+}),(B_2,g_2^{(n),+}))$.
  \HOX{We need to define the topology of $\mathcal J$ carefully}

 \subsection{Single source measurement}
 Let us now consider distorted plane  waves propagating
 on the spherical surfaces  $\L(({\hat x}_{a_\ell},\xi_{a_\ell});\rho)$, $|\ell |\leq L$.
 When $N$ is large, these waves interact and produce 3-wave and 4-wave interactions.
 These waves can be produced by sources ${f}_{a_\ell}$ supported in an arbitrarily
 small neighborhoods  $W_{a_\ell}\subset U$ of the "conic" points ${\hat x}_{a_\ell}$ of the spherical surfaces  $\L(({\hat x}_{a_\ell},\xi_{a_\ell});\rho)$.

  We see that for generic set  (i.e.\ in an open and dense set) of  manifolds in $\mathcal N$ it
holds that   there are no points ${\hat x}_{a_\ell}$  \HOX{This "generic" statement needs to be checked}
 and ${\hat x}_{a_\ell^\prime}$ with ${a_\ell}\not= {a_\ell^\prime}$ that
 can be connected with a light-like geodesic. This is due to
 the fact that this condition holds in an open set and we see that if for any $L\in \Z_+$
 we can order the points ${\hat x}_{a(\ell)}$, $\ell=1,2,\dots,L$ in causal order, such that ${\hat x}_{a(\ell)}\not \in J^+({\hat x}_{a(\ell+1)})$, and then making a small perturbation to the metric near
each point ${\hat x}_{a(\ell)}$, in causal order, we can  choose such perturbations
that the point ${\hat x}_{a(\ell)}$ is not on the surface $\cup_{j<\ell}\L^+({\hat x}_{a(j)})$. 
{Below we will assume that neighborhoods \HOX{This generic statement needs to be checked}
 $W_{a_\ell}$ of ${\hat x}_{a_\ell}$ of the  points ${\hat x}_{a_\ell}$ are so small
 that $W_{a_\ell}\subset U$ does not intersect any $\L^+({\hat x}_{a(j)})$, where
$j\not =\ell$ and $j\leq L$.

In the main text we used sources ${f}_{a_\ell}$ that are elements of  the space
$\I(Y_{a_\ell})$  and produce waves $u_{a_\ell}\in \I (N^*Y_{a_\ell},N^*K_{a_\ell})$,
where $K_{a_\ell}\subset \L^+({\hat x}_{a(\ell)})$ are 3-dimensional subsets  and 
 $Y_{a_\ell}= \L^+({\hat x}_{a(\ell)})\cap {\bf t}^{-1}(T_{a(\ell)})$ are 2-dimensional sub-manifolds.
 As the 
products of  distributions associated to two Lagrangian manifolds are difficult to analyze, we modify the construction by 
replacing $K_a$ by surfaces $K_a\subset \L^+_{g}({\hat x}_a)$ and 
  ${f}_{a_\ell}\in \I(Y_{a_\ell})$ by sources $\tilde {f}_{a_\ell}\in \I(K_{a_\ell})$,
supported in the  neighborhood  $W_{a_\ell}$ of ${\hat x}_{a_\ell}$.
As  then $\tilde {f}_{a_\ell}\in \I(N^*K_{a_\ell})$ and $N^*K_{a_\ell}$ is invariant in the future
bicharacteristic flow of $\square_{g}$, the sources $\tilde {f}_{a_\ell}$
 produce by \cite[Prop. 2.1]{GU1} waves  $ u_{a_\ell}\in \I (N^*K_{a_\ell})$,
where the principals symbols are not anymore evaluations of the principal
symbol of the source at some point but integrals of the principal symbol along a bicharacteristic.

The sources $\tilde {f}_{a_\ell}\in \I (K_{a_\ell})$  can be chosen so that there are 
 neighborhoods $V_{a_\ell}$ of the surfaces  $\L(({\hat x}_{a_\ell},\xi_{a_\ell});\rho)$  that satisfy
$W_{a_\ell}\subset V_{a_\ell}$,
the linearized waves $ u_{a_\ell}={\bf Q}_{ g}\tilde {f}_{a_\ell}$ satisfy  singsupp$( u_{a_\ell})\subset V_{a_\ell}$,
and   
\beq\label{eq: V ehto}
\hbox{$V_{a_\kappa}\cap W _{a_\ell}=\emptyset$  for all $\ell,\kappa\leq L$,
 $\ell\not =\kappa$.}
\eeq
We also assume that
$W_{a_\kappa}$  are such that $J^+_{g}(W_{a_\ell})\cap
J^-_{g}( W_{a_\ell})\subset  U$ and that the wave map coordinates
corresponding to the Euclidean metric are well defined in
this set. Note that the above properties of $W_{a_\kappa}$ can be obtained
by choosing all surfaces $K_{a_\kappa}$, $\kappa\leq L$ with a the same parameter
$s_0$ that is small enough. Below, we assume that $s_0$ is chosen in such a way. 

Consider now the non-linear interactions produced by the source
\beq\label{Fe source}
F_\e(x)=\e\left(\sum_{\ell=1}^N\tilde {f}_{a_\ell}(x)\right),
\eeq
where $\e>0$ is small.  
This source produces the wave
\beq\label{nl wave}
u(x)=u^{(4)}_\e(x)+O(\e^5),\quad u^{(4)}_\e(x)=\sum_{j=1}^4 \e^j\M^{(j)}(x),
\eeq
where we consider $O(\e^5)$ as a term which is so small that it can be considered to be negligible. \HOX{Maybe we could use $\vec \e$ with $\e_j=10^{-100-10j}$ etc and justuse the considerations of main text}
Let us next assume that we can measure  the singularities of $\M^{(j)}$ for all $j\leq 4$.
Let $S_j$ be the intersection of $U$ and singular support of  $\M^{(j)}$, $j\leq 4$.

\subsubsection{Analysis of 1st, 2nd, and 3rd order waves}
As in the source $F_\e$ we use only one parameter $\e$ instead of four parameters
$\vec\e=(\e_1,\e_2,\e_3,\e_4)$ we need to consider self-interaction terms,
that is, in our computations the permutations $\sigma:\{1,2,3,4\}\to \{1,2,3,4\}$
need to be replaced by general (non-injective) maps $\sigma:\{1,2,3,4\}\to \{1,2,3,4\}$.
This does not cause difficulties in analyzing the terms $\M^{(j)}$ for all $j\leq 3$
as because of condition (\ref{eq: V ehto}) the source terms do not cause
self-interactions with terms corresponding to different Lagrangian manifolds
$\Lambda_{a_\ell}=N^*K_{a_\ell}$ and the interactions of linearized waves
corresponding to two such manifolds produce  waves that are in
$\I(N^*K_{a_\ell},N^*K_{a_\ell,a_\kappa})$.

For generic sources singsupp$(\M^{(1)}) $ is equal to $S_1=\cup_{\ell\leq L} K_{a_\ell}\cap U$ and 
singsupp$(\M^{(2)})$  is equal to $S_2=\cup_{\ell,\kappa\leq L} K_{a_\ell,a_\kappa}\cap U
\subset S_1$. 

Let us then consider singularities of $\M^{(3)}$ using similar methods that we used to analyze the fourth
order terms. Again, observe that we need to consider also the terms for 
which the map $\sigma:\{1,2,3\}\to \{1,2,3\}$ is not bijective. For such terms,
e.g., if $\sigma(1)=r_1$, $\sigma(2)=r_2$ and $\sigma(3)=r_3$ with $r_1=r_3$,   
we see using \cite[Lemmas 1.2 and 1.3]{GU1} 
 \cite[Thm.\ 1.3.6]{Duistermaat}  
 that e.g. $v\in u_{r_1}\cdotp {\bf Q}(u_{r_1}u_{r_2})$ has
 wave front set that is subset of $N^*K_{r_1r_2}$, implying
 that wave front set of ${\bf Q}v$ is subset of  $N^*K_{r_1r_2}$.
 Using this we see that outside $S_2$ only  the terms for which map $\sigma:\{1,2,3\}\to \{1,2,3\}$
 is bijective cause singularities. Let us next show that this kind of singularities can be observed.

To this end, let us consider the term
 $\bra F_\tau,{\bf Q}(A[u_3,{\bf Q}(A[u_2,u_1])])\cet$ where
 all operators $A[v,w]$ are of the type $A_2[v,w]=g^{np}g^{mq}v_{nm}\p_p\p_q w_{jk}$,
 cf.\ (\ref{eq: tilde M1}) and (\ref{eq: tilde M2}), we obtain  formulas similar to
  (\ref{eq: def of D 2}), where the first factor in 
   (\ref{def of D}) is removed. Next we show that using this and the similar analysis we did
   for the fourth order terms with the WKB analysis, we can see that the 
   principal symbol of  $\M^{(3)}$ in the normal coordinates does not vanish
   for generic interacting distorted plane  waves. Indeed, the indicator function
   $\Theta^{(3)}_\tau$ can be obtained as a sum of terms similar to $T^{(3),\beta}_\tau$ given in
    (\ref{eq: 3rd order term}). For these terms we obtain in the Minkowski space formulas
    analogous to (\ref{first asymptotical computation}) and (\ref{t 4 formula}). Let us consider  the case when
     ${\rtext b_j}$, $j\leq 5$ are  light-like co-vectors such that
    $b^{(5)}\in \hbox{span}(b^{(1)},b^{(2)},b^{(3)})$ and vector $b^{(4)}$ is linearly independent of $b^{(1)},b^{(2)},b^{(3)}$  
 so that
    $y=A^{-1}x$ are coordinates such  that
   ${\bf p}_4=(A^{-1})^tb^{(4)}=0$. The particular example of such co-vectors that we consider, are
   \ba
& &b^{(5)}=(1,1,0,0),\quad b^{(1)}=(1,1-\frac 12\rhoepsilon_1^2+O(\rhoepsilon_1^3),\rhoepsilon_1,\rhoepsilon_1^3),\\
& &b^{(3)}=(1,-1,0,0),\quad b^{(2)}=(1,1-\frac 12\rhoepsilon_2^2+O(\rhoepsilon_2^3),\rhoepsilon_2,\rhoepsilon_2^{201})
\ea
 and $ \rhoepsilon_1=\rhoepsilon^{100}$, $\rhoepsilon_2=\rho$, and $b^{(4)}$ is some fixed vector.
 Then $b^{(5)}=\sum_{j=1}^3 a_j {\rtext b_j}$,
 where $a_1=O(1)$, $a_2\sim -a_1\rho^{99}$, and $a_3=O( \rho^{2})$ as $\rho\to 0$.   For
 $T^{(3),\beta}_\tau$ given in
    (\ref{eq: 3rd order term}), we have then 
   \HOX{WKB needs to be checked.}
   \ba
T^{(3),\beta}_\tau&=&
\bra  u_\tau, h\,\cdotp \B_3u_3\,\cdotp {\bf Q}_0(\B_2u_2\,\cdotp \B_1 u_1)\cet\\
& &\hspace{-1cm}=C\frac {\P_\beta^\prime} {\omega_{12}}
\bra u_\tau, h\,\cdotp u_3\,\cdotp v^{a-k_2+1,a-k_1+1} (\,\cdotp ;b^{(2)},b^{(1)})\cet
\\
\nonumber&&\hspace{-1cm} =C\frac {\P^\prime_\beta}{\omega_{12}} 
\int_{\R^4}
e^{i\tau(b^{(5)}\cdotp x)} 
h(x)(b^{(3)}\,\cdotp x)^{a-k_3}_+\cdotp\\
\nonumber& &\quad\quad\quad\cdotp
(b^{(2)}\,\cdotp x)^{a-k_2+1}_+
(b^{(1)}\,\cdotp x)^{a-k_1+1}_+\,dx,\\
& &\hspace{-1cm}
=\frac {C\P^\prime_\beta\det (A^\prime)}{ \omega_{12}} \int_{(\R_+)^3\times \R}
e^{i\tau {\bf p}\cdotp y} 
h(Ay){\hat x}_3^{a-k_3}{\hat x}_2^{a-k_2+1}
{\hat x}_1^{a-k_1+1}\,dy
\\
\nonumber&&\hspace{-1cm} =
 {C\det (A^\prime)\,{\P^\prime_\beta}}
\frac {(i\tau)^{-(10+3a-|\vec k_\b|)}(1+O(\tau^{-1}))}{{\bf p}_3^{(a-k_3+1)}{\bf p}_2^{(a-k_2+2)}{\bf p}_1^{(a-k_1+2)}\omega_{12}}\,,
\ea
   where ${\bf p}_j=g(b^{(5)},{\rtext b_j})=-\frac 12 \rhoepsilon_j^2+O(\rhoepsilon_j^3)$, $j\leq 3$ and $\omega_{12}=g(b^{(1)},b^{(2)})=
   -\frac 12 \rhoepsilon_2^2+O(\rhoepsilon_2^3)$. 
Moreover, with the appropriate choice of polarizations   $v_{(\ell)}^{nm}=g^{nj}g^{mk}v^{(\ell)}_{jk}$ given in  (\ref{chosen polarization}) for $j=2,3,4$, c.f. (\ref{vbb formula}), we have 
\beq\label{pre-eq: def of D prime}
\P_{\beta_1}^\prime =(v_{(3)}^{pq}b^{(1)}_pb^{(1)}_q)(v_{(2)}^{nm}b^{(1)}_nb^{(1)}_m)\D,
\eeq
 and
\ba
 \D =g_{nj}g_{mk}v_{(5)}^{nm}v_{(1)}^{jk},
 \ea 
 where we note that $v_{(5)}^{nm}$ does not need to satisfy any harmonicity type conditions.
We see that the term $T^{(3),\beta_0}_\tau$, where $k_1=k_1^{\beta_0}=4$, dominates the other terms 
$T^{(3),\beta}_\tau$,when 
$\rho\to 0$. Thus, using the real analyticity of the leading order
coefficient in the asymptotic of  $\Theta^{(3)}_\tau$, we see similarly to the main text that the principal symbol of $\M^{(3)}$ in the normal
coordinates does not vanish in the generic case on the manifolds
$\Lambda_{a_1,a_2,a_3}^{(3)}$ that are the flow out of the light-like
directions of $N^*K_{a_1,a_2,a_3}$.
     
   \HOX{We need to check here arguments!!}

   The above analysis implies that for generic manifold
   and for generic sources singsupp$(\M^{(3)}) $ is equal to 
   \ba
   S_3=U\cap (\bigcup_{(\vec x_j,\vec \xi_j)=({\hat x}_{a_\ell},\zeta_{a_\ell}),\ \ell\leq L}
    \mathcal Y((\vec x,\vec \xi),t_0,s_0)).
   \ea
   In the generic
   case we can thus assume the set $S_3$ has the above form and that the given set
  singsupp$(\M^{(3)}) $ determines it.
  
\subsubsection{Analysis of the 4th order terms}

Let us  next consider singularities of $\M^{(4)}$ outside the set $S_1\cup S_2\cup S_3$.
Outside
these surfaces we look singularities of  $\M^{(4)}$ using the indicator functions  $\Theta^{(4)}_\tau$,  constructed using gaussian
beams.
Assume that $x_5\in U\setminus (S_1\cup S_2\cup S_3)$.
As pointed out before we need to consider general (non-injective) maps $\sigma:\{1,2,3,4\}\to \{1,2,3,4\}$.
The case when $\sigma$ is a permutation is considered in the main text.
Let us thus assume that $\sigma$ is not a permutation.

We need to consider cases where the linearized waves $u_1,u_2,u_3$ are replaced by
$u_{r_1},u_{r_2},u_{r_3}$, where $r_1,r_2,r_3\in \{1,2,3,4\}$ are arbitrary,
and $u_4$, kept the same (changing it would be just
renumeration of indexes.) Also, we can assume that $\{r_1,r_2,r_3\}$
is not the set  $\{1,2,3\}$, as this case we have already analyzed.

We start by  considering the $\tilde T$-functions.
We need to consider terms similar to $\tilde T$ in the main text containing e.g.\
$\bra u_\tau,{\bf Q}(u_{r_1}u_{r_2})\cdotp {\bf Q}(u_{r_3}u_{r_4})\cet$,
when some $r_j$ are the same. Assume next that
e.g. $r_1$ and $r_3$ are the same.
We use that fact that  $v\in \I(\Lambda_1, \Lambda_2)$,
we have
WF$(v)\subset \Lambda_1\cup \Lambda_2$. Moreover,
by \cite[Thm.\ 1.3.6]{Duistermaat},
\ba
\hbox{WF}(v\,\cdotp w)\subset 
\hbox{WF}(v)\cup \hbox{WF}(w)
\cup\{(x,\xi+\eta);\ (x,\xi)\in \hbox{WF}(v),\ (x,\eta)\in \hbox{WF}(w)\}
\ea
and thus we see that
the wave front set of term ${\bf Q}(u_{r_1}u_{r_2})\cdotp {\bf Q}(u_{r_1}u_{r_4})$
is in the union of the sets $N^*K_{r_j}$,
$N^*(K_{r_1}\cap K_{r_2})$, and 
$N^*(K_{r_1}\cap K_{r_2}\cap K_{r_4})$. The elements in
the first two sets do not propagate in the bicharacteristic flow
and the last one propagates along $\mathcal Y((\vec x,\vec \xi),t_0,s_0)$.
Observe that in the generic case $\mathcal Y((\vec x,\vec \xi),t_0,s_0)\subset S_3$.
The other  $\tilde T$-functions can be analyzed similarly.
Thus in the generic case the  $\tilde T$-functions do not cause observable 
singularities in the complement of $S_3$.

Consider next the $T$-functions.
We need to consider terms similar to $T$ in the main text containing e.g.\
$\bra u_\tau,u_{r_1}\,\cdotp {\bf Q}(u_{r_2}\,\cdotp {\bf Q}(u_{r_3}\,\cdotp u_{r_4}))\cet$,
when some $r_j$ are the same. 
 Again we decompose ${\bf Q}={\bf Q}_1+{\bf Q}_2$ and consider cases $p=1,2$ as in the main text. In the case $p=2$ we need to consider
critical points of modified phase functions 
\ba
\nonumber \Psi^{\vec r}_2(z,y,\theta)
&=&\theta_1z^{r_1}+\theta_2z^{r_2}+\theta_3z^{r_3}+\theta_4y^4+\varphi(y).
\ea
Under the assumption that $x_5\not \in K_4$, we see similarly
as in the main text that there are no critical points. Now
\ba
d_{z,y,\theta}\Psi^{\vec r}_2&=&(\omega^{\vec r};
r+d_y\psi_4(y,\theta_4),
z^{r_1},z^{r_2},z^{r_3},d_{ \theta_4}\psi_4(y,\theta_4))
\ea
where $\omega^{\vec r}=(\omega^{\vec r}_j)_{j=1}^4$,
$\omega^{\vec r}_j=\sum_{n=1}^3\theta_n\delta_{n\in R(j)}$, $R(j)=\{n\in \{1,2,3\}; r_n=j\}.$
Note that $(z,\omega^{\vec r})\in \bigcap_{n\in \{r_1,r_2,r_3\}}N^*K_{n}$ and some of $r_j$
may be the same. 
We see as in the main text that as $x_5\not \in {\mathcal Y},$
alternative (A1) can not hold.
Note that if (A2) holds then $\omega^{\vec r}\in \X$, and we see as in 
the main text that this is not possible.  Thus the term with $p=2$ causes no observable singularities
near $x_5$.

In the case $p=1$, let $(z,\theta,y,\xi)$ be a critical point of
 the phase function
  \ba
 \Psi^{\vec r}_3(z,\theta,y,\xi)=\theta_1z^{r_1}+\theta_2z^{r_2}+\theta_3z^{r_3}+(y-z)\,\cdotp \xi+\theta_4y^4+\varphi(y).
 \ea
Then
 \beq\label{eq: critical points 1}
  & &\p_{\theta_j}\Psi_3=0,\ j=1,2,3\quad\hbox{yield}\quad z\in K_{{r_1}}\cap K_{r_2}\cap K_{r_3},\\
  \nonumber
& &\p_{\theta_4}\Psi_3=0\quad\hbox{yields}\quad y\in K_4,\\
  \nonumber
 & &\p_z\Psi_3=0\quad\hbox{yields}\quad\xi=\omega^{\vec r},\\
   \nonumber
& &\p_{\xi}\Psi_3=0\quad\hbox{yields}\quad y=z,\\
  \nonumber
 & &\p_{y}\Psi_3=0\quad\hbox{yields}\quad \xi=
 -\p_y\varphi(y)-\theta_4 w.
  \eeq
  The critical points we need to consider for the asymptotics satisfy also
\beq\label{eq: imag.condition copu} \\ \nonumber
\im \varphi(y)=0,\quad\hbox{so that }y\in \gamma_{x_5,\xi_5},\ \im d\varphi(y)=0,\
\re d\varphi(y)\in L^{*,+}_y\hattuM _0.\hspace{-1cm}
\eeq
Now, as we assume that $\{r_1,r_2,r_3\}$
is not the set  $\{1,2,3\}$, some of $r_j$ is equal to $4$ or two of these are the same.   For simplicity, assume that
$r_1$ is either 4 or alternatively, $r_1$ coincides with $r_2$ or $r_3$. As  $w\in N^*K_4$, these imply that 
$\p_y\varphi(y)\in N^*K_{r_2}\cap  N^*K_{r_3}\cap N^*K_4\subset \mathcal X((\vec x,\vec\xi);t_0)$.
Hence the assumption that $x_5\not \in  {\mathcal Y}$ implies that
terms with $p=1$ do not cause any singularities near $x_5$. Thus only the singularities of
$T$-type terms for which $\sigma$ is bijection,
that were analyze in the main text, are visible in $U\setminus (S_1\cup S_2\cup S_3)$.
The singularities could be analyzed  in much more careful way
if we could analyze products of $v\in \I(K_1,K_{13})$ and $w\in \I(K_2,K_{23})$
and show that such product is a sum of conormal distributions.\hiddenfootnote{HERE ARE SOME EXTRA COMPUTATIONS:
By applying method of stationary phase in $\xi$ and $z$ variable we obtain
a integral similar to (\ref{eq; T tau asympt modified}),
 \beq \label{eq; T tau asympt modified2}
& &T^{(1)}_{\tau,1,1}^{(4),\beta}
=c\tau^{4}
\int_{\R^{8}}e^{i(\theta_1y^{r_1}+\theta_2y^{r_1}+\theta_3y^{r_1}+\theta_4y^{r_1})+i\tau \varphi(y)}
c_1(y, \theta_1, \theta_2)\cdotp\\ \nonumber
& &
\cdotp \sigma_{u_3}(y, \theta_3)
q_{1,1}(y,y,\omega^{\vec r}_\beta(\theta))\sigma_{u_4}(y,\theta_4) {{\Phi}}(y,\tau )\,d\theta_1d\theta_2d\theta_3d\theta_4dy.
\eeq
where $\omega^{\vec r}_\beta$ is similar to $\omega^{\vec r}$  is defined with 
$R_\beta(j)=\{n\in \{1,2,3\}; r_n=\sigma_\beta(j)\}$ is  similar to $R(j)$.

 Similarly, to the analysis related to the function $\chi(\theta)$ that we vanish
 in a neighborhood of the set $\mathcal A_q$ given in (\ref{set Yq}),
 we can define a fucntion 
 $\tilde \chi(\theta)\in C^\infty(\R^4)$ that 
vanishes in a $\e_3$-neighborhood $\W$ (in the $g^+$ metric) of
$\tilde {\mathcal A_q}$,
\ba
\tilde {\mathcal A_q}&:=&
\bigcup_{1\leq j<k\leq 4}N_q^*K_{jk}.
\ea
The idea of the sets $N_q^*K_{jk}$ is that in the complement of  their  neighborhoods at least three of the coordinates $(\theta_1,\theta_2,\theta_3,\theta_4)$
 are comparable to $|\vec\theta|$, that is, there is $j_0$ such that $|\theta_k|\geq c_j|\vec \theta|$ for $k\not =j_0$.
 Indeed, e.g.\ $(q;(\theta_1,\theta_2,\theta_3,\theta_4))\not \in N_q^*K_{jk}$ implies 
 that $(\theta_j,\theta_k)\not =0$. As this holds for all $(j,k)$, we see
 that $\theta\not \in \W$ and $|\theta|_{\R^4}=1$, only one coordinate of  $\theta$ can vanish. This makes the symbol
$\tilde \chi(\theta)b(y,\theta)$ a sum of a product type symbols $S^{p,l}({\rbtext((-\infty,{\reftext T_0})\times N)};\R^3\times \R\setminus \{0\})$.

Let  \ba\hbox{$\tilde b_0(y,\theta)=\phi(y)\tilde \chi(\theta)b(y,\theta)$}\ea be a  product type symbol determining 
a  distribution 
$
\tilde \F^{(4),0}$ 
 that is given by the formula (\ref{f4-formula}) with
$b(y,\theta)$ being replaced by 
 $b_0(y,\theta)$.
We  can consider
 $\tilde \F^{(4),0}$ as a sum of in three terms: 
 \ba
 \tilde \F^{(4),0}=\tilde \F^{(4),0,1}+\tilde \F^{(4),0,2}+\tilde \F^{(4),0,2}\ea
  where \HOX{Check 3/2 HERE!!}
 \ba
 \tilde \F^{(4),0,1}&\in& \sum_{ \{k_1,k_2,k_3,k_4\}=\{1,2,3,4\}} \I^{p(\vec k)-4+1/2}(K_{k_2,k_3,k_4})\cap  \I(N^*K_{k_2,k_3,k_4}\setminus N^*\{q\})\\
& &\subset 
\sum_{ \{k_1,k_2,k_3,k_4\}=\{1,2,3,4\}} \I^{p(\vec k)-4+1/2-3/2+1}(N^* K_{k_2,k_3,k_4})\
,
\ea
where $p(\vec k)=p_{k_2}+p_{k_3}+p_{k_4}$ and
\ba
\tilde \F^{(4),0,2}&\in& \I^{p-4}(\{q\}\cap  \I(N^*\{q\}\setminus (\cup_{\vec k} N^*K_{k_2,k_3,k_4}))\\
& &\subset \I^{p-4-2+1}(N^*\{q\})
\ea
 and  $\tilde \F^{(4),0,2}$ is microlocally supported in some (arbitrarily) small
 neighborhood $\W$ of  $\cup (N^*K_{k_2,k_3,k_4}\cap N^*\{q\})$. 
  Let $\Lambda^{(3)}$ be the flowout of $\cup_{\vec k} N^*K_{k_2,k_3,k_4}$.
We see that $\M^{(4)}$ is Lagrangian distribution
 with known order in $I^-(x_6)\setminus  {\mathcal Y}((\vec x,\vec \xi),t_0)$, that is outside
 the set $\cap_{s_0>0}{\mathcal Y}((\vec x,\vec \xi),t_0,s_0)$ that has the Hausdorff dimension
 two. Indeed,\HOX{Check 3/2}
\ba
& &{\bf Q}_{g}\tilde \F^{(4),0,1}\in
\sum_{ \{k_1,k_2,k_3,k_4\}=\{1,2,3,4\}} \I^{p(\vec k)-4-3/2}
\I(U\setminus K_{k_2,k_3,k_4};
\Lambda^{(3)}),\\
& &{\bf Q}_{g}\tilde \F^{(4),0,2}\in 
\I^{p-5-3/2}(U\setminus\{q\};\Lambda_q),\\
  & &\hbox{WF}({\bf Q}_{g}\tilde \F^{(4),0,3})\subset  \V .\ea
where $\V\subset T^*M$ is a neighborhood of $ 
T^*\mathcal Y((\vec x,\vec \xi),t_0)$ that is
the flow out of the set $\W$ in the bicharacteristic flow.
This makes it possible to analyze $
 \bra u_\tau,\F_1^{(4)}\cet$ and see that ${\bf Q}\F_1^{(4)}$ is Lagrangian distribution
 in $U\setminus {\mathcal Y}((\vec x,\vec\xi),t_0)$.

To study the singularities of $u^{(4)}_\e(x)$, we consider next interaction of four distorted plane  waves, 
corresponding the $\vec a=(a_\ell)_{\ell=1}^4 \in \A$, $a_\ell\leq N$ that are singular on the surfaces
$K_{a_\ell}=\L(({\hat x}_{a_\ell},\xi_{a_\ell});\rho)$.
  Assume now that $\bigcap_{\ell=1}^4 \K_{a_\ell}$ contains a point $q$ and
  let $K_{a_1,a_2}=\bigcap_{\ell=1}^2 \K_{a_\ell}$ and $K_{a_1,a_2,a_3}=\bigcap_{\ell=1}^3 \K_{a_\ell}$ . 
 Then, we see that   $\Lambda_q^+$ and
  the flowout $\Lambda_{a_1,a_2,a_3}$ of $N^*K_{a_1,a_2,a_3}$ in the canonical
  relation of $\square_g$ are in $T^*U$ four dimensional Lagrangian
  manifolds that intersect only on a three dimensional set.
  Assume now that the source ${f}_{a_\ell}$ have such orders 
  that $u_{a_\ell}={\bf Q}_g{f}_{a_\ell}\in \I^{p_{a_\ell}}(K_{a_\ell})$ outside support of ${f}_{a_\ell}$.

The linearized term $\M^{(1)}$ in (\ref{nl wave}) has the wave front set that is
a subset of  $S_1=\bigcup_{a\leq N}N^*K_a$,
the term  $\M^{(2)}$  has wave front set is a subset of  $S_1\cup S_2$ where  $S_2=\bigcup_{a_1<a_2\leq N}N^*K_{a_1,a_2}$. Let us note that when analyze the term $\M^4$, we have
to consider the term $\M^{(1)}\,\cdotp\M^{(1)}$ containing terms
$u_{a_1}\,\cdotp u_{a_1}$, where $u_{a_1}\in \I^{p}(K_{a_1})$. 
Denote $b(x^1)=u_{a_1}(x^1,x^\prime)$ and assume that $\hat b(\xi_1)\sim c_0 |\xi_1|^p$ as $\xi_1\to \infty$,
we see that 

\ba
(\hat b*\hat b)(\xi_1)=\int_\R \hat b(\xi_1-t)\hat b(t)dt
\sim 
|\xi_1|^p
\int_\R c_0(1-|\xi_1|^{-1}t)^p\hat b(t)dt
\sim c_0 b(0)|\xi_1|^p,
\ea

****GENERALIZATION***

Let  $B_1\in \I^{p_1}(\{0\})$ and $B_2\in \I^{p_2}(\{0\})$,
Let us write symbol $b_1(x,\xi)$ of $B_1$ as a sum of positive homogeneous term and
a compactly $\xi$-supported term $b_1(x,\xi)=b_1^h(x,\xi)+b_1^c(x,\xi)$ and similar
lower order terms.
Then $B_1^c,B_2^c\in C^\infty$, and
 $B_1^h\,\cdotp B_2^c\in \I^{p_1}(\{0\})$,  $B_1^c\,\cdotp B_2^h\in \I^{p_2}(\{0\})$,
  $B_1^c\,\cdotp B_2^c\in C^\infty$, and symbol of  $B_1^h\,\cdotp B_2^h$ is 
  obtained via substitution $t=|\xi_1|s$,
\ba
(b_1*b_2)(x,\xi_1)&=&\int_\R b_1(x,\xi_1-t) b_2(x,t)dt\\
&=&|\xi_1|^{p_1+p_2+1}
\int_\R b_1(x,1-s) b_2(x,s)ds
\ea
and thus  $B_1^h\,\cdotp B_2^h\in \I^{p_1+p_2+1}(\{0\})$.
Hence,  $B_1\,\cdotp B_2\in \I^{r}(\{0\})$, $r=\max(p_1,p_2)$.

****GENER. 2***

Let  $B_1\in \I^{p_1}(\{0\})$ and $B_2\in \I^{p_2}(\{0\})$,
Let us estimate the symbol $b_1(x,\xi)$ of $B_1$
by
\ba
|b_1(x,\xi)|\leq b_1^h(x,\xi)+b_1^c(x,\xi)
\ea 
where $b_1^c(x,\xi)$ is a compactly $\xi$-supported  symbol and
$b_1^h(x,\xi)$ is $p_1$-homogeneous. Let $B_1^c$ be the conormal distribution
with symbol $b_1^c(x,\xi)$ and  $B_1^h=B_1-B_1^c$, and introduce similar notations for $B_2$.
Then $B_1^c,B_2^c\in C^\infty$, and
 $B_1^h\,\cdotp B_2^c\in \I^{p_1}(\{0\})$,  $B_1^c\,\cdotp B_2^h\in \I^{p_2}(\{0\})$,
  $B_1^c\,\cdotp B_2^c\in C^\infty$, and symbol $c^h$ of  $B_1^h\,\cdotp B_2^h$ can be estimated
using substitution $t=|\xi_1|s$,
\ba
|c^h(x,\xi_1)|&\leq &\int_\R b_1(x,\xi_1-t) b_2(x,t)dt\\
&\leq &|\xi_1|^{p_1+p_2+1}
\int_\R b_1(x,1-s) b_2(x,s)ds
\ea
and thus  $B_1^h\,\cdotp B_2^h\in \I^{p_1+p_2+1}(\{0\})$.
Hence,  $B_1\,\cdotp B_2\in \I^{r}(\{0\})$, $r=\max(p_1,p_2)$.

****GENER. 3***

We need also to analyze the products of 
  $B_1\in \I^{p_1,l_1}(K_1,K_{13})$ and $B_2\in \I^{p_2,l_2}(K_2,K_{13})$.
  When we analyze this in the case when $|\theta_2|>c|\theta|$ and $|\theta_3|>c|\theta|$ 
  and  $|\theta_1|>c|\theta|$, this is equivalent to analyzing product of 
  $ \I(K_{13})$ and  $ \I(K_{12})$.
  The case where only  $|\theta_2|>c|\theta|$ and $|\theta_3|>c|\theta|$ but $\theta_1$ may be
  small  is difficult. Note that it is possible that we do not need to analyze
  this to understand light-like directions on $N^*K_{123}$. Maybe we should
  first find $S_3$ from $\e^3$ terms.

***

Plan: 

1. It seems likely that in generic case we can see $S_1$, $S_2$ and $S_3$ from different $\e$-order of singularities
The difficulty why we need to consider different $\e$ order is that we need
to consider  the product of distributions 
  $v_1\in \I(K_{13})$ and  $ v_1\in\I(K_{12})$ and we can not find the order of ${\bf Q}(v_1v_2)$.

2. Outside the surfaces  $S_1$, $S_2$ and $S_3$ we can detect $S_4$ and can separated different components.
For this we use the product theorem for the wave front sets of products
of   $v_1\in \I(K_{13})$ and  $ v_1\in\I(K_{12})$ to see that the wave front set of ${\bf Q}(v_1v_2)$ is contained
in $\Lambda^{(3)}$. 

***

and thus $u_{a_1}\,\cdotp u_{a_1}\in \I^{p}(K_{a_1})$ has a non-vanishing principal symbol
at $(x,\xi)\in N^*K_{a_1}$ if and only $u_{a_1}$ has a non-vanishing principal symbol at $(x,\xi)$ and 
a non-vanishing value at $x$.

*** ABOVE: We had problem: We can not analyze terms $\M^{(2)}\cdotp \M^{(2)}$
appearing in $\M_4$. The problem can be solved by considering solutions
outside Hausdroff 2-dimensional set $\Z$.
}
\footnote{We could use Piriou,	A.: Calcul	symbolique non	 lineare pour une onde conormale simple. Ann. Inst. Four. 38, 173-188 (1988), or \cite{GU1}, Lemma 1.3. Look also
Bougrini, H.; Piriou, A.; Varenne, J. P. Propagation et interaction des symboles principaux pour les ondes conormales semi-lin'eaires. (French) [Propagation and interaction of principal symbols for semilinear conormal waves] Comm. Partial Differential Equations 23 (1998), no. 1-2, 333--370.}

\subsection{Separation of singularities from different point sources}
Summarizing, above we have seen that in the generic case  
the terms  $\M^{(j)}$, $j\leq 3$   have wave front sets which union is $S_1\cup S_2\cup S_3$ 
and,  $\M^{(4)}$  has wave front set is a subset of   $S_1\cup S_2\cup S_3\cup S_4$ where  $S_4=\bigcup_{a_1<a_2<a_3<a_4\leq N}\Lambda^+_{q(a_1,a_2,a_3,a_4)}$.
Now we notice that for a generic  (i.e.\ in an open and dense set) $((M,g),(x_0,\xi_0))\in \mathcal N$,
we see that when $\vec a$  and $\vec b$ run over any finite subset of  $\A$,
\HOX{Check the the claims on "generic manifolds"}
the sets $\Lambda^+_{q(b_1,b_2,b_3,b_4)}$ intersect
the sets $N^*K_a$, $N^*K_{a_1,a_2}$, $\Lambda_{a_1,a_2,a_3}^{(3)}$, and 
the other sets 
$\Lambda^+_{q(a_1,a_2,a_3,a_4)}$, 
on at most 3-dimensional sets. Indeed,
 by making a small perturbation to the metric in a small neighborhood of
  $q(a_1,a_2,a_3,a_4)$ the sets  $\Lambda^+_{q(a_1,a_2,a_3,a_4)}$ can be perturbed so
  that the above non-intersection condition becomes valid, and the non-intersection condition
is clearly valid in an open set of metric tensors.
Thus the set of manifolds for which the stated condition holds is  open and dense.
Let us then choose the sources $\tilde {f}_{a_\ell}$ to have generic polarizations.
We see that  in a generic set of polarizations the principal symbol
 of the 4th order interaction term $\M^{(4)}$ is non-vanishing on an open and dense set
 of $S_4$. Recall also that   in a generic situation the singular support of $\M^{(3)}$
 determines $S_3$.

 Thus in a generic situation
we can observe an open and dense subset of the surface $S_4$. Let us now consider how
we can decide which points on the union  of the surfaces 
$\Lambda^+_{q(a_1,a_2,a_3,a_4)}$ belong to the same surface.
Recall that the order of singularity (at points that do not belong to the 3-dimensional intersections) \HOX{Check 3/2}
on $\Lambda^+_{q(a_1,a_2,a_3,a_4)}$ is $(p_{a_1}+p_{a_2}+p_{a_3}+p_{a_4})-4-3/2$, see (\ref{eq: conormal}).
Let us now assume that $p_j$, $j\leq N$ are such that all the numbers
\ba
& &(p_{a_1}+p_{a_2}+p_{a_3}+p_{a_4})-1/2-n_1,
\ea
are different,
where  $n_1$ and $n_2$ are arbitrary integers (corresponding to different orders of the
classical symbols of terms) and \HOX{Gunther, do you know if the symbols of of waves $\M_j$ classical, assuming
that the sources have classical symbols, and if so, what would be a reference for this?}
the indexes $a_1,a_2,a_3,a_4\leq N$
satisfy $a_j\not =a_k$ for $j\not=k$. Let $\Sigma$ be the set of points $(x,\xi)\in T^*U$, $x\not \in S_3$ on the 
wave front set of  $\M^{(4)}$ such that the point $x$
has a neighborhood $V$ where $\M^{(4)}$  is a conormal
distribution associated smooth surface $S\subset V$ such that
$\M^{(4)}|_V\in \I^p(N^*S)$, where $p=(p_{a_1}+p_{a_2}+p_{a_3}+p_{a_4})-4-3/2$ \HOX{Check 3/2}
but $p$ can not be replaced by any smaller number. Then for
generic set of manifolds  in $\mathcal N$
the closure of the set $\Sigma$ coincides with 
 $\Lambda^+_{q(a_1,a_2,a_3,a_4)}$. \HOX{Check the construction of  $\Lambda^+_{q(a_1,a_2,a_3,a_4)}$}
This implies
that we can determine $\Lambda^+_{q(a_1,a_2,a_3,a_4)}\cap T^*U$ for  all
$\vec a=(a_1,a_2,a_3,a_4)$ and thus
we can find the whole smooth surface $\Lambda^+_{q(a_1,a_2,a_3,a_4)}\cap T^*U$.
This implies that we can find  
$G_N(J((M,g),(x_0,\xi_0)),[g],x_0,\xi_0)$. This determines the manifold
$(J((M,g),(x_0,\xi_0)),[g])$ up to a small error.

   Summarizing the above  sketch of construction
   means the following:
   \medskip
   
   {\bf Summary:} {\it For any $\e>0$ one can $L$ such that the following holds
   in a generic set of manifolds
    $((M,g),(x_0,\xi_0))\in \mathcal N$.
    Let   $({\hat x}_{a_\ell},\xi_{a_\ell})$, $|\ell |\leq L$ be constructed above and
    $F_\e$ be a source (\ref{Fe source}) that sends distorted plane  waves to directions
       $({\hat x}_{a_\ell},\xi_{a_\ell})$, $|\ell |\leq L $ with generic polarizations  $w_{a_\ell}$  and
       satisfies (\ref{eq: V ehto}). Then,
       by measuring singularities of terms 
   $\M^{(j)}$, $j\leq 4$ of the     
       the wave  $u^{(4)}_\e(x)$ on $U$
    produced by the source  $F_\e$,
   we can find the manifold  
    $J_{g}(\hat \mu(s_-),\hat \mu(s_+))$ and the conformal type 
    of the metric $[g]$ on it up to an error $\e$.}
    \medskip
   
 In other words,  when we consider the above measurement
 with single source and approximate the $O(\e^5)$ terms by zero,
 the observations of the singularities of the terms with different $\e^j$ magnitudes $j\leq 4$, 
 make it possible in a generic case to construct
a  discrete approximation of the earliest light observation point
sets, $\{\mathcal P_V(q_k);\ k\leq K\}.$
Such data, under appropriate conditions, could be used to determine a discrete approximation 
for the conformal structure of the space-time in an analogous manner used
in \cite{AKKLT,K2L}. However, making the approach described above to a detailed proof is outside
the scope of this paper  will be considered
elsewhere.}
\medskip
}}

\hiddenfootnote{
\section*{Appendix A: Reduced Einstein equation}

In this section we review known results on Einstein equations and wave maps.

\HOX{The appendices need to be shortened  in the final version of the paper.}

\subsection*{A.1.\ Summary of the used notations}


Let us recall some definitions given in Introduction, in the Subsection
\ref{subsec: Gloabal hyperbolicity}.
Let $(\hattuM ,g)$ be a $C^\infty$-smooth globally hyperbolic Lorentzian
manifold  and  $\tilde g$ be a $C^\infty$-smooth globally hyperbolic
metric on $M$ such that $g<\tilde g$. 

Recall that there is an isometry $\Phi:(M,\tilde g)\to (\R\times N,\tilde  h)$,
where $N$
is a 3-dimensional manifold and the metric $\tilde  h$ can be written as
$\tilde  h=-\beta(t,y) dt ^2+\overline h(t,y)$ where $\beta:\R\times  N\to (0,\infty)$ is a smooth function and 
$\overline h(t,\cdotp)$ is a Riemannian metric on $ N$ depending smoothly
on $t\in \R$. 
As in the main text, we identify these isometric manifolds and denote $\hattuM =\R\times N$.
%
Also, for $t\in \R$, recall that $\hattuM (t)=(-\infty,t)\times N$. We use parameters $t_1>t_0>0$
and denote $\hattuM _j=\hattuM (t_j)$, $j\in \{0,1\}$.
We use  the time-like geodesic    $\hat \mu =\mu_{g}$, $\mu_{g}:[-1,1]\to {\rbtext((-\infty,{\reftext T_0})\times N)}$ on $({\rbtext((-\infty,{\reftext T_0})\times N)},g)$
and  the set $\K_j:=J^+_{\tilde g}({p^-})\cap M_j$ with ${p^-}=
\hat \mu(s_-)\in (0,t_0)\times N$, $s_-\in (-1,1)$
and recall that  
$J^+_{\tilde g}({p^-})\cap M_j$ is compact and 
there exists $\e_0>0$ such that if $\|g-g\|_{C^0_b(\hattuM _1;g^+)}<\e_0$, then
$g|_{\K_1}<\tilde g|_{\K_1}$,
and in particular, we have  $J^+_{g}(p)\cap \hattuM _1\subset \K_1$
for all $p\in \K_1$.

Let us use local coordinates on $\hattuM _1$ and denote by
$\nabla_k=\nabla_{X_k}$  the covariant derivative with respect to the metric $g$
to the direction $X_k=\frac \p{\p x^p}$
 and by
 $\hat \nabla_k=\hat \nabla_{X_k}$ the covariant derivative with respect to the metric $g$
to the direction $X_k$.

\subsection*{A.2.\ Reduced Ricci and Einstein tensors}
Following \cite{FM} we recall that 
\beq\label{q-formula2copy}
\Ric_{\mu\nu}(g)&=& \Ric_{\mu\nu}^{(h)}(g)
+\frac 12 (g_{\mu q}\frac{\p \Gamma^q}{\p x^{\nu}}+g_{\nu q}\frac{\p \Gamma^q}{\p x^{\mu}})
\eeq
where  $\Gamma^q=g^{mn}\Gamma^q_{mn}$,
\beq\label{q-formula2copyB}
& &\hspace{-1cm}\Ric_{\mu\nu}^{(h)}(g)=
-\frac 12 g^{pq}\frac{\p^2 g_{\mu\nu}}{\p x^p\p x^q}+ P_{\mu\nu},
\\ \nonumber
& &\hspace{-2cm}P_{\mu\nu}=  
g^{ab}g_{ps}\Gamma^p_{\mu b} \Gamma^s_{\nu a}+ 
\frac 12(\frac{\p g_{\mu\nu }}{\p x^a}\Gamma^a  
+ \nonumber
g_{\nu l}  \Gamma^l _{ab}g^{a q}g^{bd}  \frac{\p g_{qd}}{\p x^\mu}+
g_{\mu l} \Gamma^l _{ab}g^{a q}g^{bd}  \frac{\p g_{qd}}{\p x^\nu}).\hspace{-2cm}
\eeq
Note that $P_{\mu\nu}$ is a polynomial of $g_{jk}$ and $g^{jk}$ and first derivatives of $g_{jk}$.
The harmonic  Einstein tensor is  
\beq\label{harmonic Ein}
\Ein^{(h)}_{jk}(g)=
\Ric_{jk}^{(h)}(g)-\frac 12 g^{pq}\Ric_{pq}^{(h)}(g)\, g_{jk}.
\eeq 
The  harmonic  Einstein tensor is extensively used to study
Einstein equations in local coordinates where one can use the Minkowski
space $\R^4$ as the background space. To do global constructions
with a background space $(M,g)$ one uses  the  reduced Einstein tensor.
The $g$-reduced Einstein tensor 
 $\Ein_{g} (g)$ and the $g$-reduced  Ricci tensor
  $\Ric_{g} (g)$  
are given 
by
\beq\label{Reduced Einstein tensor}
& &(\Ein_{g} (g))_{pq}=(\Ric_{g} (g))_{pq}-\frac 12 (g^{jk}(\Ric_{g} g)_{jk})g_{pq},\\& &
\label{Reduced Ric tensor}
(\Ric_{g} (g))_{pq}=\Ric_{pq} g-\frac 12 (g_{pn} \hat \nabla _q\hat F^n+ g_{qn} \hat \nabla _p\hat F^n)
\eeq
where $\hat F^n$ are the harmonicity functions given by
\beq\label{Harmonicity condition BBB}
\hat F^n=\Gamma^n-\hat \Gamma^n,\quad\hbox{where }
\Gamma^n=g^{jk}\Gamma^n_{jk},\quad
\hat \Gamma^n=g^{jk}\hat \Gamma^n_{jk},
\eeq
where $\Gamma^n_{jk}$ and $\hat \Gamma^n_{jk}$ are the Christoffel symbols
for $g$ and $g$, correspondingly.
Note that $\hat \Gamma^n$ depends also on $g^{jk}$.
As $\Gamma^n_{jk}-\hat \Gamma^n_{jk}$ is the difference of two connection
coefficients, it is a tensor. Thus $\hat F^n$ is tensor (actually, a vector field), implying that both
$(\Ric_{g} (g))_{jk}$ and $(\Ein_{g} (g))_{jk}$ are  2-covariant tensors.
Observe that  the $g$-reduced Einstein tensor is sum of the harmonic Einstein tensor 
and a term that is of the first order in $g$, 
\beq\label{hat g reduced einstein and reduced einstein}
(\Ein_{g} (g))_{\mu\nu}=\Ein^{(h)}_{\mu\nu} (g)+\frac 12 (g_{\mu q}\frac{\p \hat \Gamma^q}{\p x^{\nu}}+g_{\nu q}\frac{\p \hat \Gamma^q}{\p x^{\mu}}).
\eeq

%

%

%

\subsection*{A.3.\ Wave maps and reduced Einstein equations}
Let us consider the manifold $M_1=(-\infty,t_1)\times N$ with a $C^m $-smooth metric $g^{\prime}$,
$m\geq 8$,
 which
is a perturbation of the metric $g$ and satisfies the Einstein equation
\beq\label{Einstein on M_0}
\Ein(g^{\prime})=T^{\prime}\quad \hbox{on }M_1,
\eeq
or equivalently, 
\ba
\Ric(g^{\prime})=\rho^{\prime},\quad \rho^{\prime}_{jk}=T_{jk}^{\prime}-\frac 12 ((g^{\prime})^{  nm}T^{\prime}_{nm})g^{\prime}_{jk}\quad \hbox{on }M_1.
\ea
Assume also that  $g^{\prime}=g$ in the domain $A$,
where $A=M_1\setminus\K_1$ and  $\|g^{\prime}-g\|_{C^2_b(M_1,g^+)}<\e_0$,
so that $(M_1,g^{\prime})$ is globally hyperbolic. Note that then $T^{\prime}=\hat T$ in the set $A$ 
and that the metric $g^{\prime}$ coincides with $g$
in particular in the set $M^-=\R_-\times N$

We recall next the considerations of \cite{ChBook}.
Let us  consider the Cauchy problem for the wave map
$f:(M_1,g^\prime)\to (M,g)$, namely
\beq
\label{C-problem 1}& &\square_{g^{\prime},g} f=0\quad\hbox{in } M_1,\\
\label{C-problem 2}& &f=Id,\quad \hbox{in  }\R_-\times N,
\eeq
where $M_1=(-\infty,t_1)\times N\subset M$. In (\ref{C-problem 1}), 
$\square_{g^{\prime},g} f=g^\prime\,\cdotp\hat \nabla^2 f$ is the wave map operator, where $\hat \nabla$
is the covariant derivative of a map $(M_1,g^\prime)\to (M,g)$, see \cite[formula (7.32)]{ChBook}.
In  local coordinates
$X:V\to \R^4$ of $V\subset M_1$, denoted
$X(z)=(x^j(z))_{j=1}^4$ and $Y:W\to \R^4$ of $W\subset \hattuM $, denoted
$Y(z)=(y^A(z))_{A=1}^4$, 
the wave map $f:M_1\to \hattuM $ has representation $Y(f(X^{-1}(x)))=(f^A(x))_{A=1}^4$ and
the wave map operator in equation (\ref{C-problem 1}) is given by
\beq\label{wave maps2}
& &(\square_{g^{\prime},g} f)^A(x)=
(g^{\prime})^{  jk}(x)\bigg(\frac \p{\p x^j}\frac \p{\p x^k} f^A(x) -\Gamma^{{\prime} n}_{jk}(x)\frac \p{\p x^n}f^A(x)
\\
& &\quad \quad\quad\quad\quad\quad\quad\quad\nonumber
+\hat \Gamma^A_{BC}(f(x))
\,\frac \p{\p x^j}f^B(x)\,\frac \p{\p x^k}f^C(x)\bigg)\eeq
where $\hat \Gamma^A_{BC}$ denotes the Christoffel symbols of metric $g$
and  $\Gamma^{ {\prime} j}_{kl}$ are the Christoffel symbols of metric $g^{\prime}$.
When (\ref{C-problem 1})
is satisfied, we say that  $f$ is wave map with respect to the pair $(g^{\prime},g)$.
The important property of the wave maps is that, if 
$f$ is wave map with respect to the pair $(g^{\prime},g)$ and
$g=f_*g^{\prime}$ 
then, as follows from (\ref{wave maps2}), the identity map $Id:x\mapsto x$ is a wave map with respect to the pair $(g,g)$
and, the wave map equation for the identity map is equivalent to (cf.\ \cite[p.\ 162]{ChBook})
\beq\label{wave equation id}
\Gamma^n=\hat\Gamma^n,\quad \hbox{where }\Gamma^n=g^{jk}\Gamma^n_{jk},
\quad \hat \Gamma^n=g^{jk}\hat\Gamma^n_{jk}
\eeq
%
where the Christoffel symbols $\hat\Gamma^n_{jk}$ of the metric $g$ are smooth functions.

As   $g=g^{\prime}$ outside a compact set $\K_1\subset (0,t_1)\times N$,
we see that   this Cauchy problem is equivalent to the same equation
restricted to the
set $(-\infty,t_1)\times B_0$, where $B_0\subset N$ is such an open
relatively compact set that $\K_1\subset (0,t_1]\times B_0$
with the boundary condition $f=Id$ on  $(0,t_1]\times \p B_0$.
Then using results of \cite {HKM},  that can be applied for equations on manifold  as is
done  in Appendix B, and combined with the Sobolev embedding theorem, 
we see\footnote{See also:  
Thm.\ 4.2 in App.\ III  of
of \cite{ ChBook}, and its proof for the estimates
for the time on which the solution exists.} that 
 there is $\e_1>0$ such that if $\|g^\prime-g\|_{C^{m}_b(\hattuM _1;g^+)}<\e_1$,  then there  is a map  $f:M_1 \to M$ satisfying
the Cauchy problem (\ref{C-problem 1})-(\ref{C-problem 2})
and the solution depends continuously, 
 in
$C^{m-5}_b([0,t_1]\times N,g^+)$, on the metric $g^{\prime}$.
 Moreover, by the uniqueness of the wave map, we have
$f|_{M_1\setminus \K_1}=id$ so that
 $f(\K_1)\cap M_0\subset \K_0$.

 As
 the inverse function of the wave map $f$ depends continuously,
 in
$C^{m-5}_b([0,t_1]\times N,g^+)$, on the metric $g^{\prime}$
we can also assume that $\e_1$ is  so small that
$\hattuM _0\subset f(M_1)$.

Denote next $g:=f_*g^{\prime}$, $T:=f_*T^{\prime}$, and $\rho:=f_*\rho^{\prime}$ 
and define $\hat \rho=\hat T-\frac 12 (\hbox{Tr}\, \hat T)g.$ 
Then $g$ is $C^{m-6}$-smooth and the equation (\ref{Einstein on M_0}) implies 
\beq\label{Einstein on M1}
\Ein(g)=T\quad \hbox{on }\hattuM _0.
\eeq 
Since $f$ is a wave map, $g$ satisfies (\ref{wave equation id}) and thus  by the
definition of the reduced Einstein tensor,  (\ref{Reduced Einstein tensor}), we have
\ba
\Ein_{pq}(g)=
(\Ein_{g} (g))_{pq}\quad \hbox{on } \hattuM _0.
\ea
This and  (\ref{Einstein on M1}) yield the $g$-reduced Einstein equation
\beq\label{hat g reduced einstein equations}
(\Ein_{g} (g))_{pq}=T_{pq}\quad \hbox{on } \hattuM _0.
\eeq
This equation is useful for our considerations as it is a quasilinear,
hyperbolic equation on $\hattuM _0$. Recall that 
 $g$ coincides with
$g$ in $M_0\setminus \K_0$. The unique solvability of this 
Cauchy problem is studied in e.g.\ \cite[Thm.\ 4.6 and 4.13]{ ChBook} and \cite {HKM} 
and in Appendix B below.


\subsection*{A.4.\ Relation of the reduced Einstein equations and for the original Einstein equation}

The metric $g$ which solves the $g$-reduced Einstein
equation $\Ein_{g} (g)=T$ is a solution of the original
Einstein equations  $\Ein (g)=T$ if the harmonicity
functions $\hat F^n$
vanish identically. Next we recall the result that 
the harmonicity functions vanish on $\hattuM _0$
when 
\beq\label{eq: good system}
& &(\Ein_{g}(g))_{jk}=T_{jk},\quad \hbox{on }M_0,\\
\nonumber & &\nabla_pT^{pq}=0,\quad \hbox{on }M_0,\\
\nonumber & &g=g,\quad \hbox{on }M_0\setminus \K_0.
\eeq
To see this, let
us  denote $\Ein_{jk}(g)=S_{jk}$, $S^{jk}=g^{jn}g^{km}S_{nm}$,
and $T^{jk}=g^{jn}g^{km}T_{nm}$. 
Following the standard arguments,
see \cite{ ChBook}, we see from (\ref{Reduced Einstein tensor}) that in local coordinates
\ba
S_{jk}-(\Ein_{g}(g))_{jk}=\frac12(g_{jn}\hat \nabla_k \hat F^n+g_{kn}\hat \nabla_j\hat F^n-g_{jk}\hat \nabla_n \hat  F^n).
\ea
Using equations (\ref{eq: good system}), the Bianchi identity  $\nabla_pS^{pq}=0$,
and the basic property of Lorentzian connection,
$\nabla_kg^{nm}=0$,
we obtain 
\ba
0&=&
2\nabla_p(S^{pq}-T^{pq})\\
&=&\nabla_p(g^{qk}\hat  \nabla_k F^p+g^{pm}\hat\nabla_m \hat F^q-  g^{pq}\hat  \nabla_n \hat F^n)
\\
&=&g^{pm}\nabla_p \hat \nabla_m \hat F^q+(g^{qk} \nabla_p \hat \nabla_k \hat F^p-  g^{qp} \nabla_p \hat \nabla_n \hat F^n)\\
&=&g^{pm}\nabla_p \hat \nabla_m \hat F^q+W^q(\hat F)
\ea
where $\hat F=(\hat F^q)_{q=1}^4$ and
the operator 
\ba
W:(\hat F^q)_{q=1}^4\mapsto (g^{qk} (\nabla_p \hat \nabla_k \hat F^p- \nabla_k \hat \nabla_p \hat F^p))_{q=1}^4
\ea
is a linear first order differential operator which coefficients are polynomial functions of
$g_{jk}$, $g^{jk}$, $g_{jk}$, $g^{jk}$ and their first derivatives.

Thus the harmonicity functions $\hat F^q$ satisfy on $\hattuM _0$ the hyperbolic initial 
value problem 
\ba
& &g^{pm}\nabla_p \hat \nabla_m \hat F^q+W^q(\hat F)=0,\quad\hbox{on }M_0,\\
& & \hat F^q=0,\quad \hbox{on }M_0\setminus \K_0,
\ea
and as this initial 
Cauchy problem is uniquely solved by \cite[Thm.\ 4.6 and 4.13]{ ChBook}
or \cite {HKM}, we see that $\hat F^q=0$ on $\hattuM _0$. Thus
equations (\ref{eq: good system})
yield that
the Einstein equations   $\Ein(g)=T$ holds on $\hattuM _0$.

}

\hiddenfootnote{
\subsection*{Appendix B: Stability and existence of the  direct problem}

 {

Let us  start by explaining  how we can choose a  $C^\infty$-smooth metric $\tilde g$ 
such that $g<\tilde g$ and  $(M,\tilde g)$ is globally hyperbolic:
When $v(x)$ the eigenvector corresponding to the negative eigenvalue
of $g(x)$, we can choose a smooth, strictly positive function
$\eta:\hattuM \to \R_+$ such that $\tilde g^\prime:=g-\eta v\otimes v<\tilde g$. Then
$(\hattuM ,\tilde g^\prime)$ is globally hyperbolic, $\tilde g^\prime$ is smooth
and $g<\tilde g^\prime$. Thus we can replace $\tilde g$ by the smooth metric
 $\tilde g^\prime$ having the same properties that are required for $\tilde g$.
 
 Let us now return to consider existence and stability of the solutions of the
 Einstein-scalar field equations.
Let $t={\bf t}(x)$ be local time so that there is a diffeomorphism
$\Psi:M\to \R\times N$, $\Psi(x)=({\bf t}(x),Y(x))$, and
 $S(T)=\{x\in \hattuM _0; \ {\bf t}(x)=T\}$, $T\in \R$ are Cauchy surfaces.
 Let $t_0>0$. Next we identify $M$ and $\R\times N$ via the map $\Psi$ and just
 denote $M=\R\times N$.
Let us denote $M(t)=(-\infty,t)\times N$,  and $M_0=M(t_0)$.
By \cite[Cor.\ A.5.4]{BGP} the set $\K=J^+_{\tilde g}({p^-})\cap \hattuM _0$,
where $\hat p^0\in M_0$,
is compact.
Let $N_1,N_2\subset N$ be such
 open relatively compact sets with smooth boundary that
 $N_1\subset N_2$ and  $Y(J_{\tilde g}^+(\hat p^0)\cap \hattuM _0)\subset N_1$.  \HOX{We could improve explanation on $\tilde N$}
 
 To simplify citations to existing literature, 
let us define $\tilde N$ to be a compact manifold without boundary such
that $N_2$ can be considered as a subset of $\tilde N$. Using a construction based on
a suitable partition of unity, the Hopf double of the manifold $N_2$, 
and the Seeley extension of the metric tensor,  we can
endow $\tilde M=(-\infty,t_0)\times \tilde N$ with a smooth Lorentzian
metric $g^e$ (index $e$ is for "extended")
so that $\{t\}\times \tilde N$ are Cauchy surfaces of $\tilde N$
and that $g$ and $g^e$ coincide in  the set $\Psi^{-1}((-\infty,t_0)\times N_1)$ that contains the set $J^+_{g}(\hat p^0)\cap \hattuM _0$.} We extend the metric $\tilde g$ to a (possibly non-smooth) globally
hyperbolic metric $\tilde g^e$ on $\tilde M_0=(-\infty,t_0)\times \tilde N)$
such that $g^e<\tilde g^e$.

 To simplify notations below we denote $g^e=g$ and
 $\tilde g^e=\tilde g$ 
 on the whole
 $\tilde M_0$.
 Our aim is to prove the estimate (\ref{eq: Lip estim}).

Let us denote by $t={\bf t}(x)$ the local time. 
Recall that when $(g,\phi)$ is a solution of 
the scalar field-Einstein equation, we denote $u=(g-g,\phi-\hat \phi)$.
We will consider the equation for $u$, and to emphasize that the metric
depends on $u$, we denote $g=g(u)$ and assume below that  both the metric $g$ is  dominated
by $\tilde g$, that is, $ g<\tilde g$.
We use the pairs ${\bf u}(t)=(u(t,\cdotp),\p_t u(t,\cdotp))\in
H^1(\tilde N)\times L^2(\tilde N)$ and the notations
 ${\bf v}(t)=(v(t),\p_t v(t))$ etc. 
Let us consider a generalization of the system (\ref{eq: notation for hyperbolic system 1}) of the form 
\beq\label{eq: notation for hyperbolic system}
& &\square_{g(u)}u+V(x,D)u+H(u,\p u) =R(x,u,\p u)F+K,\ \ x\in \tilde M_0,\hspace{-1cm}\\
& & \nonumber \supp(u)\subset \K,
\eeq
where $\square_{g(u)}$ is the Lorentzian Laplace operator operating on the sections 
of the bundle $\B^L$ on $M_0$ and 
$\supp(F)\cup \supp(K)\subset \K$.
Note that above  $u=(g-g,\phi-\hat \phi)$ and $g(u)=g$.
Also, $F\mapsto R(x,u,\p u)F$ is a 
 linear first order differential operator which coefficients at $x$ are depending smoothly on 
 $u(x)$, $\p_j u(x)$ and the  derivatives 
of $(g,\hat \phi)$ at $x$ and
\ba V(x,D)=V^j(x)\p_j+V(x)\ea is a linear first order differential operator which coefficients 
at $x$ are depending smoothly on the  derivatives 
of $(g,\hat \phi)$  at $x$, and finally, $H(u,\p u)$ is a polynomial
of $u(x)$ and $\p_j u(x)$ which coefficients 
at $x$ are depending smoothly on the  derivatives 
of $(g,\hat \phi)$   such that
 $\p^\a_v\p^\beta_w H(v,w)|_{v=0,w=0}=0$ for
$|\a|+|\b|\leq 1$. 
By \cite[Lemma 9.7]{Ringstrom}, the equation  (\ref{eq: notation for hyperbolic system})
has at most one solution with given  $C^2$-smooth source functions $F$ and $K$.
Next we consider the existence of $u$ and its dependency on $F$ and $K$.


%

Below we use notations, c.f. (\ref{eq: notation for hyperbolic system 1}) and
(\ref{eq: notation for hyperbolic system 1b}) \ba
{\mathcal R}({\bf u},F)=R(x,u(x),\p u(x))F(x),\quad
{\mathcal H}({\bf u})=H(u(x),\p u(x)).\ea
Note that $u=0$, i.e., $g=g$ and $\phi=\hat\phi$ satisfies (\ref{eq: notation for hyperbolic system}) with $F=0$
and $K=0$.
Let us use the same notations as in  \cite {HKM} cf. also \cite[section 16]{Kato1975}, 
 to consider quasilinear wave equation on $[0,t_0]\times \tilde N$.
Let $\H^{(s)}(\tilde N)=H^s(\tilde N)\times H^{s-1}(\tilde N)$ and 
\ba
Z=\H^{(1)}(\tilde N),\quad  Y=\H^{(k+1)}(\tilde N),\quad X=\H^{(k)}(\tilde N).
\ea
The norms on these space are defined invariantly using the smooth Riemannian 
metric $h=g|_{\{0\}\times \tilde N}$ on $\tilde N$. \HOX{Maybe we should explain the 
bundles, the connection and the wave operator on section in a detailed way.}
Note that $\H^{(s)}(\tilde N)$ are in fact the Sobolev spaces of sections on the bundle $\pi:\B_K\to
\tilde N$, where $\B_K$ denotes also the  pull back bundle of $\B_K$ on $\tilde M$ in the map
$id:\{0\}\times \tilde N\to \tilde M_0$,
or on the bundle $\pi:\B_L\to \tilde N$. Below, $\nabla_h$ denotes the  standard connection of the bundle  $\B_K$
or  $\B_L$ associated to the metric $h$.

Let $k\geq 4$ be an even integer.
By definition of $H$ and $R$ we see that  there are  $0<r_0<1$ and $L_1,L_2>0$,
all depending on $g$, $\hat \phi$, $\K$, and $t_0$, such that if 
$0<r\leq r_0$ and
\beq\label{new basic conditions}
& &\|{\bf v} \|_{C([0,t_0];\H^{(k+1)}(\tilde N))}\leq r,\quad \|{\bf v}^\prime\|_{C([0,t_0];\H^{(k+1)}(\tilde N))}\leq r,\\
& &\|F\|_{C([0,t_0];H^{(k+1)}(\tilde N))}\leq r^2,\quad \|K\|_{C([0,t_0];H^{(k+1)}(\tilde N))}\leq r^2 \nonumber\\
& &\|F^\prime\|_{C([0,t_0];H^{(k+1)}(\tilde N))}\leq r^2,\quad \|K^\prime\|_{C([0,t_0];H^{(k+1)}(\tilde N))}\leq r^2 \nonumber
\eeq
then
\beq\label{new f1f2 conditions}
& &\quad \quad \|g(\cdotp;v)^{-1}\|_{C([0,t_0];H^s(\tilde N))}\leq L_1 ,\\ \nonumber
& &\|{\mathcal H} ({\bf v})\|_{C([0,t_0];H^{s-1}(\tilde N))}\leq L_2r^2,
\quad \|{\mathcal H} ({\bf v}^\prime)\|_{C([0,t_0];H^{s-1}(\tilde N))}\leq L_2r^2,\\ \nonumber
& &\|{\mathcal H} ({\bf v})-{\mathcal H} ({\bf v}^\prime)\|_{C([0,t_0];H^{s-1}(\tilde N))}\leq  L_2r\, \|{\bf v}-{\bf v}^\prime\|_{C([0,t_0];\H^{(s)}(\tilde N))}, 
\\ \nonumber
& &\|\mathcal R ({\bf v}^\prime,F^\prime)\|_{C([0,t_0];H^{s-1}(\tilde N))}\leq L_2r^2,
\quad \|\mathcal R ({\bf v},F)\|_{C([0,t_0];H^{s-1}(\tilde N))}\leq L_2r^2,\hspace{-1cm} \\ \nonumber
& &\|\mathcal R ({\bf v},F)-\mathcal R ({\bf v}^\prime, F^\prime)\|_{C([0,t_0];H^{s-1}(\tilde N))}
\\ \nonumber & &\leq  
L_2r\, \|{\bf v}-{\bf v}^\prime\|_{C([0,t_0];\H^{(s)}(\tilde N))}+L_2\|F-{F}^\prime\|_{C([0,t_0];H^{s+1}(\tilde N))\cap C^1([0,t_0];H^{s}(\tilde N))},\hspace{-2cm}
%
  \eeq
 for all $s\in [1,k+1]$. 
%

%
%

Next we write (\ref{eq: notation for hyperbolic system}) as a first order
system.  To this end,
let $\mathcal A(t,{\bf v}):\H^{(s)}(\tilde N)\to \H^{(s-1)}(\tilde N)$ be the operator 
$\mathcal A(t,{\bf v})=\mathcal A_0(t,{\bf v})+\mathcal A_1(t,{\bf v})$
where in local coordinates and in local trivialization of the bundle $\B^L$
\ba
\mathcal A_0(t,{\bf v})=-\left(\begin{array}{cc} 0 & I\\
 \frac 1{g^{00}(v)}\sum_{j,k=1}^3 g^{jk}(v)\frac \p{\p x^j}\frac\p{ \p x^k}\quad &
 \frac 1{g^{00}(v)}\sum_{m=1}^3 g^{0m}(v)\frac {\p}{\p x^m} \end{array}\right)
\ea 
with $g^{jk}(v)=g^{jk}(t,\cdotp;v)$ is a function on $\tilde N$ and 
\ba
\mathcal A_1(t,{\bf v})= \frac {-1}{g^{00}(v)}\left(\begin{array}{cc} 0 & 0\\
\sum_{j=1}^3 
 B^j(v)
   \frac \p{\p x^j}\quad &
B^{0}(v)\end{array}\right)
\ea 
where $B^j(v)$ depend on $v(t,x)$ and its first derivatives, and 
the connection coefficients (the Christoffel  symbols) corresponding  to $g(v)$. We denote $\S=(F,K)$ and 
 \ba
& &f_\S(t,{\bf v})=(f_\S^1(t,{\bf v}),f_\S^2(t,{\bf v}))\in \H^{(k)}(\tilde N),
\quad\hbox{where}\\
& & f_\S^1(t,{\bf v})=0,\quad f_\S^2(t,{\bf v})=\mathcal R({\bf v},F)(t,\cdotp)-\mathcal H({\bf v})(t,\cdotp)+K(t,\cdotp).
 \ea
 
 Note that when (\ref{new basic conditions}) are satisfied with $r<r_0$, inequalities (\ref{new f1f2 conditions}) imply that there exists $C_2>0$ so that
 \beq\label{fF bound r2}
 & &\|f_\S(t,{\bf v})\|_Y+\|f_{\S^\prime}(t,{\bf v}^\prime)\|_Y\leq C_2r^2,\\
 & &\|f_\S(t,{\bf v})-f_\S(t,{\bf v}^\prime)\|_Y\leq C_2
 r\,\|{\bf v}-{\bf v}^\prime\|_{C([0,t_0];Y)}. \nonumber
 \eeq

%
%

Let $U^{\bf v}(t,s)$ be the wave propagator corresponding to metric $g(v)$, that is,
$U^{\bf v}(t,s): {\bf h}\mapsto {\bf w}$, where ${\bf w}(t)=(w(t),\p_t w(t))$  solves
\ba
(\square_{g(v)}+V(x,D))w=0\quad\hbox{for } (t,y)\in [s,t_0]\times \tilde N,\quad\hbox{with }
 {\bf w}(s,y)={\bf h}.
 \ea 
 Let $S=(\nabla_{h}^*\nabla_{h}+1)^{k/2}:Y\to Z$ be an isomorphism.
As $k$ is an even integer, we see using multiplication estimates for Sobolev
spaces, see e.g.\
  \cite [Sec.\ 3.2, point (2)]{HKM}, that there exists $c_1>0$ (depending on
  $r_0,L_1,$ and $L^2$) so that  $\A(t,{\bf v})S-S\A(t,{\bf v})=C(t,{\bf v})$,
 where $\|C(t,{\bf v})\|_{Y\to Z}\leq c_1$ for all ${\bf v}$ satisfying
 (\ref{new basic conditions}). This yields that 
 the property (A2) in  \cite {HKM} holds, namely
 that $S\A(t,{\bf v})S^{-1}=\A(t,{\bf v})+B(t,{\bf v})$ where
$B(t,{\bf v})$ extends to a bounded operator in $Z$
 for which $\|B(t,{\bf v})\|_{Z\to Z}\leq c_1$ for all ${\bf v}$ satisfying
 (\ref{new basic conditions}).
Alternatively, to see the mapping properties of $B(t,{\bf v})$ we could use the fact that 
 $B(t,{\bf v})$ is a zeroth order pseudodifferential operator with
$H^{k}$-symbol.\hiddenfootnote{On mapping properties of such
pseudodifferential operators, see J. Marschall, Pseudodifferential operators with coefficients in Sobolev spaces. Trans. Amer. Math. Soc. 307 (1988), no. 1, 335-361.} 
 
 Thus the proof of  \cite [Lemma 2.6]{HKM} shows\hiddenfootnote{ Alternatively, as $U^{\bf v}(t,s)$ are propagators
  for linear wave equations having finite speed of wave propagation, one can prove 
the  estimates (\ref{U bounds}) by considering the wave equation
in local coordinate neighborhoods $W_j\subset W_j^\prime \subset \tilde N$,
$\Phi_j:  W_j^\prime\to \R^3$ and a partition of unity $\psi_j\in C^\infty_0(W_j)$ and
cut-off functions $\psi^\prime_j\in C^\infty_0(W_j^\prime)$.
Then  \cite [Lemma 2.6]{HKM}, used  in local coordinates, shows
that when $t_k-t_{k-1}$ is small enough then 
\ba 
\| \psi^\prime_j U^{\bf v}(t_k,t_{k-1})(\psi_j {\bf w}) \|_{Y}\leq C_{3,jk}^{\prime}e^{C_4(t-s)}\| \psi_j {\bf w} \|_{Y}.
\ea
Using  sufficiently small time steps $t_k-t_{k-1}$ and combining the estimates
for different $j$:s
together, one obtain estimates (\ref{U bounds}).}
 that 
%
there is a constant  $C_3>0$
   so that 
  \beq
  \label{U bounds}\|U^{\bf v}(t,s)\|_{Z\to Z}\leq C_3\quad\hbox{and}\quad
  \|U^{\bf v}(t,s) \|_{Y\to Y}\leq C_3\hspace{-1.5cm}
  \eeq 
  for $0\leq s<t\leq t_0$.
    By interpolation of estimates (\ref{U bounds}), we see also that 
    \beq
  \label{U bounds2}
  \|U^{\bf v}(t,s)\|_{X\to X}\leq C_3,
  \eeq  for $0\leq s<t\leq t_0$.

Let us next modify the reasoning given in \cite{Kato1975}: let  $r_1\in (0,r_0)$
be a parameter which value will be chosen later,  $C_1>0$
and 
$E$ be the space of functions ${\bf u}\in C([0,t_0];X)$ for which
\beq
\label{1 bounded}& &\|{\bf u}(t)\|_Y\leq r_1\quad\hbox{and}\\
\label{2 Lip}& &\|{\bf u}(t_1)-{\bf u}(t_2)\|_X\leq C_1|t_1-t_2|\eeq
for all $t,t_1,t_2\in [0,t_0]$.
The set $E$ is endowed by the metric of  $C([0,t_0];X)$. We 
note that  by \cite[ Lemma 7.3] 
{Kato1975},
a convex $Y$-bounded, $Y$-closed set is closed also in $X$.
Similarly, functions $G:[0,t_0]\to X$  satisfying (\ref{2 Lip}) 
form a closed subspace of $C([0,t_0];X)$. Thus $E\subset X$ is a closed set implying
that $E$ is a complete metric space.

Let
 \ba 
 & &\hspace{-1cm}W=\\ & &\hspace{-1cm}\{(F,K)\in C([0,t_0];H^{k+1}(\tilde N))^2;\ \sup_{t\in [0,t_0]}\|F(t)\|_{H^{k+1}(\tilde N)}+\|K(t)\|_{H^{k+1}(\tilde N)}< r_1\}.\hspace{-1cm}
 \ea

Following \cite[p. 44]{Kato1975}, we see that
the solution of
equation (\ref{eq: notation for hyperbolic system}) with the source $\S\in W$  is found as
a fixed point, if it exists, of the map $\Phi_\S:E\to C([0,t_0];Y)$ where 
$\Phi_\S({\bf v})={\bf u}$ is  given by
$$
{\bf u}(t)=\int_0^t U^{\bf v}(t,\tilde t)f_\S(\tilde t,{\bf v})\,d\tilde t,\quad 0\leq t\leq t_0.
$$ 
Below, we denote  ${\bf u}^{\bf v}=\Phi_\S({\bf v})$.

Since $\Phi_{\S_0}(0)=0$ where
$\S_0=(0,0)$, we see using 
the above and the inequality $\|\,\cdotp\|_X\leq \|\,\cdotp\|_Y$ that the
 function ${\bf u}^{\bf v}$ satifies
  \ba
& &  \|{\bf u}^{\bf v}\|_{C([0,t_0];Y)}\leq C_3C_2t_0 r_1^2,\\
& &  \|{\bf u}^{\bf v}(t_2)-{\bf u}^{\bf v}(t_1)\|_{X}\leq C_3C_2r_1^2|t_2-t_1|,\quad t_1,t_2\in [0,t_0].
  \ea
  When $r_1>0$ is so small that $C_3C_2(1+t_0)<r_1^{-1}$ 
  and $C_3C_2r_1^2<C_1$
  we see that $  \|\Phi_\S({\bf v})\|_{C([0,t_0];Y)}<r_1$
  and $ \|\Phi_\S({\bf v})\|_{C^{0,1}([0,t_0];X)}<C_1$.
  Hence $ \Phi_\S(E)\subset E$
 and we can consider $\Phi_\S$ as a map
 $\Phi_\S:E\to E$.
  
  As  $k>1+\frac 32$, it follows from Sobolev embedding theorem
  that $X=\H^{(k)}(\tilde N)\subset C^1(\tilde N)^2$. This yields that
  by \cite[Thm. 3]{Kato1975}, for the original reference, see Theorems III-IV
  in \cite{Kato1973},
  \ba
 & & \|(U^{\bf v}(t,s)- U^{\bf v^\prime}(t,s)){\bf h}\|_X \\
 & &\leq 
 C_3\bigg(\sup_{t^\prime \in [0,t]}\|\A(t^\prime,{\bf v})-
 \A(t^\prime,{\bf v}^\prime)\|_{Y\to X}  \|U^{\bf v}(t^\prime,0) {\bf h}\|_Y\bigg)\\
 & &\leq 
 C_3^2 \|{\bf v}-{\bf v}^\prime\|_{C([0,t_0];X)}  \|{\bf h}\|_Y.
  \ea
  Thus,
  \ba
&&\|  U^{\bf v}(t,s)f_\S(s,{\bf v})- U^{\bf v^\prime}(t,s)f_{\S}(s,{\bf v}^\prime)\|_X\\
&\leq&\|  (U^{\bf v}(t,s)- U^{\bf v^\prime}(t,s))f_{\S}(s,{\bf v})\|_X+
\|  U^{\bf v^\prime}(t,s)(f_\S(s,{\bf v})- f_{\S}(s,{\bf v}^\prime))\|_X\\
&\leq&
(1+C_3)^2C_2r_1^2 \|{\bf v}-{\bf v}^\prime\|_{C([0,t_0];X)} .
  \ea
This implies that
\ba
& & \| \Phi_\S({\bf v})-\Phi_{\S}({\bf v}^\prime)\|_{C([0,t_0];X)}
\leq t_0(1+C_3)^2C_2r_1^2 \|{\bf v}-{\bf v}^\prime\|_{C([0,t_0];X))}.
\ea

Assume next that  $r_1>0$ is so small that we have also
\ba
t_0(1+C_3)^2C_2r_1^2<\frac 12.
\ea
 cf.\
  Thm.\ I in \cite{HKM} (or
  (9.15)  and (10.3)-(10.5) in \cite{Kato1975}). For $\S\in W$ this 
 implies that $\Phi_\S:E\to E$ is a contraction with a contraction constant 
 $C_L\leq \frac 12$,
  and thus 
 $\Phi_\S$ has a unique fixed point ${\bf u}$  in the space $E\subset C^{0,1}([0,t_0];X)$.
 
 Moreover, elementary considerations related to fixed point
 of the map $\Phi_\S$ show that ${\bf u}$  in  $C([0,t_0];X)$ depends in $E\subset C([0,t_0];X)$
 Lipschitz-continuously on
 $\S\in  W\subset C([0,t_0];H^{k+1}(\tilde N))^2$. Indeed, if $\|\S-\S^\prime\|_{C([0,t_0];H^{k+1})^2}<\e$,
 we see that
  \beq\label{fF bound r2 B}
 \|f_\S(t,{\bf v})-f_{\S^\prime}(t,{\bf v})\|_Y\leq C_2\e,\quad t\in [0,t_0],
 \eeq
 and when  (\ref{fF bound r2}) and (\ref{U bounds}) are satisfied 
 with $r=r_1$,
 we have
  \ba
& & \| \Phi_\S({\bf v})-\Phi_{\S^\prime}({\bf v}^\prime)\|_{C([0,t_0];Y)}
\leq C_3C_2t_0 r_1^2.
\ea
 Hence 
 \ba
 \|\Phi_\S({\bf v})-\Phi_{\S^\prime}({\bf v})\|_{C([0,t_0];Y)}\leq  t_0C_3C_2\e.
 \ea
 This and standard estimates for fixed points, yield that when $\e$ is small enough 
 the fixed point ${\bf u}^\prime$ of the map  $\Phi_{\S^\prime}:E\to E$ 
 corresponding to the source ${\S^\prime}$ and
  the fixed point ${\bf u}$ of the map $\Phi_{\S}:E\to E$ 
  corresponding to the source ${\S}$ satisfy
  \beq\label{stability}
  \|{\bf u}-{\bf u}^\prime\|_{C([0,t_0];X)}\leq  \frac 1{1-C_L}t_0C_3C_2\e.
  \eeq
{
Thus the solution ${\bf u}$ 
  depends in $C([0,t_0];X)$ Lipschitz continuously on
 $\S\in  C([0,t_0];H^{k+1}(\tilde N))^2$ (see also \cite[Sect.\ 16]{Kato1975}, 
 and \cite{Ringstrom})}. 
In fact,  for analogous systems it is possible to show that $u$ is in $C([0,t_0];Y)$,
but one can not obtain  
 Lipschitz or H\"older  stability for $u$ in the $Y$-norm, see  \cite{Kato1975}, Remark 7.2. 
 
 Finally, we note that the fixed point ${\bf u}$ of $\Phi_\S$ can be found as a limit
 $ {\bf u}=\lim_{n\to \infty}{\bf u}_n$ in $C([0,t_0];X)$, where ${\bf u}_0=0$ and  ${\bf u}_n=\Phi_\S({\bf u}_{n-1})$.
Denote ${\bf u}_n=(g_n-g,\phi_n-\hat \phi)$. We see that if
 $\supp({\bf u}_{n-1})\subset J_{g}(\supp(\S))$ 
then also $\supp(g_{n-1}-g)\subset J_{g}(\supp(\S))$.
Hence  for all $x\in M_0\setminus J_{g}(\supp(\S))$ we see that $J_{g_{n-1}}^-(x)
\cap J_{g}(\supp(\S))=\emptyset.$ Then, using the definition of the map $\Phi_\S$ we see that
 $\supp({\bf u}_n)\subset J_{g}(\supp(\S))$. Using induction we see that this holds
for all $n$ and hence we see that the solution ${\bf u}$ satisfies
 \beq\label{eq: support condition for u}
 \supp({\bf u})\subset J_{g}(\supp(\S)).
 \eeq 
 \hiddenfootnote{REMOVE MATERIAL IN THIS FOOTNOTE OR MOVE IT ELSEWHERE:
  Let us next use the above considerations to study system 
   (\ref{eq: notation for hyperbolic system 1}) on $\R\times N$.
   For this end, we denote below the background metric of  $\tilde M=\R\times \tilde N$
   by $\tilde g$ and the background metric of  $\hattuM _0=\R\times  N$
   by $g$.
 Consider a source 
\ba
\S\in {C^1([0,T];H^{s_0-1}(\tilde N))\cap
C^0([0,t_0];H^{s_0}(\tilde N))}
\ea  that is supported in a compact set $\K$ and
small enough in the norm of this space and let $\tilde u$ 
be the solution of the  
system (\ref{eq: notation for hyperbolic system 1}) on $\R\times \tilde N$.
Assume that  $u$ is some $C^2$-solution of  
(\ref{eq: notation for hyperbolic system 1}) on $\R\times N$
with the same source $\S$. Our aim is to show that it
coincides with $\tilde u$ in its support.
%
For this end, let $\tilde g^\prime$ be a metric tensor which
has the same eigenvectors as $\tilde g$ and the same
eigenvalues except that  the unique negative eigenvalue of $\tilde g$
 is multiplied by $(1-h_1)$, $h_1>0$. 
Assume that $h_1>0$ is so small
that $J^+_{\hattuM _0,\tilde g^\prime}(p^+)\subset  \hattuM _0\subset N_0\times (-\infty,t_0]$
and assume that $\e$ is so small that $g(\tilde u)<\tilde g^\prime$ for all 
sources $\S$ with $\|\S\|_{C^1([0,t_0];H^{s_0-1}(\tilde N))\cap
C^0([0,t_0];H^{s_0}(\tilde N))}\leq \e$. Then
 $J^+_{\tilde M,g(\tilde u)}(p^+)\subset  (-\infty,t_0]\times N_0$.
This  proves the estimate (\ref{eq: Lip estim}).

Next us we consider solution $u$ on $\hattuM _0$ corresponding to source $\S$.
Let  $T(h_1,\S)$
be the supremum of all $T^\prime\leq t_0$ for which
$J^+_{g(u)}(p^+)\cap \hattuM _0(T^\prime)\subset N_1\times (-\infty,T^\prime]$.
Assume that $ T(h_1,\S)<t_0$.
Then for $T^\prime= T(h_1,\S)$ we see that
$u$ can be continued by zero from $ (-\infty,T^\prime]\times N_1$
to a function on $(-\infty,T^\prime]\times \tilde N$ that 
satisfies the system (\ref{eq: notation for hyperbolic system 1}). Thus
$u$ coincides in $ (-\infty,T^\prime]\times N_1$
with $\tilde u$ and vanishes outside this set.
Since  $J^+_{\tilde M,g(\tilde u)}(p^+)\subset  (-\infty,t_0]\times N_0$, this implies that
also $u$ satisfies 
$J^+_{g(u)}(p^+)\cap \hattuM _0(T^\prime)\subset (-\infty,T^\prime]\times N_0$.
As $u$ is $C^2$-smooth and $N_0\subset \subset N_1$
we see that there is $T^\prime_1> T(h_1,\S)$ so 
that $J^+_{g(u)}(p^+)\cap \hattuM _0(T^\prime_1)\subset  (-\infty,T^\prime_1]
\times N_1$ that is in contradiction with definition of $ T(h_1,\S)$. Thus
we have to have  $T(h_1,\S)=t_0$.


This shows that when the norm of $\S$ is smaller than the above chosen $\e$
then $J^+_{g(u)}(p^+)\cap \hattuM _0\subset N_0\times (-\infty,t_0]$.
Then $u$ vanishes outside $N_1\times (-\infty,t_0]$ by
the support condition in (\ref{eq: notation for hyperbolic system 1}).
Thus $u$ and $\tilde u$ are both supported on  $N_1\times (-\infty,t_0]$ 
and coincide there. As $\tilde u$ is unique, 
this shows that for a sufficiently small sources $ \S$
 supported on  $\K$   the solution $u$ of
(\ref{eq: notation for hyperbolic system 1})  in $M_0$
exists and is unique.}

\bigskip

}

\extension{

\section*{Appendix A: Reduced Einstein equation}

In this section we review known results on the Einstein equations and wave maps.

\subsection*{A.1.\ Summary of the used notations}


Let us recall some definitions given in Introduction.
Let $(\hattuM ,g)$ be a $C^\infty$-smooth globally hyperbolic Lorentzian
manifold  and  $\tilde g$ be a $C^\infty$-smooth globally hyperbolic
metric on $M$ such that $g<\tilde g$.
Let us  start by explaining  how one can construct a  $C^\infty$-smooth metric $\tilde g$ 
such that $g<\tilde g$ and  $(M,\tilde g)$ is globally hyperbolic:
When $v(x)$ is an eigenvector corresponding to the negative eigenvalue
of $g(x)$, we can choose a smooth, strictly positive function
$\eta:\hattuM \to \R_+$ such that $$\tilde g^\prime:=g-\eta v\otimes v<\tilde g.$$ Then
$(\hattuM ,\tilde g^\prime)$ is globally hyperbolic, $\tilde g^\prime$ is smooth
and $g<\tilde g^\prime$. Thus we can replace $\tilde g$ by the smooth metric
 $\tilde g^\prime$ having the same properties that are required for $\tilde g$.

Recall that there is an isometry $\Phi:(M,\tilde g)\to (\R\times N,\tilde  h)$,
where $N$
is a 3-dimensional manifold and the metric $\tilde  h$ can be written as
$\tilde  h=-\beta(t,y) dt ^2+\overline h(t,y)$ where $\beta:\R\times  N\to (0,\infty)$ is a smooth function and 
$\overline h(t,\cdotp)$ is a Riemannian metric on $ N$ depending smoothly
on $t\in \R$. 
As in the main text we identify these isometric manifolds and denote $\hattuM =\R\times N$.
%
Also, for $t\in \R$, recall that $\hattuM (t)=(-\infty,t)\times N$. We use parameters $t_1>t_0>0$
and denote $\hattuM _j=\hattuM (t_j)$, $j\in \{0,1\}$.
\modified{We use  the time-like geodesic    $\hat \mu =\mu_{g}$, $\mu_{g}:[-1,1]\to M_0$ on $(M_0,g)$
and  the set $\K_j:=J^+_{\tilde g}(\hat \mu(-1))\cap M_j$ with $\hat \mu(-1)=
\in (-\infty,t_0)\times N$. Then 
$J^+_{\tilde g}(\hat \mu(-1))\cap M_j$ is compact. Also,
there exists $\e_0>0$ such that if 
$g$ is a Lorentz metric in $M_1$ such that 
$\|g-g\|_{C^0_b(\hattuM _1;g^+)}<\e_0$, then
$g|_{\K_1}<\tilde g|_{\K_1}$. In particular, this implies that we have  $J^+_{g}(p)\cap \hattuM _1\subset \K_1$
for all $p\in \K_1$. Later, we use this property to deduce that
when $g$ satisfies the $g$-reduced Einstein equations 
in $M_1$, with a source that is supported in $ \K_1$ and has
small enough norm is a suitable space, then $g$ coincides with $g$ in
$M_1\setminus \K_1$ and satisfies $g<\tilde g$.}

Let us use local coordinates on $\hattuM _1$ and denote by
$\nabla_k=\nabla_{X_k}$  the covariant derivative with respect to the metric $g$
in the direction $X_k=\frac \p{\p x^k}$
 and by
 $\hat \nabla_k=\hat \nabla_{X_k}$ the covariant derivative with respect to the metric $g$
to the direction $X_k$.

\subsection*{A.2.\ Reduced Ricci and Einstein tensors}
Following \cite{FM} we recall that 
\beq\label{q-formula2copy AAA}
\Ric_{\mu\nu}(g)&=& \Ric_{\mu\nu}^{(h)}(g)
+\frac 12 (g_{\mu q}\frac{\p \Gamma^q}{\p x^{\nu}}+g_{\nu q}\frac{\p \Gamma^q}{\p x^{\mu}})
\eeq
where  $\Gamma^q=g^{mn}\Gamma^q_{mn}$,
\beq\label{q-formula2copyBa}
& &\hspace{-1cm}\Ric_{\mu\nu}^{(h)}(g)=
-\frac 12 g^{pq}\frac{\p^2 g_{\mu\nu}}{\p x^p\p x^q}+ P_{\mu\nu},
\\ \nonumber
& &\hspace{-2cm}P_{\mu\nu}=  
g^{ab}g_{ps}\Gamma^p_{\mu b} \Gamma^s_{\nu a}+ 
\frac 12(\frac{\p g_{\mu\nu }}{\p x^a}\Gamma^a  
+ \nonumber
g_{\nu l}  \Gamma^l _{ab}g^{a q}g^{bd}  \frac{\p g_{qd}}{\p x^\mu}+
g_{\mu l} \Gamma^l _{ab}g^{a q}g^{bd}  \frac{\p g_{qd}}{\p x^\nu}).\hspace{-2cm}
\eeq
Note that $P_{\mu\nu}$ is a polynomial of $g_{jk}$ and $g^{jk}$ and first derivatives of $g_{jk}$.
The harmonic  Einstein tensor is  
\beq\label{harmonic Ein}
\Ein^{(h)}_{jk}(g)=
\Ric_{jk}^{(h)}(g)-\frac 12 g^{pq}\Ric_{pq}^{(h)}(g)\, g_{jk}.
\eeq 
The  harmonic  Einstein tensor is extensively used to study
the Einstein equations in local coordinates where one can use the Minkowski
space $\R^4$ as the background space. To do global constructions
with a background space $(M,g)$ one uses  the  reduced Einstein tensor.
The $g$-reduced Einstein tensor 
 $\Ein_{g} (g)$ and the $g$-reduced  Ricci tensor
  $\Ric_{g} (g)$  
are given 
by
\beq\label{Reduced Einstein tensora}
& &(\Ein_{g} (g))_{pq}=(\Ric_{g} (g))_{pq}-\frac 12 (g^{jk}(\Ric_{g} g)_{jk})g_{pq},\\& &
\label{Reduced Ric tensora}
(\Ric_{g} (g))_{pq}=\Ric_{pq} (g)-\frac 12 (g_{pn} \hat \nabla _q\hat F^n+ g_{qn} \hat \nabla _p\hat F^n)
\eeq
where $\hat F^n$ are the harmonicity functions given by
\beq\label{Harmonicity condition CCC}
\hat F^n=\Gamma^n-\hat \Gamma^n,\quad\hbox{where }
\Gamma^n=g^{jk}\Gamma^n_{jk},\quad
\hat \Gamma^n=g^{jk}\hat \Gamma^n_{jk},
\eeq
where $\Gamma^n_{jk}$ and $\hat \Gamma^n_{jk}$ are the Christoffel symbols
for $g$ and $g$, correspondingly.
Note that $\hat \Gamma^n$ depends also on $g^{jk}$.
As $\Gamma^n_{jk}-\hat \Gamma^n_{jk}$ is the difference of two connection
coefficients, it is a tensor. Thus $\hat F^n$ is tensor (actually, a vector field), implying that both
$(\Ric_{g} (g))_{jk}$ and $(\Ein_{g} (g))_{jk}$ are  2-covariant tensors.
A direct calculation shows that  the $g$-reduced Einstein tensor is the sum of the harmonic Einstein tensor 
and a term that is a zeroth order in $g$, 
\beq\label{hat g reduced einstein and reduced einstein}
(\Ein_{g} (g))_{\mu\nu}=\Ein^{(h)}_{\mu\nu} (g)+\frac 12 (g_{\mu q}\frac{\p \hat \Gamma^q}{\p x^{\nu}}+g_{\nu q}\frac{\p \hat \Gamma^q}{\p x^{\mu}}).
\eeq
We also use the wave operator  
\beq\label{eq: wave ope}\\
\nonumber
\square_g\phi =
\sum_{p,q=1}^4(-\det(g(x)))^{-\frac 12}\frac  \p{\p x^p}
\left ((-\det(g(x)))^{\frac12}
g^{pq}(x)\frac \p{\p x^q}\phi(x)\right),
\eeq
can be written  for 
as 
\beq\label{eq: wave ope BB}\square_g \phi =g^{jk}\p_j\p_k \phi-g^{pq} \Gamma^n_{pq} \p_n \phi=
g^{jk}\p_j\p_k \phi- \Gamma^n \p_n \phi.
\eeq

%

%

%

\subsection*{A.3.\ Wave maps and reduced Einstein equations}
Let us consider the manifold $M_1=(-\infty,t_1)\times N$ with a $C^m $-smooth metric $g^{\prime}$,
$m\geq 8$,
 which
is a perturbation of the metric $g$ and satisfies the Einstein equation
\beq\label{Einstein on M_0}
\Ein(g^{\prime})=T^{\prime}\quad \hbox{on }M_1,
\eeq
or equivalently, 
\ba
\Ric(g^{\prime})=\rho^{\prime},\quad \rho^{\prime}_{jk}=T_{jk}^{\prime}-\frac 12 ((g^{\prime})^{  nm}T^{\prime}_{nm})g^{\prime}_{jk}\quad \hbox{on }M_1.
\ea
Assume also that  $g^{\prime}=g$ in the domain $A$,
where $A=M_1\setminus\K_1$ and  $\|g^{\prime}-g\|_{C^2_b(M_1,g^+)}<\e_0$,
so that $(M_1,g^{\prime})$ is globally hyperbolic. Note that then $T^{\prime}=\hat T$ in the set $A$ 
and that the metric $g^{\prime}$ coincides with $g$
in particular in the set $M^-=\R_-\times N$

We recall next the considerations of \cite{ChBook}.
Let us  consider the Cauchy problem for the wave map
$f:(M_1,g^\prime)\to (M,g)$, namely
\beq
\label{C-problem 1 AAA}& &\square_{g^{\prime},g} f=0\quad\hbox{in } M_1,\\
\label{C-problem 2 AAA}& &f=Id,\quad \hbox{in  }\R_-\times N,
\eeq
where $M_1=(-\infty,t_1)\times N\subset M$. In (\ref{C-problem 1 AAA}), 
$\square_{g^{\prime},g} f=g^\prime\,\cdotp\hat \nabla^2 f$ is the wave map operator, where $\hat \nabla$
is the covariant derivative of a map $(M_1,g^\prime)\to (M,g)$, see \cite[Ch. VI, formula (7.32)]{ChBook}.
In  local coordinates
$X:V\to \R^4$ of $V\subset M_1$, denoted by
$X(z)=(x^j(z))_{j=1}^4$ and $Y:W\to \R^4$ of $W\subset \hattuM $, denoted by
$Y(z)=(y^A(z))_{A=1}^4$, 
the wave map $f:M_1\to \hattuM $ has the representation $Y(f(X^{-1}(x)))=(f^A(x))_{A=1}^4$ and
the wave map operator in equation (\ref{C-problem 1 AAA}) is given by
\beq\label{wave maps2}
& &(\square_{g^{\prime},g} f)^A(x)= 
(g^{\prime})^{  jk}(x)\bigg(\frac \p{\p x^j}\frac \p{\p x^k} f^A(x) -\Gamma^{{\prime} n}_{jk}(x)\frac \p{\p x^n}f^A(x)
\\
& &\quad \quad\quad\quad\quad\quad\quad\quad\nonumber
+\hat \Gamma^A_{BC}(f(x))
\,\frac \p{\p x^j}f^B(x)\,\frac \p{\p x^k}f^C(x)\bigg)\eeq
where $\hat \Gamma^A_{BC}$ denotes the Christoffel symbols of metric $g$
and  $\Gamma^{ {\prime} j}_{kl}$ are the Christoffel symbols of metric $g^{\prime}$.
When (\ref{C-problem 1 AAA})
is satisfied, we say that  $f$ is a wave map with respect to the pair $(g^{\prime},g)$.

\extension{It follows from \cite[App.\ III, Thm. 4.2 and sec.\ 4.2.2] {ChBook},
that if $g^\prime \in C^m(M_0)$, $m\geq 5$ 
is sufficiently close to $g$
in $C^m(M_0)$, then (\ref{C-problem 1 pre})-(\ref{C-problem 2 pre}) has a unique 
solution $f\in C^0([0,t_1];H^{m-1}(N))\cap  C^1([0,t_1];H^{m-2}(N))$.
This comes from the fact that the Christoffel symbols of $g^\prime$
are in $ C^{m-1}(M_0)$.
Moreover, when $m$ is even, using 
\cite[Thm.\ 7]{Kato1975} 
for $f$ and $\p_t^p f$,  we see that
the solution $f$ is in $f\in \cap_{p=0}^{m-1}C^p([0,t_1];H^{m-1-p}(N))\subset C^{m-3}(M_0)$ and
$f$ depends in $C^{m-3}(M_0)$ continuously on $g^\prime\in C^m(M_0)$.
We note that these smoothness results for $f$ are not optimal.}  

The wave map operator
$\square_{g^{\prime},g} $ is a coordinate invariant operator.
The important property of the wave maps is that, if 
$f$ is wave map with respect to the pair $(g^{\prime},g)$ and
$g=f_*g^{\prime}$ 
then, as follows from (\ref{wave maps2}), the identity map $Id:x\mapsto x$ is a wave map with respect to the pair $(g,g)$
and, the wave map equation for the identity map is equivalent to (cf.\ \cite[p.\ 162]{ChBook})
\beq\label{wave equation id}
\Gamma^n=\hat\Gamma^n,\quad \hbox{where }\Gamma^n=g^{jk}\Gamma^n_{jk},
\quad \hat \Gamma^n=g^{jk}\hat\Gamma^n_{jk}
\eeq
%
where the Christoffel symbols $\hat\Gamma^n_{jk}$ of the metric $g$ are smooth functions.

Since  $g=g^{\prime}$ outside a compact set $\K_1\subset (0,t_1)\times N$,
we see that   this Cauchy problem is equivalent to the same equation
restricted to the
set $(-\infty,t_1)\times B_0$, where $B_0\subset N$ is an open
relatively compact set such that $\K_1\subset (0,t_1]\times B_0$
with the boundary condition $f=Id$ on  $(0,t_1]\times \p B_0$.
%
 Moreover, by the uniqueness of the wave map, we have
$f|_{M_1\setminus \K_1}=id$ so that
 $f(\K_1)\cap M_0\subset \K_0$.

 As 
 the inverse function of the wave map $f$ depends continuously,
 in
$C^{m-3}_b([0,t_1]\times N,g^+)$, on the metric $g^{\prime} \in C^m(M_0)$
we can also assume that $\e_1$ is  so small that
$\hattuM _0\subset f(M_1)$.

Denote next $g:=f_*g^{\prime}$, $T:=f_*T^{\prime}$, and $\rho:=f_*\rho^{\prime}$ 
and define $\hat \rho=\hat T-\frac 12 (\hbox{Tr}\, \hat T)g.$ 
Then $g$ is $C^{m-6}$-smooth and the equation (\ref{Einstein on M_0}) implies 
\beq\label{Einstein on M1}
\Ein(g)=T\quad \hbox{on }\hattuM _0.
\eeq 
Since $f$ is a wave map and $g=f_*g^\prime$, we have that
the identity map is a $(g,g)$-wave map and thus $g$ satisfies (\ref{wave equation id}) and thus  by the
definition of the reduced Einstein tensor,  (\ref{Reduced Ric tensor})-(\ref{Reduced Einstein tensor}), we have
\ba
\Ein_{pq}(g)=
(\Ein_{g} (g))_{pq}\quad \hbox{on } \hattuM _0.
\ea
This and  (\ref{Einstein on M1}) yield the $g$-reduced Einstein equation
\beq\label{hat g reduced einstein equations}
(\Ein_{g} (g))_{pq}=T_{pq}\quad \hbox{on } \hattuM _0.
\eeq
This equation is useful for our considerations as it is a quasilinear,
hyperbolic equation on $\hattuM _0$. Recall that 
 $g$ coincides with
$g$ in $M_0\setminus \K_0$. The unique solvability of this 
Cauchy problem is studied in e.g.\ \cite[Thm.\ 4.6 and 4.13]{ ChBook}, \cite {HKM} 
and  Appendix B below.


\subsection*{A.4.\ Relation of the reduced Einstein equations and the original Einstein equation}

The metric $g$ which solves the $g$-reduced Einstein
equation $\Ein_{g} (g)=T$ is a solution of the original
Einstein equations  $\Ein (g)=T$ if the harmonicity
functions $\hat F^n$
vanish identically. Next we recall the result that 
the harmonicity functions vanish on $\hattuM _0$
when 
\beq\label{eq: good system}
& &(\Ein_{g}(g))_{jk}=T_{jk},\quad \hbox{on }M_0,\\
\nonumber & &\nabla_pT^{pq}=0,\quad \hbox{on }M_0,\\
\nonumber & &g=g,\quad \hbox{on }M_0\setminus \K_0.
\eeq
To see this, let
us  denote $\Ein_{jk}(g)=S_{jk}$, $S^{jk}=g^{jn}g^{km}S_{nm}$,
and $T^{jk}=g^{jn}g^{km}T_{nm}$. 
Following standard arguments,
see \cite{ ChBook}, we see from (\ref{Reduced Einstein tensor}) that in local coordinates
\ba
S_{jk}-(\Ein_{g}(g))_{jk}=\frac12(g_{jn}\hat \nabla_k \hat F^n+g_{kn}\hat \nabla_j\hat F^n-g_{jk}\hat \nabla_n \hat  F^n).
\ea
Using equations (\ref{eq: good system}), the Bianchi identity  $\nabla_pS^{pq}=0$,
and the basic property of Lorentzian connection,
$\nabla_kg^{nm}=0$,
we obtain 
\ba
0&=&
2\nabla_p(S^{pq}-T^{pq})\\
&=&\nabla_p(g^{qk}\hat  \nabla_k F^p+g^{pm}\hat\nabla_m \hat F^q-  g^{pq}\hat  \nabla_n \hat F^n)
\\
&=&g^{pm}\nabla_p \hat \nabla_m \hat F^q+(g^{qp} \nabla_n \hat \nabla_p \hat F^n-  g^{qp} \nabla_p \hat \nabla_n \hat F^n)\\
&=&g^{pm}\nabla_p \hat \nabla_m \hat F^q+W^q(\hat F)
\ea
where $\hat F=(\hat F^q)_{q=1}^4$ and
the operator 
\ba
W:(\hat F^q)_{q=1}^4\mapsto (g^{qk} (\nabla_p \hat \nabla_k \hat F^p- \nabla_k \hat \nabla_p \hat F^p))_{q=1}^4
\ea
is a linear first order differential operator which coefficients are polynomial functions of
$g_{jk}$, $g^{jk}$, $g_{jk}$, $g^{jk}$ and their first derivatives.

Thus the harmonicity functions $\hat F^q$ satisfy on $\hattuM _0$ the hyperbolic initial 
value problem 
\ba
& &g^{pm}\nabla_p \hat \nabla_m \hat F^q+W^q(\hat F)=0,\quad\hbox{on }M_0,\\
& & \hat F^q=0,\quad \hbox{on }M_0\setminus \K_0,
\ea
and as this initial 
Cauchy problem is uniquely solvable by \cite[Thm.\ 4.6 and 4.13]{ ChBook}
or \cite {HKM}, we see that $\hat F^q=0$ on $\hattuM _0$. Thus
equations (\ref{eq: good system})
yield
that the Einstein equations  $\Ein(g)=T$ hold on $\hattuM _0$.

We note that in the $(g,g)$-wave map coordinates, where $\hat F^q=0$, the
wave operator (\ref{eq: wave ope BB}) has the form
\beq\label{eq: wave operator in wave gauge}
\square_g \phi=g^{jk}\p_j\p_k \phi-g^{pq}\hat \Gamma^n_{pq} \p_n \phi.
\eeq
Thus, the scalar field equation $\square_g \phi-m^2\phi=0$
does not involve derivatives of $g$.

 \subsection*{A.5.\ Linearized Einstein-scalar field equations} 
 
 Next we consider  the linearized equation are obtained as the derivatives
 of the solutions of the non-linear, generalized Einstein-matter field
equations 
(\ref{eq: adaptive model with no source}). 

 
{Observe that if  a family $\mathcal F_\e=(\mathcal F^1_\e,\mathcal F^2_\e)$
of sources and a family  $(g_\e,\phi_\e)$ of functions 
are solutions of  the non-linear reduced Einstein-scalar field  
equations 
(\ref{eq: adaptive model with no source})
that depend smoothly on $\e\in [0,\e_0)$ in $C^{15}(M_0)$
and satisfy $\mathcal F_\e|_{\e=0}=0$, $(g_\e,\phi_\e)|_{\e=0}=(g,\hat \phi)$,
then $\p_\e  \mathcal F_\e|_{\e=0}=f=(f^1,f^2)$, and  $\p_\e  (g_\e,\phi_\e)|_{\e=0}=(\dot g,\dot \phi)$,
satisfy
the linearized version of 
the equation (\ref{eq: adaptive model with no source}) that has the form (in local coordiantes),
\beq\label{linearized eq: adaptive model with no source BB}
& &\square_{g} \dot g_{jk}+A_{jk}(\dot g,\dot \phi,\p \dot g,\p \dot \phi)=f^1_{jk},\quad
\hbox{in }M_0,\\
\nonumber
& &\square_{g}\dot \phi_\ell +
B_\ell(\dot g,\dot \phi,\p \dot g,\p \dot \phi)
=f^2_\ell,\quad \ell=1,2,3,\dots,L.
\eeq
Here $A_{jk}$ and $B_\ell$ are first order linear differential
operators whose coefficients depend on $g$ and $\hat \phi$.
Let us write these equations in more explicit form. We see that
the linearized reduced Einstein tensor is in local coordinates of the form 
\ba
e_{pq}(\dot g)&:=&\p_\e (\Ein_{g} g_\e)_{pq}|_{\e =0}\\
&=&
-\frac 12 g^{jk}\hat \nabla_j\hat \nabla_k  \dot g_{pq}
+\frac 14(  g^{nm} g^{jk}\hat \nabla_j\hat \nabla_k  \dot g_{nm}) g_{pq}
\\
& &+A_{pq}^{abn}\hat \nabla_n\dot g_{ab}+B_{pq}^{ab}\dot g_{ab},
\ea
where $A_{pq}^{abn}(x)$ and $B_{pq}^{ab}(x)$ depend on  $g_{jk}$
and its derivatives (these terms can be computed explicitly using  (\ref{q-formula2copyBa})).
The linearized scalar field stress-energy tensor is the linear first order differential operator 
\ba
& &t^{(1)}_{pq}(\dot g)+t^{(2)}_{pq}(\dot \phi):=\p_\e ({\bf T}_{jk}(g_\e,\phi_\e)
)
|_{\e =0}\\
&=&
\sum_{\ell=1}^L\bigg(
\p_j\hat \phi_\ell \,\p_k\dot \phi_\ell +
\p_j\dot \phi_\ell \,\p_k\hat  \phi_\ell \\
& &\quad
-\frac 12 \dot g_{jk} \hat  g^{pq}\p_p \hat  \phi_\ell \,\p_q \hat \phi_\ell
-\frac 12 g_{jk} \dot g^{pq}\p_p \hat  \phi_\ell \,\p_q \hat \phi_\ell\\
& &\quad
-\frac 12 g_{jk} g^{pq}\p_p\dot \phi_\ell \,\p_q \hat \phi_\ell
-\frac 12 g_{jk} \hat  g^{pq}\p_p \hat  \phi_\ell \,\p_q\dot  \phi_\ell\\
& &\quad
- m^2\hat  \phi_\ell \dot \phi_\ellg_{jk}
-\frac 12 m^2\hat \phi_\ell^2\dot g_{jk}\bigg).
\ea
Thus when $\mathcal F_\e=(\mathcal F^1_\e,\mathcal F^2_\e)$
is a family 
of sources and   $(g_\e,\phi_\e)$ a family of functions 
that satisfy  the non-linear reduced Einstein-scalar field  
equations 
(\ref{eq: adaptive model with no source}),
the  $\e$-derivatives $\dot u=(\dot g,\dot \phi)$ 
and $\p_\e  \mathcal F_\e|_{\e=0}=f=(f^1,f^2)$
satisfy  in local $(g,g)$-wave coordinates
\beq\label{eq: final linearization}
& &
e_{pq}(\dot g)-t^{(1)}_{pq}(\dot g)-t^{(2)}_{pq}(\dot \phi)=f^1_{pq},\\
\nonumber
& &\square_{g}\dot \phi_\ell
-g^{nm}g^{kj}(\p_n\p_j\hat \phi_\ell+ \hat \Gamma_{nj}^p\p_p
\hat \phi_\ell)\dot g_{mk}-m^2\phi_\ell=f^2_\ell.
\eeq
We call this the linearized Einstein-scalar field equation. 
 }

\subsection*{A.6.\ Linearization of the conservation law.}

Assume that $u=(g,\phi)$ and $\F=(\F^1,\F^2)$  
satisfy equation (\ref{eq: adaptive model with no source}).
Then
  the conservation law (\ref {conservation law0})
 gives for all $j=1,2,3,4$ equations (see \cite[Sect. 6.4.1]{ChBook}) 
 \ba
0&=&\frac 12\nabla_p^g (g^{pk}T_{jk})\\
&=&\frac 12\nabla_p^g (g^{pk}({\bf T}_{jk}(g,\phi)+\F^1_{jk}))\\
&=& \sum_{\ell=1}^L(g^{pk} \nabla^g_p\p_k\phi_\ell )\,\p_j\phi_\ell 
-(m^2_\ell \phi_\ell \p_p \phi_\ell) \delta_j^p +\frac 12 \nabla^g_p (g^{pk}\F^1_{jk})
\\
&=& \sum_{\ell=1}^L(g^{pk} \nabla^g_p\p_k\phi_\ell -m^2_\ell \phi_\ell ) \,\p_j \phi_\ell +\frac 12 \nabla^g_p (g^{pk}\F^1_{jk}).
\ea
This yields by (\ref{eq: adaptive model with no source})
\beq
\label{mod. conservation law aaa}
\frac 12 g^{pk} \nabla^g_p \F^1_{jk}+\sum_{\ell=1}^L \F^2_\ell \,\p_j\phi_\ell=0,
\quad j=1,2,3,4.
\eeq
Next assume that $u_\e=(g_\e,\phi_\e)$ and $\F_\e$
satisfy equation (\ref{eq: adaptive model with no source}) and
$C^1$-smooth functions of $\e\in (-\e_0,\e_0)$ taking values in $H^1(N)$-tensor fields, and $(g_\e,\phi_\e)|_{\e=0}=(g,\hat \phi)$
and  $\F_\e|_{\e=0}=0$. 
 Denote $(f^1,f^2)=\p_\e\F_\e|_{\e=0}$. 
 Then by taking $\e$-derivative of (\ref{mod. conservation law aaa}) at $\e=0$
 we get
\beq
\label{linearized conservation law aaa}
\frac 12 g^{pk} \hat \nabla_p f^1_{jk}+\sum_{\ell=1}^L f^2_\ell \,\p_j\hat \phi_\ell=0,
\quad j=1,2,3,4.
\eeq
We call this the linearized conservation law.

\section*{Appendix B: Stability and existence of the  direct problem}

 Let us next consider existence and stability of the solutions of the
 Einstein-scalar field equations.
Let $t={\bf t}(x)$ be local time so that there is a diffeomorphism
$\Psi:M\to \R\times N$, $\Psi(x)=({\bf t}(x),Y(x))$, and
 $$S(T)=\{x\in \hattuM _0; \ {\bf t}(x)=T\},\quad T\in \R$$ are Cauchy surfaces.
 Let $t_0>0$. Next we identify $M$ and $\R\times N$ via the map $\Psi$ and just
 denote $M=\R\times N$.
Let us denote $$M(t)=(-\infty,t)\times N,\quad M_0=M(t_0).$$
By \cite[Cor.\ A.5.4]{BGP} the set $\K=J^+_{\tilde g}(\hat p^-)\cap \hattuM _0$,
where $\hat p^0\in M_0$,
is compact.
Let $N_1,N_2\subset N$ be such
 open relatively compact sets with smooth boundary so that
 $N_1\subset N_2$ and  $Y(J_{\tilde g}^+(\hat p^0)\cap \hattuM _0)\subset N_1$.  
 
 To simplify citations to the existing literature, 
let us define $\tilde N$ to be a compact manifold without boundary such
that $N_2$ can be considered as a subset of $\tilde N$. Using a construction based on
a suitable partition of unity, the Hopf double of the manifold $N_2$, 
and the Seeley extension of the metric tensor,  we can
endow $\tilde M=(-\infty,t_0)\times \tilde N$ with a smooth Lorentzian
metric $g^e$ (index $e$ is for "extended")
so that $\{t\}\times \tilde N$ are Cauchy surfaces of $\tilde N$
and that $g$ and $g^e$ coincide in  the set $\Psi^{-1}((-\infty,t_0)\times N_1)$ that contains the set $J^+_{g}(\hat p^0)\cap \hattuM _0$. We extend the metric $\tilde g$ to a (possibly non-smooth) globally
hyperbolic metric $\tilde g^e$ on $\tilde M_0=(-\infty,t_0)\times \tilde N)$
such that $g^e<\tilde g^e$.

 To simplify notations below we denote $g^e=g$ and
 $\tilde g^e=\tilde g$ 
 on the whole
 $\tilde M_0$.
 Our aim is to prove the estimate (\ref{eq: Lip estim}).

Let us denote by $t={\bf t}(x)$ the local time. 
Recall that when $(g,\phi)$ is a solution of 
the scalar field-Einstein equation, we denote $u=(g-g,\phi-\hat \phi)$.
We will consider the equation for $u$, and to emphasize that the metric
depends on $u$, we denote $g=g(u)$ and assume below that  both the metric $g$ is  dominated
by $\tilde g$, that is, $ g<\tilde g$.
We use the pairs $${\bf u}(t)=(u(t,\cdotp),\p_t u(t,\cdotp))\in
H^1(\tilde N)\times L^2(\tilde N)$$ the notations
 ${\bf v}(t)=(v(t),\p_t v(t))$ etc. 
We consider sources $\F$ and $F$, cf. Appendix C.
Let us consider a generalization of the system (\ref{eq: notation for hyperbolic system 1}) of the form 
\beq\label{eq: notation for hyperbolic system}
& &\square_{g(u)}u+V(x,D)u+H(u,\p u) =R(x,u,\p u,F,\p F)+\F,\ \ x\in \tilde M_0,\hspace{-1cm}\\
& & \nonumber \supp(u)\subset \K,
\eeq
where $\square_{g(u)}$ is the Lorentzian Laplace operator operating on the sections 
of the bundle $\B^L$ on $M_0$ and 
$\supp(F)\cup \supp(\F)\subset \K$.
Note that above  $u=(g-g,\phi-\hat \phi)$ and $g(u)=g$.
Also, $F\mapsto R(x,u,\p u,F,\p F)$ is a 
 non-linear first order differential operator where $R$ is a smooth function
 that depends on all variables $x$,
 $u(x)$, $\p_j u(x)$, $F(x),$ and $\p_jF(x)$ smoothly and also on the  derivatives 
of $(g,\hat \phi)$ at $x$. Let
\ba 
V(x,D)=V^j(x)\p_j+V(x)
\ea 
be a linear first order differential operator whose coefficients 
at $x$ are depending smoothly on the  derivatives 
of $(g,\hat \phi)$  at $x$, and finally, $H(u,\p u)$ is a polynomial
of $u(x)$ and $\p_j u(x)$ which coefficients 
at $x$ are depending smoothly on the  derivatives 
of $(g,\hat \phi)$   such that
 $$\p^\a_v\p^\beta_w H(v,w)|_{v=0,w=0}=0\hbox { for }
|\a|+|\b|\leq 1.$$ 
We note that also the non-linear scalar equation (\ref{eq: wave-eq general}) could 
be analyzed in the same framework  by doing simplifications in the above, that is, by 
defining $g(u)$ to be independent of $g$ and $H(u,\p u)$ to be the term $au^2$.

By \cite[Lemma 9.7]{Ringstrom}, the equation  (\ref{eq: notation for hyperbolic system})
has at most one solution with given  $C^2$-smooth source functions $F$ and $\F$.
Next we consider the existence of $u$ and its dependency on $F$ and $\F$.


%

Below we use notations, c.f. (\ref{eq: notation for hyperbolic system 1}) and
(\ref{eq: notation for hyperbolic system 1b}) and Appendix C, \ba
{\mathcal R}({\bf u},F)=R(x,u(x),\p u(x),F(x),\p F(x)),\quad
{\mathcal H}({\bf u})=H(u(x),\p u(x)).\ea
Note that $u=0$, i.e., $g=g$ and $\phi=\hat\phi$ satisfies (\ref{eq: notation for hyperbolic system}) with $F=0$
and $\F=0$.
Let us use the same notations as in  \cite {HKM} cf. also \cite[section 16]{Kato1975}, 
 to consider quasilinear wave equation on $[0,t_0]\times \tilde N$.
Let $\H^{(s)}(\tilde N)=H^s(\tilde N)\times H^{s-1}(\tilde N)$ and 
\ba
Z=\H^{(1)}(\tilde N),\quad  Y=\H^{(k+1)}(\tilde N),\quad X=\H^{(k)}(\tilde N).
\ea
The norms on these space are defined invariantly using the smooth Riemannian 
metric $h=g|_{\{0\}\times \tilde N}$ on $\tilde N$.
Note that $\H^{(s)}(\tilde N)$ are in fact the Sobolev spaces of sections on the bundle $\pi:\B_K\to
\tilde N$, where $\B_K$ denotes also the  pull back bundle of $\B_K$ on $\tilde M$ in the map
$id:\{0\}\times \tilde N\to \tilde M_0$,
or on the bundle $\pi:\B_L\to \tilde N$. Below, $\nabla_h$ denotes the  standard connection of the bundle  $\B_K$
or  $\B_L$ associated to the metric $h$.

Let $k\geq 4$ be an even integer.
By definition of $H$ and $R$ we see that  there are  $0<r_0<1$ and $L_1,L_2>0$,
all depending on $g$, $\hat \phi$, $\K$, and $t_0$, such that if 
$0<r\leq r_0$ and
\beq\label{new basic conditions}
& &\|{\bf v} \|_{C([0,t_0];\H^{(k+1)}(\tilde N))}\leq r,\quad \|{\bf v}^\prime\|_{C([0,t_0];\H^{(k+1)}(\tilde N))}\leq r,\\
& &\|F\|_{C([0,t_0];H^{(k+1)}(\tilde N))}\leq r^2,\quad \|\F\|_{C([0,t_0];H^{(k)}(\tilde N))}\leq r^2 \nonumber\\
& &\|F^\prime\|_{C([0,t_0];H^{(k+1)}(\tilde N))}\leq r^2,\quad \|\F^\prime\|_{C([0,t_0];H^{(k)}(\tilde N))}\leq r^2 \nonumber
\eeq
then
\beq\label{new f1f2 conditions}
& &\quad \quad \|g(\cdotp;v)^{-1}\|_{C([0,t_0];H^s(\tilde N))}\leq L_1 ,\\ \nonumber
& &\|{\mathcal H} ({\bf v})\|_{C([0,t_0];H^{s-1}(\tilde N))}\leq L_2r^2,
\quad \|{\mathcal H} ({\bf v}^\prime)\|_{C([0,t_0];H^{s-1}(\tilde N))}\leq L_2r^2,\\ \nonumber
& &\|{\mathcal H} ({\bf v})-{\mathcal H} ({\bf v}^\prime)\|_{C([0,t_0];H^{s-1}(\tilde N))}\leq  L_2r\, \|{\bf v}-{\bf v}^\prime\|_{C([0,t_0];\H^{(s)}(\tilde N))}, 
\\ \nonumber
& &\|\mathcal R ({\bf v}^\prime,F^\prime)\|_{C([0,t_0];H^{s-1}(\tilde N))}\leq L_2r^2,
\quad \|\mathcal R ({\bf v},F)\|_{C([0,t_0];H^{s-1}(\tilde N))}\leq L_2r^2,\hspace{-1cm} \\ \nonumber
& &\|\mathcal R ({\bf v},F)-\mathcal R ({\bf v}^\prime, F^\prime)\|_{C([0,t_0];H^{s-1}(\tilde N))}
\\ \nonumber & &\leq  
L_2r\, \|{\bf v}-{\bf v}^\prime\|_{C([0,t_0];\H^{(s)}(\tilde N))}+L_2\|F-{F}^\prime\|_{C([0,t_0];H^{s+1}(\tilde N))\cap C^1([0,t_0];H^{s}(\tilde N))},\hspace{-2cm}
%
  \eeq
 for all $s\in [1,k+1]$. 
%

%
%

Next we write (\ref{eq: notation for hyperbolic system}) as a first order
system.  To this end,
let $\mathcal A(t,{\bf v}):\H^{(s)}(\tilde N)\to \H^{(s-1)}(\tilde N)$ be the operator 
$$\mathcal A(t,{\bf v})=\mathcal A_0(t,{\bf v})+\mathcal A_1(t,{\bf v})$$
where in local coordinates and in local trivialization of the bundle $\B^L$
\ba
\mathcal A_0(t,{\bf v})=-\left(\begin{array}{cc} 0 & I\\
 \frac 1{g^{00}(v)}\sum_{j,k=1}^3 g^{jk}(v)\frac \p{\p x^j}\frac\p{ \p x^k}\quad &
 \frac 1{g^{00}(v)}\sum_{m=1}^3 g^{0m}(v)\frac {\p}{\p x^m} \end{array}\right)
\ea 
with $g^{jk}(v)=g^{jk}(t,\cdotp;v)$ is a function on $\tilde N$ and 
\ba
\mathcal A_1(t,{\bf v})= \frac {-1}{g^{00}(v)}\left(\begin{array}{cc} 0 & 0\\
\sum_{j=1}^3 
 B^j(v)
   \frac \p{\p x^j}\quad &
B^{0}(v)\end{array}\right)
\ea 
where $B^j(v)$ depend on $v(t,x)$ and its first derivatives, and 
the connection coefficients (the Christoffel  symbols) corresponding  to $g(v)$. We denote $\S=(F,\F)$ and 
 \ba
& &f_\S(t,{\bf v})=(f_\S^1(t,{\bf v}),f_\S^2(t,{\bf v}))\in \H^{(k)}(\tilde N),
\quad\hbox{where}\\
& & f_\S^1(t,{\bf v})=0,\quad f_\S^2(t,{\bf v})=\mathcal R({\bf v},F)(t,\cdotp)-\mathcal H({\bf v})(t,\cdotp)+\F(t,\cdotp).
 \ea
 
 Note that when (\ref{new basic conditions}) are satisfied with $r<r_0$, inequalities (\ref{new f1f2 conditions}) imply that there exists $C_2>0$ so that
 \beq\label{fF bound r2}
 & &\|f_\S(t,{\bf v})\|_Y+\|f_{\S^\prime}(t,{\bf v}^\prime)\|_Y\leq C_2r^2,\\
 & &\|f_\S(t,{\bf v})-f_\S(t,{\bf v}^\prime)\|_Y\leq C_2
 r\,\|{\bf v}-{\bf v}^\prime\|_{C([0,t_0];Y)}. \nonumber
 \eeq

%
%

Let $U^{\bf v}(t,s)$ be the wave propagator corresponding to metric $g(v)$, that is,
$U^{\bf v}(t,s): {\bf h}\mapsto {\bf w}$, where ${\bf w}(t)=(w(t),\p_t w(t))$  solves
\ba
(\square_{g(v)}+V(x,D))w=0\quad\hbox{for } (t,y)\in [s,t_0]\times \tilde N,\quad\hbox{with }
 {\bf w}(s,y)={\bf h}.
 \ea 
 Let $$S=(\nabla_{h}^*\nabla_{h}+1)^{k/2}:Y\to Z$$ be an isomorphism.
As  $k$ is an even integer, we see using multiplication estimates for Sobolev
spaces, see e.g.\
  \cite [Sec.\ 3.2, point (2)]{HKM}, that there exists $c_1>0$ (depending on
  $r_0,L_1,$ and $L^2$) so that  $$\A(t,{\bf v})S-S\A(t,{\bf v})=C(t,{\bf v}),$$
 where $\|C(t,{\bf v})\|_{Y\to Z}\leq c_1$ for all ${\bf v}$ satisfying
 (\ref{new basic conditions}). This yields that 
 the property (A2) in  \cite {HKM} holds, namely
 that $$S\A(t,{\bf v})S^{-1}=\A(t,{\bf v})+B(t,{\bf v})$$ where
$B(t,{\bf v})$ extends to a bounded operator in $Z$
 for which $$\|B(t,{\bf v})\|_{Z\to Z}\leq c_1$$ for all ${\bf v}$ satisfying
 (\ref{new basic conditions}).
Alternatively, to see the mapping properties of $B(t,{\bf v})$ we could use the fact that 
 $B(t,{\bf v})$ is a zeroth order pseudodifferential operator with
$H^{k}$-symbol.\hiddenfootnote{On mapping properties of such
pseudodifferential operators, see J. Marschall, Pseudodifferential operators with coefficients in Sobolev spaces. Trans. Amer. Math. Soc. 307 (1988), no. 1, 335-361.} 
 
 Thus the proof of  \cite [Lemma 2.6]{HKM} shows\hiddenfootnote{ Alternatively, as $U^{\bf v}(t,s)$ are propagators
  for linear wave equations having finite speed of wave propagation, one can prove 
the  estimates (\ref{U bounds}) by considering the wave equation
in local coordinate neighborhoods $W_j\subset W_j^\prime \subset \tilde N$,
$\Phi_j:  W_j^\prime\to \R^3$ and a partition of unity $\psi_j\in C^\infty_0(W_j)$ and
cut-off functions $\psi^\prime_j\in C^\infty_0(W_j^\prime)$.
Then  \cite [Lemma 2.6]{HKM}, used  in local coordinates, shows
that when $t_k-t_{k-1}$ is small enough then 
\ba 
\| \psi^\prime_j U^{\bf v}(t_k,t_{k-1})(\psi_j {\bf w}) \|_{Y}\leq C_{3,jk}^{\prime}e^{C_4(t-s)}\| \psi_j {\bf w} \|_{Y}.
\ea
Using  sufficiently small time steps $t_k-t_{k-1}$ and combining the estimates
for different $j$:s
together, one obtain estimates (\ref{U bounds}).}
 that 
%
there is a constant  $C_3>0$
   so that 
  \beq
  \label{U bounds}\|U^{\bf v}(t,s)\|_{Z\to Z}\leq C_3\quad\hbox{and}\quad
  \|U^{\bf v}(t,s) \|_{Y\to Y}\leq C_3\hspace{-1.5cm}
  \eeq 
  for $0\leq s<t\leq t_0$.
    By interpolation of estimates (\ref{U bounds}), we see also that 
    \beq
  \label{U bounds2}
  \|U^{\bf v}(t,s)\|_{X\to X}\leq C_3,
  \eeq  for $0\leq s<t\leq t_0$.

Let us next modify the reasoning given in \cite{Kato1975}: let  $r_1\in (0,r_0)$
be a parameter which value will be chosen later,  $C_1>0$
and 
$E$ be the space of functions ${\bf u}\in C([0,t_0];X)$ for which
\beq
\label{1 bounded}& &\|{\bf u}(t)\|_Y\leq r_1\quad\hbox{and}\\
\label{2 Lip}& &\|{\bf u}(t_1)-{\bf u}(t_2)\|_X\leq C_1|t_1-t_2|\eeq
for all $t,t_1,t_2\in [0,t_0]$.
The set $E$ is endowed by the metric of  $C([0,t_0];X)$. We 
note that  by \cite[ Lemma 7.3] 
{Kato1975},
a convex $Y$-bounded, $Y$-closed set is closed also in $X$.
Similarly, functions $G:[0,t_0]\to X$  satisfying (\ref{2 Lip}) 
form a closed subspace of $C([0,t_0];X)$. Thus $E\subset X$ is a closed set implying
that $E$ is a complete metric space.

Let
 \ba 
 & &\hspace{-1cm}W=\{(F,\F)\in C([0,t_0];H^{k+1}(\tilde N))
 \times C([0,t_0];H^{k}(\tilde N));
 \\
 & &\quad\quad
 \ \sup_{t\in [0,t_0]}\|F(t)\|_{H^{k+1}(\tilde N)}+\|\F(t)\|_{H^{k}(\tilde N)}< r_1\}.\hspace{-1cm}
 \ea

Following \cite[p. 44]{Kato1975}, we see that
the solution of
equation (\ref{eq: notation for hyperbolic system}) with the source $\S\in W$  is found as
a fixed point, if it exists, of the map $\Phi_\S:E\to C([0,t_0];Y)$ where 
$\Phi_\S({\bf v})={\bf u}$ is  given by
$$
{\bf u}(t)=\int_0^t U^{\bf v}(t,\tilde t)f_\S(\tilde t,{\bf v})\,d\tilde t,\quad 0\leq t\leq t_0.
$$ 
Below, we denote  ${\bf u}^{\bf v}=\Phi_\S({\bf v})$.

As  $\Phi_{\S_0}(0)=0$ where
$\S_0=(0,0)$, we see using 
the above and the inequality $\|\,\cdotp\|_X\leq \|\,\cdotp\|_Y$ that the
 function ${\bf u}^{\bf v}$ satifies
  \ba
& &  \|{\bf u}^{\bf v}\|_{C([0,t_0];Y)}\leq C_3C_2t_0 r_1^2,\\
& &  \|{\bf u}^{\bf v}(t_2)-{\bf u}^{\bf v}(t_1)\|_{X}\leq C_3C_2r_1^2|t_2-t_1|,\quad t_1,t_2\in [0,t_0].
  \ea
  When $r_1>0$ is so small that $C_3C_2(1+t_0)<r_1^{-1}$ 
  and $C_3C_2r_1^2<C_1$
  we see that 
  \ba
  & & \|\Phi_\S({\bf v})\|_{C([0,t_0];Y)}<r_1,\\
  & & \|\Phi_\S({\bf v})\|_{C^{0,1}([0,t_0];X)}<C_1.
  \ea
  Hence $ \Phi_\S(E)\subset E$
 and we can consider $\Phi_\S$ as a map
 $\Phi_\S:E\to E$.
  
  As  $k>1+\frac 32$, it follows from Sobolev embedding theorem
  that $X=\H^{(k)}(\tilde N)\subset C^1(\tilde N)^2$. This yields that
  by \cite[Thm. 3]{Kato1975}, for the original reference, see Theorems III-IV
  in \cite{Kato1973},
  \beq\label{eq: stab.}
 & & \|(U^{\bf v}(t,s)- U^{\bf v^\prime}(t,s)){\bf h}\|_X \\
 & &\leq  \nonumber
 C_3\bigg(\sup_{t^\prime \in [0,t]}\|\A(t^\prime,{\bf v})-
 \A(t^\prime,{\bf v}^\prime)\|_{Y\to X}  \|U^{\bf v}(t^\prime,0) {\bf h}\|_Y\bigg)\\
 & &\leq  \nonumber
 C_3^2 \|{\bf v}-{\bf v}^\prime\|_{C([0,t_0];X)}  \|{\bf h}\|_Y.
  \eeq
  Thus,
  \ba
&&\|  U^{\bf v}(t,s)f_\S(s,{\bf v})- U^{\bf v^\prime}(t,s)f_{\S}(s,{\bf v}^\prime)\|_X\\
&\leq&\|  (U^{\bf v}(t,s)- U^{\bf v^\prime}(t,s))f_{\S}(s,{\bf v})\|_X+
\|  U^{\bf v^\prime}(t,s)(f_\S(s,{\bf v})- f_{\S}(s,{\bf v}^\prime))\|_X\\
&\leq&
(1+C_3)^2C_2r_1^2 \|{\bf v}-{\bf v}^\prime\|_{C([0,t_0];X)} .
  \ea
This implies that
\ba
& & \| \Phi_\S({\bf v})-\Phi_{\S}({\bf v}^\prime)\|_{C([0,t_0];X)}
\leq t_0(1+C_3)^2C_2r_1^2 \|{\bf v}-{\bf v}^\prime\|_{C([0,t_0];X))}.
\ea

Assume next that  $r_1>0$ is so small that we have also
\ba
t_0(1+C_3)^2C_2r_1^2<\frac 12.
\ea
 cf.\
  Thm.\ I in \cite{HKM} (or
  (9.15)  and (10.3)-(10.5) in \cite{Kato1975}). For $\S\in W$ this 
 implies that $$\Phi_\S:E\to E$$ is a contraction with a contraction constant 
 $C_L\leq \frac 12$,
  and thus 
 $\Phi_\S$ has a unique fixed point ${\bf u}$  in the space $E\subset C^{0,1}([0,t_0];X)$.
 
 Moreover, elementary considerations related to fixed points
 of the map $\Phi_\S$ show that ${\bf u}$  in  $C([0,t_0];X)$ depends in $E\subset C([0,t_0];X)$
 Lipschitz-continuously on
 $\S\in  W\subset C([0,t_0];H^{k+1}(\tilde N)\times H^{k}(\tilde N))$. Indeed, if 
 $$\|\S-\S^\prime\|_{C([0,t_0];H^{k+1}(\tilde N)\times H^{k}(\tilde N))}<\e,$$
 we see that
  \beq\label{fF bound r2 B}
 \|f_\S(t,{\bf v})-f_{\S^\prime}(t,{\bf v})\|_Y\leq C_2\e,\quad t\in [0,t_0],
 \eeq
 and when  (\ref{fF bound r2}) and (\ref{U bounds}) are satisfied 
 with $r=r_1$,
 we have
  \ba
& & \| \Phi_\S({\bf v})-\Phi_{\S^\prime}({\bf v}^\prime)\|_{C([0,t_0];Y)}
\leq C_3C_2t_0 r_1^2.
\ea
 Hence 
 \ba
 \|\Phi_\S({\bf v})-\Phi_{\S^\prime}({\bf v})\|_{C([0,t_0];Y)}\leq  t_0C_3C_2\e.
 \ea
 This and standard estimates for fixed points, yield that when $\e$ is small enough 
 the fixed point ${\bf u}^\prime$ of the map  $\Phi_{\S^\prime}:E\to E$ 
 corresponding to the source ${\S^\prime}$ and
  the fixed point ${\bf u}$ of the map $\Phi_{\S}:E\to E$ 
  corresponding to the source ${\S}$ satisfy
  \beq\label{stability}
  \|{\bf u}-{\bf u}^\prime\|_{C([0,t_0];X)}\leq  \frac 1{1-C_L}t_0C_3C_2\e.
  \eeq
{
Thus the solution ${\bf u}$ 
  depends in $C([0,t_0];X)$ Lipschitz continuously on
 $\S\in  C([0,t_0];H^{k+1}(\tilde N)\times H^{k}(\tilde N))$ (see also \cite[Sect.\ 16]{Kato1975}, 
 and \cite{Ringstrom})}.
In fact,  for analogous systems it is possible to show that $u$ is in $C([0,t_0];Y)$,
but one can not obtain  
 Lipschitz or H\"older  stability for $u$ in the $Y$-norm, see  \cite{Kato1975}, Remark 7.2. 
 
 Finally, we note that the fixed point ${\bf u}$ of $\Phi_\S$ can be found as a limit
 $ {\bf u}=\lim_{n\to \infty}{\bf u}_n$ in $C([0,t_0];X)$, where ${\bf u}_0=0$ and  ${\bf u}_n=\Phi_\S({\bf u}_{n-1})$.
Denote ${\bf u}_n=(g_n-g,\phi_n-\hat \phi)$. We see that if
 $$\supp({\bf u}_{n-1})\subset J_{g}(\supp(\S))$$ 
then also $$\supp(g_{n-1}-g)\subset J_{g}(\supp(\S)).$$
Hence  for all $x\in M_0\setminus J_{g}(\supp(\S))$ we see that $$J_{g_{n-1}}^-(x)
\cap J_{g}(\supp(\S))=\emptyset.$$ Then, using the definition of the map $\Phi_\S$ we see that
 $\supp({\bf u}_n)\subset J_{g}(\supp(\S))$. Using induction we see that this holds
for all $n$ and hence we have that the solution ${\bf u}$ satisfies
 \beq\label{eq: support condition for u}
 \supp({\bf u})\subset J_{g}(\supp(\S)).
 \eeq 
 \hiddenfootnote{REMOVE MATERIAL IN THIS FOOTNOTE OR MOVE IT ELSEWHERE:
  Let us next use the above considerations to study system 
   (\ref{eq: notation for hyperbolic system 1}) on $\R\times N$.
   For this end, we denote below the background metric of  $\tilde M=\R\times \tilde N$
   by $\tilde g$ and the background metric of  $\hattuM _0=\R\times  N$
   by $g$.
 Consider a source 
\ba
\S\in {C^1([0,T];H^{s_0-1}(\tilde N))\cap
C^0([0,t_0];H^{s_0}(\tilde N))}
\ea  that is supported in a compact set $\K$ and
small enough in the norm of this space and let $\tilde u$ 
be the solution of the  
system (\ref{eq: notation for hyperbolic system 1}) on $\R\times \tilde N$.
Assume that  $u$ is some $C^2$-solution of  
(\ref{eq: notation for hyperbolic system 1}) on $\R\times N$
with the same source $\S$. Our aim is to show that it
coincides with $\tilde u$ in its support.
%
For this end, let $\tilde g^\prime$ be a metric tensor which
has the same eigenvectors as $\tilde g$ and the same
eigenvalues except that  the unique negative eigenvalue of $\tilde g$
 is multiplied by $(1-h_1)$, $h_1>0$. 
Assume that $h_1>0$ is so small
that $J^+_{\hattuM _0,\tilde g^\prime}(p^+)\subset  \hattuM _0\subset N_0\times (-\infty,t_0]$
and assume that $\e$ is so small that $g(\tilde u)<\tilde g^\prime$ for all 
sources $\S$ with $\|\S\|_{C^1([0,t_0];H^{s_0-1}(\tilde N))\cap
C^0([0,t_0];H^{s_0}(\tilde N))}\leq \e$. Then
 $J^+_{\tilde M,g(\tilde u)}(p^+)\subset  (-\infty,t_0]\times N_0$.
This  proves the estimate (\ref{eq: Lip estim}).

Next us we consider solution $u$ on $\hattuM _0$ corresponding to source $\S$.
Let  $T(h_1,\S)$
be the supremum of all $T^\prime\leq t_0$ for which
$J^+_{g(u)}(p^+)\cap \hattuM _0(T^\prime)\subset N_1\times (-\infty,T^\prime]$.
Assume that $ T(h_1,\S)<t_0$.
Then for $T^\prime= T(h_1,\S)$ we see that
$u$ can be continued by zero from $ (-\infty,T^\prime]\times N_1$
to a function on $(-\infty,T^\prime]\times \tilde N$ that 
satisfies the system (\ref{eq: notation for hyperbolic system 1}). Thus
$u$ coincides in $ (-\infty,T^\prime]\times N_1$
with $\tilde u$ and vanishes outside this set.
Since  $J^+_{\tilde M,g(\tilde u)}(p^+)\subset  (-\infty,t_0]\times N_0$, this implies that
also $u$ satisfies 
$J^+_{g(u)}(p^+)\cap \hattuM _0(T^\prime)\subset (-\infty,T^\prime]\times N_0$.
As  $u$ is $C^2$-smooth and $N_0\subset \subset N_1$
we see that there is $T^\prime_1> T(h_1,\S)$ so 
that $J^+_{g(u)}(p^+)\cap \hattuM _0(T^\prime_1)\subset  (-\infty,T^\prime_1]
\times N_1$ that is in contradiction with definition of $ T(h_1,\S)$. Thus
we have to have  $T(h_1,\S)=t_0$.


This shows that when the norm of $\S$ is smaller than the above chosen $\e$
then $J^+_{g(u)}(p^+)\cap \hattuM _0\subset N_0\times (-\infty,t_0]$.
Then $u$ vanishes outside $N_1\times (-\infty,t_0]$ by
the support condition in (\ref{eq: notation for hyperbolic system 1}).
Thus $u$ and $\tilde u$ are both supported on  $N_1\times (-\infty,t_0]$ 
and coincide there. Since $\tilde u$ is unique, 
this shows that for a sufficiently small sources $ \S$
 supported on  $\K$   the solution $u$ of
(\ref{eq: notation for hyperbolic system 1})  in $M_0$
exists and is unique.}

\bigskip

\section*{Appendix C: 
An example  satisfying microlocal linearization stability}


\subsection*{C.1.\ Formulation of the direct problem} 
Let us define some physical fields  and introduce a
 model as a system of partial differential
 equations. Later we will motivate this system by discussion
 of the corresponding Lagrangians, but we postpone this
 discussion to  Appendix C.3 as it is not completely rigorous.

We assume that there are $C^\infty$-background fields $g$,  $\hat \phi$,
 on $M$. 
%

We consider 
a Lorentzian metric  $g$ on $\hattuM _0$ and $\phi=(\phi_\ell)_{\ell=1}^L$ where $\phi_\ell$
are scalar fields on $\hattuM _0=(-\infty,t_0)\times N$.

{Let $P=P_{jk}(x)dx^jdx^k$ be a symmetric tensor on $\hattuM _0$, corresponding below
to a direct perturbation to the stress energy tensor, and 
$Q=(Q_{\ell}(x))_{{\ell}=1}^K$ where  $Q_{\ell}(x)$ are real-valued functions on $\hattuM _0$, where $K\geq L+1$.
We denote by ${\mathcal V}( \phi_\ell;S_\ell)$
the potential functions of the fields $\phi_\ell$, 
\beq\label{eq: Slava-aaa}
{\mathcal V}( \phi_\ell;S_\ell)=\frac 12m^2\bigg(\phi_\ell +\frac 1{m^2}S_\ell\bigg)^2.
\eeq 
These potentials depend on the source variables $S_\ell$.
The way how  $S_\ell$,
called below the adaptive source functions, depend on other 
fields is explained later. We assume that there are smooth background
fields $\hat P$ and $\hat Q$. For a while we consider
the case when $\hat P=0$ and $\hat Q=0$, and discuss later
the generalization to non-vanishing background fields.

Using the $\phi$ and $P$ fields, we define the stress-energy tensor
 \beq\label{eq: T tensor1}
\hspace{-5mm} T_{jk}=\sum_{\ell=1}^L(\p_j\phi_\ell \,\p_k\phi_\ell 
-\frac 12 g_{jk}g^{pq}\p_p\phi_\ell \,\p_q\phi_\ell-\mathcal V( \phi_\ell;S_\ell)g_{jk})
+P_{jk}.\hspace{-14mm}
\eeq


We assume that $P$ and $Q$ are supported on
$\K=J^+_{\tilde g}(\hat p^-)\cap \hattuM _0$. 
Let  us represent the stress energy tensor (\ref{eq: T tensor1})
 in the form
  \ba
T_{jk}
&=&P_{jk}+Zg_{jk}+{\bf T}_{jk}(g,\phi),\quad Z=-(\sum_{\ell=1}^L S_\ell\phi_\ell+
\frac 1{2m^2}S_\ell^2),
\\
{\bf T}_{jk}(g,\phi)&=&\sum_{\ell=1}^L(\p_j\phi_\ell \,\p_k\phi_\ell 
-\frac 12 g_{jk}g^{pq}\p_p\phi_\ell \,\p_q\phi_\ell-
\frac 12 m^2\phi_\ell^2g_{jk}),
 \ea
where we call $Z$ the stress energy density caused by the sources $S_\ell$.

\generalizations{We will also add in our considerations linear, first order differential operators
\beq\label{B-interaction}
B_\ell(\phi)=\sum_{\kappa=1}^L a^\kappa_\ell\phi_\kappa(x)
+ \sum_{\kappa,\alpha,\beta =1}^Lb^{\kappa,\a,\beta}_\ell \phi_\kappa (x)
\phi_\alpha (x)\phi_\beta (x)\hspace{-1cm}
\eeq modeling interaction of matter fields,
 where $a^{\kappa}_\ell,b^{\kappa,\a,\beta}_\ell\in \R$ are constants. 
\HOX{ We have  added here some first order terms
to the wave equation that connect $\phi$-variables together. WE NEED TO ADD 
TERMS IN STRESS ENERGY TENSOR AND IN THE APPENDIX B AND THE MAIN THEOREM}}

{Now we are ready to formulate the direct problem for 
the adaptive Einstein-scalar field equations. Let $g$ and $\phi$ satisfy
\beq\label{eq: adaptive model}& &\Ein_{g}(g) =P_{jk}+Zg_{jk}+{\bf T}_{jk}(g,\phi),\quad Z=-(\sum_{\ell=1}^L S_\ell\phi_\ell+
\frac 1{2m^2}S_\ell^2),
\\ \nonumber
& &\square_g\phi_\ell -\mathcal V^{\prime}( \phi_\ell;S_\ell) =0
\generalizations{+B_\ell(\phi)}
\quad
\hbox{in }\hattuM _0,\quad \ell=1,2,3,\dots,L,
\\ \nonumber
& &S_\ell={\mathcal S}_\ell(g,\phi, \nabla
 \phi,Q,\nabla Q,P,\nabla^g P),\quad
\hbox{in }\hattuM _0,
 \\ \nonumber
& & g=g,\quad \phi_\ell=\hat \phi_\ell,\quad
\hbox{in }\hattuM _0\setminus \K.
\eeq
Above, $\mathcal V^{\prime}(\phi ;s)=\p_\phi \mathcal V(\phi;s)$ so that
$
\mathcal V^{\prime}( \phi_\ell;S_\ell)=m^2\phi_\ell+S_\ell.
$
}
We assume that  the background fields $g$,  $\hat \phi$,
 satisfy these equations with  $\hat Q=0$ and  $\hat P=0$.}

We consider  here $P=(P_{jk})_{j,k=1}^4$
and $Q=(Q_\ell)_{\ell=1}^K$   as   fields that we can control and call those
the controlled source fields. 
 Local existence of the solution for small sources $P$ and $Q$ is considered
in Appendix B. 

%
%
To obtain a physically meaningful model,
we need to consider how  the adaptive source functions $\mathcal S_\ell$ should be chosen 
so that the physical conservation law in relativity
\beq\label{conservation law CC}
\nabla_k(g^{kp} T_{pq})=0
\eeq
 is satisfied. Here $\nabla=\nabla^g$ is the connection corresponding to the metric $g$.

  We
note that the conservation law 
is a necessary condition for the equation (\ref{eq: adaptive model})
to have solutions for which
$\Ein_{g}(g)=\Ein(g)$, i.e., that the solutions of
(\ref{eq: adaptive model}) are solutions of  the Einstein field equations.

{The functions ${\mathcal S}_\ell(g,\phi, \nabla
 \phi,Q,\nabla Q,P,\nabla^g P)$ model the devices that we use to perform active
 measurements. Thus, even though the Assumption S below may appear quite technical,
 it can be viewed as the instructions on  how to build a device that can be used 
 to measure the structure of the spacetime far away. Outside the support
 of the measurement device (that contain the union of the supports of $Q$ and $P$) we have just assumed that the standard coupled Einstein-scalar field
 equations hold, c.f. (\ref{S-vanish condition}). 
We can consider them in the form
\ba
& &{\mathcal S}_\ell(g,\phi, \nabla
 \phi,Q,\nabla Q,P,\nabla^g P)=Q_\ell+{\mathcal S}^{2nd}_\ell(g,\phi, \nabla
 \phi,Q,\nabla Q,P,\nabla^g P)
\ea
for $\ell=1,2,\dots,L$
where $Q_\ell$ are the primary sources and ${\mathcal S}^{2nd}_\ell$,
that depend also on $Q_\ell$ with $\ell=L+1,L+2,\dots,K$,
corresponds to the response of the measurement device that forces
the conservation law to be valid.

The solution $(g,\phi)$ of (\ref{eq: adaptive model}) is a solution of the equations
(\ref{eq: adaptive model with no source}) when we denote
 \ba
& &\F^1_{jk}=P_{jk}+Zg_{jk},\\
& &\F^2_\ell=\mathcal V^{\prime}( \phi_\ell;S_\ell)-\mathcal V^{\prime}( \phi_\ell;0)=S_\ell.
\ea
Our next goal is to construct suitable adaptive source
functions ${\mathcal S}_\ell$ and consider what kind of sources $\F^1$ and $\F^2$
of the above form can be obtained by varying $P$ and $Q$.

We will consider  adaptive source
functions ${\mathcal S}_\ell$  satisfying the following conditions:

 {
\medskip 

\noindent{\bf Assumption S}:

\noindent
The adaptive source functions ${\mathcal S}_\ell(g,\phi, \nabla
 \phi,Q,\nabla Q,P,\nabla^g P)$ have the following properties: 
\smallskip

{(i) Denoting  $c=\nabla
 \phi$, $C=\nabla^g P$, and $H=\nabla Q$
 we assume 
  that ${\mathcal S}_\ell(g,\phi, c,Q,H,P,C)$ are smooth {non-linear} functions, 
  of the pointwise values $g_{jk}(x),\phi(x), \nabla
 \phi(x),$ $Q(x),\nabla Q(x),P(x),$ and $\nabla^g P(x)$,
  defined
 near $(g,\phi, c,Q,H,P,C)=(g,\hat \phi,\nabla \hat \phi ,0,0,0,0)$, that
satisfy
\beq\label{S-vanish condition}
 {\mathcal S}_\ell(g,\phi,c,0,0,0,0)=0.
\eeq
 
\smallskip


 \smallskip

(ii)  
We assume that  ${\mathcal S}_\ell$ is independent of $P(x)$ and the dependency of ${\mathcal S}$ on  $\nabla^g P$ 
and $\nabla Q$ is only due to the dependency
in the term $g^{pk}\nabla^g_p( P_{jk}+Zg_{jk})=
g^{pk}\nabla^g_p P_{jk}+\nabla^g_jQ_K$, associated to
the divergence of the perturbation of $T$,
that is,
there exist functions $\tilde {\mathcal S}_\ell$ so that
\ba
{\mathcal S}_\ell(g,\phi, c,Q,H,P,C)=\tilde {\mathcal S}_\ell(g,\phi,c,Q,
R),\quad R=(g^{pk}\nabla^g_p ( P_{jk}+Q_Kg_{jk}))_{j=1}^4.\ea
%
}
 

 \bigskip


Below, denote $\hat R=g^{pk}\hat\nabla_p \hat P_{jk}+\hat\nabla_j\hat Q_K.$
Note that we still are considering the  case when $\hat Q=0$ and $\hat P=0$
so that $\hat R=0$, too. This implies
 that for the background fields that adaptive source functions $\mathcal S_\ell$ vanish.

To simplify notations, we also denote below
$\tilde {\mathcal S}_\ell$ just by ${\mathcal S}_\ell$ and indicate the function
which we use by the used variables in these functions.

Below we will denote $Q=(Q^{\prime},Q_K)$, $Q^{\prime}=(Q_\ell)_{\ell=1}^{K-1}$.
There are examples when the background fields $(g,\hat \phi)$ and 
the adaptive source functions ${\mathcal S}_\ell$ exists and satisfy the Assumption S. 
  This is shown later in the case the following condition is valid for the 
     background fields:
 \medskip

 {\bf Condition A}:
Assume that at any $x\in \overline U$ there is
a permutation $\sigma:\{1,2,\dots,L\}\to \{1,2,\dots,L\}$, denoted $\sigma_x$, such that the 
$5\times 5$ matrix $[ B_{jk}^\sigma(\hat \phi(x),\nabla \hat \phi(x))]_{j,k\leq 5}$
is invertible, where
  \ba
[ B_{jk}^\sigma(\phi(x),\nabla \phi(x))]_{k,j\leq 5}=\left[\begin{array}{c}
(\,\p_j  \phi_{\sigma(\ell)}(x))_{\ell\leq 5,\ j\leq 4}\\
(\phi_{\sigma(\ell)}(x))_{\ell\leq 5}\end{array}\right].
 \ea
 Below, for a permutation $\sigma:\{1,2,\dots,L\}\to \{1,2,\dots,L\}$ we denote by $U_{g,\sigma}$
 the open set of points $x\in U$ for which $[ B_{jk}^\sigma(\hat \phi(x),\nabla \hat \phi(x))]_{j,k\leq 5}$
is invertible. So, we assume that the sets $U_{g,\sigma}$,
$\sigma\in \Sigma(L)$ is an open covering of $U$.
 \medskip

 Our next aim is to prove the following:
  \medskip

 \begin{theorem} \label{thm: Good S functions}
Let $L\geq 5$ and assume that $\hat Q=0$ and $\hat P=0$ so that
 $\hat R=0$.
Moreover, assume that 
 Condition A is valid. Then for all permutations $\sigma:\{1,2,\dots,L\}\to \{1,2,\dots,L\}$ there exists  functions $  {\mathcal S}_{\ell,\sigma}$ satisfying Assumption S
  such that
%
\smallskip

(i) For all
$x\in   U_{g,\sigma}$ 
 the differential of $$  {\mathcal S}_\sigma(g,\hat \phi,  \nabla
 \hat \phi,Q,R)=(  {\mathcal S}_{\ell,\sigma}(g,\hat \phi, \nabla
 \hat \phi,Q,R))_{\ell=1}^L$$ with respect to $Q$ and $R$, that is, the map
  \beq\label{surjective 1}
D_{Q,R}  {\mathcal S}_\sigma(g(x),\hat \phi(x),  \nabla
 \hat \phi(x),Q,R)|_{Q=\hat Q(x), R=\hat R(x)}:\R^{K+4}\to \R^L\hspace{-1cm}
 \eeq
 is surjective.
\smallskip
 
(ii)  
The adaptive source functions ${\mathcal S}_\sigma$ are such that
 for
 $(Q_\ell)_{\ell=1}^K$
and  $(P_{jk})$ that are sufficiently close to $\hat Q=0$ and $\hat P=0$ in the $C^{{3}}_b(M_0)$-topology 
and supported in $  U_{g,\sigma}$
the equations (\ref{eq: adaptive model}) with source functions  ${\mathcal S}_\sigma$ have
a unique solution $(g,\phi)$ and 
the conservation
 law (\ref{conservation law CC}) is valid. 
 \smallskip

 (iii) Under the same assumptions as in (ii), when  $(g,\phi)$ is a solution  of (\ref{eq: adaptive model}) 
 with the controlled source functions $P$ and $Q$, we have
  $Q_K=Z$. This means that  the physical field $Z$  can be directly controlled. 
 
 \end{theorem}

}}

\noindent
{\bf Proof.} As  one can enumerate the $\ell$-indexes 
of the fields $\phi_\ell$ as one wishes,
it is enough to prove the claim with one $\sigma$. We
consider below the case when $\sigma=Id$.

 Consider 
  a symmetric (0,2)-tensor $P$ and a scalar functions $Q_\ell$
  that are $C^3$-smooth and compactly supported in  $U_{g,\sigma}$.
Let $[P_{jk}(x)]_{j,k=1}^4$ be  the coefficients of 
  $P$ in local coordinates and $Q(x)=(Q_\ell(x))_{\ell=1}^L$. 
 

To obtain adaptive required adaptive source functions, let us start implications of the conservation law (\ref {conservation law CC}).
To this end, consider $C^2$-smooth functions 
 $S_\ell(x)$ on $U_{g,\sigma}$.

Note that since $[\nabla_p,\nabla_n]=[\p_p,\p_n]=0$ (see \cite[Sect. III.6.4.1]{ChBook}), 
 \ba
& &\nabla_ p (g^{pj} {\bf T}_{jk}(g,\phi)=
\\ & &\sum_{\ell=1}^L \nabla_ p
g^{pj}  (\p_j\phi_\ell \,\p_k\phi_\ell 
-\frac 12 g_{jk}g^{nm}\p_n\phi_\ell \,\p_m\phi_\ell-
\frac 12 m^2\phi_\ell^2g_{jk})
\\
&=&
\sum_{\ell=1}^L (g^{pj}  \nabla_ p \p_j\phi_\ell) \,\p_k\phi_\ell
+\sum_{\ell=1}^L  (g^{pj}  \p_j\phi_\ell \, \nabla_ p \p_k\phi_\ell)\\
& &-
 \frac 12
 \sum_{\ell=1}^L \delta^p_k 
 (g^{nm}(\nabla_ p \p_n\phi_\ell) \,\p_m\phi_\ell+
 g^{nm}\p_n\phi_\ell \,(\nabla_ p  \p_m\phi_\ell))
  - \sum_{\ell=1}^L m^2  \delta^p_k   \phi_\ell \p_p \phi_\ell
\\
&=&
\sum_{\ell=1}^L (g^{pj}  \nabla_ p \p_j\phi_\ell) \,\p_k\phi_\ell
+\sum_{\ell=1}^L  (g^{pj}  \p_j\phi_\ell \, (\nabla_ p \p_k\phi_\ell))\\
& &-
 \frac 12
 \sum_{\ell=1}^L 
 (g^{nm}\p_m\phi_\ell\,(\nabla_ n \p_k\phi_\ell)+
 g^{nm}\p_n\phi_\ell \,(\nabla_ m  \p_k\phi_\ell))
 - \sum_{\ell=1}^L m^2   \phi_\ell \p_k \phi_\ell
\\
&=&
\sum_{\ell=1}^L (g^{pj}  \nabla_ p \p_j\phi_\ell -m^2   \phi_\ell) \,\p_k\phi_\ell. 
\ea
Thus conservation law (\ref {conservation law CC})
gives for all $j=1,2,3,4$ equations
 \ba
0&=&\nabla_p^g (g^{pk}T_{jk})\\
&=&\nabla_p^g (g^{pk}({\bf T}_{jk}(g,\phi)+P_{jk}+Zg_{jk}))\\
&=& \sum_{\ell=1}^L\bigg((g^{pk} \nabla^g_p\p_k\phi_\ell )\,\p_j\phi_\ell 
-(m^2_\ell \phi_\ell \p_p \phi_\ell) \delta_j^p \\
& &\quad \quad- \nabla^g_p (g^{pk}g_{jk}
(S_\ell\phi_\ell{\color{black}+\frac 1{2m^2}S_\ell^2})+g^{pk}P_{jk})\bigg)
\\
&=& \sum_{\ell=1}^L\bigg(
(g^{pk} \nabla^g_p\p_k\phi_\ell -m^2_\ell \phi_\ell ) \,\p_j \phi_\ell 
\\ & &\quad \quad-\nabla^g_p (g^{pk}g_{jk}(S_\ell\phi_\ell
{\color{black}+\frac 1{2m^2}S_\ell^2})
+g^{pk}P_{jk})\bigg)
\\
&=&\sum_{\ell=1}^L S_\ell \,\p_j\phi_\ell- \nabla^g_p (g^{pk}g_{jk}
(S_\ell \phi_\ell{\color{black}+\frac 1{2m^2}S_\ell^2}))+ g^{pk} \nabla^g_p P_{jk}\\
&=&\left (\sum_{\ell=1}^L S_\ell \,\p_j\phi_\ell\right)-\p_j\left(
\sum_{\ell=1}^L S_\ell \phi_\ell {\color{black}+\frac 1{2m^2}S_\ell^2}\right)+g^{pk}\nabla_p^g P_{jk}.
\ea
Summarizing, the conservation law yields 
\beq
\label{mod. conservation law a}
\left (\sum_{\ell=1}^L S_\ell \,\p_j\phi_\ell\right)-\p_j\left(
\sum_{\ell=1}^L S_\ell \phi_\ell +\frac 1{2m^2}S_\ell^2\right)+g^{pk}\nabla_p^g P_{jk}=0,\hspace{-1cm}
\eeq
for $j=1,2,3,4$.

Recall that the field $Z$ has the definition
\beq\label{meson number}
\sum_{\ell=1}^L S_\ell \phi_\ell+\frac 1{2m^2}S_\ell^2=-Z.
\eeq

Then, the  conservation law  (\ref {conservation law CC}) holds if we have
\beq\label{meson number2}& &
\sum_{\ell=1}^L S_\ell \,\p_j\phi_\ell
=- g^{pk}\nabla^g_p  {\color{black}V}_{jk},\quad
 {\color{black}V}_{jk}=(P_{jk}+g_{jk}Z),
\eeq
for $j=1,2,3,4.$

Equations (\ref{meson number}) and (\ref{meson number2}) give together  five point-wise equations
 for the  functions $S_1,\dots,S_L$.

%
%
%
Recall that we consider here the case when $\sigma=Id$.
By Condition A, at any $x\in U_{g,\sigma}$ that the 
$5\times 5$ matrix $( B_{jk}^\sigma(\hat \phi(x),\nabla \hat \phi(x)))_{j,k\leq 5}$
is invertible, where
  \ba
( B_{jk}^\sigma(\phi(x),\nabla \phi(x)))_{j,k\leq 5}=\left(\begin{array}{c}
( \,\p_j  \phi_{\cell}(x))_{j\leq 4,\ \ell\leq 5}\\
(   \phi_{\cell}(x)) _{\ell\leq 5}\end{array}\right).
 \ea

 We  consider
a $\R^K$ valued function  $Q(x)=(Q^\prime(x),Q_K(x))$,
where  $$Q^\prime=(Q_\ell)_{\ell=1}^{K-1}.$$

Also, below $ {\color{black}R}_j=g^{pk}\nabla^g_p  {\color{black}V} _{jk},$ $ {\color{black}V}_{jk}=P_{jk}+g_{jk}Z$ and  we require that identity
\beq
& &Q_K=Z
\label{S R Z equations 3bbb}
\eeq
holds.

Motivated by equations (\ref{meson number}), (\ref{meson number2}),
and (\ref{S R Z equations 3bbb}),
our next aim is to consider a point $x\in U_{g,\sigma}$,
and construct  functions
${{\rbtext {\mathcal S}^{(1)}}_{\sigma,\ell}}(\phi,\nabla \phi, Q^\prime,Q_K,R,g)$, $\ell=1,2,\dots,L$ that satisfy
\beq\label{S R Z equations 1bbb}
& &\sum_{\ell=1}^5 {{\rbtext {\mathcal S}^{(1)}}_{\sigma,\ell}}(\phi,\nabla \phi, Q^\prime,Q_K,R,g)
 \,\p_j\phi_{\cell}=- R_{j} -
\sum_{\ell=6}^{L} Q_{\sigma,\ell} \,\p_j\phi_{\cell},\hspace{-1cm} \\
& &\sum_{\ell=1}^5 {{\rbtext {\mathcal S}^{(1)}}_{\sigma,\ell}}(\phi,\nabla \phi, Q^\prime,Q_K,R,g)\, \phi_{\cell}=-\bigg(Q_K+\sum_{\ell=6}^L Q_{\sigma,\ell} \phi_{\cell}
+\label{S R Z equations 2bbb}\\ \nonumber\\ & &\quad
+\sum_{\ell=1}^L\frac 1{2m^2} {{\rbtext {\mathcal S}^{(1)}}_{\sigma,\ell}}(\phi,\nabla \phi, Q^\prime,Q_K,R,g)^2\bigg)
.\hspace{-1cm}
\eeq
%
%
Let
\ba
& &(Y_\sigma(\phi,\nabla \phi))(x)=\psi(x)(B^\sigma(\phi,\nabla \phi))^{-1},\quad \hbox{for
$x\in U_{g,\sigma}$,}
\\
& &(Y_\sigma(\phi,\nabla \phi))(x)=0,\quad \hbox{for
$x\not \in U_{g,\sigma}$,}
\ea 
where $\psi\in C^\infty_0(U_{g,\sigma})$ has value 1 in $\supp(Q)\cup\supp(P)$.

Then we define 
${{\rbtext {\mathcal S}^{(1)}}_{\sigma,\ell}}={{\rbtext {\mathcal S}^{(1)}}_{\sigma,\ell}}(g,\phi,\nabla \phi, Q^\prime,Q_K,R)$, $\ell=1,2,\dots,L,$ to be the solution of the system 
\beq\label{S sigma formulas}
& &(S_{\sigma,\ell})_{\ell\leq 5}= Y_\sigma(\phi,\nabla \phi) \left(\begin{array}{c}
(-  {\color{black}R}_j -
\sum_{\ell=6}^L Q_{\sigma,\ell} \,\p_j\phi_{\cell})_{j\leq 4}
\\
-Q_K-\sum_{\ell=6}^L Q_{\sigma,\ell}\phi_{\cell}-\sum_{\ell=1}^L\frac 1{2m^2}S_{\sigma,\ell}^2 \end{array}\right)
\hspace{-1cm}\\ \nonumber
%
\hspace{-1cm}\\
& & \nonumber
(S_{\sigma,\ell})_{\ell\geq 6}= (Q_{\cell})_{\ell\geq 6}.
\eeq  
{When $Q$ and $R$ are sufficiently small,
 this equation can be solved point-wisely, 
 at each point $x\in U_{g,\sigma}$, using iteration by the Banach fixed point theorem.}


Let
 \ba
( K_{jk}^\sigma(\phi(x),\nabla \phi(x)))_{j,k\leq 5}=\left(\begin{array}{c}
( \,\p_j  \phi_{\cell}(x))_{j\leq 4,\ 6\leq \ell\leq L}\\
(   \phi_{\cell}(x)) _{6\leq \ell\leq L}\end{array}\right).
 \ea
Then
we see that
 the differential of 
$ {\mathcal S}_\sigma=( {\mathcal S}_{\sigma,\ell})_{\ell=1}^L$
with respect to $(Q^{\prime},Q_{K}, {R})$ at $(Q,R)=(0,0)$, that is,
 \beq\label{eq: deriv1}\\
 \nonumber
& &D_{Q^{\prime}, Q_{K}, {R}} {\mathcal S}_\sigma(g,\hat \phi, \nabla
 \hat \phi,Q^{\prime},Q_{K}, {R})|_{Q=0,R=0}:\R^{K+4}\to \R^L,\\
 & &\nonumber (Q^{\prime},Q_{K}, R)
\mapsto 
 -\left(\begin{array}{cc}
Y_\sigma(\hat \phi,\nabla\hat\phi)   &Y_\sigma(\hat\phi,\nabla\hat\phi)K(\hat\phi,\nabla\hat\phi) 
\\
0 &I_\sigma  \end{array}\right)
 \left(\begin{array}{c}
\left(\begin{array}{c}
  R\\
 Q_K \end{array}\right) \\
Q^\prime\\
\end{array}\right),
 \eeq
 is surjective, where $I_\sigma=[\delta_{k,j+5}]_{k\leq K-1,\ j\leq L-5,}\in \R^{(K-1)\times (L-5)}$. Hence (i) is valid.




%

%
%
%

By their construction, the functions ${\mathcal S}_\sigma=( {\mathcal S}_{\sigma,\ell})_{\ell=1}^L$ satisfy the equations
(\ref{meson number}) and (\ref{meson number2})
for all $x\in U_{g,\sigma}$ and also equation  (\ref{S R Z equations 3bbb}) holds.
\hiddenfootnote{ 
We make the following note related 
the case considered in the main part of the paper when $Q,P,R\in \I^m(Y)$, where $Y$ is 2-dimensional sufrace,
and we  need to use the principal symbol ${\bf r}$ of $R$
as an independent variable:
Let us consider also a point $x_0\in \hattuM _0$ and $\eta$
be a light-like covector choose
coordinates so that $g=\diag(-1,1,1,1)$ and $\eta=(1,1,0,0).$
When $c=(c^k)_{k=1}^4\in \R^4$ and $P^{jk}=C^{jk}(x\,\cdotp \eta)^a_+$, where $C^{jk}$ is such a symmetric
matrix that
$C^{11}=c_1$, $
C^{12}=C^{21}=\frac  12 c_2$,
 $C^{13}=C^{31}=c_3$, and  $C^{14}=C^{41}=c_4$ and other $C^{jk}$ are zeros. 
Then we have
\ba
\nabla^g_j (C^{jk}(x\,\cdotp \eta)_+^a)&=&\eta_jP^{jk}=(C^{1k}+ C^{2k})(x\,\cdotp \eta)_+^{a-1}=c_k(x\,\cdotp \eta)_+^{a-1}.
\ea
As we can always choose coordinates
that $g$ and $\eta$ have at a given point the above forms,
we can obtain arbitrary vector ${\bf r}$, as
principal symbol of $R$, by considering
$Q^{\prime},Q_{L+1},P\in \I^m(Y)$, where $Y$ is 2-submanifold, with principal symbols of
$\tilde {\bf p}$ and $\tilde {\bf z}$ satisfing equations corresponding to
equations
$g^{jk}\hat \nabla_j ({\bf p}_{kn}+ {\bf z}g_{kn})\in \I^m(Y)$, 
and the sub-principal symbols of ${\bf p}_{kn}$ and ${\bf z}$
vary arbitrarily.} 
Hence (iii) is valid.

 Above, the equation (\ref{meson number2}) is valid by construction of the
 functions $( {\mathcal S}_\ell)_{\ell=1}^L$. Thus the conservation law is valid.
 This proves (ii). \hfill \Box \medskip



%

%

Note that as  the adaptive source functions
 $\mathcal S_\ell$  were constructed in Theorem \ref{thm: Good S functions}
 using the inverse function theorem, the results of Theorem \ref{thm: Good S functions} are valid also if
 $\hat Q$  and $\hat P$ are sufficiently small non-vanishing fields
 and $g$ and $\hat \phi$ satisfy the Einstein scalar field equations  (\ref{eq: adaptive model})
 with these background fields. Next we return to the case when $\hat P=0$  and $\hat Q=0$.

 \subsection*{C.2.\  Microlocal linearization stabililty}
 Below we consider the case when $\hat P=0$  and $\hat Q=0$ and use the adaptive source functions
 $\mathcal S_\ell$  constructed in Theorem \ref{thm: Good S functions} and
 its proof.
 
 Assume that $Y\subset \hattuM _0$
is a 2-dimensional space-like submanifold and consider local coordinates defined 
in  $V\subset \hattuM _0$. Moreover, assume that 
in these local coordinates $Y\cap V\subset \{x\in \R^4;\ x^jb_j=0,\  x^jb^\prime_j=0\}$,
where $b^\prime_j\in \R$  and let ${\bf p}\in \I^{n}(Y)$, $n\leq n_0=-17$,  be defined by 
\beq\label{eq: b b-prime BBB}
{\bf p}_{jk}(x^1,x^2,x^3,x^4)=
\re \int_{\R^2}e^{i(\theta_1b_m+\theta_2b^\prime_m)x^m}\symbolP_{jk}(x,\theta_1,\theta_2)\,
d\theta_1d\theta_2.\hspace{-2cm}
\eeq
Here, we assume that 
 $\symbolP_{jk}(x,\theta)$, $\theta=(\theta_1,\theta_2)$ 
 are  classical symbols and we denote their principal symbols by
 $\sigma_p({\bf p}_{jk})(x,\theta)$. When 
$x\in Y$ and $\xi=(\theta_1b_m+\theta_2b^\prime_m)dx^m$ so that $(x,\xi)\in N^*Y$,
we denote the value of the principal symbol $\sigma_p({\bf p})$ at $(x,\theta_1,\theta_2)$ by 
 $\tilde \symbolP ^{(a)}_{jk}(x,\xi)$, that is,  $ \tilde \symbolP ^{(a)}_{jk} (x,\xi)=\sigma_p({\bf p}_{jk})(x,\theta_1,\theta_2)$,
 and say that it is the principal symbol of ${\bf p}_{jk}$ at $(x,\xi)$, associated to the phase function $\phi(x,\theta_1,\theta_2)=(\theta_1b_m+\theta_2b^\prime_m)x^m$.
The above defined principal symbols can be defined invariantly, see \cite{GuU1}.

We assume that also ${\bf q}^\prime,{\bf z}\in \I^{n}(Y)$ have 
representations (\ref{eq: b b-prime BBB}) with classical symbols. Below
we consider symbols in local coordinates.
{Let
 us denote the principal symbols of ${\bf p},{\bf q^{\prime}},{\bf z}\in \I^{n}(Y)$ 
 by  $\tilde \symbolP ^{(a)}(x,\xi),$ $\tilde \symbolQ_{1}^{(a)}(x,\xi)$,
 $\tilde \symbolQ_1^{(a)}(x,\xi)$, respectively
  and let $\tilde \symbolP ^{(b)}(x,\xi)$ and $\tilde \symbolQ_1^{(b)}(x,\xi)$ denote the sub-principal
symbols of ${\bf p}$ and ${\bf z}$, correspondingly, at $(x,\xi)\in N^*Y$.


We will below consider what happens when  $({\bf p}_{jk}+{\bf z}g_{jk})\in \I^{n}(Y)$
 satisfies 
\beq\label{eq: leading order lineariz. conse. law}
g^{lk} \nabla_l^{g} ({\bf p}_{jk}+{\bf z}g_{jk})\in \I^{n}(Y),\quad j=1,2,3,4.
\eeq
Note that a priori this function is only in $\I^{n+1}(Y)$, so the
assumption (\ref{eq: leading order lineariz. conse. law}) means that
$g^{lk} \nabla_l^{g} ({\bf p}_{jk}+{\bf z}g_{jk})$
is one degree smoother than it should be a priori.

When (\ref{eq: leading order lineariz. conse. law}) is valid, we say
that  {\it  the leading order of singularity of the wave satisfies
the linearized conservation
law}. This  corresponds to the assumption that the principal
symbol of the sum of divergence of the first two terms appearing in the stress energy tensor
on the right hand side of (\ref{eq: adaptive model}) vanishes.
}


{By \cite{GuU1},  
the identity  (\ref{eq: leading order lineariz. conse. law}) is equivalent to the vanishing
of the principal symbol on $N^*Y$,} that is,
\beq\label{eq: lineariz. conse. law symbols mLS}
& &g^{lk}\xi_l(\tilde \symbolP ^{(a)}_{kj} (x,\xi)
+g_{kj}(x)\tilde \symbolQ_1^{(a)}(x,\xi))=0,\ \hbox{for 
$j\leq 4$ and $\xi\in N^*_x Y$}.
\eeq
We say that this is the {\it linearized conservation law for principal symbol of $R$}.

  Let us consider source fields
that have the form $Q^{\prime}_\e=((Q_\e)_\ell)_{\ell=1}^{K-1}=\e{\bf q^{\prime}}$, $(Q_\e)_K=\e{\bf z}$ and 
 $P_\e=\e {\bf p}$.
We denote 
${\bf q}=({\bf q^{\prime}},{\bf z})$. We assume that 
${\bf q^{\prime}}$, ${\bf z}$, and ${\bf p}$
%
%
are supported in $\hat V\subset \subset U$. 

Let $u_\e=(g_\e,\phi_\e)$ be the  solution of
 (\ref{eq: adaptive model}) with source  $P_\e$  and $Q_\e$.
 Then $u_\e$ depends $C^4$-smoothly on $\e$ and
$(g_\e,\phi_\e)|_{\e=0}=(g,\hat \phi)$.
Denote $\p_\e  (g_\e,\phi_\e)|_{\e=0}=(\dot g,\dot \phi)$.
When $\e_0$ is small enough,
$P_\e$ and $Q_\e$ are supported in $U_{g_\e}$ for all $\e\in (0,\e_0)$.

Let $$R_\e=
g_\e ^{pk}\nabla_p^{g_\e} ((P_\e)_{jk}+ g^\e_{jk}(Q_\e)_K)$$ and 
$$
({S_\e})_\ell=
 {\mathcal S}_\ell(g_\e,\phi_\e, \nabla
 \phi_\e ,Q_\e ^{\prime},(Q_\e)_{K},R_\e),
$$
where
 $\mathcal S_\ell$  are   the adaptive source functions constructed in Theorem \ref{thm: Good S functions} and
 its proof.

Then $S_\e|_{\e=0}=0$ and $\p_\e  S_\e|_{\e=0}={\dot S}$
satisfy
\beq\label{eq: formula for sources 1}\\
\nonumber
{\dot S}_\ell=
D_{Q',Q_K, R} {\mathcal S}_\ell(g,\hat \phi, \nabla
 \hat \phi,Q^{\prime},Q_{K},R)\bigg|_{Q'=0,Q_K=0,R=0}
 \left(\begin{array}{c}
{\bf q'}
\\
{\bf z}
\\
{\bf r}
 \end{array}\right),
 \eeq
where 
${\bf r}=
g^{pk}\hat \nabla_p ({\bf p}_{jk}+g_{jk}{\bf z})$.

Functions $\dot u=(\dot g,\dot \phi)$ satisfy
the  linearized Einstein-scalar field equation (\ref{eq: final linearization}).
The  linearized Einstein-scalar field equation (\ref{eq: final linearization})
is
\ba
& &e_{pq}(\dot g)-t^{(1)}_{pq}(\dot g)-t^{(2)}_{pq}(\dot \phi)={f}^1_{pq}
\\
& &\square_{g}\dot \phi^\ell
-g^{nm}g^{kj}(\p_n\p_j\hat \phi_\ell+ \hat \Gamma_{nj}^p\p_p
\hat \phi_\ell)\dot g_{mk}
-m^2\dot \phi^\ell= {f}^2_\ell,
\ea
where 
\beq\label{eq: formula for sources 2}
& &{f}^1_{pq}=
{\bf p}_{pq}-(\sum_{\ell=1}^L \dot S_\ell \hat \phi_\ell)g_{pq},\quad
\hbox{and }
-(\sum_{\ell=1}^L \dot S_\ell \hat \phi_\ell)={\bf z},\\
\nonumber 
& & {f}^2_\ell=\dot S_\ell. \nonumber 
\eeq
By Theorem \ref{thm: Good S functions} (ii),
  $u_\e=(g_\e,\phi_\e)$ satisfy  the conservation
 law (\ref{conservation law CC}).
This yields that  $\dot u=(\dot g,\dot \phi)$ satisfies the
  linearized Einstein-scalar field equation (\ref{eq: final linearization})
 and  the linearized conservation law (\ref{eq: lineariz. conse. law PRE}) is valid, too.
The linearized  conservation law (\ref{eq: lineariz. conse. law PRE}) gives,
by the considerations before (\ref {eq: formula for sources 2}),
\beq\label{lin cons 2B}
& &\hspace{-1.5cm}\left (\sum_{\ell=1}^L {f}^2_\ell \, \p_j\hat \phi_\ell\right)
+ g^{pk}\hat \nabla_p {f}^1_{kj}=0.\hspace{-1cm}
\eeq

Below,
we use the adaptive source functions
 $\mathcal S_\ell$  constructed in Theorem \ref{thm: Good S functions} and
 its proof.
 
%

{We see that 
\beq\label{eq: source f}
{f}=F(x;{\bf p},{\bf q})=(F^{(1)}(x;{\bf p},{\bf q}),F^{(2)}(x;{\bf p},{\bf q}))
\eeq
has by formulas (\ref{eq: formula for sources 1}) and (\ref{eq: formula for sources 2})
and Assumption S 
the form
\beq
\label{lin wave eq source A1a} F^{(1)}_{jk}(x;{\bf p},{\bf q})={\bf p}_{jk}+{\bf z}(x) g_{jk}(x)\eeq
  and
\beq
\label{lin wave eq source B}
F^{(2)}(x;{\bf p},{\bf q})
&=& M_{(2)} {\bf q^{\prime}}+ L_{(2)}{\bf z}+N^j_{(2)} g^{lk}\hat \nabla_l ( {\bf p}_{jk}+{\bf z}g_{jk}),
  \eeq
  where 
  \ba
  & &M_{(2)}=M_{(2)}(\hat \phi(x),\hat \nabla  \hat \phi(x),g(x)),\\
  & & L_{(2)}=L_{(2)}(\hat \phi(x),\hat \nabla  \hat \phi(x),g(x)),\\
  & &N^j_{(2)}=N^j_{(2)}(\hat \phi(x),\hat \nabla  \hat \phi(x),g(x))
  \ea
  are, in local coordinates, matrices whose elements  are smooth functions 
of  $\hat \phi(x),\hat \nabla  \hat \phi(x),$ and $g(x)$.
By Thm.\ \ref{thm: Good S functions} (i), the union of the image spaces of  the matrices 
$M_{(2)}(x)$ and
$L_{(2)}(x) $,  
and $N^j_{(2)}(x) $,  $j=1,2,3,4$,
span the space $\R^L$  for all $x\in U_{g,\sigma}$.
}

Consider $n\in \Z$, $t_0,s_0>0$, $Y=Y(x_0,\zeta_0;t_0,s_0)$, $K=K(x_0,\zeta_0;t_0,s_0)$
and $ (x,\xi)\in N^*Y$ (to recall the definitions of these notations, see formula (\ref{associated submanifold})
and definitions below it).
Let  $\mathcal Z=\mathcal Z (x,\xi)$ be the set of the values of 
the principal symbol $\tilde f(x,\xi)=(\tilde f_1(x,\xi),\tilde f_2(x,\xi))$,
at $(x,\xi)$,  
of the source ${f}=(f_1(x),f_2(x))\in \I^{n}(Y)$ that
 satisfy  the linearized conservation law for principal symbols (\ref{mu linearized conservation law for symbols}). 
 
 We use the following auxiliary result:


\begin{lemma}\label{lem: Lagrangian} 
Assume the the Condition A is satisfied and $\hat Q=0$  and $\hat P=0$.
Let $k_0\geq 8$,  $s_1\geq k_0+5$, and 
$Y\subset W_{g}$
be a 2-dimensional space-like submanifold and $y\in Y$, $\xi\in N^*_yY$,
and let $\W$ be a conic neighborhood of $(y,\xi)$ in $T^*M$. Also, let $y\in U_{g,\sigma}$
with some permutation $\sigma$.
Let us consider an open, relatively compact local coordinate neighborhood $V\subset  U_{g,\sigma}$ of $y$ such that
in the coordinates $X:V\to \R^4$, $X^j(x)=x^j$,
we have $X(Y\cap V)\subset \{x\in \R^4;\ x^jb^1_j=0,\  x^jb^2_j=0\}$.
Let  $n_1\in \Z_+$  be sufficiently large and  $n\leq -n_1$.
 Let us consider
 ${\bf p},{\bf q^{\prime}},{\bf z}\in \I^{n}(Y)$,  supported in $V$, 
 that  have classical symbols 
 with principal symbols $\tilde \symbolP ^{(a)}(x,\xi),$ $\tilde \symbolQ_{1}^{(a)}(x,\xi)$,
 $\tilde \symbolQ_1^{(a)}(x,\xi)$, correspondingly, at $(x,\xi)\in N^*Y$.
 Moreover, assume that the principal symbols of ${\bf p}$ and ${\bf z}$ 
 satisfy the linearized conservation law 
 for the principal symbols, that is, (\ref{eq: lineariz. conse. law symbols mLS}),  at 
 all $N^*Y\cap N^*K$ and assume that they vanish outside the conic neighborhood $\W$ of
 $(y,\xi)$ in $T^*M$. Let ${f}=(f^1,f^2) \in\I^{n}(Y)$ be given by
(\ref{eq: formula for sources 1}) and (\ref{eq: formula for sources 2}).

Then  the principal symbol $\tilde f(y,\xi)=(\tilde f_1(y,\xi),\tilde f_2(y,\xi))$ 
of the source ${f}$ at $(y,\xi)$ is the set $\mathcal Z=\mathcal Z(y,\xi)$.
 Moreover,  by varying ${\bf p},{\bf q^{\prime}},{\bf z}$  so that the linearized conservation law  (\ref{eq: lineariz. conse. law symbols mLS})  for principal symbols 
is satisfied, {the principal symbol  $\tilde f(y,\xi)$} at $(y,\xi)$
 achieves all values in the $(L+6)$ dimensional space $\mathcal Z$.

\end{lemma}

\noindent 
{\bf Proof.} 
{
Let us use local coordinates $X:V\to \R^4$ where $V\subset M_0$ is a neighborhood of $x$.
In these coordinates, let $\tilde \symbolP ^{(b)}(x,\xi)$ and $\tilde \symbolQ_{2}^{(b)}(x,\xi)$ denote the sub-principal
symbols of ${\bf p}$ and ${\bf z}$, respectively, at $(x,\xi)\in N^*Y$. 
 Moreover, let $\tilde \symbolP ^{(c)}_j(x,\xi)=\frac \p {\p x^j}\tilde \symbolP ^{(a)}(x,\xi)$ and $\tilde d^{(c)}_j(x,\xi)=\frac \p {\p x^j}\tilde d^{(a)}_2(x,\xi)$, $j=1,2,3,4$ be the $x$-derivatives of the principal symbols
and let us denote \ba
\tilde \symbolP ^{(c)}(x,\xi)=(\tilde \symbolP ^{(c)}_j(x,\xi))_{j=1}^4,\quad 
\tilde d^{(c)}(x,\xi)=(\tilde d^{(c)}_{j}(x,\xi))_{j=1}^4.
\ea


 Let ${f}=({f}_1,{f}_2)=F(x;{\bf p},{\bf q})$ 
 be defined by (\ref{lin wave eq source A1a}) and
 (\ref{lin wave eq source B}).
{When the principal symbols of  ${\bf p},{\bf q^{\prime}},{\bf z}\in \I^{n}(Y)$ 
 are such that  the linearized conservation law  (\ref{eq: lineariz. conse. law symbols mLS})  for principal symbols 
is satisfied, we see that  ${f}\in \I^{n}(Y)$ has the principal symbol
  $\tilde f(x,\xi)=(\tilde f_{1}(x,\xi),\tilde f_{2}(x,\xi))$ at $(x,\xi)$, given by  
  \ba
  & &\tilde  f_{1}(x,\xi) = s_1(x,\xi),\\
\nonumber 
& & \tilde  f_{2}(x,\xi)= s_2(x,\xi),
\ea
where   
     \beq\label{lin wave eq source symbols}
& &s_{1}(x,\xi) = (\tilde \symbolP ^{(a)}+g\tilde \symbolQ_1^{(a)})(x,\xi),\\
\nonumber 
& & s_{2}(x,\xi)= \bigg(
 M_{(2)}(x)\tilde \symbolQ_{1}^{(a)} + J_{(2)}(\tilde \symbolP ^{(c)}+g\tilde d^{(c)}) +\\
 \nonumber
 & &\quad +
L_{(2)}\tilde \symbolQ_1^{(a)} +N^j_{(2)} \, g^{lk}\xi_l (\tilde \symbolP_{1}^{(b)}+g\tilde \symbolQ_1^{(b)})_{jk}\bigg)(x,\xi).
  \eeq
Here, roughly speaking, the $J_{(2)}$ term appears when the $\nabla$-derivatives in $R$ hit to the symbols of the conormal distributions having the form
(\ref{eq: b b-prime BBB}).  
We emphasize that here the symbols $s_1(x,\xi)$  and $s_2(x,\xi)$ are 
well defined objects (in fixed local coordinates) also when   the linearized conservation law  (\ref{eq: lineariz. conse. law symbols mLS})  for principal symbols is not valid.  When (\ref{eq: lineariz. conse. law symbols mLS})
 is valid, ${f}\in \I^{n}(Y)$ and $s_1(x,\xi)$  and $s_1(x,\xi)$
coincide with the principal symbols of ${f}_1$ and ${f}_2$.}

   Observe that  the map $(c^{(b)}_{jk})\mapsto (g^{lk}\xi_lc^{(b)}_{jk})_{j=1}^4$, 
 defined as
   $\hbox{Symm}(\R^{4\times 4})\to \R^4$, is surjective.
  Denote 
  \ba
& &  \tilde m^{(a)}=(\tilde \symbolP ^{(a)}+g\tilde \symbolQ_1^{(a)})(x,\xi),
 \\
 & &\tilde m^{(b)}=(\tilde \symbolP ^{(b)}+g\tilde \symbolQ_1^{(b)})(x,\xi),
\\
 & &\tilde m^{(c)}=(\tilde \symbolP ^{(c)}+g\tilde d^{(c)})(x,\xi).
 \ea
 As noted above,
by {(\ref{eq: deriv1})},   the union of the image spaces of  the matrices 
$M_{(2)}(x)$ and
$L_{(2)}(x) $,  
and $N^j_{(2)}(x) $,  $j=1,2,3,4$,
span the space $\R^L$  for all $x\in U$.  
Hence the map 
  \ba
  {\bf A}:(\tilde m^{(a)},\tilde m^{(b)},\tilde m^{(c)},
  \tilde \symbolQ_{1}, \tilde \symbolQ_1^{(a)})|_{(x,\xi)}\mapsto
{(s_1(x,\xi),s_2(x,\xi)),}
  \ea 
  given by (\ref{lin wave eq source symbols}), considered as a 
  map ${\bf A}:\mathcal Y=(\hbox{Symm}(\R^{4\times 4}))^{1+1+4}\times \R^{K}\times \R\to
  \hbox{Symm}(\R^{4\times 4})\times \R^L$,
  is surjective. 
  Let $\mathcal X$ be  the set of
  elements $(\tilde m^{(a)}(x,\xi),\tilde m^{(b)}(x,\xi),\tilde m^{(c)}(x,\xi),
  \tilde \symbolQ_{1}^{(a)}(x,\xi),$ $\tilde \symbolQ_1^{(a)}(x,\xi))\in \mathcal Y$
  where 
   $\tilde m^{(a)}(x,\xi)=(\tilde \symbolP ^{(a)}+g\tilde \symbolQ_1^{(a)})(x,\xi)$ is such that the  pair  $(\tilde \symbolP ^{(a)}(x,\xi),\tilde \symbolQ_1^{(a)}(x,\xi))$  satisfies
 the linearized conservation law for principal symbols, see (\ref{eq: lineariz. conse. law symbols mLS}). 
 Then    $\mathcal X$ has codimension 4 in $Y$ and we see
 that the image  ${\bf A}(\mathcal X)$ has  in $ \hbox{Symm}(\R^{4\times 4})\times \R^L$ co-dimension less or equal
 to 4.
 
By (\ref{lin cons 2B})}) and the considerations above it, we have that $\bf f$ satisfies the linearized
 conservation law 
(\ref{eq: lineariz. conse. law PRE}). This implies that  its
principal symbol ${\bf A}((\tilde m^{(a)}(x,\xi),\tilde m^{(b)}(x,\xi),$ $ \tilde m^{(c)}(x,\xi),
  \tilde \symbolQ_{1}^{(a)}(x,\xi),$ $\tilde \symbolQ_1^{(a)}(x,\xi)))$ 
has to satisfy the linearized conservation law for principal symbols 
(\ref{mu linearized conservation law for symbols})
and  hence ${\bf A}(\mathcal X)\subset \mathcal Z$. 
As  $\mathcal Z$ has codimension 4, this and the above
prove that    ${\bf A}(\mathcal X)=\mathcal Z$.
\hfill \Box \medskip

Now we are ready to prove the microlocal stability result for 
the Einstein-scalar field equation (\ref{eq: adaptive model with no source}).
Note that the claim of the following theorem does not 
involve the  adaptive source functions constructed in Theorem \ref{thm: Good S functions}
as these functions are needed only as an auxiliary tool in the proof.

 \begin{theorem}\label{microloc linearization stability thm Einstein} ($\mu$-LS i.e., Microlocal linearization stability) 
Let $k_0\geq 8$,  $s_1\geq k_0+5$, and 
$Y\subset W_{g}$
be a 2-dimensional space-like submanifold, and $y\in Y$
and let $\W$ be a conic neighborhood of $(y,\eta)$ in $T^*M$. 
Also, 
consider in an open local coordinate neighborhood $V\subset W_{g}$ of $y\in Y$ such that
in the coordinates $X:V\to \R^4$, $X^j(x)=x^j$,
we have $X(Y\cap V)\subset \{x\in \R^4;\ x^jb^1_j=0,\  x^jb^2_j=0\}$.
Then there there is  $n_1\in \Z_+$ such that if $(y,\eta)\in N^*Y$ is light-like, $n\leq -n_1$ and  $\tilde c_{jk}(x,\theta)$ and $\tilde d_{\ell}(x,\theta)$
 are positively $n$-homogeneous in the $\theta$-variable in the domain $ \{(x,\theta)\in V\times \R^2;\ |\theta|>1\}$  
and satisfy 
\beq\label{mu linearized conservation law for symbols BB}
& &g^{lk}(\theta_1  b^1_l
+\theta_2  b^2_l)
\tilde c_{kj} (x,\theta_1,\theta_2)
=0,\quad |\theta|>1,\quad j=1,2,3,4,
\eeq
then there are $f^1\in \I^{n}(Y)$
and $f^2\in \I^{n}(Y)$
 supported in $V$ such that the principal symbols of these distributions vanish outside $\W$ and are equal to 
$\tilde c_{jk}(y,\eta)$  and $\tilde d_{\ell}(y,\eta)$  at $(y,\eta)$, respectively.
Moreover, 
 $f=(f^1,f^2)$  satisfies the linearized conservation law (\ref{eq: lineariz. conse. law PRE})
and  
 $f$  has the LS-property (\ref{eq: LSp}) 
 in $C^{s_1}(M_0)$
  with a family  $\F_\e$, $\e\in [0,\e_0)$
  such that all functions   $\F_\e$,$\e\in [0,\e_0)$ are supported in $V$.
 \end{theorem}

\noindent {\bf Proof.}
Let $\sigma\in \Sigma(K)$ be such that $y\in U_{g,\sigma}$.
Let
 ${\bf p}$ and ${\bf q}$  be the functions constructed in
 Lemma \ref{lem: Lagrangian} such that
they are supported in $V$ and their principal symbols vanish outside $\W$.  We can assume 
 that these functions are supported in $W_0= V\cap U_{g,\sigma}$. 
Let  $P_\e=\e {\bf p}$ and $Q_\e=\e {\bf q}$   be sources depending on $\e\in \R$
and $u_\e=(g_\e,\phi_\e)$ be the  solution of
 (\ref{eq: adaptive model}) with the sources  $P_\e$  and $Q_\e$.
 Also, let
\ba
& &\F_\e^1=P_\e+Z_\e g_\e ,\quad Z_\e=-(\sum_{\ell=1}^L S^\e_\ell\phi^\e_\ell+
\frac 1{2m^2}(S^\e_\ell)^2),\\
& &(\F_\e^2)_\ell=S^\e_\ell,
\ea
where 
\ba
& &S^\e_\ell= \mathcal S_\ell(g_\e,\phi_\e,\nabla \phi_\e,Q_\e,\nabla Q_\e,P_\e,
\nabla^{g_\e} P_\e),
\ea
where 
 $\mathcal S_\ell$ are the adaptive source functions constructed in Theorem \ref{thm: Good S functions} and
 its proof.

By (\ref{S-vanish condition}) also $S^\e_\ell$ and 
the family $\F_\e$, $\e\in [0,\e_0]$ of non-linear sources are supported in  $V_0$
and we have shown that  $u_\e=(g_\e,\phi_\e)$ and $\F_\e$
satisfy
the reduced Einstein-scalar field equation (\ref{eq: adaptive model with no source})
and  the conservation
 law (\ref{conservation law CC}).
This proves Theorem \ref{microloc linearization stability thm Einstein}.
\hfill \Box \medskip

  {\extension 
  
  \subsection*{C.3. A special case when the whole metric is determined}
  Finally, we consider the case when  $\hat Q$ and $\hat P$ are non-zero, and we define
\ba
{\cal D}^{mod}(g,\hat \phi,\e)=\{[(U_g,g|_{U_g},\phi|_{U_g},F|_{U_g})]&;&(g,\phi,F)\hbox{ are smooth 
solutions,}\\
& & 
\hspace{-6cm}\hbox{ of  
(\ref{eq: adaptive model}) with $F=(P-\hat P,Q-\hat Q)$, }F\in C^{19}_0(W_{g};\B^K),
\\
& &\hspace{-6cm}\hbox{$J_g^+(\supp(F))\cap J_g^-(\supp(F))\subset 
W_g, $ $\mathcal N^{({ {17}})}(F)<\e$, $\mathcal N^{({{17}})}_{g}(g)<\e$}\}.
\ea
Next we consider the case when $\hat Q^{(j)}$ and  $\hat P^{(j)}$ are not assumed to be zero,
see Fig.\ A11. 

\begin{corollary}\label{coro of main thm original Einstein}
Assume that $(\hattuM^{(1)} ,g^{(1)})$ and $(\hattuM^{(2)} ,g^{(2)})$ are  globally hyperbolic
manifolds and $\hat \phi^{(j)}$, $\hat Q^{(j)}$,  and $\hat P^{(j)}$ are background 
fields satisfying (\ref{eq: adaptive model}).
Also, assume that there are neighborhoods $U_{g^{(j)}}$, $j=1,2$ of time-like geodesics $\mu_j\subset \hattuM^{(j)}$ where the  Condition $\mu$-LS is valid,
and points $p_j^-=\mu_j(s_-)$ and $p_j^+=\mu_j(s_+)$.
Moreover,
 assume  
also that for $j=1,2$  there are sets $W_j\subset M_j$ such that 
$\hat \phi^{(j)}$, $\hat Q^{(j)}$ and $\hat P^{(j)}$ are zero (and thus
the metric tensors $g_j$ have vanishing Ricci curvature)
in $W_j$ and that  $I_{g_j}(\hat p^-_j,\hat p^+_j)\subset W_j \cup U$. 
If there is $\e>0$  such that 
\ba
{\cal D}^{mod}(g^{(1)},\hat \phi^{(1)},\e)={\cal D}^{mod}(g^{(2)},\hat \phi^{(2)},\e)
\ea
then the metric $\Psi^*g_2$ is isometric to $g_1$ in $I_{g_1}(\hat p^-_1,\hat p^+_1)$.
\end{corollary}
\medskip

\begin{center}\label{Fig-4}

\psfrag{1}{}
\psfrag{2}{}
\psfrag{3}{}
\psfrag{4}{$W_{g}$}
\psfrag{5}{$\mu_{a}$}
\includegraphics[width=4.5cm]{vacuum2.eps}
\end{center}
{\it FIGURE A11: A schematic figure where the space-time is represented as the  2-dimensional set $\R^{1+1}$ on the setting in Corollary \ref{coro of main thm original Einstein}.
The set $J_{g}(\hat p^-,\hat p^+)$, i.e., the diamond with the black boundary, is contained in the union of the blue set $U$
and the set $W$. The set $W$ is in the figure the area outside of the green curves. 
The sources are controlled in the set $U\setminus W$ and  the set $W$ consists of vacuum.
}

\medskip

\noindent
{\bf Proof.}\  (of Corollary \ref{coro of main thm original Einstein})  
{In the above proof of  Theorem \ref{alternative main thm Einstein}, we used the  
 assumption that $\hat Q=0$ and $\hat P=0$ to obtain  
 equations (\ref{w1-3 solutions}). In the   
setting of Cor.\ \ref{coro of main thm original Einstein}  
 where the background source fields  
 $\hat Q$ and $\hat P$  
  are not zero,  
we need to assume   
in the computations related to sources (\ref{eq: f vec e sources})  
that there are  
neighborhoods $V_j$ of the geodesics $\gamma_j$ that satisfy  
$\supp({f}_j)\subset V_j$,  
the linearized waves $u_j={\bf Q}_{g}{f}_j$ satisfy  singsupp$(u_j)\subset V_j$,  
and     
$V_i\cap V_j\cap (\supp(\hat Q)\cup \supp(\hat P))=\emptyset$ for $i\not =j$.}  
To this end, we have first consider  measurements for the linearized waves   
and check for given $(\vec x,\vec\xi)$ that no two geodesic $\gamma_{x_j,\xi_j}(\R_+)$  
do intersect at $U$ and restrict all considerations for such $(\vec x,\vec\xi)$.  
Notice that such $(\vec x,\vec\xi)$ form an open and dense set in   
$(TU)^4$. If then the sources are supported in balls so small
that each ball is in some set $U_{g,\sigma}$ and the width $\hat s$ of the used  
distorted plane  waves is chosen to be small enough, we see that condition   
$V_i\cap V_j\cap (\supp(\hat Q)\cup \supp(\hat P))=\emptyset$ is satisfied.  
  
The above restriction causes only minor modifications in the above
proof and thus, mutatis mutandis, we see that we can determine  
the conformal  
type of the  metric in relatively compact sets 
$I_{g}(\hat\mu(s^\prime),\hat\mu(s^{\prime\prime}))\setminus (\supp(\hat Q)\cup \supp(\hat P))$ for all $s_-<s^\prime<s^{\prime\prime}<s_+$.  
By glueing these  manifolds and $U$ together, we find the conformal  
type of the  metric in $I_{g}(\hat p^-,\hat p^+)$.  
After this the claim follows similarly to the  
 proof of  Theorem \ref{alternative main thm Einstein} using 
 Corollary 1.3 of \cite{Paper-part-2}.  
\hfill \Box \medskip

In the setting of Corollary \ref{coro of main thm original Einstein} the set
 $W$ is such $I(\hat p^-,\hat p^+)\cap (M\setminus W)\subset U$. This means that if we restrict to the domain $I(\hat p^-,\hat p^+)$
then we have the Vacuum Einstein equations in the unknown domain $I(\hat p^-,\hat p^+)\setminus U$
and have matter only in the domain $U$ where we implement our measurement (c.f.\ a space
ship going around in a system of black holes, see \cite{ChM}). This could be considered as an ``Inverse problem
for the vacuum Einstein equations''. 

  }

\subsection*{C.4. Motivation using Lagrangian formulation}

To motivate the system (\ref{eq: adaptive model}) of partial differential equations, we 
give in this subsection a non-rigorous discussion.

Following  \cite[Ch. III, Sect. 6.4, 7.1, 7.2, 7.3]{ChBook}
and  \cite[p. 36]{AnderssonComer}
we start by considering  the 
Lagrangians, associated to gravity, scalar fields $\phi=(\phi_\ell)_{\ell=1}^L$ and non-interacting fluid fields, that is, the number density four-currents ${\bf n}=({\bf n}_\kappa(x))_{\kappa=1}^J$
(where each ${\bf n}_\kappa$ is a vector field,  see \cite[p. 33]{AnderssonComer}).
We consider also products of  vector fields ${\bf n}_\kappa(x)$ 
and $\frac 12$-density $|\det(g)|^{1/2}$ that denote ${\bf p}_\kappa$,
\beq\label{half density}
{\bf p}_\kappa^j(x)\frac\p{\p x^j}={\bf n}^j_\kappa(x)|\det(g)|^{1/2}\frac\p{\p x^j}.
\eeq
see  \cite[p. 53]{Dirac}. Also, $\rho=(-g_{jk}{\bf n}^j_\kappa{\bf n}^k_\kappa)^{1/2}$ 
corresponds to the energy density of the fluid.
Below, we use the variation of density with respect to the metric,
\beq\label{eee}
\frac {\delta}
{\delta g_{jk}}( 
\sum_{\kappa=1}^J  (-g_{nm}{\bf p}^n_\kappa{\bf p}^m_\kappa)^{1/2})
&=&- \sum_{\kappa=1}^J   \frac 12 (-g_{nm}{\bf p}^n_\kappa{\bf p}^m_\kappa)^{-1/2}
{\bf p}_\kappa^j
{\bf p}_\kappa^k
\\ \nonumber
&=& - \sum_{\kappa=1}^J   \frac 12 \rho\,{\bf n}_\kappa^j
{\bf n}_\kappa^k\, |\det(g)|^{1/2}. 
\eeq
Due to this, we denote  
\beq\label{P formula}\\
\nonumber
P&=& \sum_{\kappa=1}^J   \frac 12 \rho\,{\bf n}^\kappa_j
{\bf n}^\kappa_k dx^j\otimes dx^k,\quad\hbox{where }
{\bf n}^\kappa_k =g_{ki}  {\bf n}_\kappa^i=g_{ki}  {\bf p}_\kappa^i|\det(g)|^{-1/2}.
\eeq

Below, we consider a model for $g$, $\phi$, and ${\bf p}$.
We also add in to the model a Lagrangian  associated with
some scalar valued source fields $S=(S_\ell)_{\ell=1}^L$  and
$Q=(Q_k)_{k=1}^K$. We consider the action corresponding to
the coupled Lagrangian\hiddenfootnote{
An alternative Lagrangian would be the one where the term
$P_{pq}g^{pq}$ is replaced by $\frac 12 P_{pq}P_{jk}g^{pj}g^{qk}$.
Then in the stress energy tensor $P_{jk}$ is replaced by
${\mathcal T}^{(1)}_{jk}= F_{jp}F_{kq}g^{qp}$, i.e. ${\mathcal T}^{(1)}=\hat F\hat G\hat F.$
By [Duistermaat, Prop. 1.3.2], 
its derivative with respect to $P_{jk}$, that is,
$((D_F{\mathcal T}^{(1)}_{jk}|_{\hat F})^{ab})$ (COMPUTE THIS)
is a surjective pointwise map is surjective at points $x\in M$ if
$\hat F_{qp}$ is positive definite symmetric (and thus invertible) matrix (CHECK THIS!!).
}
%
\ba
& &\mathcal A=\int_{M}\bigg(L_{grav}(x)+L_{fields}(x)+L_{source}(x)\bigg)\,dV_g(x),\\
& &L_{grav}=\frac 12 R(g),\\
& &L_{fields}=\sum_{\ell=1}^L\bigg(-\frac 12 g^{jk}\p_j\phi_\ell \,\p_k\phi_\ell 
-{\mathcal V}( \phi_\ell;S_\ell)\bigg)
+  \\
& &\quad\quad\quad\quad\quad+\sum_{\kappa=1}^J\bigg(- \frac 12 (-g_{jk}{\bf p}^j_\kappa {\bf p}^k_\kappa)^{\frac 12}
\bigg)|\det (g)|^{-\frac 12},
\\
& &L_{source}=\e \mathcal H_\e(g,S,Q, {\bf p},\phi),
\ea
where $R(g)$ is the scalar curvature, $dV_g=(-\det (g))^{1/2}dx$ is the volume form on $(M,g)$, 
\beq\label{eq: Slava-a}
{\mathcal V}( \phi_\ell;S_\ell)=\frac 12m^2\bigg(\phi_\ell +\frac 1{m^2}S_\ell\bigg)^2
\eeq are energy potentials
of the scalar fields $\phi_\ell$ that depend on $S_\ell$,
 and  $\mathcal H_\e(g,S,Q, {\bf p},\phi)$ is a function modeling the measurement
 device we use.
 We assume that $\mathcal H_\e$ is bounded and   its derivatives with respect
 to $S,Q,{\bf p}$  are very large (like of order $O((\e)^{-2})$) and its derivatives
with respect of $g$ and $\phi$ are bounded when $\e>0$ is small.
We note that the  above Lagrangian for the fluid fields is the sum of the single fluid Lagrangians.
where for all fluids the master function  $\Lambda(s)=s^{1/2}$, that is, the energy density
of each fluid is given by
$\rho=\Lambda(-g_{jk}{\bf n}^j{\bf n}^k)$. On fluid Lagrangians, see the discussions in
 \cite[p. 33-37]{AnderssonComer},  \cite[Ch. III, Sect. 8]{ChBook}, \cite[p. 53]{Dirac}, and
 \cite{Taub} and \cite[p.\ 196]{Felice}.

%
When  we compute the critical points of the Lagrangian $L$ and 
neglect the  $O(\e)$-terms, the equation  $\frac {\delta  \mathcal A}{\delta g}=0$,
together with formulas (\ref{eee}) and  (\ref{P formula}), 
give the Einstein equations  with a stress-energy tensor $T_{jk}$
defined in (\ref{eq: T tensor1}).
The equation  $\frac {\delta  \mathcal A}{\delta \phi}=0$
gives the wave equations with sources $S_\ell$. We assume that $O(\e^{-1})$ order equations
obtained from 
the equation $(\frac {\delta \mathcal A}{\delta S},\frac {\delta  \mathcal A}{\delta Q},
\frac {\delta  \mathcal A}{\delta {\bf p}})=0$ 
fix the values of the scalar functions $Q$ and 
the fields ${\bf p}^\kappa$, $\kappa=1,2,\dots,J$, and moreover, 
  yield for the sources $S=(S_\ell)_{\ell=1}^L$ equations of the form
\beq\label{S! formula}
S_\ell=\mathcal S_\ell(g,\phi, \nabla
 \phi,Q,\nabla Q,P,\nabla^g P)
 \eeq
where $P$ is given by (\ref{P formula}).
Let us aslo write (\ref{S! formula}) using different notations, as
\ba
S_\ell=Q_\ell+\mathcal S_\ell^{2nd}(g,\phi, \nabla
 \phi,Q,\nabla Q,P,\nabla^g P).
 \ea 
 Summarizing, we have obtained, up to the above used approximations,
the model  (\ref{eq: adaptive model}). However, note that above the field $P$  is not directly controlled 
but instead, we control ${\bf p}$ and the value of the field $P$ is determined by
the solution $n$ and  formula  (\ref{P formula}). In this sense $P$  is not controlled,
but an observed field. 

 Above, the function
 $\mathcal H_\e$ models the way the measurement device works. 
Due to this we will assume that $\mathcal H_\e$ and thus functions ${\mathcal S}_\ell$ may
be quite complicated.   The interpretation of the above is that in  each measurement event we use a device that
fixes the values of the scalar functions $Q$, ${\bf p}$, and gives the equations
for $S^{2nd}$ that 
 tell how the sources of the $\phi$-fields adapt to these changes
 so that 
the physical conservation laws are satisfied.

}

\canberemoved{
\noindent
{{\bf Remark 5.1.}  {\gtext The proof of 
Theorem  \ref{main thm3 new} 
can be used to do an approximative reconstruction of the set $I(p^-,p^+)=I^+(p^-)\cap I^-(p^+)$ and its conformal structure with only one
measurement under certain geometric conditions, that is, one can use one suitably constructed source $f$,
measure $L_U(f)$, and construct an approximation
of the conformal class of the metric $g$ in $I(p^-,p^+)$. 
In the proof above we showed that it is possible
to use the non-linearity to create an artificial point source at 
a point $q\in I(p^-,p^+)$.
Using the same method it is possible to create  an arbitrary finite number
of artificial point sources with a single source $f$. For this, we can use 
a source $f$ that is  
a sum of a large number of modified point sources in $U$ that produce 
spherical waves having different order of smoothnesses. The interaction
of such spherical waves produces even a larger
number of artificial point sources in $I(p^-,p^+)$ that all have a different order of smoothness.}
 As a very rough analogy,
we can produce in $I(p^-,p^+)$ an arbitrarily dense collection
of artificial points sources that all have different ``colors''.
Using this we
can show using similar methods to those in \cite{AKKLT,HLO} that
the measurement with one  source
$f$ determines,
in a suitable  compact class of Lorentzian manifolds 
having no conjugate points,
an approximation (in a suitable sense) of the conformal class of the
 manifold $(J(p^-,p^+),g)$. The details of this construction will
 be considered elsewhere.}
}

\medskip

\noindent{\bf Acknowledgements.} The authors express their gratitude to
 MSRI, the Newton  
Institute, the Fields Institute and  
the Mittag-Leffler Institute, where parts of this work have been done.

YK was partly supported by EPSRC and the AXA professorship.
ML was partly supported by the Finnish Centre of
Excellence in Inverse Problems Research 2012-2017 and an Academy Professorship.
GU was partly supported  by NSF, a Clay Senior Award at MSRI,  a Chancellor Professorship at UC Berkeley,  
a Rothschild Distinguished Visiting Fellowship at the Newton Institute, the Fondation de Sciences Math\'ematiques de Paris,
FiDiPro professorship,
and a Simons Fellowship.



 \end{document}